\documentclass[12pt,a4paper]{article}
\usepackage{amsmath, amssymb, amsthm}
\usepackage{graphicx,amssymb,amsfonts,amsbsy,amsmath,amsthm}
\usepackage{mathtools}
\usepackage{hyperref}
\usepackage[utf8]{inputenc}
\usepackage{floatrow}
\usepackage{float}
\usepackage{subfig}
\usepackage{amssymb}
\usepackage{amsthm}
\usepackage{epstopdf}
\usepackage{cite}
\usepackage{anysize}
\usepackage{color,soul}
\usepackage{graphicx}
\usepackage{tabularx}
\usepackage{multirow}
\usepackage{float} 
\floatstyle{plaintop}
\restylefloat{table}
\usepackage{soul}
\usepackage{ulem}
\usepackage{authblk} 
\hypersetup{colorlinks,linkcolor={blue},citecolor={blue},urlcolor={black}}

\marginsize{1cm}{1cm}{1cm}{1cm} 

\input epsf.sty
\usepackage{lineno}
\usepackage[english]{babel}

\numberwithin{equation}{section}
\newtheorem{theorem}{Theorem}[section]
\newtheorem{lemma}[theorem]{Lemma}

\newtheorem{proposition}[theorem]{Proposition}

\newtheorem{remark}[theorem]{Remark}

\title{Nonsmooth bifurcations and delay-induced stability switching in a Filippov system with dynamic resource capacity}       
\author[1,2]{Sarita Bugalia}  
\author[3]{Manoj Kumar}  
\author[2]{Jai Prakash Tripathi 
	\thanks{Corresponding author: \texttt{jtripathi85@gmail.com}}} \affil[1]{School of Animal and Comparative Biomedical Sciences, 
	College of Agriculture, Life and Environmental Sciences, 
	The University of Arizona, Tucson, Arizona 85721, USA} \affil[2]{Department of Mathematics, 
	Central University of Rajasthan, 
	Bandar Sindri, Kishangarh, Ajmer 305817, Rajasthan, India} 
\affil[3]{School of Engineering and Science, 
	Indian Institute of Technology Madras, Zanzibar Campus, 
	Bweleo, Zanzibar, Tanzania} 

\date{}

\begin{document}
	
	\maketitle
	
	\begin{abstract}
		Threshold-dependent switching coupled with delayed feedback and 
		dynamically evolving resource constraints can give rise to distinct 
		nonlinear phenomena. Motivated by infectious disease transmission 
		under healthcare resource constraints, we develop and analyze a 
		time-delayed Filippov system coupling disease transmission, 
		threshold-based hospitalization, and dynamic healthcare capacity. 
		The per capita hospitalization rate decreases discontinuously at 
		the threshold, and any sliding region on the switching manifold 
		$\Sigma=\{I=I_0\}$ is necessarily repelling, with its width 
		determined by the evolving healthcare capacity. We establish conditions for a nondegenerate subthreshold fold in the smooth extension of the above-threshold subsystem and show that discontinuity-induced bistability can occur through the coexistence of an admissible stable endemic equilibrium and the locally stable disease-free equilibrium when the basic reproduction number $\mathcal{R}_0$ is below unity. 
		The admissible positive-capacity pseudo-equilibrium is globally 
		asymptotically stable in the algebraically extended reduced 
		tangential dynamics but plays a saddle-type organizing role in 
		the full Filippov system. The capacity-response delay can induce 
		repeated stability switching of the endemic equilibrium. A 
		worst-case reproduction number $\mathcal{R}_0^{\max}$ provides a 
		sufficient condition for disease elimination across both 
		hospitalization regimes. In the delay-free setting, we establish 
		a transversal codimension-two Hopf-boundary interaction and 
		show that the Hopf-generated periodic orbit possesses a local 
		first-contact curve quadratically tangent to the Hopf locus. 
		Numerical investigations illustrate crossing periodic orbits near 
		quadratic tangencies, boundary-equilibrium transitions, and 
		repeated delay-induced stability switching. An illustrative 
		physics-informed neural network calibration to cumulative 
		COVID-19 cases from Kenya provides a data-informed 
		parameterization. These results demonstrate how dynamic resource 
		availability and delayed feedback influence switching-induced 
		dynamics, equilibrium stability, and disease persistence. 
	\end{abstract} 
	
	
	
	\section{Introduction}
	
	Piecewise-smooth dynamical systems provide a natural framework for
	modeling processes governed by threshold-dependent mechanisms. In
	Filippov systems, the governing vector field is piecewise smooth and
	changes when the solution reaches a prescribed switching manifold
	\cite{FILI2013,BERN2008}. Unlike smooth models, Filippov systems can
	exhibit non-smooth phenomena such as sliding modes,
	boundary-equilibrium transitions, grazing events, and other
	discontinuity-induced bifurcations
	\cite{Kuznetsov2003,Kowalczyk2005,Kowalczyk2006,Colombo2012}.
	Time delays arising from information transmission, response
	mechanisms, or adjustment processes can further influence the
	stability and bifurcation structure of dynamical systems
	\cite{Hale1993,Sieber2006}.
	
	An important mathematical problem arises when discontinuous
	switching is coupled with a dynamically evolving auxiliary variable
	whose evolution depends on the past states of the system. Such a
	variable may influence both the stability of the smooth subsystems
	and the geometry of the sliding and crossing regions on the switching
	manifold. In the present work, we investigate this interaction in a
	delayed Filippov system with dynamic resource capacity, using
	infectious disease transmission under healthcare resource constraints
	as the primary motivating application. In particular, the switching
	mechanism depends only on the current infection level, while the
	time delay enters through the continuous healthcare-capacity
	equation.
	
	Emerging and re-emerging contagious diseases, such as Ebola virus
	disease \cite{GOOD2014,WHOEBOLA2014}, severe acute respiratory
	syndrome (SARS) \cite{DONN2002,GUM2002}, Zika virus
	\cite{AMAN2016}, the 2009 H1N1 influenza epidemic
	\cite{FRAS2009,TANG2010,XIAO2015}, foot-and-mouth disease
	\cite{SUN2022}, and COVID-19, have become significant sources of
	morbidity and mortality during crises. The severity of an epidemic
	places immense burdens on healthcare organizations, especially
	in areas with limited medical resources. An effective strategy
	during an epidemic is to marshal as many medical resources as
	possible, including healthcare personnel, medicine, and hospital
	beds. For instance, during the initial phase of the COVID-19
	outbreak in Wuhan, China, the construction of square-cabin
	hospitals contributed to strengthening the healthcare response
	\cite{JLI2020}. Thus, investing in medical resources can profoundly
	influence the management and spread of contagious diseases.
	
	The impact of limited medical resources on disease control has
	been extensively studied using compartmental epidemic models
	\cite{SHAN2014,SHAN2016,YU2023,ZHAO2020,LIGH2017}.
	Shan et al.~\cite{SHAN2014,SHAN2016} proposed SIR-type models
	in which the treatment rate is described by a continuous function
	of hospital-bed availability, and demonstrated rich bifurcation
	phenomena, including backward, saddle-node, Hopf, and cusp-type
	Bogdanov-Takens bifurcations. Yu et al.~\cite{YU2023}
	developed an SEIHR model in which the hospitalization rate
	depends on available hospital beds and showed that the system
	can undergo backward and Hopf bifurcations, as well as a
	Bogdanov-Takens bifurcation. Similarly, Zhao et al.~\cite{ZHAO2020}
	analyzed an epidemic model with limited hospital beds and
	identified backward and Hopf bifurcations together with a
	Bogdanov-Takens bifurcation. Collectively, these studies highlight
	the important role of medical resource constraints in shaping
	epidemic dynamics and informing disease control strategies.
	
	In epidemic modeling, threshold-dependent switching arises when
	public-health interventions, treatment allocation, or resource-dependent
	control measures are activated or modified as disease prevalence
	crosses a prescribed threshold. Such intervention policies have been
	widely investigated in compartmental epidemic models
	\cite{WANG2013,CHONG2016,WANG2014,CHEN2018,DENG2021,
		WANG2019,WANG2025}. Despite these advances, relatively few studies
	have examined the interplay between medical resource constraints
	and threshold-based treatment policies within Filippov epidemic
	models. Wang et al.~\cite{WANG2018} proposed a model with a constant
	treatment rate below a critical infection level and capacity-dependent
	treatment above it, identifying boundary saddle-node, boundary node
	or focus, grazing, and sliding homoclinic bifurcations.
	Mu et al.~\cite{MU2019} investigated a Filippov model for avian
	influenza A(H7N9), showing that equilibrium stability depends
	sensitively on the tolerance threshold. Qin et al.~\cite{QIN2016}
	analyzed selective treatment strategies under medical resource
	constraints and identified boundary saddle, boundary node,
	crossing, buckling, and grazing bifurcations. 
	
	However, in reality, medical resources deployed to curb infections
	are often limited, and intervention strategies may need to vary
	depending on both infection levels and time. Green \cite{GREEN2002}
	estimated bed unavailability in obstetrics and intensive care
	units using data from New York State and showed that hospitals
	must plan capacity based on timely placement rather than target
	occupancy. Moreover, expanding healthcare capacity is often
	time-consuming, costly, and labor-intensive, leading to a delay
	between the availability of infection information and the
	implementation of corresponding healthcare-capacity adjustments.
	From a dynamical-systems perspective, the combination of
	non-smooth threshold switching, dynamically evolving healthcare
	capacity, and time delay introduces analytical challenges that
	are not captured by either smooth delay models or Filippov
	systems with fixed resource availability.
	
	The goal of this work is to determine how threshold-based hospitalization and delayed healthcare-capacity response jointly affect equilibrium stability, admissibility, sliding geometry, and disease persistence. We analyze the two smooth subsystems and the induced Filippov dynamics, identifying a subthreshold fold, discontinuity-induced bistability, and repeated delay-induced stability switching. We establish the structurally repelling nature of sliding, characterize the tangential stability of pseudo-equilibria, and derive a regime-independent sufficient condition for disease elimination.

	A further contribution is the analytical characterization of a transversal codimension-two Hopf-boundary interaction in the delay-free system and the quadratic first-contact geometry of its Hopf-generated periodic-orbit family. The nondegeneracy and transversality conditions are verified numerically at the identified interaction point. Supporting simulations and an illustrative PINN-based calibration to cumulative COVID-19 data from Kenya complement the analysis.
	
	The remainder of the paper is organized as follows. 
	Section~\ref{proposedmodel} introduces the model.
	Section~\ref{SUBDYN} analyzes the two smooth subsystems and
	their equilibrium and bifurcation structure.
	Section~\ref{SLIDDYN} studies the non-smooth dynamics on the
	switching manifold, including pseudo-equilibria, tangency points,
	boundary equilibria, and the codimension-two Hopf-boundary
	interaction. 
	Section~\ref{Est_param} presents the data-informed parameter
	calibration using physics-informed neural networks.
	Section~\ref{NumSim} presents numerical simulations.
	Finally, Section~\ref{CONCLUSION} provides a discussion,
	biological interpretation, and concluding remarks.
	Additional mathematical proofs, technical calculations, and
	supporting numerical results are provided in the supplementary
	material. 
	
	\section{A delayed Filippov model with dynamic healthcare capacity}\label{proposedmodel}
	
	At time $t>0$, the total human population, $N(t)$, is divided into five different compartments: susceptible individuals, $S(t)$, who can contract the disease; exposed individuals, $E(t)$, who have contracted the disease but are typically not infectious \cite{WANG2022}; reported infected individuals who are infectious, $I(t)$; hospitalized individuals, $H(t)$; and recovered individuals, $R(t)$. Hence, 
	\begin{equation}
		N(t) = S(t) + E(t) + I(t) + H(t) + R(t). 
	\end{equation}
	
	In the model, susceptible individuals are recruited through birth or immigration and through waning immunity from the recovered class, and are removed by infection and natural mortality. Exposed individuals arise through disease transmission and progress to the infectious class at a constant rate, while also experiencing natural mortality. The infectious population increases through progression from the exposed class and decreases due to natural death, disease-induced mortality, recovery, and hospitalization. Hospitalized individuals are admitted from the infectious class and leave through recovery or mortality. Recovered individuals accumulate through recovery from both infectious and hospitalized classes and lose immunity at a constant rate, returning to the susceptible class. We include disease-induced mortality to account for infections in which disease-related deaths play an important role. 
	
	In the absence of infection-driven expansion, we assume that the positive effective healthcare capacity adjusts logistically toward a baseline capacity level $K$, with intrinsic adjustment rate $r$. During an epidemic, additional capacity may be mobilized in response to the reported infection burden, with the per capita expansion rate increasing and saturating as the reported infection level increases. To represent this saturation effect, we use the Holling type II functional response. In this model, $\lambda$ represents the maximal per capita infection-responsive expansion rate of healthcare  capacity, and $\omega$ denotes the half-saturation constant, that is, the infection level at which half of the maximal infection-responsive expansion rate $\lambda$ is reached. Thus, in the absence of delay, the rate of change in healthcare capacity can be described as follows: 
	\begin{equation}
		\frac{dB}{dt} = rB\left(1-\frac{B}{K} \right) + \frac{\lambda I B}{\omega + I}.
	\end{equation}
	
	Administrative, logistical, and implementation processes may introduce a time lag between reported infection information and the corresponding adjustment of healthcare capacity. In the proposed model, we assume that the growth of healthcare  capacity at time $t$ is influenced by the previously reported number of infected individuals, represented by $I(t-\tau)$, where $\tau\geq0$ denotes the healthcare-capacity response delay. Given these assumptions, we formulate the rate of change in healthcare capacity as follows:
	\begin{equation}
		\frac{dB}{dt} = rB\left(1-\frac{B}{K} \right) + \frac{\lambda I(t-\tau) B}{\omega + I(t-\tau)}.
	\end{equation}
	
	The rate at which infected individuals can be hospitalized depends on both
	the infection burden and the available healthcare capacity. We assume that
	when the number of reported infected individuals remains below or equal to
	a prescribed threshold $I_0>0$, the maximum per capita hospitalization rate
	$\mu_1$ is applied, reflecting relatively unconstrained healthcare
	resources. When the infection burden exceeds $I_0$, the hospitalization
	rate becomes capacity dependent, increasing with available healthcare
	capacity and, for positive healthcare capacity, decreasing as the infection
	burden grows relative to that capacity. 
	
	Here, $B(t)$ represents effective mobilizable healthcare capacity expressed
	on an individual-equivalent scale, making $B$ commensurate with the
	infectious population $I$ in the ratio $B/(B+I)$. Let $\mu_0$ and $\mu_1$ denote the minimum and maximum per capita
	hospitalization rates, with $0<\mu_0<\mu_1$, corresponding to
	strongly constrained and relatively sufficient healthcare capacity,
	respectively. We define
	\begin{equation}\label{treatrate}
		T(I,B) =
		\left\{
		\begin{array}{ll}
			\mu_1, & I\leq I_0, \\[2mm]
			\displaystyle
			\mu_0+(\mu_1-\mu_0)\frac{B}{B+I},
			& I>I_0.
		\end{array}
		\right.
	\end{equation}
	
	At $I=I_0$, we assign $T(I_0,B)=\mu_1$ according to the
	lower-infection regime. This pointwise convention does not affect the
	Filippov dynamics on the switching manifold, which are determined from
	the two one-sided vector fields through the Filippov set-valued
	regularization. For every finite $B\geq0$, 
	\[
	T(I_0^+,B)
	=
	\mu_0+(\mu_1-\mu_0)\frac{B}{B+I_0}
	<
	\mu_1
	=
	T(I_0^-,B),
	\]
	so the hospitalization law exhibits a downward jump as the infection
	level crosses the threshold from $I<I_0$ to $I>I_0$. 
	
	\begin{figure}[H]
		\centering
		\includegraphics[scale=0.35]{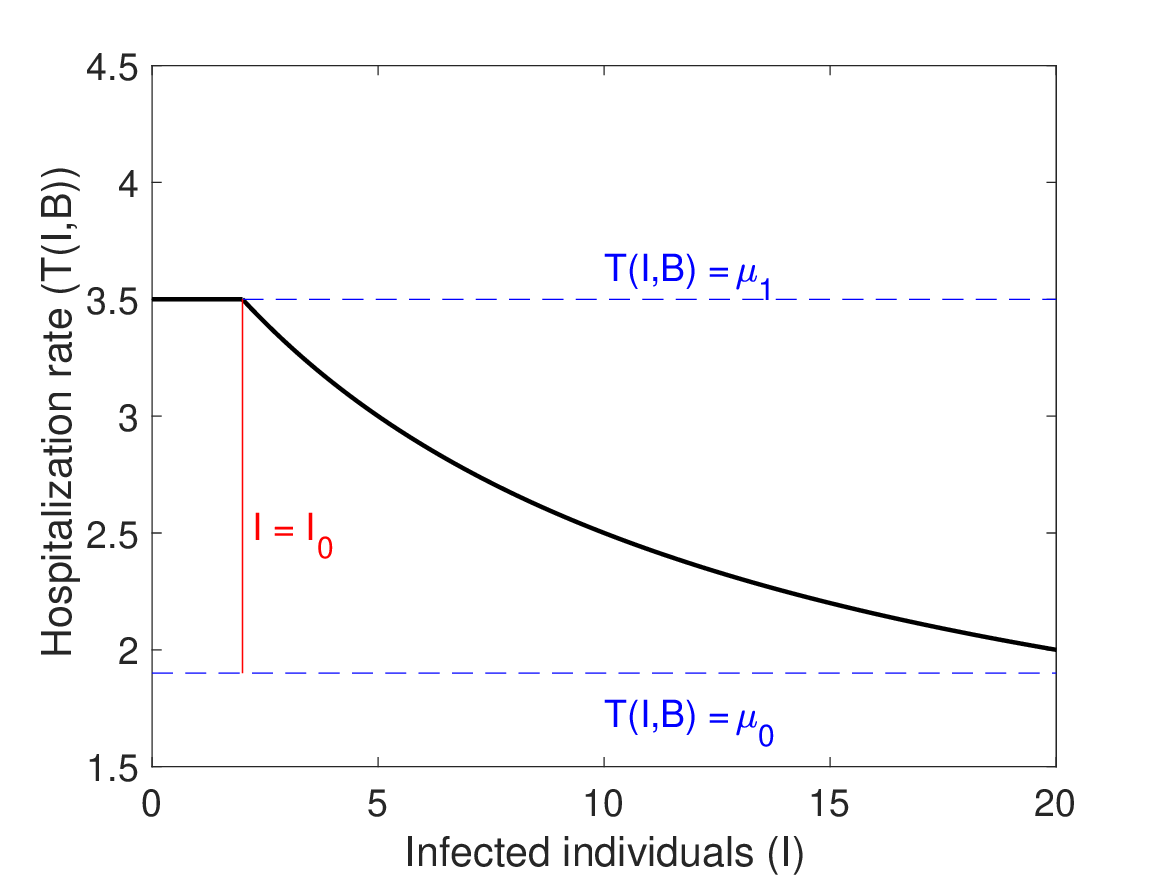}
		\caption{The curve of the hospitalization rate $T(I,B)$ defined in \eqref{treatrate}.}
		\label{traet1}
	\end{figure}
	
	The hospitalization function \eqref{treatrate} is illustrated in
	Fig.~\ref{traet1}. This formulation is closely related to the
	threshold-based treatment functions proposed in
	\cite{LIGH2017,WANG2018}, with the key distinction that the healthcare
	capacity $B(t)$ evolves dynamically. 
	
	Based on this formulation, we obtain the following properties of $T(I,B)$:
	\begin{itemize}
		\item $T(I,B)>0$ for $I\ge0$ and $B\ge0$, with $T(0,B)=\mu_1$ and, for $I>I_0$, $T(I,0)=\mu_0$.
		\item For $I>I_0$, $\partial T(I,B)/\partial I\leq0$, with strict inequality whenever $B>0$, and $\lim_{I\to\infty}T(I,B)=\mu_0>0$ for fixed $B\geq0$.
		\item For $I>I_0$, $\partial T(I,B)/\partial B>0$, with $\lim_{B\to\infty}T(I,B)=\mu_1$ and $\lim_{B\to0^+}T(I,B)=\mu_0$.
	\end{itemize}
	
	Given the above formulation, the disease transmission dynamics incorporating healthcare capacity and a response delay $\tau$ are described by the following delay differential equations: 
	
	\begin{equation}\label{ddemodel}
		\begin{aligned}
			\frac{dS}{dt} &= \Lambda - \beta SI - \mu S + \xi R, \\
			\frac{dE}{dt} &= \beta SI - (\mu + a)E, \\
			\frac{dI}{dt} &= aE - (\mu + \delta_1 + \gamma_1)I - T(I,B)I, \\
			\frac{dH}{dt} &= T(I,B)I - (\mu + \delta_2 + \gamma_2)H, \\
			\frac{dR}{dt} &= \gamma_1 I + \gamma_2 H - (\mu + \xi)R, \\
			\frac{dB}{dt} &= rB\left(1-\frac{B}{K}\right) + \frac{\lambda I(t-\tau)B}{\omega + I(t-\tau)},
		\end{aligned}
	\end{equation}
	
	with initial conditions
	$
	S(\theta)=\phi_1(\theta),
	E(\theta)=\phi_2(\theta), 
	I(\theta)=\phi_3(\theta),
	H(\theta)=\phi_4(\theta), 
	R(\theta)=\phi_5(\theta),
	B(\theta)=\phi_6(\theta),
	$
	where $\theta\in[-\tau,0]$,
	$
	(\phi_1,\phi_2,\phi_3,\phi_4,\phi_5,\phi_6)
	\in C([-\tau,0],\mathbb{R}_+^6),
	$ 
	and 
	$C([-\tau,0],\mathbb{R}_+^6)$ denotes the space of continuous non-negative history functions from $[-\tau,0]$ to $\mathbb{R}_+^6$, equipped with the supremum norm. Biological descriptions of the parameters in model \eqref{ddemodel} are provided in Table~\ref{table_1}.
	
	\begin{table}[H]
		\centering
		\caption{Biological interpretations of model parameters.}
		\label{table_1}
		\begin{tabular}{|p{2cm}|p{13cm}|}
			\hline
			Parameters & Biological interpretations \\
			\hline
			$\Lambda$ & Recruitment rate into the susceptible population \\
			$\beta$ & Disease transmission coefficient \\
			$1/\mu$ & Mean lifetime in the absence of disease-induced mortality \\
			$1/a$ & Mean latent period in the absence of natural mortality \\ 
			$\delta_1$ ($\delta_2$) & Disease-induced death rate of infectious (hospitalized) individuals \\
			$\gamma_1$ ($\gamma_2$) & Recovery rate of infectious (hospitalized) individuals \\
			$1/\xi$ & Mean duration of immunity in the absence of natural mortality \\ 
			$r$ & Intrinsic adjustment rate of effective healthcare capacity \\
			$\lambda$ & Maximal per capita infection-responsive expansion rate
			of effective healthcare capacity \\ 
			$K$ & Baseline effective healthcare capacity in the absence of
			infection-driven expansion \\ 
			$\omega$ & Half-saturation constant for infection-responsive capacity expansion \\ 
			$\mu_1$ ($\mu_0$) & Maximum (minimum) per capita hospitalization rate \\ 
			$I_0$ & Hospitalization threshold \\
			$\tau$ & Healthcare-capacity response delay to reported infection information \\  
			\hline
		\end{tabular}
	\end{table}
	
	We assume that
	$\Lambda,\beta,\mu,a,\gamma_1,\gamma_2,\xi,r,K,\lambda,\omega,I_0>0$,
	$\delta_1,\delta_2\geq0$, $0<\mu_0<\mu_1$, and $\tau\geq0$,
	unless otherwise stated. Here, $\mathbb{R}_+=[0,\infty)$.
	
	System \eqref{ddemodel} with \eqref{treatrate} defines a piecewise-smooth
	delay system \cite{BERN2008}, interpreted in the Filippov sense
	\cite{FILI2013,JIAO2022}. Switching depends on the current state $I(t)$,
	while delay enters only through the healthcare-capacity equation.
	
	We let $X=(S,E,I,H,R,B)$ and define
	\begin{align}
		F_1(X(t),X(t-\tau)) &=
		\begin{pmatrix}
			\Lambda-\beta SI-\mu S+\xi R\\
			\beta SI-(\mu+a)E\\
			aE-(\mu+\delta_1+\gamma_1)I-\mu_1 I\\
			\mu_1 I-(\mu+\delta_2+\gamma_2)H\\
			\gamma_1 I+\gamma_2 H-(\mu+\xi)R\\
			rB\left(1-\frac{B}{K}\right)+\frac{\lambda I(t-\tau)B}{\omega+I(t-\tau)}
		\end{pmatrix},
		\label{sub1}
		\\
		F_2(X(t),X(t-\tau)) &=
		\begin{pmatrix}
			\Lambda-\beta SI-\mu S+\xi R\\
			\beta SI-(\mu+a)E\\
			aE-(\mu+\delta_1+\gamma_1)I-
			\left(\mu_0+(\mu_1-\mu_0)\frac{B}{B+I}\right)I\\
			\left(\mu_0+(\mu_1-\mu_0)\frac{B}{B+I}\right)I-(\mu+\delta_2+\gamma_2)H\\
			\gamma_1 I+\gamma_2 H-(\mu+\xi)R\\
			rB\left(1-\frac{B}{K}\right)+\frac{\lambda I(t-\tau)B}{\omega+I(t-\tau)}
		\end{pmatrix}.
		\label{sub2}
	\end{align}
	
	Then system \eqref{ddemodel} with \eqref{treatrate} is interpreted as the Filippov differential inclusion 
	\begin{equation}\label{filisystem}
		\frac{dX}{dt}\in
		\mathcal{F}(X(t),X(t-\tau))
		=
		\begin{cases}
			F_1(X(t),X(t-\tau)), & I<I_0,\\[2mm]
			\operatorname{co}\left\{
			F_1(X(t),X(t-\tau)),
			F_2(X(t),X(t-\tau))
			\right\}, & I=I_0,\\[2mm]
			F_2(X(t),X(t-\tau)), & I>I_0,
		\end{cases}
	\end{equation}
	where $\operatorname{co}\{\cdot\}$ denotes the closed convex hull of the two one-sided vector fields. Thus, a Filippov solution is a continuous function $X(t)$ agreeing with the prescribed history on $[-\tau,0]$, locally absolutely continuous for $t\geq0$ on its interval of existence, and satisfying \eqref{filisystem} for almost every $t>0$ in that interval.
	
	For system \eqref{filisystem}, the hypersurface
	$
	I=I_0
	$ 
	defines the switching manifold
	\begin{equation*}
		\Sigma
		=
		\left\{
		X\in\mathbb{R}_+^6
		\mid
		I=I_0
		\right\},
	\end{equation*}
	which separates the state space into the two smooth regions
	\begin{align*}
		\Sigma_1
		&=
		\left\{
		X\in\mathbb{R}_+^6
		\mid
		I<I_0
		\right\},\\
		\Sigma_2
		&=
		\left\{
		X\in\mathbb{R}_+^6
		\mid
		I>I_0
		\right\}.
	\end{align*}
	For $\tau>0$, switching in history space is determined by $I(t)$,
	the infectious component at $\theta=0$.
	The dynamics in $\Sigma_1$ and $\Sigma_2$ are governed by the
	corresponding one-sided smooth vector fields $F_{\Sigma_1}=F_1$ and
	$F_{\Sigma_2}=F_2$, respectively. For convenience, we refer to these
	vector fields as subsystems $F_{\Sigma_i}$, $i=1,2$, and use their 
	continuous one-sided extensions when evaluating them on the switching 
	manifold $\Sigma$. 
	
	The one-sided vector fields are continuous up to $\Sigma$ and locally
	bounded on their respective regions, with local extensions to open
	neighborhoods of the non-negative state space. Their Filippov
	regularization is upper semicontinuous with nonempty, compact, convex
	values. For $\tau>0$, the method of steps reduces local existence to
	that for a nonautonomous Filippov inclusion on each delay interval;
	for $\tau=0$, standard Filippov existence theory applies directly
	\cite{FILI2013}. Forward uniqueness is not assumed on the sliding
	region. The following theorem establishes non-negativity, boundedness,
	and global continuation for every maximal Filippov solution.
	
	\begin{theorem} \label{thm_well_posed}
		Let $(S(t),E(t),I(t),H(t),R(t),B(t))$ be any maximal Filippov solution of system \eqref{ddemodel} corresponding to a non-negative continuous initial history. Then the solution is defined for all $t\geq0$ and remains non-negative. Moreover, 
		\[
		N(t)\le
		\max\left\{N(0),\frac{\Lambda}{\mu}\right\},
		\]
		and
		\[
		B(t)\le
		\max\left\{
		B(0),
		K\left(1+\frac{\lambda}{r}\right)
		\right\},
		\quad t\ge0.
		\]
		If $B(0)>0$, then $B(t)>0$ for all $t\ge0$. In particular, the set
		\begin{equation}
			\Omega=\left\{(S,E,I,H,R,B)\in\mathbb{R}_+^6 \,\Big|\,
			0\le S+E+I+H+R\le\frac{\Lambda}{\mu},
			\; 0\le B\le K\left(1+\frac{\lambda}{r}\right)\right\}
		\end{equation}
		is positively invariant in the current-state space for every Filippov
		solution corresponding to a non-negative continuous initial history
		whose initial state belongs to $\Omega$.
	\end{theorem}
	
	The proof of Theorem~\ref{thm_well_posed} is provided in Supplementary Section~S2.
	
	\section{Rigorous analysis of smooth subsystems} \label{SUBDYN}
	We analyze the equilibria, stability and bifurcations of the smooth subsystems $F_{\Sigma_i}$, $i=1,2$, obtained by freezing the switching rule in the respective hospitalization regimes. 
	
	\subsection{Dynamics of subsystem  $F_{\Sigma_1}$}\label{subdyn1} 
	\subsubsection{Disease-free equilibria and basic reproduction number}  
	Subsystem $F_{\Sigma_1}$ always admits two disease-free equilibria denoted by $D_{01}$ and $D_{02},$ where 
	\begin{align*}
		D_{01} = \Big(\frac{\Lambda}{\mu}, 0, 0, 0, 0, 0\Big), \quad D_{02} = \Big(\frac{\Lambda}{\mu}, 0, 0, 0, 0, K\Big),
	\end{align*}
	respectively. 
	
	We compute the basic reproduction number, denoted by $\mathcal R_0$,
	using the next-generation matrix method \cite{PVAN2002}. Although the
	full model contains a discrete delay, the delay enters only through the
	healthcare-capacity equation and does not enter the linearized infected
	subsystem at $D_{02}$. Hence, $\mathcal R_0$ can be obtained from the
	standard next-generation matrix of the infected subsystem. 
	
	The infected compartments are $(E,I)$. The next-generation matrix calculation is provided in Supplementary Section~S5.1. A direct calculation yields
	\begin{equation}
		\mathcal R_0
		= \frac{\beta S^0  a}{(\mu+a)(\mu+\delta_1+\gamma_1+\mu_1)}
		= \frac{\beta \Lambda a}{\mu(\mu+a)(\mu+\delta_1+\gamma_1+\mu_1)}.
	\end{equation}
	
	The basic reproduction number $\mathcal R_0$ represents the average number of secondary infections generated by a typical newly infected individual introduced into a wholly susceptible population at the disease-free equilibrium. It can be interpreted as the product of the transmission rate $\beta$, the number of susceptible individuals at the disease-free state $\Lambda/\mu$, the fraction $a/(\mu+a)$ of exposed individuals who progress to the infectious class, and the average infectious period  $1/(\mu+\delta_1+\gamma_1+\mu_1)$. 
	
	\subsubsection{Local and global stability of disease-free equilibria}
	
	The local stability of $D_{01}$ is determined by the characteristic equation of the linearization of subsystem $F_{\Sigma_1}$ at $D_{01}$. The corresponding characteristic equation factorizes as 
	\begin{equation}\label{eq_chara_dfe1}
		\begin{aligned}
			(s + \mu)
			(s - r)
			(s + \mu + \xi)
			(s + \mu  + \delta_2 + \gamma_2)
			& \left[
			s^2
			+ (2\mu + a 
			+ \delta_1 + \gamma_1 + \mu_1)s \right. \\ 
			& \left. + (\mu + a)(\mu + \delta_1 + \gamma_1 + \mu_1)(1-\mathcal R_0)
			\right]
			= 0. 
		\end{aligned}
	\end{equation}
	Since $r>0$, the characteristic equation admits a strictly positive root $s=r$. Therefore, the disease-free equilibrium $D_{01}$ is always unstable. 
	
	Similarly, the local stability of the disease-free equilibrium $D_{02}$ is determined by the characteristic equation of the linearization of subsystem $F_{\Sigma_1}$ at $D_{02}$. The corresponding characteristic equation factorizes as
	\begin{equation}\label{eq_chara_dfe2}
		\begin{aligned}
			(s + \mu)
			(s + r)
			(s + \mu + \xi)
			(s + \mu  + \delta_2 + \gamma_2)
			& \left[
			s^2
			+ (2\mu + a 
			+ \delta_1 + \gamma_1 + \mu_1)s \right. \\ 
			& \left. + (\mu + a)(\mu + \delta_1 + \gamma_1 + \mu_1)(1-\mathcal R_0)
			\right]
			= 0. 
		\end{aligned}
	\end{equation}
	
	The first four factors give strictly negative real eigenvalues. The remaining two eigenvalues are determined by the quadratic factor. The trace associated with this quadratic is 
	\[
	-\bigl(2\mu+a+\delta_1+\gamma_1+\mu_1\bigr)<0,
	\]
	and its determinant is
	\[
	(\mu+a)(\mu+\delta_1+\gamma_1+\mu_1)(1-\mathcal R_0).
	\] 
	
	Therefore, when $\mathcal R_0 < 1$, the determinant is positive and all eigenvalues
	have negative real parts, implying that $D_{02}$ is locally asymptotically stable.
	When $\mathcal R_0 > 1$, the determinant becomes negative, and $D_{02}$ loses 
	stability. 
	
	The preceding analysis yields the following global stability result. 
	\begin{theorem}\label{thm:dfe_global}
		If $\mathcal R_0<1$, then the disease-free equilibrium $D_{02}$ is
		globally asymptotically stable for the smooth extension of subsystem $F_{\Sigma_1}$ with
		initial conditions in $\Omega$ satisfying $B(0)>0$.
	\end{theorem} 
	
	The proof is provided in Supplementary Section~S3.

	The instability of $D_{01}$ associated with the eigenvalue $r$ does not
	indicate instability with respect to infection. Rather, the positive eigenvalue $r$ is associated with
	the healthcare-capacity variable, implying that a small positive
	capacity perturbation grows away from $B=0$ under the logistic growth
	mechanism. Thus, the epidemiologically relevant disease-free equilibrium
	is $D_{02}$, which represents disease elimination with positive
	long-term healthcare capacity. When $\mathcal R_0<1$, small infectious
	perturbations cannot invade $D_{02}$, whereas for $\mathcal R_0>1$,
	$D_{02}$ is unstable with respect to infectious perturbations. 
	
	\subsubsection{Transcritical bifurcation at the disease-free equilibrium $D_{02}$} 
	We observe that the characteristic equation of subsystem $F_{\Sigma_1}$ has a simple zero root when $\mathcal R_0 = 1$ (see Eq. \eqref{eq_chara_dfe2}), while all other eigenvalues 
	have strictly negative real parts. This indicates the possibility of a transcritical bifurcation at the disease-free equilibrium $D_{02}$. To rigorously establish the direction and local 
	dynamics of this bifurcation, we apply to the autonomous epidemiological subsystem the framework of Castillo-Chavez and Song \cite{CAS2004}, which is based on center manifold theory \cite{GUCK1983}. The following theorem establishes the occurrence of a forward transcritical bifurcation at $\mathcal R_0=1$. The proof is deferred to Supplementary Section~S4.  
	
	\begin{theorem}\label{thm_transcritical}
		The subsystem $F_{\Sigma_1}$ undergoes a forward transcritical bifurcation at the disease-free equilibrium $D_{02}$ when $\mathcal R_0 = 1$. 
	\end{theorem}
	
	\subsubsection{Existence and stability of endemic equilibria}
	
	Setting the right-hand side of subsystem $F_{\Sigma_1}$ to zero yields a
	unique positive-capacity endemic equilibrium 
	\[
	D^*=(S^*,E^*,I^*,H^*,R^*,B^*),
	\]
	which exists whenever $\mathcal{R}_0>1$ in the smooth extension of $F_{\Sigma_1}$. Its components are given by 
	
	{\small
		\begin{align*}
			S^*
			&=
			\frac{\Lambda}{\mu+\beta I^*}
			+
			\frac{\xi}
			{(\mu+\beta I^*)(\mu+\xi)}
			\left[
			\gamma_1+
			\frac{\gamma_2\mu_1}
			{\mu+\delta_2+\gamma_2}
			\right]I^*,
			\quad
			E^*
			=
			\frac{(\mu+\delta_1+\gamma_1+\mu_1)I^*}{a},
			\\
			I^*
			&=
			\frac{
				\mu(\mu+a)(\mu+\delta_1+\gamma_1+\mu_1)
				(\mu+\xi)(\mu+\delta_2+\gamma_2)
				(\mathcal R_0-1)
			}{
				\beta\Big(
				\mu(a+\mu+\xi)
				\big(
				\gamma_1(\mu+\gamma_2+\delta_2)+\gamma_2\mu_1
				\big)
				+
				(a+\mu)(\mu+\xi)
				\big(
				\gamma_2(\mu+\delta_1)
				+
				(\mu+\delta_2)(\mu+\delta_1+\mu_1)
				\big)
				\Big)
			},
			\\
			H^*
			&=
			\frac{\mu_1}{\mu+\delta_2+\gamma_2}I^*,
			\quad
			R^*
			=
			\frac{1}{\mu+\xi}
			\left[
			\gamma_1+
			\frac{\gamma_2\mu_1}
			{\mu+\delta_2+\gamma_2}
			\right]I^*,
			\quad
			B^*
			=
			\frac{K}{r}
			\left(
			r+\frac{\lambda I^*}{\omega+I^*}
			\right).
		\end{align*}
	}
	
	For $D^*$ to be an admissible equilibrium of $F_{\Sigma_1}$, it must
	also satisfy $I^*<I_0$. To analyze the local stability of the endemic equilibrium $D^*$, we 
	derive the characteristic equation associated with the linearization of 
	subsystem $F_{\Sigma_1}$ at $D^*$. The corresponding characteristic 
	equation is
	\begin{equation}\label{chddee}
		(s-k_{66})
		\left(
		s^5+l_{11}s^4+l_{12}s^3+l_{13}s^2+l_{14}s+l_{15}
		\right)
		=0,
	\end{equation}
	where the coefficients
	$l_{11},l_{12},l_{13},l_{14},l_{15}$ and $k_{66}$ are explicitly
	given in Supplementary Section~S5.
	
	The root $s=k_{66}$ corresponds to the healthcare-capacity dynamics.
	Using the expression for $B^*$, we obtain 
	\[
	k_{66}
	=
	-\left(
	r+\frac{\lambda I^*}{\omega+I^*}
	\right)
	<0,
	\]
	and hence this root is always negative. The remaining five roots are
	determined by the fifth-degree polynomial. To characterize their
	stability, we introduce the Routh-Hurwitz determinants $H_1,\ldots,H_5$, given in Supplementary Section~S5.2.
	
	By the Routh-Hurwitz criterion, all characteristic roots have strictly
	negative real parts if and only if
	\[
	H_1>0,\quad
	H_2>0,\quad
	H_3>0,\quad
	H_4>0,\quad
	H_5>0.
	\]
	Under these conditions, the endemic equilibrium $D^*$ is locally asymptotically stable.
	
	The following result provides a sufficient Lyapunov matrix criterion
	for global convergence to the endemic equilibrium within the smooth
	regime $I<I_0$.
	
	\begin{theorem}\label{thm_global_stab}
		Suppose that $\mathcal R_0>1$ and $I^*<I_0$. Assume that there exist a symmetric positive-definite matrix $P=P^\top\succ0$ and a constant $\eta>0$ such that 
		\[
		PQ(0)+Q(0)^\top P\preceq0,
		\quad
		PQ(I_0)+Q(I_0)^\top P\preceq-\eta I_5,
		\]
		where $Q(I)$ is given in Supplementary Section~S7. Then $D^*$ is globally asymptotically stable relative to nonnegative histories with $E(0)+I(0)>0$ and $B(0)>0$ whose trajectories remain in $I<I_0$.
	\end{theorem}
	
	The proof, including an analytical demonstration that the matrix criterion is feasible for a sufficiently small positive value of $\beta I_0$, is provided in Supplementary Section~S7. The proof uses the divergence of the infection integral for infection-positive trajectories. The criterion is conditional; feasibility is not asserted for every parameter set used below.
	
	\begin{remark}
		The healthcare-capacity response delay $\tau$ does not enter the local
		stability conditions above or the global convergence criterion in
		Theorem~\ref{thm_global_stab}. In the regime $I<I_0$, the
		hospitalization rate is constant, $T(I,B)=\mu_1$, and the
		healthcare-capacity variable $B$ does not feed back into the
		epidemiological subsystem $(S,E,I,H,R)$. Consequently, $\tau$ may
		affect the transient evolution of $B(t)$, but it does not alter the
		disease-invasion threshold or the stability conditions of the
		epidemiological equilibria while the trajectory remains in
		$I<I_0$. 
	\end{remark} 
	
	\subsection{Dynamics of subsystem $F_{\Sigma_2}$} \label{subdyn2} 
	\subsubsection{Existence and stability of equilibria} 
	
	We first consider the equilibria of the smooth extension of
	$F_{\Sigma_2}$. The disease-free state $D_{02}$ formally satisfies the corresponding
	equilibrium equations but is not admissible in $F_{\Sigma_2}$, whose
	domain is $I>I_0$. The state $D_{01}$ is a singular point of the
	extended hospitalization function. In the full Filippov system, the epidemiologically
	relevant disease-free equilibrium $D_{02}$ lies in the below-threshold
	regime $F_{\Sigma_1}$, where its stability is determined by
	$\mathcal R_0$. Thus, $\mathcal R_0$ is the disease-invasion threshold
	associated with the below-threshold disease-free equilibrium, rather than
	a reproduction number obtained from an admissible disease-free equilibrium
	of $F_{\Sigma_2}$. The positive-capacity endemic equilibria are denoted by
	\[
	D_1^*=(S_1^*,E_1^*,I_1^*,H_1^*,R_1^*,B_1^*),
	\]
	where 
	
	{\footnotesize   
		\begin{align*}
			S_1^* &= \frac{\Lambda}{\mu + \beta I_1^*} + \frac{\xi}{(\mu + \beta I_1^*)(\mu + \xi)}\Big[\gamma_1  + \frac{\gamma_2}{\mu + \delta_2 + \gamma_2}\Big(\mu_0 + (\mu_1-\mu_0)\frac{B_1^*}{B_1^* + I_1^*} \Big) \Big]I_1^*,\\
			E_1^* &= \frac{\mu + \delta_1 + \gamma_1}{a}I_1^* + \frac{1}{a}\Big(\mu_0 + (\mu_1-\mu_0)\frac{B_1^*}{B_1^* + I_1^*} \Big)I_1^*, \quad H_1^* = \frac{1}{\mu + \delta_2 + \gamma_2}\Big(\mu_0 + (\mu_1-\mu_0)\frac{B_1^*}{B_1^* + I_1^*} \Big)I_1^*,\\
			R_1^* &= \frac{1}{\mu + \xi}\Big[\gamma_1 + \frac{\gamma_2}{\mu + \delta_2 + \gamma_2}\Big(\mu_0 + (\mu_1-\mu_0)\frac{B_1^*}{B_1^* + I_1^*} \Big) \Big]I_1^*, \quad
			B_1^* = \frac{K}{r}\Big(r+\frac{\lambda I_1^*}{\omega + I_1^*}\Big).
	\end{align*}} 
	
	Here, $I_1^*$ is a positive root of the following equation: 
	\begin{equation}\label{eq12}
		G_1 I^3 + G_2 I^2 + G_3 I + G_4 = 0,
	\end{equation}
	where the coefficients $G_1, G_2, G_3,$ and $G_4$ are provided in Supplementary Section~S5. 
	
	Using Descartes' rule of signs, Eq.~\eqref{eq12} can have up to three
	positive roots, depending on the signs of the coefficients $G_i$,
	$i=1,2,3,4$. More precisely, the number of positive roots is equal to
	the number of sign changes or is less than it by an even integer. The
	possible numbers of positive roots, and hence candidate endemic
	equilibria of the smooth extension of $F_{\Sigma_2}$, are summarized
	in Table~\ref{end_subsystem2}. The $\mathcal R_0$ column indicates the
	corresponding side of the disease-invasion threshold
	$\mathcal R_0=1$. 
	
	\begin{table}[H] 
		\centering
		\caption{Possible numbers of positive roots of the endemic-equilibrium equation for the smooth extension of  $F_{\Sigma_2}$ based on Descartes' rule of signs.}  
		\label{end_subsystem2}
		\begin{tabular}{cccccc} 
			\hline
			$G_1$ & $G_2$ & $G_3$ & $G_4$ & $\mathcal R_0$ & Possible number of positive roots \\
			\hline
			$+$ &$+$ & $+$ & $+$ & $<1$ & 0 \\
			$+$&$+$ & $-$ & $+$ & $<1$ & 0 or 2\\
			$+$&$-$ & $-$ & $+$ & $<1$ & 0 or 2 \\
			$+$& $-$ & $+$ & $+$ & $<1$ & 0 or 2 \\
			$+$ &$+$ & $+$ & $-$ & $>1$ & 1 \\
			$+$&$+$ & $-$ & $-$ & $>1$ & 1 \\
			$+$&$-$ & $-$ & $-$ & $>1$ & 1 \\
			$+$& $-$ & $+$ & $-$ & $>1$ & 1 or 3 \\
			\hline
		\end{tabular}
	\end{table}
	
	The positive roots identified in Table~\ref{end_subsystem2} correspond
	to biologically feasible endemic equilibria only when the associated
	equilibrium components are positive. Furthermore, an endemic equilibrium
	is admissible in subsystem $F_{\Sigma_2}$ only when
	$
	I_1^*>I_0.
	$

	\subsubsection{Nondegenerate fold of endemic equilibria in the smooth extension of $F_{\Sigma_2}$} 
	
	\begin{theorem}\label{thm:fold} 
		Assume that endemic equilibria of the smooth extension of  $F_{\Sigma_2}$ correspond to the 
		positive roots of 
		\begin{equation}\label{eq:G}
			G(I,\beta)=G_1(\beta) I^3 + G_2(\beta) I^2 + G_3(\beta) I + G_4(\beta)=0, 
			\quad I>0,
		\end{equation}
		where $G_1(\beta)>0$ and $G \in C^2(I,\beta)$, i.e., $G$ is twice continuously differentiable in $(I,\beta)$. 
		Suppose there exists a point $(I_{\rm SN},\beta_{\rm SN})$, with $I_{\rm SN}>0$, such that 
		\begin{equation}\label{eq:foldconds}
			G(I_{\rm SN},\beta_{\rm SN})=0, 
			\quad 
			\partial_I G(I_{\rm SN},\beta_{\rm SN})=0,
		\end{equation}
		and the nondegeneracy conditions
		\begin{equation}\label{eq:foldnondeg}
			\partial_{II}G(I_{\rm SN},\beta_{\rm SN})\neq 0, 
			\quad 
			\partial_\beta G(I_{\rm SN},\beta_{\rm SN})\neq 0
		\end{equation}
		hold. Then the endemic-equilibrium equation of the smooth extension of $F_{\Sigma_2}$ undergoes a nondegenerate fold at $(I_{\rm SN},\beta_{\rm SN})$. In particular, there exists a neighborhood of $(I_{\rm SN},\beta_{\rm SN})$ in which equation \eqref{eq:G} has exactly two positive solutions on one side of $\beta_{\rm SN}$ and no positive solutions on the other side.
	\end{theorem}
	
	\begin{proof}
		From condition \eqref{eq:foldconds}, $I_{\rm SN}$ is a double root of
		$G(I,\beta_{\rm SN})$. The conditions
		$\partial_{II}G(I_{\rm SN},\beta_{\rm SN})\neq0$ and
		$\partial_\beta G(I_{\rm SN},\beta_{\rm SN})\neq0$ ensure,
		respectively, the nondegeneracy of the double root and its unfolding
		with respect to $\beta$. Hence, by the standard fold theorem, the
		local solution set of $G(I,\beta)=0$ near
		$(I_{\rm SN},\beta_{\rm SN})$ is smoothly equivalent to
		\[
		\eta^2=c(\beta-\beta_{\rm SN}),
		\quad c\neq0.
		\]
		Thus, two positive equilibrium roots exist on one side of
		$\beta_{\rm SN}$ and none on the other.
	\end{proof}
	
	Suppose further that
	\[
	\mathcal R_0(\beta_{\rm SN})<1.
	\]
	By continuity of $\mathcal R_0(\beta)$, the two positive roots created
	on the appropriate side of the fold occur locally at parameter values
	for which $\mathcal R_0<1$. Therefore, the bifurcation is more precisely
	interpreted as a subthreshold fold of the endemic-equilibrium branches
	of the smooth extension of $F_{\Sigma_2}$. It should not, by itself, be
	interpreted as a classical backward bifurcation, since the disease-free
	equilibrium at which $\mathcal R_0$ is defined belongs to the
	below-threshold regime $F_{\Sigma_1}$ and is not an admissible
	disease-free equilibrium of $F_{\Sigma_2}$. The positive roots are
	admissible in the full Filippov system only when they are biologically
	feasible and satisfy $I_1^*>I_0$.
	
	The fold conditions alone do not determine the stability of the endemic
	branches. If an admissible positive branch is locally
	asymptotically stable, then for
	$\mathcal R_0<1$ the stable endemic equilibrium of the above-threshold
	regime can coexist with the locally stable disease-free equilibrium
	$D_{02}$ of the below-threshold regime. Such coexistence constitutes
	discontinuity-induced subthreshold bistability in the full Filippov
	system. Consequently, $\mathcal R_0<1$ guarantees local stability of the
	disease-free equilibrium against small infectious perturbations but does
	not necessarily exclude an admissible stable endemic state in the
	above-threshold regime. We will further examine the fold structure, branch admissibility, and
	associated bistability in the numerical simulation section. 
	
	For the stability analysis of an admissible endemic equilibrium of subsystem $F_{\Sigma_2}$, we first consider the case when the time delay is zero, i.e., $\tau=0$, in the following section. 
	
	\subsubsection{Local stability analysis and bifurcations for $\tau=0$}\label{hopf1} 
	For the delay-free subsystem $F_{\Sigma_2}$, we denote the endemic equilibrium $D_1^*=(S_1^*,E_1^*,I_1^*,H_1^*,R_1^*,B_1^*)$ by $D^*=(S^*,E^*,I^*,H^*,R^*,B^*)$ for convenience. The characteristic equation of its Jacobian matrix is
	\begin{equation} \label{chara_end} 
		C(s) = s^6 + f_{11} s^5 + f_{12} s^4 + f_{13} s^3 + f_{14} s^2 + f_{15} s +  f_{16}  = 0,
	\end{equation}  
	where  $f_{11}$,  $f_{12}$,  $f_{13}$,  $f_{14}$,  $f_{15}$, and $f_{16}$ are provided in Supplementary Section~S5. 
	
	We further define the Routh-Hurwitz determinants $H_1,\ldots,H_6$ in Supplementary Section~S5.2.
	By the Routh-Hurwitz criterion, all characteristic roots have strictly negative real parts if and only if $H_1 > 0, H_2 > 0, H_3 > 0, H_4 > 0, H_5 > 0,$ and $H_6 > 0$. Under these conditions, $D^*$ is locally asymptotically stable; if a characteristic root has positive real part, $D^*$ is unstable.
	
	Next, we analyze the conditions for a Hopf bifurcation in subsystem
	$F_{\Sigma_2}$, choosing $\lambda$, the per capita growth rate of
	infection-driven healthcare capacity, as the bifurcation parameter. 
	We derive the critical value $\lambda=\lambda^{*}$ at which the endemic
	equilibrium loses stability through a Hopf bifurcation. We assume that there exists a smooth branch of endemic equilibria $D(\lambda)$ of subsystem $F_{\Sigma_2}$ in a neighborhood of $\lambda=\lambda^{*}$, with $D(\lambda^{*})=D^{*}$. By Liu criterion \cite{Liu1994}, if  
	\[ 
	C(s, \lambda^*) = s^6 + f_{11}(\lambda^*) s^5 + f_{12}(\lambda^*) s^4 + f_{13}(\lambda^*) s^3 + f_{14}(\lambda^*) s^2 + f_{15} (\lambda^*) s + f_{16}(\lambda^*)
	\] 
	is the characteristic equation at $(D^*, \lambda^*)$, then a simple imaginary pair crosses transversally under the following conditions:
	\begin{itemize}
		\item[(i)] $f_{16} (\lambda^*)>0, H_1 (\lambda^*) > 0, H_2 (\lambda^*) > 0, H_3 (\lambda^*) > 0, H_4 (\lambda^*) > 0,$ and $H_5 (\lambda^*) = 0.$
		\item[(ii)] $\frac{d}{d \lambda}(H_5 (\lambda^*)) \neq 0,$
	\end{itemize}
	where $H_1 (\lambda^*), H_2 (\lambda^*), H_3 (\lambda^*), H_4 (\lambda^*),$ and $H_5 (\lambda^*)$ are the Hurwitz determinants at the bifurcation parameter $\lambda^*$. 
	
	The detailed transversality calculation is provided in Supplementary Section~S6.
	
	\begin{theorem}
		Let $D^*$ be the positive equilibrium of subsystem $F_{\Sigma_2}$. If the conditions derived above are satisfied at $\lambda=\lambda^*$, then a simple pair of imaginary eigenvalues crosses the imaginary axis transversally. If, additionally, the first Lyapunov coefficient is nonzero, then $D^*$ undergoes a nondegenerate Hopf bifurcation as $\lambda$ passes through $\lambda^*$. 
	\end{theorem} 
	
	\subsubsection{Local stability and Hopf bifurcation when $\tau > 0$}\label{Sec_Hopf_tau} 
	
	In this section, we discuss the conditions of local stability and Hopf bifurcation around the equilibrium $D^*$ when $\tau >0$. The characteristic equation of the subsystem $F_{\Sigma_2}$ at any equilibrium $D^*$ is given by 
	\begin{equation}
		\begin{aligned} \label{chdde2} 
			s^6 + g_{11}s^5 + g_{12}s^4 + g_{13}s^3 + g_{14}s^2 & + g_{15}s + g_{16} \\
			& + (g_{17}s^4 + g_{18}s^3 + g_{19}s^2 + g_{20}s + g_{21})e^{-s\tau} = 0, 
		\end{aligned} 
	\end{equation} 
	where $g_{11}, g_{12}, g_{13}, g_{14}, g_{15}, g_{16}, g_{17}, g_{18}, g_{19}, g_{20}$ and $g_{21}$ are provided in Supplementary Section~S5.  
	
	Let $s=i\nu$ $(\nu > 0)$ be a root of equation \eqref{chdde2}. Substituting $s=i\nu$ in equation \eqref{chdde2} and separating the imaginary and real parts, we obtain
	the two equations given in Supplementary Section~S5.3.
	
	Squaring and adding the above equations, we further obtain 
	\begin{equation}\label{eq15}
		\nu^{12} + h_{11}\nu^{10} + h_{12}\nu^{8} + h_{13}\nu^6 + h_{14}\nu^4 + h_{15}\nu^2 + h_{16} = 0, 
	\end{equation}
	where $h_{11}, h_{12}, h_{13}, h_{14}, h_{15}$ and $h_{16}$ are provided in Supplementary Section~S5. 
	
	Let $\nu^2 = \zeta$. Then equation \eqref{eq15} becomes
	\begin{equation}\label{eq16}
		\zeta^{6} + h_{11}\zeta^{5} + h_{12}\zeta^4 + h_{13}\zeta^3 + h_{14}\zeta^2 + h_{15}\zeta + h_{16} = 0.
	\end{equation}
	
	Suppose that equation \eqref{eq16} has $m\geq1$ positive roots
	$
	\zeta_j>0, j=1,2,\ldots,m,
	$ 
	and define
	$
	\nu_j=\sqrt{\zeta_j}.
	$ 
	For each $\nu_j$, let
	\begin{align*}
		F_{1j} &= g_{20}\nu_j-g_{18}\nu_j^3,\\
		F_{2j} &= g_{21}+g_{17}\nu_j^4-g_{19}\nu_j^2,\\
		F_{3j} &= \nu_j^6-g_{12}\nu_j^4+g_{14}\nu_j^2-g_{16},\\
		F_{4j} &= g_{13}\nu_j^3-g_{11}\nu_j^5-g_{15}\nu_j.
	\end{align*}
	Assuming $F_{1j}^2+F_{2j}^2>0$, the corresponding sine and cosine expressions follow from Eq.~(S5.7) (see Supplementary Section~S5.3).
	
	Let
	\begin{equation}\label{eq:theta_j}
		\theta_j
		=
		\operatorname{atan2}
		\left(
		\frac{F_{1j}F_{3j}-F_{2j}F_{4j}}
		{F_{1j}^2+F_{2j}^2},
		\frac{F_{2j}F_{3j}+F_{1j}F_{4j}}
		{F_{1j}^2+F_{2j}^2}
		\right),
	\end{equation}
	where $\operatorname{atan2}$ is taken modulo $2\pi$, so that $\theta_j\in[0,2\pi)$. The corresponding family of critical delays is therefore given by
	\begin{equation}\label{eq:critical_delays}
		\tau_j^{(k)}
		=
		\frac{\theta_j+2k\pi}{\nu_j},
		\quad
		j=1,2,\ldots,m,\quad
		k=0,1,2,\ldots,
	\end{equation}
	where only positive values of $\tau_j^{(k)}$ are considered.
	
	Further, differentiating equation \eqref{chdde2} with respect to $\tau$ and substituting
	$s=i\nu_j$ and $\tau=\tau_j^{(k)}$, we obtain
	\begin{equation}\label{eq:transversality}
		\operatorname{Re}
		\left(
		\frac{ds}{d\tau}
		\right)^{-1}_{s=i\nu_j,\,
			\tau=\tau_j^{(k)}}
		=
		\frac{
			A_{11}A_{13}+A_{12}A_{14}
		}{
			A_{13}^2+A_{14}^2
		},
	\end{equation}
	where $A_{11}, A_{12}, A_{13}$ and $A_{14}$ are provided in Supplementary Section~S5 and are evaluated at
	$s=i\nu_j$ and $\tau=\tau_j^{(k)}$. 
	
	\begin{lemma}
		If
		$
		A_{11}A_{13}+A_{12}A_{14}\neq0,
		$ 
		at a critical delay $\tau=\tau_j^{(k)}$, then
		\[
		\operatorname{Re}
		\left(
		\frac{ds}{d\tau}
		\right)_{s=i\nu_j,\,
			\tau=\tau_j^{(k)}}
		\neq0,
		\]
		and hence the pair of characteristic roots crosses the imaginary axis transversally. In particular, the sign of
		$
		A_{11}A_{13}+A_{12}A_{14}
		$ 
		determines the direction of crossing.
	\end{lemma}
	
	From the above discussion, we conclude the following result. 
	
	\begin{theorem}
		Suppose that equation \eqref{eq16} has positive roots
		$\zeta_j$, $j=1,2,\ldots,m$, and let the corresponding critical delays
		$\tau_j^{(k)}$ be defined by \eqref{eq:critical_delays}. Assume that the endemic equilibrium $D^*$ has all characteristic roots with strictly negative real parts at $\tau=0$ and that each relevant pair of purely imaginary characteristic roots $\pm i\nu_j$ is simple, with no simultaneous crossings of distinct pairs and $F_{1j}^2+F_{2j}^2>0$ for every positive frequency. If
		$ 
		A_{11}A_{13}+A_{12}A_{14}\neq0
		$ 
		at $\tau=\tau_j^{(k)}$, then the characteristic roots cross the imaginary axis transversally at that critical delay. Consequently, the local stability of $D^*$ can change only when $\tau$ passes through one of the critical values $\tau_j^{(k)}$. A positive crossing direction increases the number of characteristic roots with positive real part by two, whereas a negative crossing direction decreases it by two. Therefore, depending on the ordering and crossing directions of the
		critical delays, the equilibrium can undergo multiple stability switches
		as $\tau$ varies. Under the additional standard nonlinear nondegeneracy conditions for Hopf bifurcation, a nondegenerate Hopf bifurcation occurs at a critical delay $\tau=\tau_j^{(k)}$.
	\end{theorem}
	
	The transversality condition determines the direction in which the characteristic roots cross the imaginary axis, but it does not determine whether the corresponding Hopf bifurcation is supercritical or subcritical. Determination of the criticality and orbital stability of the bifurcating periodic solutions requires a nonlinear normal-form calculation, such as the evaluation of the first Lyapunov coefficient, and is not pursued here.
	
	\begin{remark}
		The Hopf bifurcation indicates that, once infections
		exceed the threshold $I_0$ and hospitalization becomes dependent on
		available healthcare capacity, the endemic equilibrium loses stability and
		oscillatory epidemic dynamics emerges. In the absence of delay, these
		oscillations arise from the nonlinear feedback between infection prevalence
		and capacity-constrained hospitalization. In the delayed system, the use of
		past infection information in healthcare capacity adjustment can further change the 
		stability of the endemic equilibrium, inducing or suppressing oscillations as
		the delay varies. The existence of several critical-delay families can therefore produce repeated destabilization and restabilization of the endemic equilibrium, consistent with the stability-switching behavior observed numerically.
	\end{remark} 
	
	\section{Non-smooth dynamics on the switching manifold $\Sigma = \left\{I = I_0 \right\}$} \label{SLIDDYN} 
	
	In this section, we explore the non-smooth dynamics on the switching manifold $\Sigma = \left\{I = I_0 \right\}$ of the Filippov system \eqref{ddemodel}, and define the local trajectories through the boundary $\Sigma$. These local trajectories initiate from the subregions $\Sigma_1$ and $\Sigma_2$ and are defined by the subsystems $F_{\Sigma_1}$ and $F_{\Sigma_2}$. To extend the definition of a trajectory to the discontinuous boundary $\Sigma$, it is necessary to distinguish the different local configurations according to the directions of the two one-sided vector fields \cite{FILI2013,GUAR2010}, as follows: 
	\begin{itemize}
		\item Attracting sliding region $\Sigma_s$, where the vector fields on both sides point toward the switching manifold.
		\item Repelling (escaping) sliding region $\Sigma_e$, where the vector fields on both sides point away from the switching manifold.
		\item Crossing region $\Sigma_c$, where the two vector fields cross the switching manifold from one side to the other. 
	\end{itemize} 
	
	\subsection{Sliding and crossing regions}
	
	Firstly, denote
	$
	L(X)=I-I_0
	$ 
	as the switching function, with 
	$ 
	\nabla L=(0,0,1,0,0,0).
	$ 
	On either side of $\Sigma$, the normal components of the corresponding vector fields are
	\begin{align*}
		\mathcal{L}_1(X)
		&=
		\left\langle \nabla L,F_{\Sigma_1}\right\rangle
		=
		aE-(\mu+\delta_1+\gamma_1+\mu_1)I,\\
		\mathcal{L}_2(X)
		&=
		\left\langle \nabla L,F_{\Sigma_2}\right\rangle 
		=
		aE-
		\Big(
		\mu+\delta_1+\gamma_1+\mu_0
		+(\mu_1-\mu_0)\frac{B}{B+I}
		\Big)I.
	\end{align*}
	
	Following Filippov's convention, sliding or escaping motion can occur at points on $\Sigma$ where the two normal components have opposite signs, i.e.,
	$ 
	\mathcal{L}_1(X)\mathcal{L}_2(X)<0.
	$ 
	On the manifold $\Sigma$, we set $I=I_0$. Solving
	$ 
	\mathcal{L}_1=0
	$ and
	$ 
	\mathcal{L}_2=0
	$ 
	for $E$ gives the two critical values
	\begin{align*}
		E_{c1}
		&=
		\frac{(\mu+\delta_1+\gamma_1+\mu_1)I_0}{a},\\
		E_{c2}(B)
		&=
		\frac{
			\Big(
			\mu+\delta_1+\gamma_1+\mu_0
			+(\mu_1-\mu_0)\frac{B}{B+I_0}
			\Big)I_0
		}{a}.
	\end{align*}
	
	Since $0<\mu_0<\mu_1$ and $I_0>0$, for every finite $B\geq0$ we have
	\[
	\mu_0+(\mu_1-\mu_0)\frac{B}{B+I_0}
	<
	\mu_1.
	\]
	Thus, the effective per capita hospitalization rate immediately above the
	switching threshold is always smaller than the hospitalization rate
	$\mu_1$ below the threshold. Consequently,
	\begin{equation}\label{Ec_order}
		E_{c2}(B)<E_{c1}.
	\end{equation}
	
	Because $\mathcal{L}_1$ and $\mathcal{L}_2$ are linear and strictly increasing in $E$ for fixed $I_0$ and $B$, their sign patterns on $\Sigma$ are
	\[
	\begin{cases}
		\mathcal{L}_1<0,\quad \mathcal{L}_2<0,
		& E<E_{c2}(B),\\[4pt]
		\mathcal{L}_1<0,\quad \mathcal{L}_2=0,
		& E=E_{c2}(B),\\[4pt]
		\mathcal{L}_1<0,\quad \mathcal{L}_2>0,
		& E_{c2}(B)<E<E_{c1},\\[4pt]
		\mathcal{L}_1=0,\quad \mathcal{L}_2>0,
		& E=E_{c1},\\[4pt]
		\mathcal{L}_1>0,\quad \mathcal{L}_2>0,
		& E>E_{c1}.
	\end{cases}
	\]
	
	The ordering \eqref{Ec_order} has an important structural consequence for
	the nonsmooth dynamics.
	
	\begin{theorem}\label{thm_repelling_sliding}
		Assume that $0<\mu_0<\mu_1$. The Filippov  sliding condition 
		$
		\mathcal{L}_1(X)\mathcal{L}_2(X) < 0
		$
		holds if and only if
		$
		E_{c2}(B)<E<E_{c1}.
		$
		Moreover, throughout this interval,
		$
		\mathcal{L}_1(X)<0<\mathcal{L}_2(X).
		$
		Therefore, whenever sliding occurs, it is necessarily repelling
		(escaping), and the system admits the repelling  sliding region 
		\begin{equation}\label{rep_region}
			\Sigma_e=
			\left\{
			(S,E,I,H,R,B)\in\mathbb{R}_+^6
			\,\middle|\,
			E_{c2}(B)<E<E_{c1},\ I=I_0
			\right\}. 
		\end{equation}
		In particular, the hospitalization function considered here
		structurally excludes an attracting sliding region for every finite 
		$B\geq0$. 
	\end{theorem}
	
	\begin{proof}
		From \eqref{Ec_order}, 
		$
		E_{c2}(B)<E_{c1}.
		$
		For
		$
		E_{c2}(B)<E<E_{c1},
		$
		we have
		$
		\mathcal{L}_1(X)<0<\mathcal{L}_2(X).
		$
		Hence, $F_{\Sigma_1}$ points from the switching manifold toward
		$\Sigma_1$, while $F_{\Sigma_2}$ points toward $\Sigma_2$.
		Thus, both one-sided vector fields point away from $\Sigma$, yielding
		a repelling (escaping) sliding configuration.
		
		Since $E_{c2}(B)<E_{c1}$ for every finite $B\geq0$, the reverse
		configuration
		$
		\mathcal{L}_1(X)>0>\mathcal{L}_2(X),
		$
		required for attracting sliding, cannot occur.
	\end{proof}
	
	Theorem~\ref{thm_repelling_sliding} shows that the downward jump in the
	per capita hospitalization rate at $I=I_0$ gives the switching manifold
	an intrinsically repelling geometry. Consequently, the threshold cannot
	act as an attracting regulation surface, and the pseudo-equilibria,
	tangencies, and boundary-equilibrium transitions considered below should
	be interpreted within this repelling threshold structure.
	
	For $E<E_{c2}(B)$ or $E>E_{c1}$, the normal components have
	the same sign and trajectories locally cross the switching manifold.
	Accordingly,
	\[
	\Sigma_c=\Sigma_c^-\cup\Sigma_c^+,
	\]
	where
	\[
	\Sigma_c^-=
	\left\{
	X\in\Sigma:\ E<E_{c2}(B)
	\right\},
	\quad
	\Sigma_c^+=
	\left\{
	X\in\Sigma:\ E>E_{c1}
	\right\}.
	\]
	For $E<E_{c2}(B)$, both normal components are negative, whereas for
	$E>E_{c1}$, both are positive. The sets
	\[
	E=E_{c2}(B)
	\quad\text{and}\quad
	E=E_{c1}
	\]
	form the tangency boundaries separating the crossing and repelling
	sliding regions. \\ 
	
	\noindent
	\textbf{Interpretation of healthcare capacity on the repelling region:}
	The lower threshold $E_{c2}(B)$ depends explicitly on the healthcare
	capacity $B(t)$, so the lower boundary of the repelling region is
	state-dependent. Since
	\[
	\frac{\partial E_{c2}}{\partial B}
	=
	\frac{(\mu_1-\mu_0)I_0^2}
	{a(B+I_0)^2}
	>0,
	\]
	an increase in healthcare capacity raises $E_{c2}(B)$ and narrows the
	repelling sliding interval
	$ 
	\big(E_{c2}(B),E_{c1}\big),
	$ 
	whereas a decrease in healthcare capacity lowers $E_{c2}(B)$ and widens
	the interval. 
	
	\subsection{Filippov sliding vector field}   
	On the repelling sliding region $\Sigma_e$, the associated Filippov tangential vector field is obtained from the convex combination 
	\[
	F_{\Sigma_e}
	=
	\alpha F_{\Sigma_1}
	+
	(1-\alpha)F_{\Sigma_2},
	\]
	where $\alpha\in[0,1]$ is chosen so that the normal component of the convex combination vanishes 
	\[
	\left\langle
	\nabla L,F_{\Sigma_e}
	\right\rangle
	=
	\alpha\mathcal{L}_1
	+
	(1-\alpha)\mathcal{L}_2
	=
	0.
	\]
	Therefore,
	\begin{equation}\label{alpha_slide}
		\alpha
		=
		\frac{\mathcal{L}_2}
		{\mathcal{L}_2-\mathcal{L}_1}.
	\end{equation}
	Using the expressions for $\mathcal{L}_1$ and $\mathcal{L}_2$, we obtain
	\begin{align*}
		\alpha
		&=
		\frac{
			aE-(\mu+\delta_1+\gamma_1)I
			-\Big(
			\mu_0+(\mu_1-\mu_0)\frac{B}{B+I}
			\Big)I
		}{
			\Big(
			\mu_1-\mu_0
			-(\mu_1-\mu_0)\frac{B}{B+I}
			\Big)I
		},
	\end{align*}
	for
	$ 
	E\in[E_{c2}(B),E_{c1}], 
	I=I_0.
	$ 
	For points in the interior of $\Sigma_e$, we have 
	$ 
	0<\alpha<1.
	$ 
	
	Using the Filippov tangency condition
	$ 
	\left\langle
	\nabla L,F_{\Sigma_e}
	\right\rangle=0,
	$ 
	the effective hospitalization rate along the tangential sliding vector field satisfies 
	\[
	0
	=
	aE-(\mu+\delta_1+\gamma_1)I_0
	-
	T_{\mathrm{slide}}(B,E)I_0.
	\]
	Therefore,
	\begin{equation}\label{T_slide}
		T_{\mathrm{slide}}(B,E)
		=
		\frac{
			aE-(\mu+\delta_1+\gamma_1)I_0
		}{I_0}.
	\end{equation}
	Substituting $T_{\mathrm{slide}}$ into the $H$-equation gives
	\[
	\frac{dH}{dt}
	=
	T_{\mathrm{slide}}(B,E)I_0
	-
	(\mu+\delta_2+\gamma_2)H
	\]
	and hence
	\[
	\frac{dH}{dt}
	=
	aE-(\mu+\delta_1+\gamma_1)I_0
	-
	(\mu+\delta_2+\gamma_2)H.
	\]
	
	For the delayed model, the current state satisfies $I(t)=I_0$ on $\Sigma_e$, whereas the delayed quantity $I(t-\tau)$ generally retains information from the preceding history. Therefore, the reduced  tangential dynamics on $\Sigma_e$ are 
	\begin{equation}\label{Rep_sys_delay} 
		\begin{aligned}
			\frac{dS}{dt}
			&=
			\Lambda-\beta SI_0-\mu S+\xi R,\\
			\frac{dE}{dt}
			&=
			\beta SI_0-(\mu+a)E,\\
			\frac{dH}{dt}
			&=
			aE-(\mu+\delta_1+\gamma_1)I_0
			-(\mu+\delta_2+\gamma_2)H,\\
			\frac{dR}{dt}
			&=
			\gamma_1I_0+\gamma_2H-(\mu+\xi)R,\\
			\frac{dB}{dt}
			&=
			rB\left(1-\frac{B}{K}\right)
			+
			\frac{\lambda I(t-\tau)B}
			{\omega+I(t-\tau)}.
		\end{aligned}
	\end{equation}
	
	In the delay-free case $\tau=0$, or at a sliding pseudo-equilibrium for which 
	$ 
	I(t)=I(t-\tau)=I_0,
	$ 
	system \eqref{Rep_sys_delay} reduces to
	\begin{equation}\label{Rep_sys} 
		\begin{aligned}
			F_{\Sigma_e}(S,E,H,R,B)
			=
			\begin{pmatrix}
				\Lambda-\beta SI_0-\mu S+\xi R\\[4pt]
				\beta SI_0-(\mu+a)E\\[4pt]
				aE-(\mu+\delta_1+\gamma_1)I_0
				-(\mu+\delta_2+\gamma_2)H\\[4pt]
				\gamma_1I_0+\gamma_2H-(\mu+\xi)R\\[4pt]
				rB\left(1-\frac{B}{K}\right)
				+\frac{\lambda I_0B}{\omega+I_0}
			\end{pmatrix}.
		\end{aligned}
	\end{equation}
	
	Because $\Sigma_e$ is a repelling sliding region, the vector field \eqref{Rep_sys_delay} describes the reduced tangential dynamics associated with the Filippov convexification on $\Sigma_e$. Stability with respect to this reduced vector field should therefore be distinguished from asymptotic stability in the full forward-time Filippov system, since perturbations transverse to $\Sigma_e$ are directed away from the switching manifold. 
	
	\subsubsection{Existence and stability of pseudo-equilibria}
	\label{pseudo_equilibria}
	
	A pseudo-equilibrium is a steady state of the sliding dynamics on the switching manifold $I=I_0$, rather than an equilibrium of either smooth subsystem. An admissible pseudo-equilibrium must lie in the repelling
	sliding region $\Sigma_e$ established in Theorem~\ref{thm_repelling_sliding}. 
	
	When $B=0$, the sliding vector field $F_{\Sigma_e}$ admits an algebraic
	pseudo-equilibrium, denoted by
	\[
	D^{0s}
	=
	(S^{s},E^{s},I_0,H^{s},R^{s},0).
	\]
	
	When $B>0$, the sliding vector field $F_{\Sigma_e}$ admits an algebraic
	pseudo-equilibrium, denoted by
	\[
	D^{s}
	=
	(S^{s},E^{s},I_0,H^{s},R^{s},B^{s}),
	\]
	where
	
	\begin{align*}
		S^s &= \frac{(\mu +a)}{\beta}\Big(\frac{(\mu + \delta_1 + \gamma_1)}{a} - \frac{\gamma_1}{a \gamma_2}(\mu + \delta_2 + \gamma_2) \Big) 
		+ \frac{(\mu +a)}{\beta I_0}\Big(\frac{(\mu + \delta_2 + \gamma_2)(\mu + \xi)}{a \gamma_2} \Big)R^s, \\
		E^s &= \Big(\frac{(\mu + \delta_1 + \gamma_1)}{a} - \frac{\gamma_1}{a \gamma_2}(\mu + \delta_2 + \gamma_2) \Big) I_0 
		+ \Big(\frac{(\mu + \delta_2 + \gamma_2)(\mu + \xi)}{a \gamma_2} \Big)R^s, \\
		H^s &= \frac{(\mu + \xi)}{\gamma_2}R^s - \frac{\gamma_1}{\gamma_2}I_0, \quad 
		R^s = \frac{\Lambda - \frac{(\mu +a)}{\beta}\Big(\frac{(\mu + \delta_1 + \gamma_1)}{a} - \frac{\gamma_1}{a \gamma_2}(\mu + \delta_2 + \gamma_2) \Big) (\beta I_0 + \mu)}
		{\Big( \frac{(\mu +a)}{\beta I_0}\Big(\frac{(\mu + \delta_2 + \gamma_2)(\mu + \xi)}{a \gamma_2} \Big)\Big) (\beta I_0 + \mu) - \xi}, \\
		B^s &= \frac{K}{r}\Big(r + \frac{\lambda I_0}{\omega + I_0} \Big).
	\end{align*}
	
	The above states are admissible pseudo-equilibria only when their components are biologically feasible and they lie in the repelling sliding region. In particular, $D^{0s}$ is admissible if
	\[
	E_{c2}(0)<E^s<E_{c1},
	\]
	whereas $D^s$ is admissible if
	\begin{equation}\label{admissibility_Ds}
		E_{c2}(B^s)<E^s<E_{c1}.
	\end{equation}
	It is worth noting that the location of $D^s$ does not depend explicitly on $\mu_0$ or $\mu_1$. These parameters instead determine the boundaries $E_{c1}$ and $E_{c2}(B^s)$ and therefore determine whether the pseudo-equilibrium is admissible within $\Sigma_e$.
	
	For the stability of pseudo-equilibria, we obtain the following result. 
	
	\begin{theorem}\label{thm:pseudo_stability}
		Assume that $D^s$ is biologically feasible and satisfies the admissibility condition \eqref{admissibility_Ds}. Then $D^s$ is globally asymptotically stable for the algebraically extended reduced tangential system \eqref{Rep_sys} in the region $B>0$. Consequently, every trajectory of the reduced sliding dynamics with $B(0)>0$ that remains in the admissible sliding region $\Sigma_e$ converges to $D^s$. If $D^{0s}$ is biologically feasible and admissible within $\Sigma_e$, then $D^{0s}$ is unstable with respect to the reduced sliding dynamics. 
	\end{theorem}
	
	The proof is provided in Supplementary Section~S11. 
	
	The above result shows that $D^s$ has a distinct geometric role in the Filippov system. Within the algebraically extended tangential dynamics, trajectories with positive healthcare capacity converge to $D^s$. However, since $\Sigma_e$ is a repelling sliding region, this tangential attraction does not imply attraction in the full forward-time Filippov system. Perturbations transverse to the switching manifold are directed away from $\Sigma_e$ and can drive the epidemic into either of the two hospitalization regimes. Thus, when admissible, $D^s$ has a saddle-type organizing role: it is attracting along the reduced tangential dynamics but repelling in the transverse Filippov direction. In contrast, the instability of $D^{0s}$ shows that the zero-capacity  pseudo-equilibrium is unstable even with respect to the reduced sliding dynamics.
	
	\subsubsection{Tangency points}  
	
	Tangency points occur where the normal component of one of the smooth vector fields vanishes. Using the critical values $E_{c1}$ and $E_{c2}(B)$ derived above, the tangency sets are
	\[
	\Sigma^t_1=\{X\in\Sigma:E=E_{c1}\},\quad
	\Sigma^t_2=\{X\in\Sigma:E=E_{c2}(B)\}.
	\]
	These sets bound the repelling sliding region $\Sigma_e$ and represent generic quadratic tangencies under the nondegeneracy conditions established in Supplementary Section~S8.
	
	\subsubsection{Boundary equilibria} 
	
	A boundary equilibrium is an equilibrium of one of the smooth subsystems located on the switching manifold $I=I_0$. Under a regular crossing condition, it marks a critical parameter value at which an equilibrium branch crosses the switching manifold and changes from an admissible equilibrium to a virtual equilibrium, or vice versa. Its existence therefore indicates a possible transition between the below-threshold and capacity-constrained hospitalization regimes. Thus, a boundary equilibrium $D^b$ associated with subsystem $F_{\Sigma_i}$ satisfies 
	\[
	F_{\Sigma_i}(S,E,I_0,H,R,B)=0,
	\quad i=1,2.
	\]
	
	At an equilibrium with $I=I_0$, the healthcare-capacity equation yields two possible equilibrium capacity levels 
	\[
	B^b_1=0,
	\quad
	B^b_2=
	K\left(1+\frac{\lambda I_0}{r(\omega+I_0)}\right).
	\]
	However, these conditions on $B$ alone are not sufficient for the existence of a boundary equilibrium. The normal component of the corresponding smooth vector field must also vanish.
	
	For $F_{\Sigma_1}$, the $I$-component of the equilibrium condition is equivalent to $E=E_{c1}$, whereas for $F_{\Sigma_2}$ it is equivalent to $E=E_{c2}(B)$.
	
	Therefore, a pseudo-equilibrium can coincide with a boundary equilibrium
	only when one of its admissibility inequalities becomes an equality.
	
	In particular, $D^{0s}$ may coincide with a boundary equilibrium when
	\[
	E^s=E_{c1}
	\quad \text{or} \quad
	E^s=E_{c2}(0),
	\]
	whereas $D^s$ may coincide with a boundary equilibrium when
	\[
	E^s=E_{c1}
	\quad \text{or} \quad
	E^s=E_{c2}(B^s).
	\]
	These equalities impose additional scalar conditions on the model
	parameters and, under the usual nondegeneracy condition, define
	codimension-one boundary-equilibrium loci in parameter space. Away from
	these loci, $D^{0s}$ and $D^s$ remain equilibria of the algebraically
	extended reduced tangential dynamics but are not equilibria of either
	smooth subsystem. As parameters vary, a pseudo-equilibrium may therefore enter or leave
	the admissible sliding region through one of its tangency boundaries,
	providing the connection between the pseudo-equilibria described above
	and the boundary-equilibrium transitions considered subsequently. 
	
	\subsubsection{Regular and virtual equilibria} 
	
	Regular equilibria are admissible steady states that lie within the regime under which they are defined and can therefore represent feasible long-term epidemic outcomes. Virtual equilibria satisfy the equilibrium equations of a subsystem but lie outside its corresponding region of validity. Hence, they are not directly realizable under that regime, although they may influence changes in admissibility and nearby switching behavior. 
	
	The disease-free equilibria $D_{01}$ and $D_{02}$ are independent of the threshold $I_0$. For the subsystem with $I<I_0$, an endemic equilibrium $D^*$ is regular if $I^*<I_0$, and virtual if $I^*>I_0$. For the subsystem with $I>I_0$, an endemic equilibrium $D_1^*$ is regular if $I_1^*>I_0$, and virtual if $I_1^*<I_0$. Crossing the threshold $I=I_0$ generates a non-smooth boundary-collision rather than a smooth bifurcation. 
	
	\subsubsection{Codimension-two Hopf-boundary interaction}
	\label{subsec:codim2_HB} 
	
	We investigate the interaction between a Hopf bifurcation of the smooth
	extension of subsystem $F_{\Sigma_2}$ and a boundary-equilibrium
	transition of the full Filippov system. We consider the maximum
	hospitalization rate $\mu_1$ and the baseline healthcare capacity $K$
	as bifurcation parameters and define
	\[
	p=(\mu_1,K)\in\mathbb{R}_{+}^{2}.
	\]
	
	Codimension-one boundary and sliding bifurcations, as well as
	two-parameter degeneracies in Filippov systems, have been studied in
	general settings \cite{BERN2008,Kuznetsov2003,Kowalczyk2005,Kowalczyk2006,
		Guardia2011,Dercole2011,Simpson2018}. In particular, boundary-Hopf
	interactions and the associated cycle-grazing loci have been
	characterized for planar Filippov systems \cite{Dercole2011}.
	Here, our contribution is the explicit characterization of the
	Hopf-boundary interaction in a system with dynamically evolving
	resource capacity. We derive an analytical boundary-equilibrium locus
	and establish conditions for its transversal intersection with the Hopf
	locus. The quadratic first-contact scaling itself follows from the
	standard square-root amplitude relation near a nondegenerate Hopf
	bifurcation. Our analysis identifies how this local geometry arises
	within the capacity-dependent repelling sliding structure of the
	proposed Filippov system.
	
	\medskip
	\noindent
	\textbf{Analytical boundary-equilibrium locus:}
	Suppose that a positive endemic-equilibrium branch of the smooth
	extension of $F_{\Sigma_2}$ reaches the switching manifold, so that 
	$
	I^*=I_0.
	$ 
	For a positive healthcare-capacity equilibrium $B^*>0$, the healthcare-capacity 
	equation gives
	\[
	0=
	rB^*
	\left(1-\frac{B^*}{K}\right)
	+
	\frac{\lambda I_0B^*}{\omega+I_0}.
	\]
	Hence,
	\begin{equation}
		B^*=B_b(K)
		=
		K\left(
		1+\frac{\lambda I_0}{r(\omega+I_0)}
		\right)
		=:Q_0K,
		\quad
		Q_0
		=
		1+\frac{\lambda I_0}{r(\omega+I_0)}>0.
		\label{eq:B_boundary}
	\end{equation}
	
	Let
	\[
	T_b
	=
	\mu_0+
	(\mu_1-\mu_0)
	\frac{B_b}{B_b+I_0}
	\]
	denote the hospitalization rate associated with the smooth extension of
	$F_{\Sigma_2}$ at a positive-capacity boundary equilibrium. Setting the 
	$E$, $I$, $H$, and $R$ equations equal to zero gives
	\begin{equation}
		E_b
		=
		\frac{(\mu+\delta_1+\gamma_1+T_b)I_0}{a},
		\quad
		S_b
		=
		\frac{(\mu+a)(\mu+\delta_1+\gamma_1+T_b)}
		{a\beta},
		\label{eq:boundary_SE}
	\end{equation}
	and
	\begin{equation}
		H_b
		=
		\frac{T_bI_0}{\mu+\delta_2+\gamma_2},
		\quad
		R_b
		=
		\frac{I_0}{\mu+\xi}
		\left(
		\gamma_1+
		\frac{\gamma_2T_b}{\mu+\delta_2+\gamma_2}
		\right).
		\label{eq:boundary_HR}
	\end{equation}
	
	Substituting \eqref{eq:boundary_SE} and
	\eqref{eq:boundary_HR} into the susceptible equation yields
	\[
	0
	=
	\Lambda
	-
	(\beta I_0+\mu)
	\frac{(\mu+a)(\mu+\delta_1+\gamma_1+T_b)}
	{a\beta}
	+
	\frac{\xi I_0}{\mu+\xi}
	\left(
	\gamma_1+
	\frac{\gamma_2T_b}{\mu+\delta_2+\gamma_2}
	\right).
	\]
	Therefore,
	\begin{equation}
		T_b
		=
		\frac{
			\Lambda-A(\mu+\delta_1+\gamma_1)+D\gamma_1
		}{
			A-\dfrac{D\gamma_2}{\mu+\delta_2+\gamma_2}
		},
		\label{eq:T_boundary}
	\end{equation}
	provided
	\[
	A-\frac{D\gamma_2}{\mu+\delta_2+\gamma_2}\neq0,
	\]
	where
	\[
	A
	=
	\frac{(\beta I_0+\mu)(\mu+a)}{a\beta},
	\quad
	D
	=
	\frac{\xi I_0}{\mu+\xi}.
	\]
	
	For fixed values of the remaining model parameters, $T_b$ is independent
	of $(\mu_1,K)$. On the other hand, the hospitalization law gives
	\[
	T_b
	=
	\mu_0+
	(\mu_1-\mu_0)
	\frac{Q_0K}{Q_0K+I_0}.
	\]
	Solving for $\mu_1$ gives the positive-capacity  boundary-equilibrium locus
	\begin{equation}
		\mu_1^{\,b}(K)
		=
		T_b+
		\frac{(T_b-\mu_0)I_0}{Q_0K}.
		\label{eq:boundary_curve}
	\end{equation}
	Thus, within the biologically feasible parameter region, a positive
	endemic-equilibrium branch of the smooth extension of $F_{\Sigma_2}$
	reaches the switching manifold precisely when
	\[
	\mu_1=\mu_1^{\,b}(K),
	\]
	provided that the corresponding equilibrium components are positive and
	$\mu_1^{\,b}(K)>\mu_0$.
	
	Define
	\begin{equation}
		\mathcal{B}(\mu_1,K)
		=
		\mu_1-\mu_1^{\,b}(K).
		\label{eq:boundary_function}
	\end{equation}
	Then $\mathcal{B}=0$ defines the positive-capacity  boundary-equilibrium 
	curve. Since
	\[
	\frac{\partial\mathcal{B}}{\partial\mu_1}=1,
	\]
	one has
	\[
	\nabla_{(\mu_1,K)}\mathcal{B}\neq0,
	\]
	and therefore the boundary-equilibrium locus is a regular
	codimension-one curve. Moreover,
	\[
	\frac{d\mu_1^{\,b}}{dK}
	=
	-\frac{(T_b-\mu_0)I_0}{Q_0K^2},
	\]
	so the boundary-equilibrium curve depends nontrivially on $K$  whenever
	$T_b\neq\mu_0$.
	
	\medskip
	\noindent
	\textbf{Hopf locus of the smooth subsystem $F_{\Sigma_2}$:}
	Let
	$
	D^*(p)
	=
	\bigl(
	S^*(p),E^*(p),I^*(p),H^*(p),R^*(p),B^*(p)
	\bigr)
	$ 
	denote a local endemic-equilibrium branch of the smooth extension of
	$F_{\Sigma_2}$. As derived in Section~\ref{subdyn2}, its infectious
	component satisfies
	\begin{equation}
		G(I;p)
		=
		G_1(p)I^3+G_2(p)I^2+G_3(p)I+G_4(p)=0.
		\label{eq:G_codim2}
	\end{equation}
	Suppose that at $p=p^*$,
	\begin{equation}
		G_I(I_0;p^*)\neq0.
		\label{eq:simple_equilibrium_condition}
	\end{equation}
	Then the implicit function theorem yields a unique smooth local
	equilibrium branch $I^*=I^*(p)$ through $I^*(p^*)=I_0$, and all
	components of $D^*(p)$ depend smoothly on $p$ near $p^*$.
	
	The characteristic polynomial of the Jacobian matrix at $D^*(p)$,
	see Eq.~\eqref{chara_end}, is
	\begin{equation}
		C(s;p)
		=
		s^6+f_{11}(p)s^5+f_{12}(p)s^4
		+f_{13}(p)s^3+f_{14}(p)s^2
		+f_{15}(p)s+f_{16}(p),
		\label{eq:char_codim2}
	\end{equation}
	where the coefficients $f_{1j}$ are given in
	Supplementary Section~S5.
	
	Let $H_j(p)$, $j=1,\ldots,6$, denote the Hurwitz determinants
	introduced in Section~\ref{hopf1}. At a Hopf point $p^*$ we require
	\begin{equation}
		f_{16}>0,\quad
		H_j>0,\quad j=1,2,3,4,
		\quad
		H_5=0,
		\label{eq:Hopf_Hurwitz_conditions}
	\end{equation}
	together with a simple complex-conjugate critical pair
	\[
	s_{1,2}(p^*)=\pm i\nu^*,
	\quad
	\nu^*>0,
	\]
	while all remaining eigenvalues have strictly negative real parts.
	
	Because the critical pair is simple, there exists locally a smooth
	critical eigenvalue branch
	\[
	s_H(p)=\alpha(p)+i\nu(p),
	\quad
	s_H(p^*)=i\nu^*.
	\]
	We define
	\begin{equation}
		\mathcal{H}(\mu_1,K)
		=
		\alpha(\mu_1,K)
		=
		\Re s_H(\mu_1,K).
		\label{eq:Hopf_function}
	\end{equation}
	The local Hopf locus is therefore given by
	\[
	\mathcal{H}(\mu_1,K)=0.
	\]
	
	\begin{theorem} 
		\label{thm:codim2_HB}
		Consider the delay-free Filippov system with
		$p=(\mu_1,K)$ as the bifurcation parameter vector.
		Suppose that there exists a biologically admissible point 
		\[
		p^*=(\mu_1^*,K^*)
		\]
		satisfying the following conditions: 
		
		\begin{enumerate}
			\item[(H1)]
			A positive endemic-equilibrium branch of the smooth extension of
			$F_{\Sigma_2}$ reaches the switching manifold:
			$
			I^*(p^*)=I_0,
			$ equivalently
			$
			\mathcal{B}(p^*)=0.
			$ 
			
			\item[(H2)]
			The equilibrium root is simple:
			$
			G_I(I_0;p^*)\neq0.
			$ 
			
			\item[(H3)]
			The Hopf conditions
			$
			f_{16}>0,$
			$
			H_j>0,  j=1,2,3,4,
			H_5=0
			$ 
			hold, and the critical spectrum consists of a simple pair
			$
			s_{1,2}(p^*)=\pm i\nu^*, 
			\nu^*>0,
			$ 
			while all remaining eigenvalues have strictly negative real parts.
			
			\item[(H4)]
			The boundary-equilibrium and Hopf conditions are independent:
			\begin{equation}
				\det
				\begin{pmatrix}
					\dfrac{\partial\mathcal{B}}{\partial\mu_1}
					&
					\dfrac{\partial\mathcal{B}}{\partial K}
					\\[3mm]
					\dfrac{\partial\mathcal{H}}{\partial\mu_1}
					&
					\dfrac{\partial\mathcal{H}}{\partial K}
				\end{pmatrix}_{p=p^*}
				\neq0.
				\label{eq:codim2_det}
			\end{equation}
		\end{enumerate}
		
		Then $p^*$ is a transversal codimension-two Hopf-boundary
		interaction point: the codimension-one Hopf locus of the smooth
		extension of $F_{\Sigma_2}$ and the codimension-one
		boundary-equilibrium locus of the Filippov system intersect
		transversally at $p^*$.
		
		If, additionally,
		$
		\ell_1(p^*)\neq0,
		$ 
		then the Hopf bifurcation of the smooth extension of
		$F_{\Sigma_2}$ is nondegenerate.
	\end{theorem}
	
	\begin{proof}
		Condition (H2) and the implicit function theorem provide a smooth local
		endemic-equilibrium branch $D^*(p)$ through $p^*$. Since
		$\partial\mathcal B/\partial\mu_1=1$, the boundary-equilibrium condition
		$\mathcal B=0$ defines a regular codimension-one curve.
		
		Condition \eqref{eq:codim2_det} implies
		$\nabla\mathcal H(p^*)\neq0$, so the Hopf condition
		$\mathcal H=0$ also defines a regular codimension-one curve.
		Moreover,
		\[
		\det D_{(\mu_1,K)}(\mathcal B,\mathcal H)(p^*)\neq0
		\]
		implies that
		$\nabla\mathcal B(p^*)$ and $\nabla\mathcal H(p^*)$ are linearly
		independent. Hence, the two curves intersect transversally at an
		isolated point $p^*$, and their simultaneous occurrence has
		codimension two. If $\ell_1(p^*)\neq0$, the standard Hopf
		nondegeneracy condition is satisfied for the smooth subsystem.
		
		It remains to verify the nonsmooth character of the boundary collision.
		Let
		$
		L(X)=I-I_0$, with $
		\nabla L=(0,0,1,0,0,0)^{\top}.
		$
		Because $D^*(p^*)$ is an equilibrium of the smooth extension of
		$F_{\Sigma_2}$ on $\Sigma$,
		\[
		\left\langle
		\nabla L,
		F_{\Sigma_2}(D^*(p^*))
		\right\rangle=0.
		\]
		Writing $T_2$ for the above-threshold hospitalization law, the infectious equilibrium equation gives
		\[
		aE^*
		=
		\left(
		\mu+\delta_1+\gamma_1+T_2(I_0,B^*)
		\right)I_0.
		\]
		Consequently,
		\begin{align*}
			\left\langle
			\nabla L,
			F_{\Sigma_1}(D^*(p^*))
			\right\rangle
			&=
			aE^*-(\mu+\delta_1+\gamma_1+\mu_1^*)I_0\\
			&=
			\left(T_2(I_0,B^*)-\mu_1^*\right)I_0.
		\end{align*}
		Since
		\[
		T_2(I_0,B^*)
		=
		\mu_0+
		(\mu_1^*-\mu_0)
		\frac{B^*}{B^*+I_0},
		\]
		we obtain
		\begin{equation}
			\left\langle
			\nabla L,
			F_{\Sigma_1}(D^*(p^*))
			\right\rangle
			=
			-
			\frac{(\mu_1^*-\mu_0)I_0^2}
			{B^*+I_0}
			<0.
			\label{eq:normal_F1_boundary}
		\end{equation}
		Thus, the opposite vector field has a nonzero normal component at the
		collision point. The event is therefore a boundary-equilibrium
		collision rather than a simultaneous tangency of the two smooth vector
		fields.
	\end{proof}
	
	\noindent
	\textbf{Local interaction of the Hopf-generated periodic orbit with the
		switching manifold:}
	Theorem~\ref{thm:codim2_HB} establishes that the endemic equilibrium
	simultaneously reaches the switching manifold and satisfies the Hopf
	condition at $p=p^*$. To determine whether the small periodic orbit
	generated by this Hopf bifurcation also reaches the switching manifold,
	define
	\begin{equation}
		\Delta(p)
		=
		I^*(p)-I_0.
		\label{eq:Delta_HB}
	\end{equation}
	Thus,
	\[
	\Delta(p^*)=0.
	\]
	Locally, $\Delta=0$ and $\mathcal B=0$ describe the same boundary-equilibrium locus, while $\Delta$ measures the signed displacement of the smooth endemic equilibrium from the switching manifold in the infectious direction. 
	
	\begin{proposition} 
		\label{prop:HB_grazing}
		Let $p^*$ be the transversal codimension-two Hopf-boundary point
		established in Theorem~\ref{thm:codim2_HB}. Assume, in addition, that
		
		\begin{enumerate}
			\item[(G1)]
			the smooth Hopf bifurcation is nondegenerate:
			$ 
			\ell_1(p^*)\neq0;
			$ 
			
			\item[(G2)]
			the equilibrium branch crosses the switching manifold regularly:
			$ 
			\nabla_p\Delta(p^*)\neq0;
			$ 
			
			\item[(G3)]
			if $q\in\mathbb C^6$ is a right eigenvector of
			$J_2(p^*)$ associated with $i\nu^*$, then its infectious component
			is nonzero:
			$ 
			q_I\neq0.
			$ 
		\end{enumerate}
		
		Then the small-amplitude periodic-orbit family generated by the smooth Hopf bifurcation has a nonzero component transverse to
		\[
		\Sigma=\{I=I_0\}.
		\]
		On the side $\Delta>0$, where the equilibrium lies in the admissible
		region of $F_{\Sigma_2}$, there exists a local first-contact curve of
		this periodic orbit with $\Sigma$ emanating from $p^*$.
		
		In suitable regular local unfolding coordinates $(\Delta,\eta)$, where
		$\eta=0$ denotes the Hopf locus and the orientation of $\eta$ is chosen
		so that the small periodic orbit exists for $\eta<0$, the first-contact
		curve satisfies
		\begin{equation}
			\eta
			=
			-C_g\Delta^2
			+
			o(\Delta^2),
			\quad
			C_g>0.
			\label{eq:HB_grazing_curve}
		\end{equation}
		Hence, the first-contact locus is quadratically tangent to the Hopf
		locus at $p^*$. 
	\end{proposition}
	
	\begin{proof}
		By Theorem~\ref{thm:codim2_HB}, $J_2(p^*)$ possesses a simple critical
		pair
		$
		\pm i\nu^*, 
		\nu^*>0,
		$
		with all remaining eigenvalues having strictly negative real parts.
		Therefore, the smooth extension of $F_{\Sigma_2}$ possesses a
		two-dimensional local center manifold. 
		
		Since $\Delta=0$ and $\mathcal B=0$ locally define the same
		boundary-equilibrium curve, conditions
		\[
		\nabla\Delta(p^*)\neq0,
		\quad
		\nabla\mathcal B(p^*)\neq0
		\]
		imply that their gradients are nonzero scalar multiples of one another
		at $p^*$. Condition (H4) therefore implies
		\[
		\det D_p(\Delta,\mathcal H)(p^*)\neq0.
		\]
		Consequently, $\Delta$ together with a suitably oriented Hopf coordinate
		$\eta$ forms a regular local coordinate system near $p^*$.
		
		In these coordinates, the parameter-dependent Hopf normal form on the
		center manifold can be written as
		\begin{equation}
			\dot z
			=
			\bigl(\sigma\eta+i\nu(\Delta,\eta)\bigr)z
			+
			c(\Delta,\eta)z|z|^2
			+
			\mathcal O(|z|^4),
			\label{eq:HB_Hopf_normalform}
		\end{equation}
		where
		\[
		\nu(0,0)=\nu^*,
		\quad
		\Re c(0,0)\neq0,
		\quad
		\sigma=\operatorname{sgn}(\Re c(0,0))
		\]
		by condition (G1). Here $\eta=\sigma\mathcal H(p)$, so that the small periodic
		orbit exists for $\eta<0$ and its amplitude satisfies
		\begin{equation}
			\rho
			=
			C_H\sqrt{-\eta}
			+
			o\!\left(\sqrt{|\eta|}\right),
			\quad
			C_H>0,
			\label{eq:HB_Hopf_amplitude}
		\end{equation}
		to leading order near $(\Delta,\eta)=(0,0)$.
		
		Recall that
		\[
		L(X)=I-I_0.
		\]
		For a state on the center manifold near the equilibrium branch,
		\[
		X=D^*(p)+zq+\bar z\,\bar q+\mathcal O(|z|^2),
		\]
		and therefore
		\begin{equation}
			L(X)
			=
			\Delta(p)
			+
			2\Re(q_Iz)
			+
			\mathcal O(|z|^2).
			\label{eq:L_center_manifold}
		\end{equation}
		Because $q_I\neq0$, the critical Hopf mode has a nonzero component
		normal to $\Sigma$.
		
		Along the small periodic orbit,
		\[
		z(t)
		=
		\rho e^{i\theta(t)}
		+
		\mathcal O(\rho^2),
		\]
		so
		\[
		L(t)
		=
		\Delta
		+
		2\rho\Re\!\left(q_Ie^{i\theta(t)}\right)
		+
		\mathcal O(\rho^2).
		\]
		Hence,
		\begin{equation}
			L_{\min}
			=
			\Delta
			-
			2|q_I|\rho
			+
			\mathcal O(\rho^2).
			\label{eq:Lmin_HB}
		\end{equation}
		
		For an orbit lying entirely in the smooth region $F_{\Sigma_2}$,
		$L_{\min}>0$. Its first contact with $\Sigma$ occurs when
		$L_{\min}=0$, giving
		\begin{equation}
			\Delta
			=
			2|q_I|\rho
			+
			\mathcal O(\rho^2).
			\label{eq:Delta_rho_HB}
		\end{equation}
		In particular, the relevant first-contact branch satisfies
		$\Delta>0$.
		
		Substituting \eqref{eq:HB_Hopf_amplitude} into
		\eqref{eq:Delta_rho_HB} yields
		\[
		\Delta
		=
		2|q_I|C_H\sqrt{-\eta}
		+
		o\!\left(\sqrt{|\eta|}\right).
		\]
		Therefore,
		\[
		\eta
		=
		-\frac{1}{4|q_I|^2C_H^2}\Delta^2
		+
		o(\Delta^2),
		\]
		and hence
		\[
		C_g
		=
		\frac{1}{4|q_I|^2C_H^2}>0.
		\]
		This proves \eqref{eq:HB_grazing_curve}.
		
		Moreover, at the phase at which $L$ attains its minimum,
		$\dot L=0$, while the leading harmonic term gives 
		\[
		\ddot L
		=
		2|q_I|\rho(\nu^*)^2
		+
		\mathcal O(\rho^2)
		>0
		\]
		for sufficiently small $\rho>0$. Thus, the first contact of the
		Hopf-generated orbit with $\Sigma$ is a nondegenerate quadratic
		tangency at leading order. 
	\end{proof}
	
	Proposition~\ref{prop:HB_grazing} provides the local dynamical
	interpretation of the codimension-two Hopf-boundary point. The
	boundary-equilibrium and Hopf loci intersect transversally at $p^*$,
	whereas the first-contact locus of the Hopf-generated periodic orbit is
	quadratically tangent to the Hopf locus. The proposition describes only
	the local approach and first contact with the switching manifold.
	Dynamics after contact, including possible crossing or interaction with
	the repelling sliding region, require a separate nonsmooth return-map
	or discontinuity-map analysis and are not classified here.
	
	\subsection{Non-smooth dynamics of the delayed model on the switching manifold} 
	
	The time delay $\tau>0$ enters only through the 
	healthcare-capacity equation and does not affect the 
	normal components of the vector fields associated with 
	$L(I)=I-I_0$. Consequently, the switching manifold 
	$\Sigma=\{I=I_0\}$, sliding conditions, tangency curves, 
	and repelling sliding region $\Sigma_e$ remain identical 
	to those of the delay-free system.
	
	Restricting the delayed system to $\Sigma_e$ yields 
	system~\eqref{Rep_sys_delay}. Along sliding motion, 
	$I(t)=I_0$, while the delayed argument $I(t-\tau)$ 
	retains the pre-sliding history during the first delay 
	interval. Thus, the delay affects only the $B$-dynamics, 
	leaving the $(S,E,H,R)$-subsystem unchanged. 
	
	\subsubsection{Effect of delay on sliding pseudo-equilibria}
	
	At either pseudo-equilibrium $D^{0s}$ or $D^{s}$, the equilibrium
	history is constant and satisfies
	$
	I(t)=I(t-\tau)=I_0.
	$ 
	Consequently, the delay does not alter the algebraic location of the
	sliding pseudo-equilibria.
	
	\begin{proposition}\label{prop:sliding_delay_stability}
		For every $\tau\geq0$ and for histories contained in the switching
		manifold, the healthcare-capacity response delay does not alter the
		stability properties of the admissible pseudo-equilibria with respect to the
		algebraically extended reduced tangential dynamics. In particular, $D^{0s}$ remains unstable in
		the $B$-direction, whereas $D^{s}$ retains the global asymptotic
		stability for $B>0$ established in Theorem~\ref{thm:pseudo_stability}.
	\end{proposition}
	
	The proof is provided in Supplementary Section~S12.
	
	If a trajectory reaches $\Sigma$ with a preceding history not entirely
	contained in the switching manifold, then $I(t-\tau)$ may initially
	differ from $I_0$ and act as a transient forcing in the $B$-equation.
	Provided the trajectory remains on $\Sigma$, this effect disappears
	after one delay interval. Thus, the healthcare-capacity response delay may influence the
	short-term capacity adjustment after entry to the switching manifold,
	but it does not alter the stability properties of the pseudo-equilibria
	established above. Moreover, since $\Sigma_e$ is repelling by
	Theorem~\ref{thm_repelling_sliding}, these tangential stability
	properties do not imply attraction in the full Filippov system. 
	
	\subsection{Uniform disease elimination under a worst-case threshold}  
	
	Define the worst-case threshold 
	\[
	\mathcal R_0^{\max}
	=
	\frac{\beta(\Lambda/\mu)\,a}
	{(\mu+a)\bigl(\mu+\delta_1+\gamma_1+\mu_0\bigr)}.
	\]
	This quantity is obtained by combining the asymptotic upper bound
	$\limsup_{t\to\infty}S(t)\leq\Lambda/\mu$ with the minimum hospitalization rate
	$T(I,B)\geq\mu_0$. It therefore provides a sufficient, regime-independent
	condition for decay of the infected subsystem across both smooth regions and
	on the switching manifold. It is not the basic reproduction number obtained
	by linearization at the disease-free equilibrium.
	
	Recall that the basic reproduction number associated with invasion near the
	disease-free equilibrium in the below-threshold regime is
	\[
	\mathcal R_0
	=
	\frac{\beta(\Lambda/\mu)\,a}
	{(\mu+a)\bigl(\mu+\delta_1+\gamma_1+\mu_1\bigr)}.
	\]
	Since $0<\mu_0<\mu_1$, it follows that
	\[
	\mathcal R_0<\mathcal R_0^{\max}.
	\]
	Thus, $\mathcal R_0<1$ characterizes local disease invasion near the
	disease-free equilibrium, whereas $\mathcal R_0^{\max}<1$ is a stronger
	condition that guarantees disease elimination uniformly across the switching
	regimes.
	
	\begin{theorem}
		\label{thm:global_dfe_worstcase}
		Suppose that $\mathcal R_0^{\max}<1$. Then, for every admissible
		nonnegative initial history and every delay $\tau\geq0$: 
		\begin{itemize}
			\item if $B(0)=0$, then $B(t)=0$ for all $t\geq0$, and the solution
			converges to
			$ 
			D_{01}=(\Lambda/\mu,0,0,0,0,0);
			$ 
			\item if $B(0)>0$, then $B(t)\to K$ and the solution converges to
			$
			D_{02}=(\Lambda/\mu,0,0,0,0,K).
			$ 
		\end{itemize}
	\end{theorem}
	
	The proof is provided in Supplementary Section~S13. 
	
	Since
	$ 
	\mathcal R_0<\mathcal R_0^{\max},
	$
	three parameter regimes can be distinguished. If
	$\mathcal R_0^{\max}<1$, uniform elimination follows from
	Theorem~\ref{thm:global_dfe_worstcase}. If
	$ 
	\mathcal R_0<1\leq\mathcal R_0^{\max},
	$ 
	the disease-free equilibrium remains locally stable to small infectious
	perturbations, but uniform elimination across both hospitalization
	regimes is not guaranteed. This intermediate regime may therefore
	permit threshold- and initial-condition-dependent dynamics, including
	coexistence with an admissible endemic state when such a state exists.
	If $\mathcal R_0>1$, the disease-free equilibrium is unstable to
	infection invasion.
	
	\section{A data-informed parameterization using physics-informed neural networks} 
	\label{Est_param}
	
	We provide a data-informed illustration of model calibration using
	cumulative reported COVID-19 cases from Kenya. During the COVID-19
	pandemic, healthcare systems in many countries experienced substantial
	resource constraints, particularly during major outbreak waves
	\cite{AJB2021}. We therefore use Kenya as a case study and obtain
	cumulative reported case data from the World Health Organization (WHO)
	\cite{WHO_Data}. The resulting data set consists of 236 weekly
	cumulative-case observations. The purpose of this analysis is to provide
	a data-informed parameterization of selected model parameters rather than
	a formal statistical identification of all epidemiological and
	healthcare-capacity parameters.
	
	\begin{table}[H]
		\centering
		\caption{Fixed and calibrated parameter values used in the 
			PINN-based fitting to cumulative reported COVID-19 cases for Kenya.}
		\label{param_val}
		\begin{tabular}{|p{2cm}|p{8cm}|p{2cm}|}
			\hline
			Parameters & Numerical values & References \\
			\hline
			
			$a$ & $1/2$ week$^{-1}$ & \cite{WANG2022} \\
			
			$\delta_1$ & $0.0119$ week$^{-1}$ & \cite{ZHOU2020} \\
			
			$\delta_2$ & $0.0077$ week$^{-1}$ & \cite{ZHOU2020} \\
			
			$\gamma_1$ & $1/2$ week$^{-1}$ & \cite{WANG2022,NIT2021} \\
			
			$\gamma_2$ & $1.6758$ week$^{-1}$ & \cite{ZHOU2020} \\
			
			$\xi$ & $0.0389$ week$^{-1}$ & \cite{BUG2023,CHOW2022} \\
			
			$r$ & $1.3139$ week$^{-1}$ & \cite{ZHOU2020} \\
			
			$\omega$ & $50$ individuals & Assumed \\
			
			$\mu_1$ & $0.77$ week$^{-1}$ & \cite{Bugalia2020} \\
			
			$\mu_0$ & $0.0182$ week$^{-1}$ & \cite{AJB2021} \\
			
			$\mu$ & $0.00028$ week$^{-1}$ & \cite{KenyaLifeExpectancy} \\ 
			
			$K$ & $1000$ individual-equivalent capacity units & Assumed \\
			
			$I_0$ & $4500$ individuals & Assumed \\
			
			$\Lambda$ & $1.4621\times10^{4}$ individuals week$^{-1}$ 
			& Computed \\
			
			$\beta$ & $1.2648\times10^{-8}$ individual$^{-1}$ week$^{-1}$ 
			& Calibrated \\
			
			$\lambda$ & $0.8583$ week$^{-1}$ & Calibrated \\
			
			$\tau$ & $0.9065$ weeks & Calibrated \\
			
			\hline
		\end{tabular}
	\end{table}
	
	\medskip
	\noindent \textbf{Parameter calibration:} 
	We treat the transmission rate $\beta$, the infection-responsive 
	healthcare-capacity expansion rate $\lambda$, and the capacity-response 
	delay $\tau$ as trainable parameters together with the neural-network 
	weights, and optimize them using the Adam optimizer. The remaining model 
	parameters are fixed using values reported in the literature or assumed 
	values chosen to represent biologically relevant time scales. All rate 
	parameters and the time variable are expressed consistently on a weekly 
	time scale. 
	
	The normalization, network architecture, loss function, and implementation 
	details are provided in Supplementary Section~S9.
	
	Because the observed data used in the PINN consist only of cumulative 
	reported cases, whereas $S(t)$, $E(t)$, $I(t)$, $H(t)$, $R(t)$, and $B(t)$ 
	are not directly observed in the present calibration, the calibrated values 
	of $\beta$, $\lambda$, and $\tau$ should be interpreted as effective 
	parameter values that provide a model trajectory consistent with the 
	cumulative-case data and the imposed dynamical constraints. In particular, 
	the present analysis does not establish structural or practical 
	identifiability of these parameters. 
	
\begin{figure}[H]
		\centering
		(a)\includegraphics[width=0.3\textwidth]
		{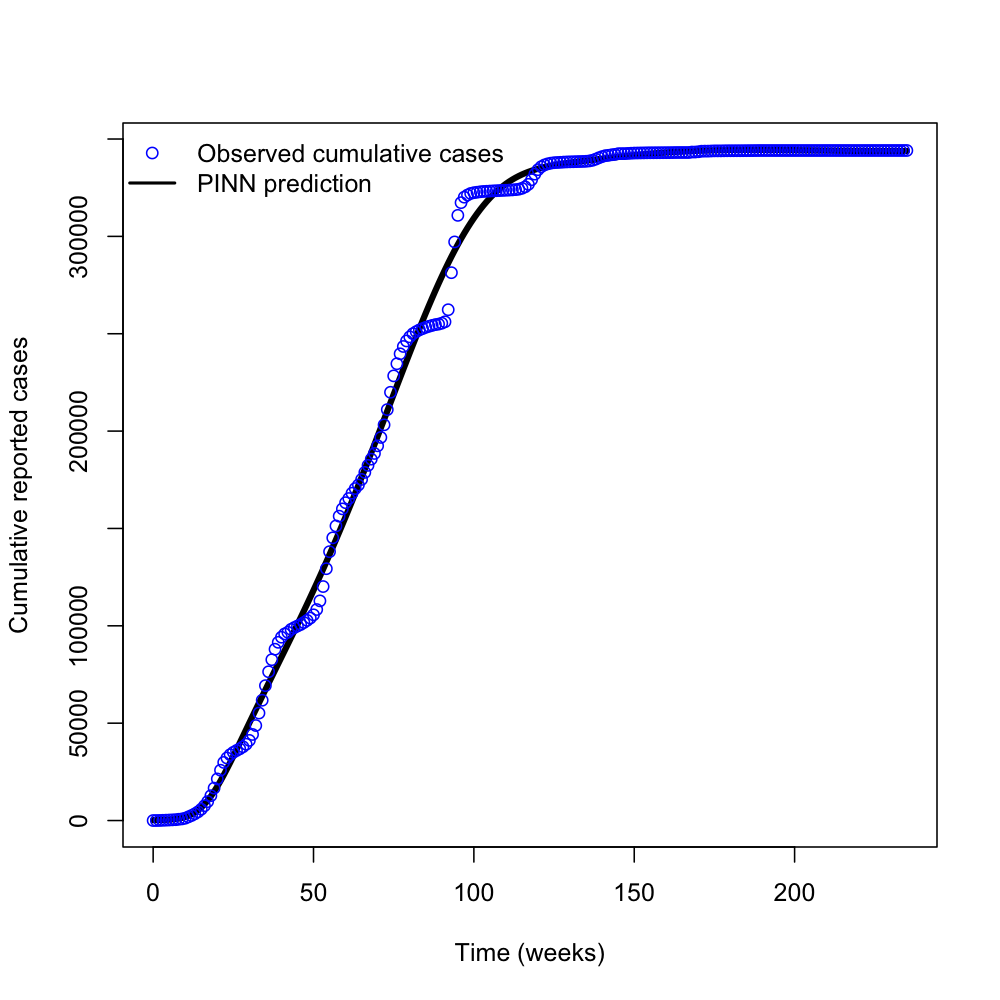}
		(b)\includegraphics[width=0.3\textwidth]
		{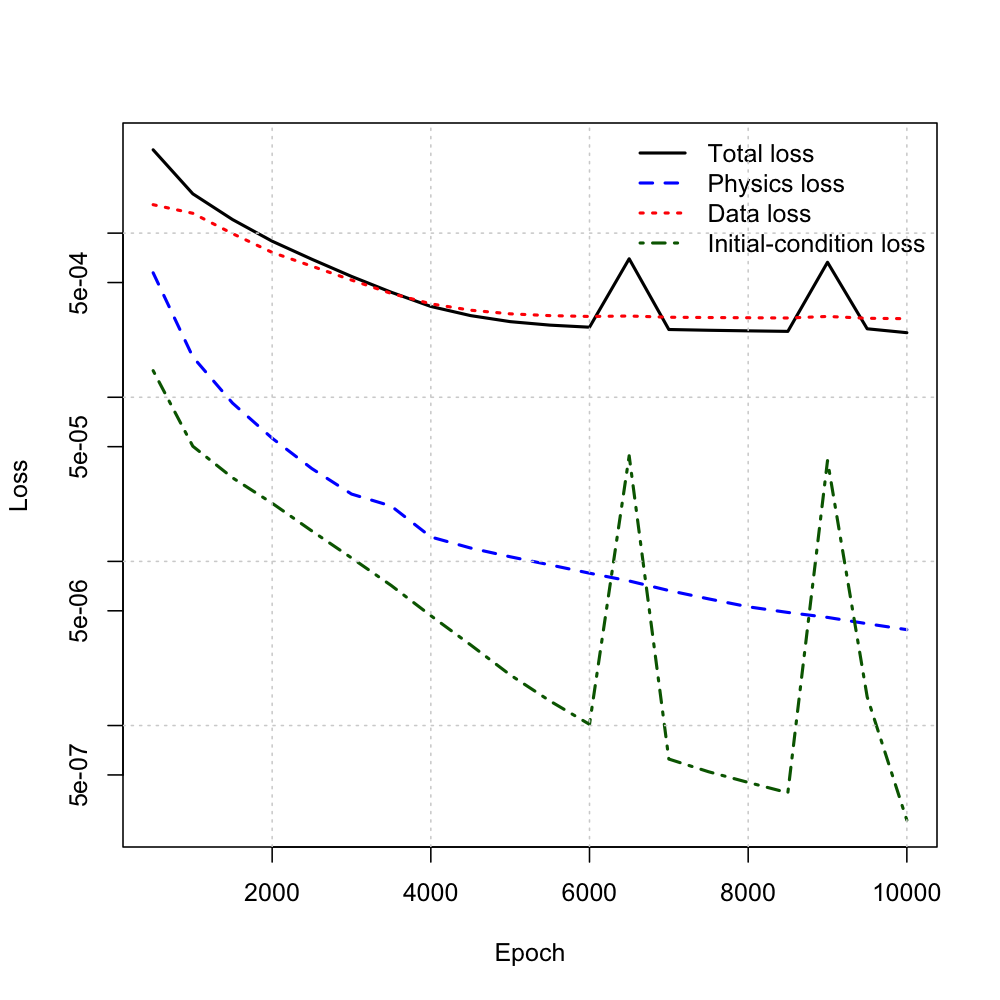}
		\caption{
			PINN-based fitting of model \eqref{ddemodel} to cumulative reported
			COVID-19 cases for Kenya. The left panel (a) shows the observed cumulative
			cases (open circles) and the PINN-predicted cumulative-case trajectory
			(solid line) after rescaling to the original population scale. The right
			panel (b) shows the total, physics-informed, data-fitting, and
			initial-condition losses over 10,000 training epochs. The PINN
			reproduces the overall cumulative epidemic trajectory while approximately 
			satisfying the residual constraints imposed by the governing model 
			equations.
		}
		\label{Kenya_fit}
	\end{figure}

	\medskip
	\noindent \textbf{Training and evaluation:}
	Training is carried out for 10{,}000 epochs using 236 observed weekly 
	cumulative-case observations and the Adam optimizer with a learning rate 
	of $10^{-3}$. The smallest total loss attained during training is 
	$
	\mathcal{L}_{\mathrm{total}}^{\mathrm{best}}
	=
	2.474008\times10^{-4}.
	$ 
	The goodness of fit is evaluated using the root-mean-square error (RMSE), 
	mean absolute error (MAE), normalized root-mean-square error (NRMSE), 
	and the coefficient of determination ($R^2$). For the Kenya data, the 
	calibrated model gives 
	\[
	\mathrm{RMSE}=5971.5250, \quad
	\mathrm{MAE}=3387.1916, \quad
	\mathrm{NRMSE}=0.023867,\quad
	\text{and} \quad
	R^2=0.997663.
	\]
	The normalized root-mean-square error is approximately $2.39\%$,
	indicating close descriptive agreement between the predicted and observed
	cumulative-case trajectory. Fig.~\ref{Kenya_fit} shows the observed and
	predicted cumulative cases together with the evolution of the training
	loss components. The PINN reproduces the overall cumulative epidemic
	trajectory while smoothing some short-term variations in the reported data.
	
	The calibrated parameter values yield $\mathcal R_0\approx0.515$ and $\mathcal R_0^{\max}\approx1.245$. Thus, although the positive-capacity disease-free equilibrium is locally asymptotically stable, the sufficient condition for disease elimination across both hospitalization regimes is not satisfied.
	
	The subsequent numerical bifurcation analysis explores broader parameter 
	regimes, and the corresponding parameter sets are used to investigate 
	the possible dynamical behavior of the model rather than being interpreted 
	as calibrated epidemiological estimates.
	
	\section{Numerical simulations}\label{NumSim} 
	We illustrate the analytical results by varying selected parameters in the smooth subsystems and the full Filippov system. All numerical computations and simulations in this section are performed in MATLAB.
	
	\subsection{Numerical simulations for subsystem $F_{\Sigma_1}$} 
\begin{figure}[H]
		\centering
		(a)\includegraphics[scale=0.3]{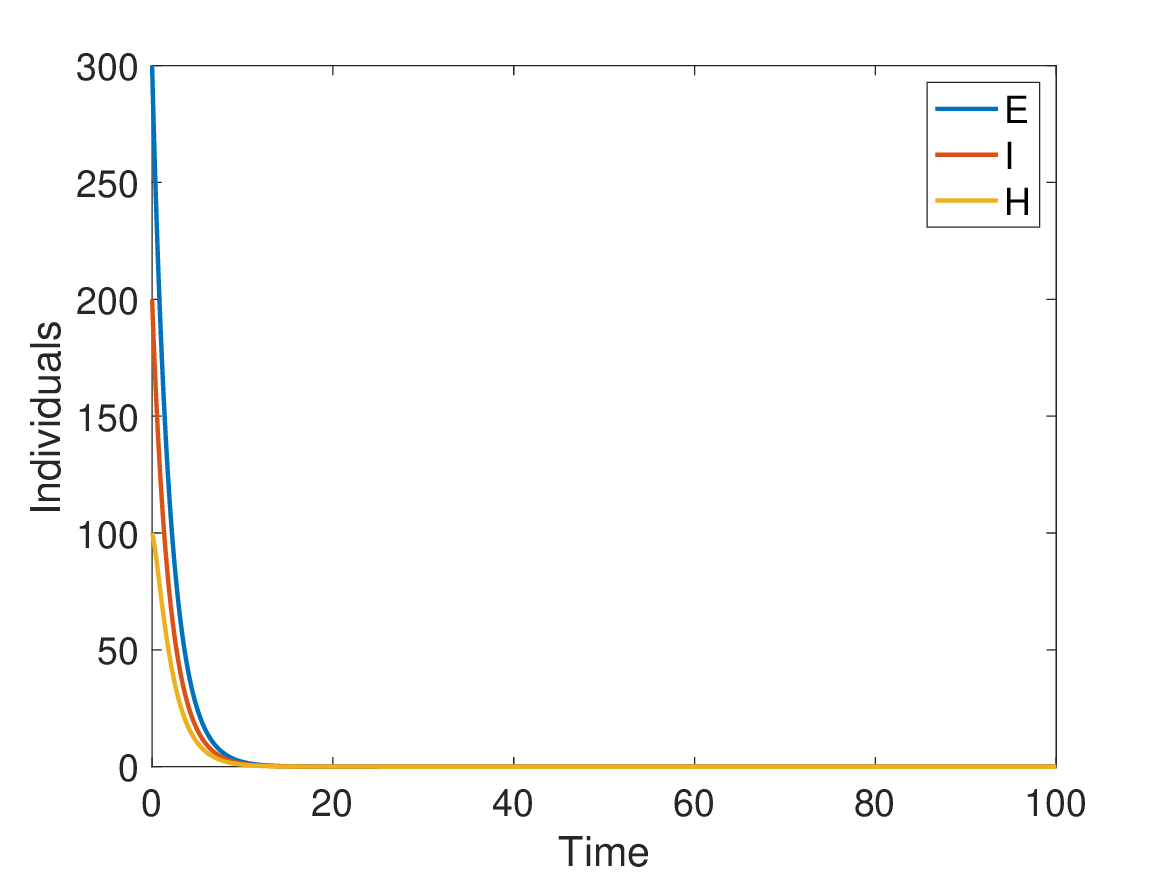} 
		(b)\includegraphics[scale=0.3]{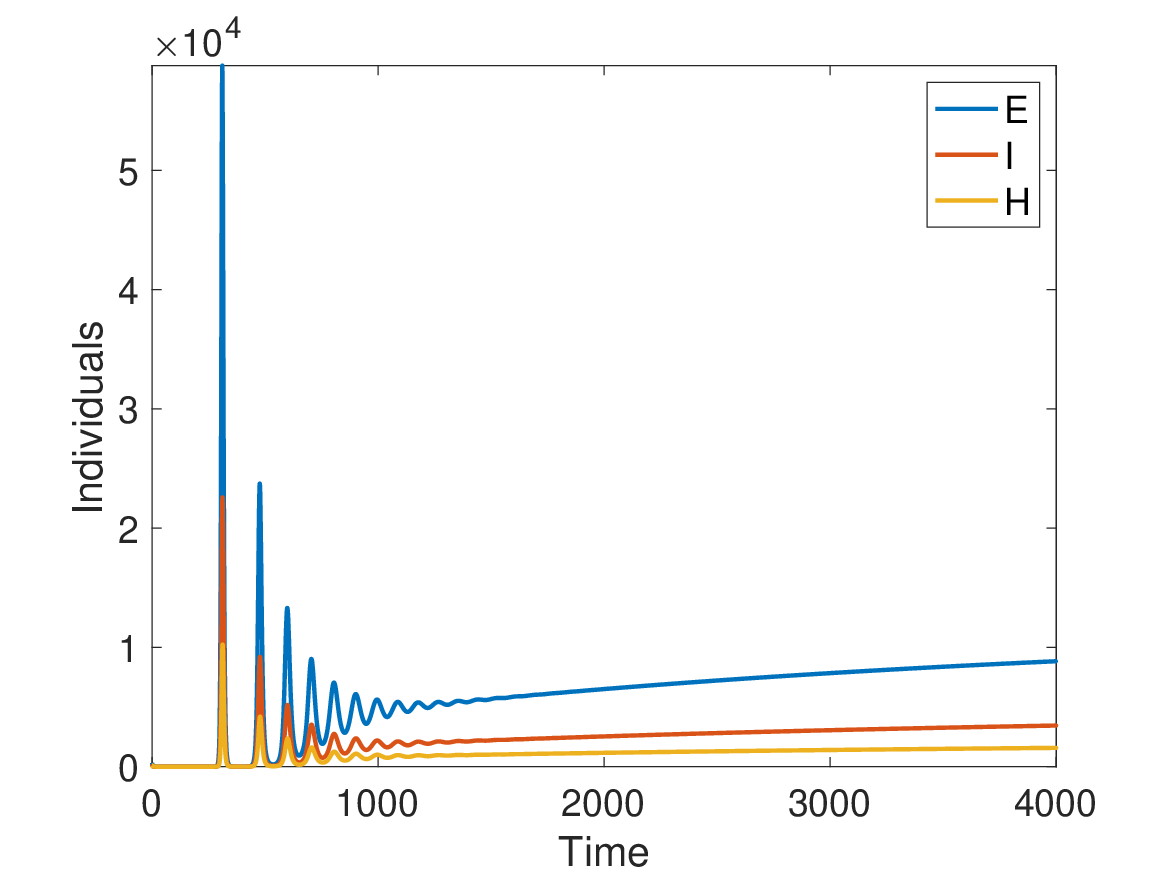} 
		\caption{Dynamics of the delay-free subsystem $F_{\Sigma_1}$. 
			(a) The disease-free equilibrium ($D_{02}$) is asymptotically stable for 
			$\beta = 0.00000021$ and $\mathcal R_0 = 0.82 < 1$. 
			(b) A unique endemic equilibrium ($D^*$) is asymptotically stable for 
			$\beta = 0.000007$ and $\mathcal R_0 = 27.28 > 1$.}
		\label{subsys1} 
	\end{figure}

	We numerically illustrate the analytical results for subsystem $F_{\Sigma_1}$, both in the absence and presence of time delay. Unless otherwise stated, the parameters are fixed as $\Lambda = 1400$ and $\xi = 0.001$, with all remaining values taken from  Table~\ref{param_val}. Under these settings, Fig.~\ref{subsys1} shows that the  disease-free equilibrium $D_{02}$ is locally asymptotically stable when $\mathcal R_0 < 1$, whereas a unique positive-capacity endemic equilibrium $D^*$ exists and is asymptotically stable for the illustrated parameter set with $\mathcal R_0 > 1$. This suggests that, below the hospitalization threshold, the treatment response is sufficient to preserve the classical threshold behavior, i.e. reducing $\mathcal R_0$ below unity eliminates infection, whereas values above unity lead to steady endemic persistence. 
	
	Additional simulations for positive delay are provided in Supplementary Section~S10.1.
	
	Fig.~\ref{fig:forwardbif} further illustrates the bifurcation structure of subsystem $F_{\Sigma_1}$ with respect to the basic reproduction number $\mathcal R_0$. For $\mathcal R_0 \le 1$, only disease-free equilibria exist and are represented by the zero-infected branch. As $\mathcal R_0$ crosses unity, the disease-free equilibrium loses stability and a unique endemic equilibrium emerges continuously and increases monotonically with $\mathcal R_0$. The monotone increase of the endemic infected population with $\mathcal R_0$ indicates that higher transmission potential leads to a larger persistent disease burden. This behavior is a characteristic of a forward (transcritical) bifurcation and confirms that, within subsystem $F_{\Sigma_1}$, reducing $\mathcal R_0$ below one is sufficient to ensure disease elimination. 
	
	\begin{figure}[H]
		\centering
		\includegraphics[scale=0.3]{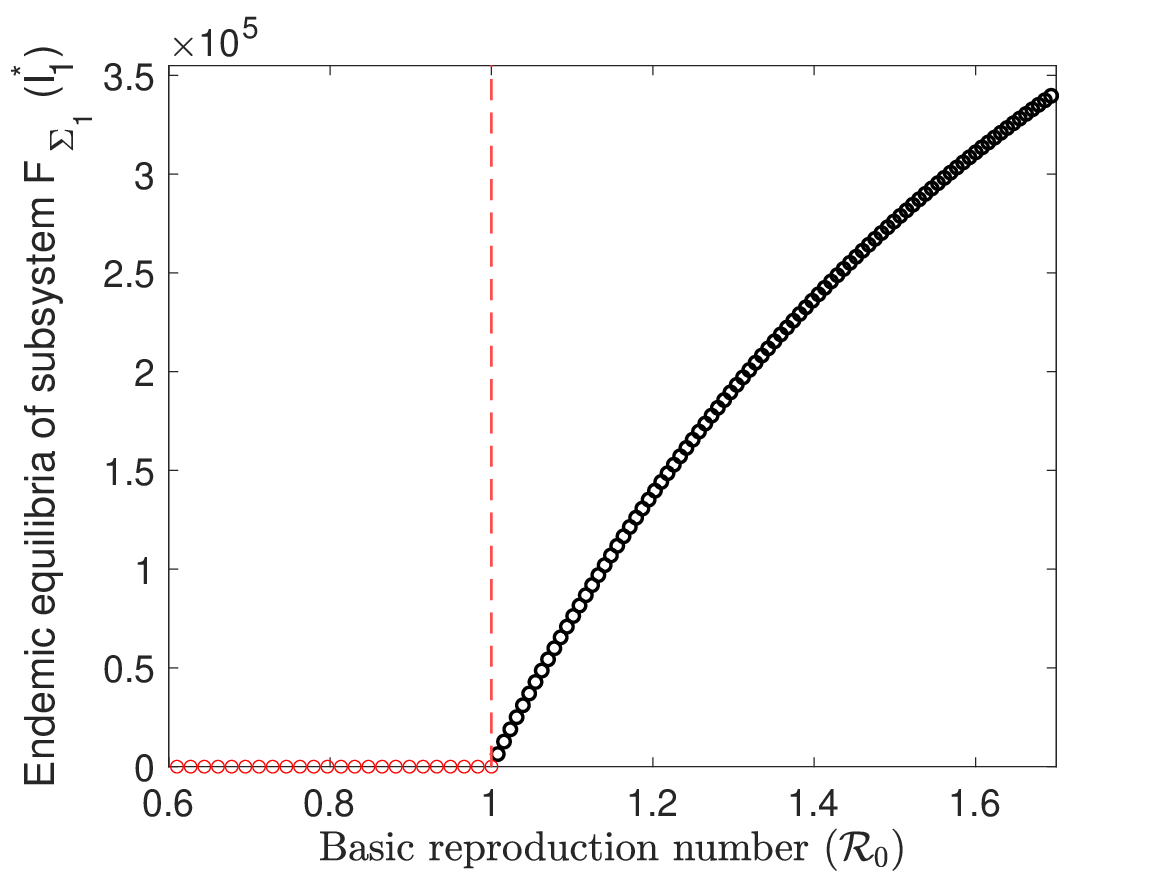}
		\caption{
			Equilibria of subsystem $F_{\Sigma_1}$ as a function of the basic reproduction
			number $\mathcal R_0$.
			Black circles indicate the endemic equilibrium for $\mathcal R_0>1$,
			while red circles denote the displayed disease-free branch $(I^*=0)$ for 
			$0\le \mathcal R_0\le 1$.
			The vertical red dashed line marks the threshold  $\mathcal R_0=1$.}
		\label{fig:forwardbif}
	\end{figure}
	
	\subsection{Numerical simulations for subsystem  $F_{\Sigma_2}$}
	
	The disease-free equilibrium simulations for the smooth extension are provided in Supplementary Section~S10.2.
	
\begin{figure}[H]
		\centering
		(a)\includegraphics[scale=0.3]{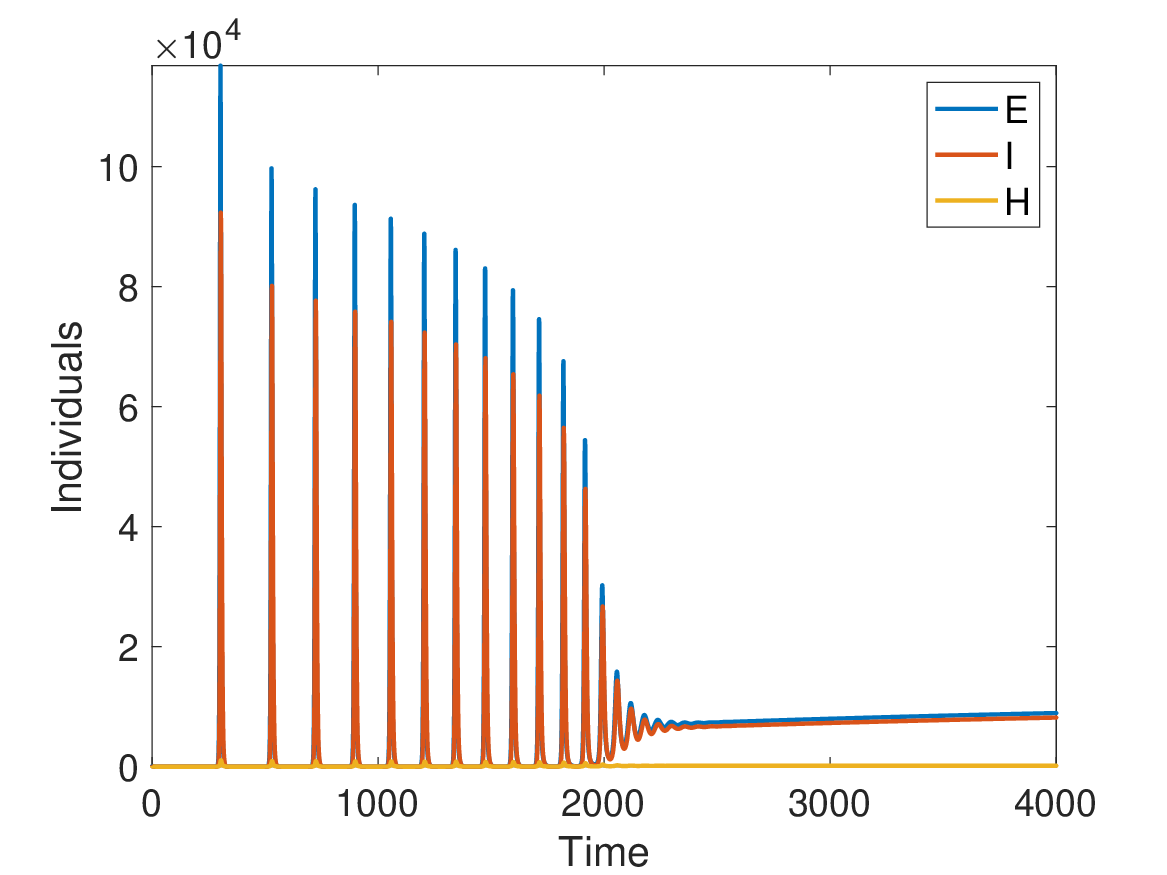}
		(b)\includegraphics[scale=0.3]{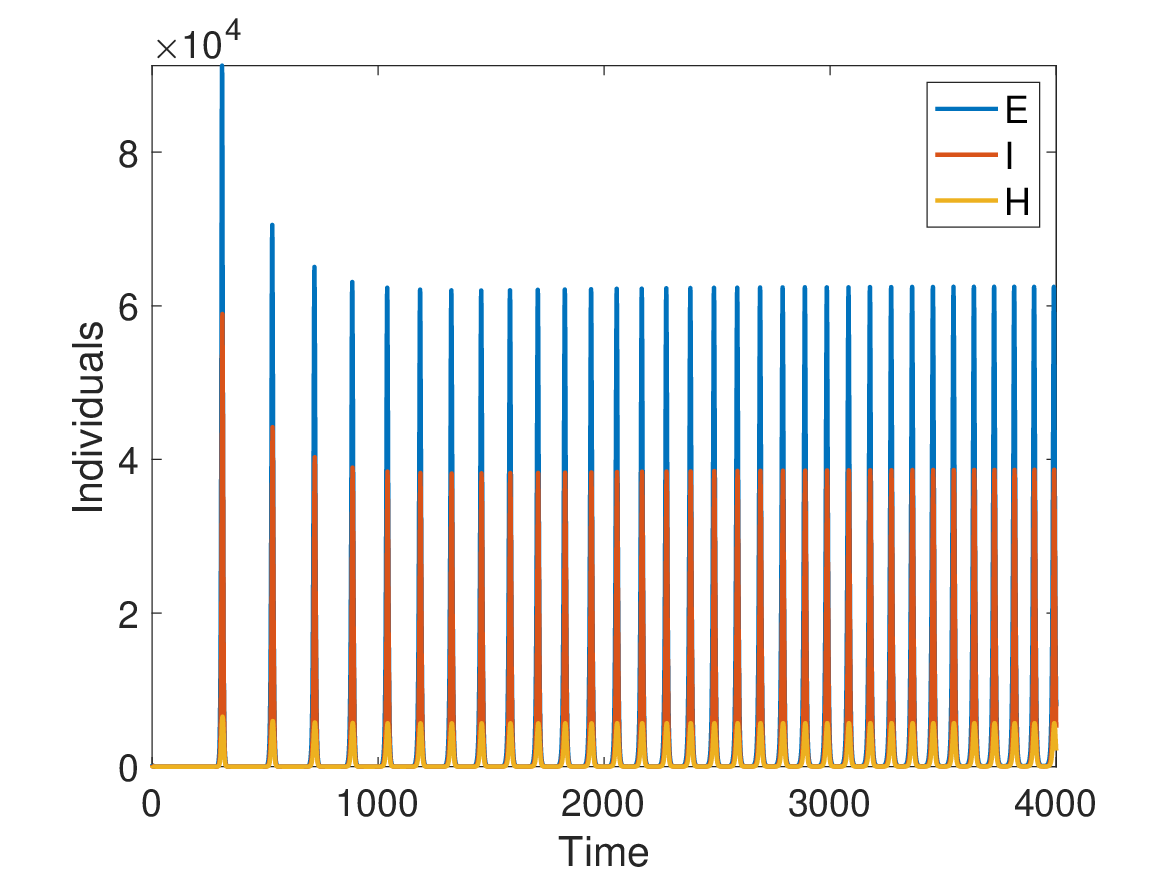}
		(c)\includegraphics[scale=0.3]{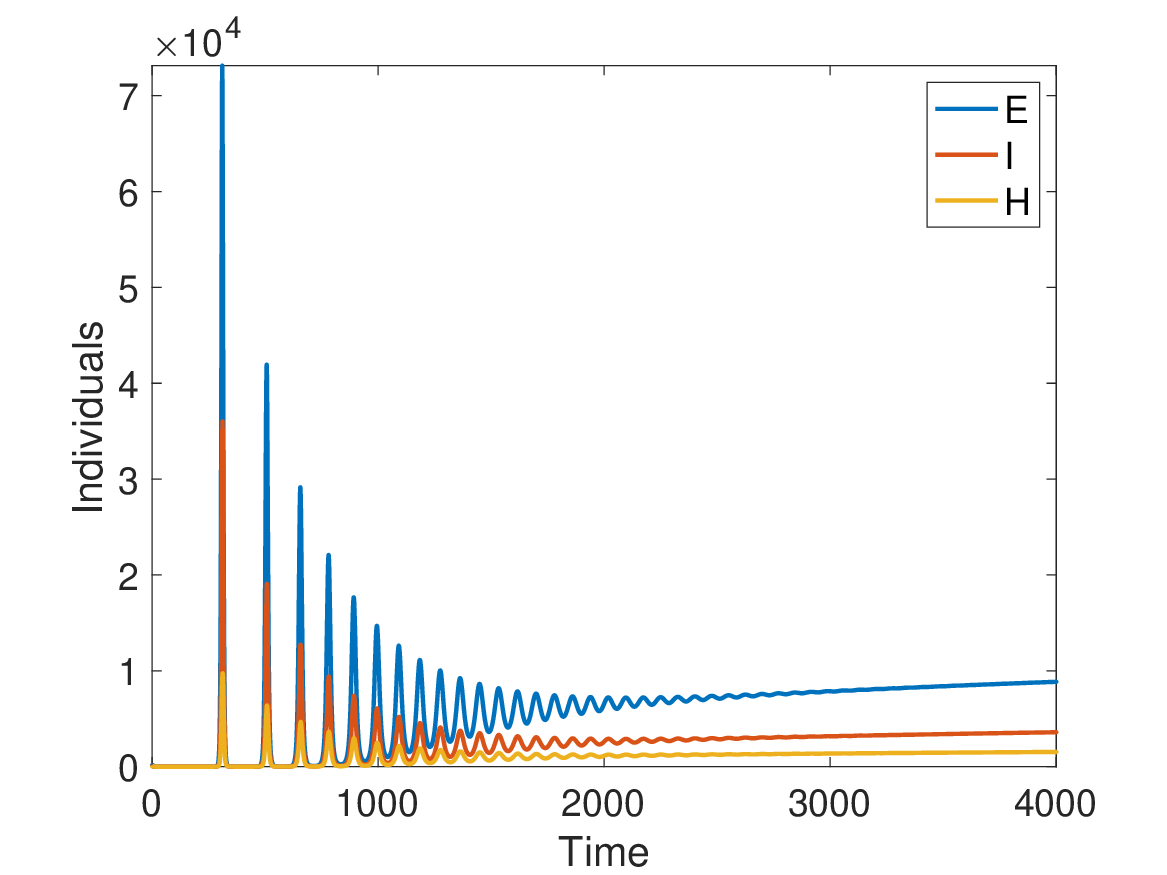}
		\caption{Dynamics of the non-delayed subsystem $F_{\Sigma_2}$ for (a) $K = 100$, (b) $K = 10{,}000$, and (c) $K = 30{,}000$, illustrating changes in stability of the endemic equilibrium when $\mathcal R_0 = 27.2820 > 1$.}\label{subsys2}
	\end{figure}

	Next, we examine the dynamics of the non-delayed subsystem
	$F_{\Sigma_2}$ for $\mathcal R_0>1$ by varying the baseline effective
	healthcare capacity $K$. Fig.~\ref{subsys2} illustrates a nonmonotone
	dependence of the endemic-equilibrium stability on $K$. For $K=100$,
	$D_1^*$ is locally asymptotically stable
	(Fig.~\ref{subsys2}(a)), whereas at $K=10000$, it is unstable
	and sustained oscillations emerge, consistent with a Hopf bifurcation
	(Fig.~\ref{subsys2}(b)). For larger values of $K$, the endemic
	equilibrium regains stability (Fig.~\ref{subsys2}(c)). Thus, increasing
	$K$ first destabilizes the endemic state and generates recurrent
	epidemic oscillations, while further increases can restore convergence
	to a steady endemic level. For these simulations, we fixed
	$\beta=7\times10^{-6}$, $\Lambda=1400$, and $\xi=0.001$, with the
	remaining parameters as in Table~\ref{param_val}. 
	
	\begin{figure}[H]
		\centering
		\includegraphics[scale=0.3]{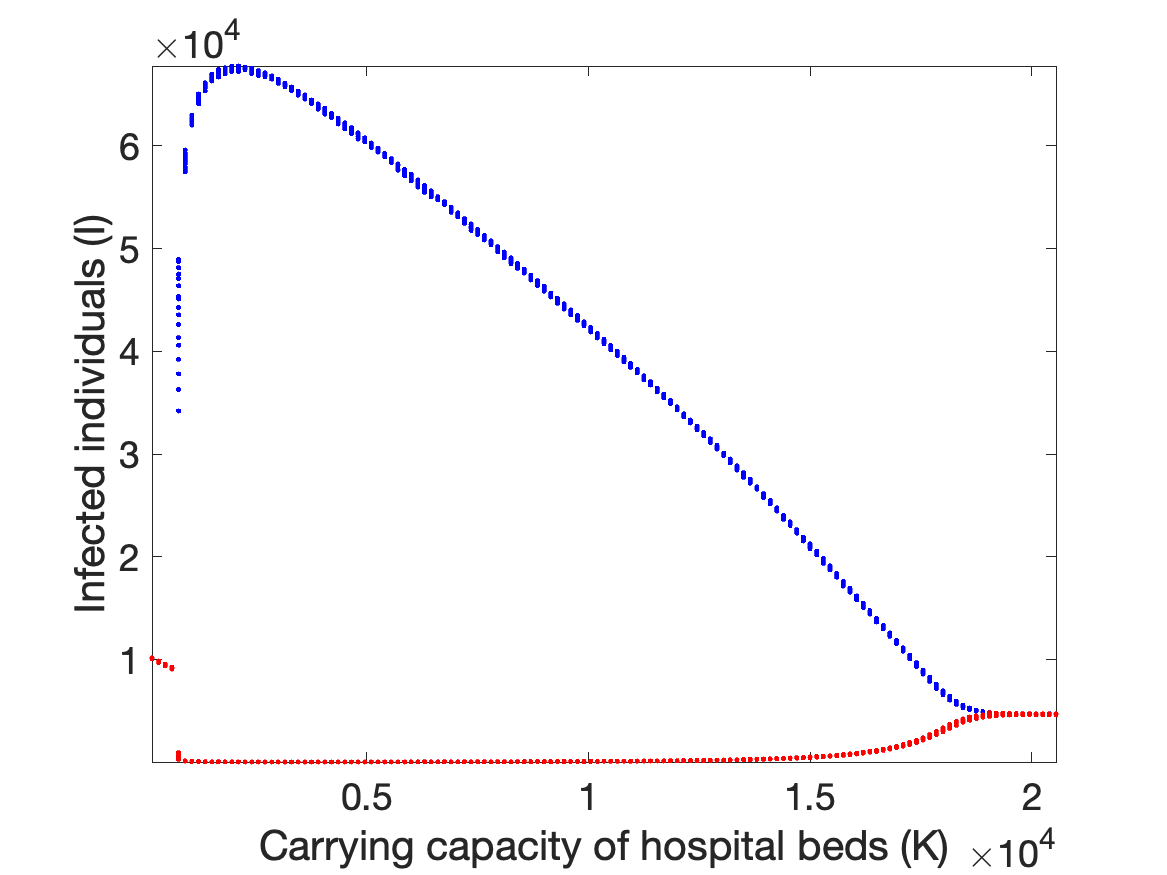} 
		\caption{Hopf bifurcation (endemic bubble) of the non-delayed subsystem $F_{\Sigma_2}$ with respect to $K$, showing the emergence of periodic solutions. The blue and red curves denote the upper and lower envelopes of the resulting limit cycle, respectively.}
		\label{bifdiag} 
	\end{figure}
	
	This nonmonotone behavior is summarized in the bifurcation diagram
	Fig.~\ref{bifdiag}, where $D_1^*$ is stable for small and large values
	of $K$ and unstable over an intermediate range supporting sustained
	oscillations. For each value of $K$, the system was integrated over a
	sufficiently long time interval, the initial transient was discarded,
	and the post-transient extrema of $I(t)$ were used to construct the
	upper and lower oscillatory envelopes. The parameter values are the
	same as those used in Fig.~\ref{subsys2}.
	
	A mechanistic interpretation of this nonmonotone dependence on $K$
	follows from the feedback between infection, healthcare capacity, and
	hospitalization. For $I>I_0$, variations in $B$ modify
	\[
	T(I,B)=\mu_0+(\mu_1-\mu_0)\frac{B}{B+I}.
	\]
	At intermediate values of $K$, this feedback can destabilize the
	endemic equilibrium, whereas for sufficiently large $K$,
	$T(I,B)$ approaches $\mu_1$ and becomes less sensitive to further
	changes in $B$, providing a qualitative explanation for the
	restabilization observed in Fig.~\ref{bifdiag}.
	
	\begin{figure}[H]
		\centering
		(a)\includegraphics[scale=0.3]{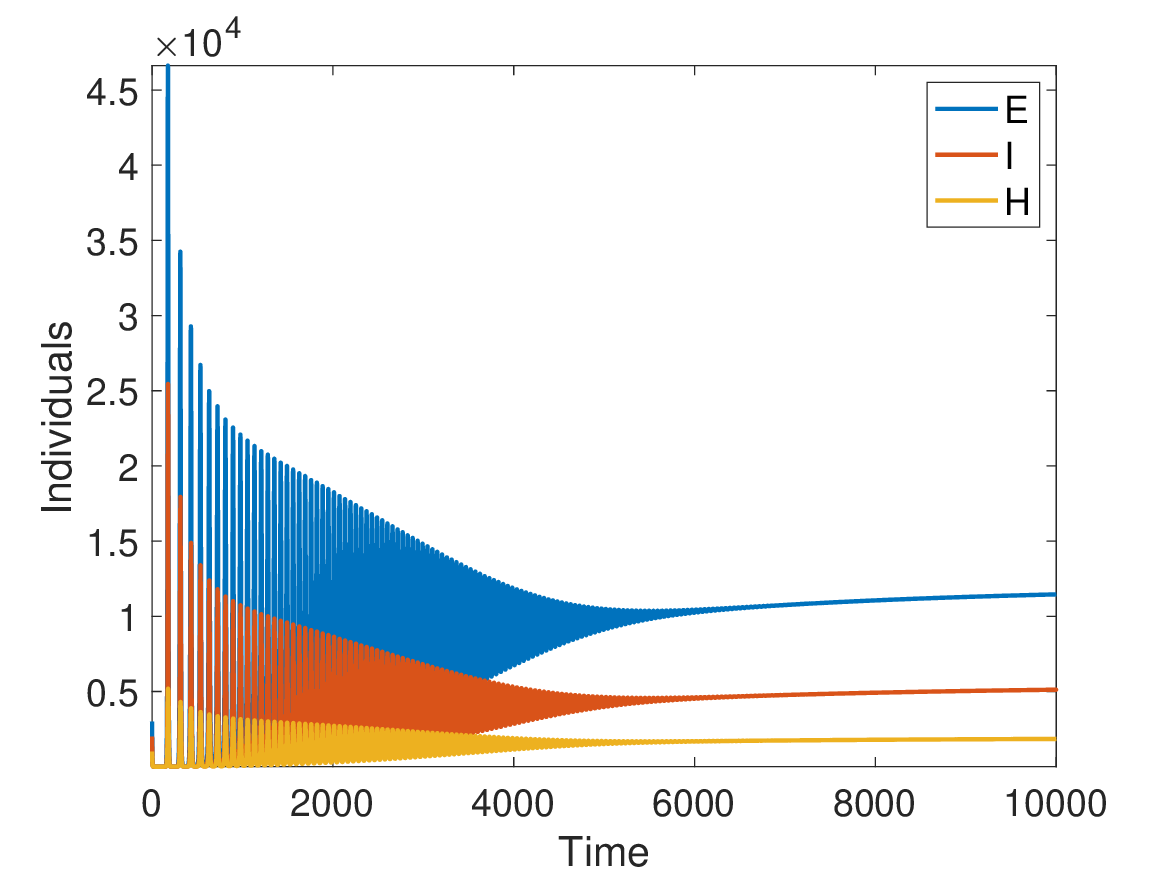}
		(b)\includegraphics[scale=0.3]{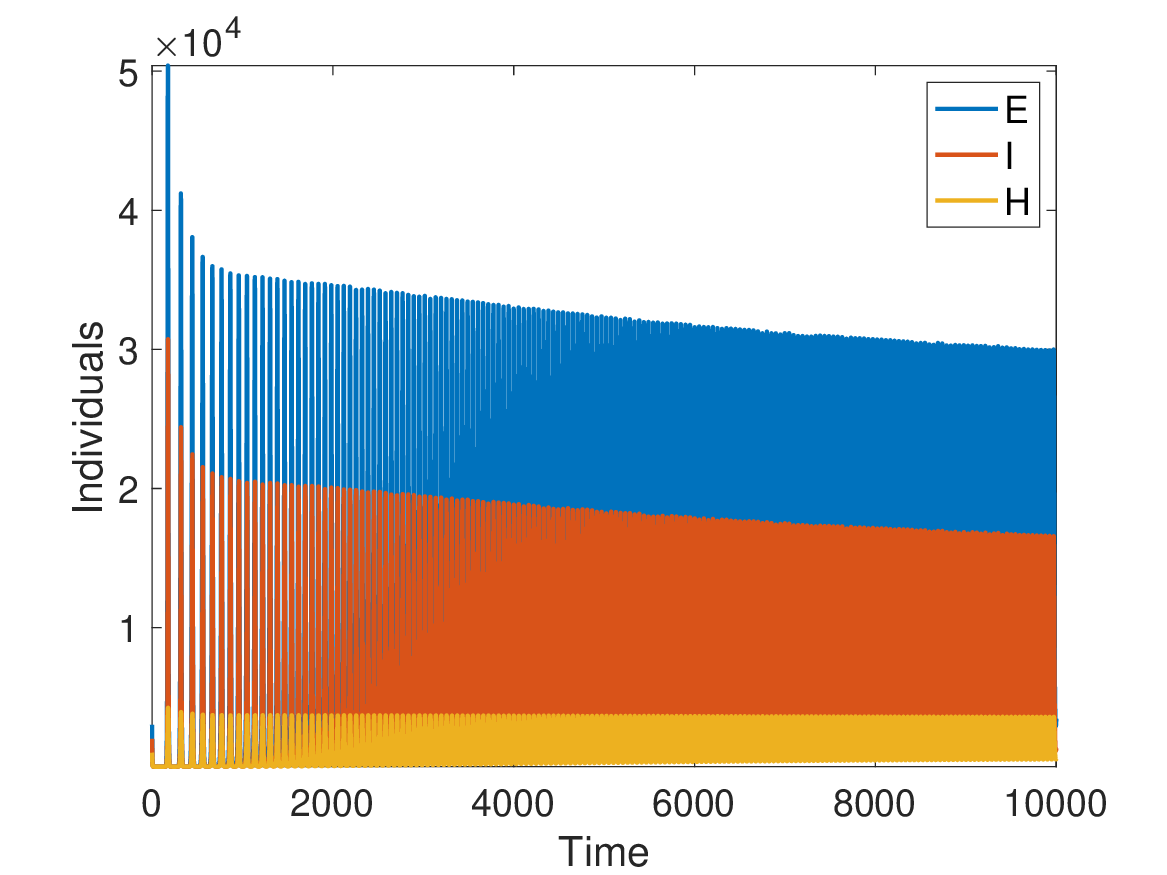}
		(c)\includegraphics[scale=0.3]{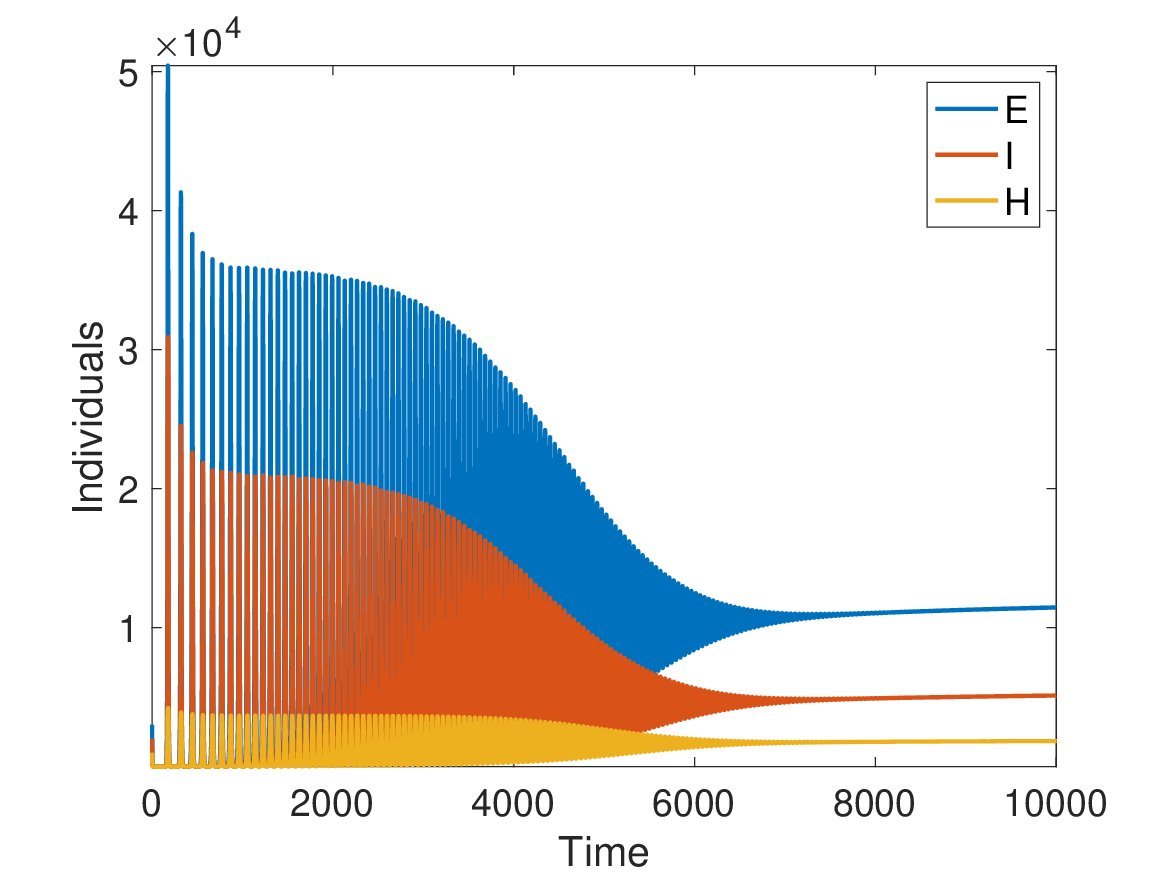}
		(d)\includegraphics[scale=0.3]{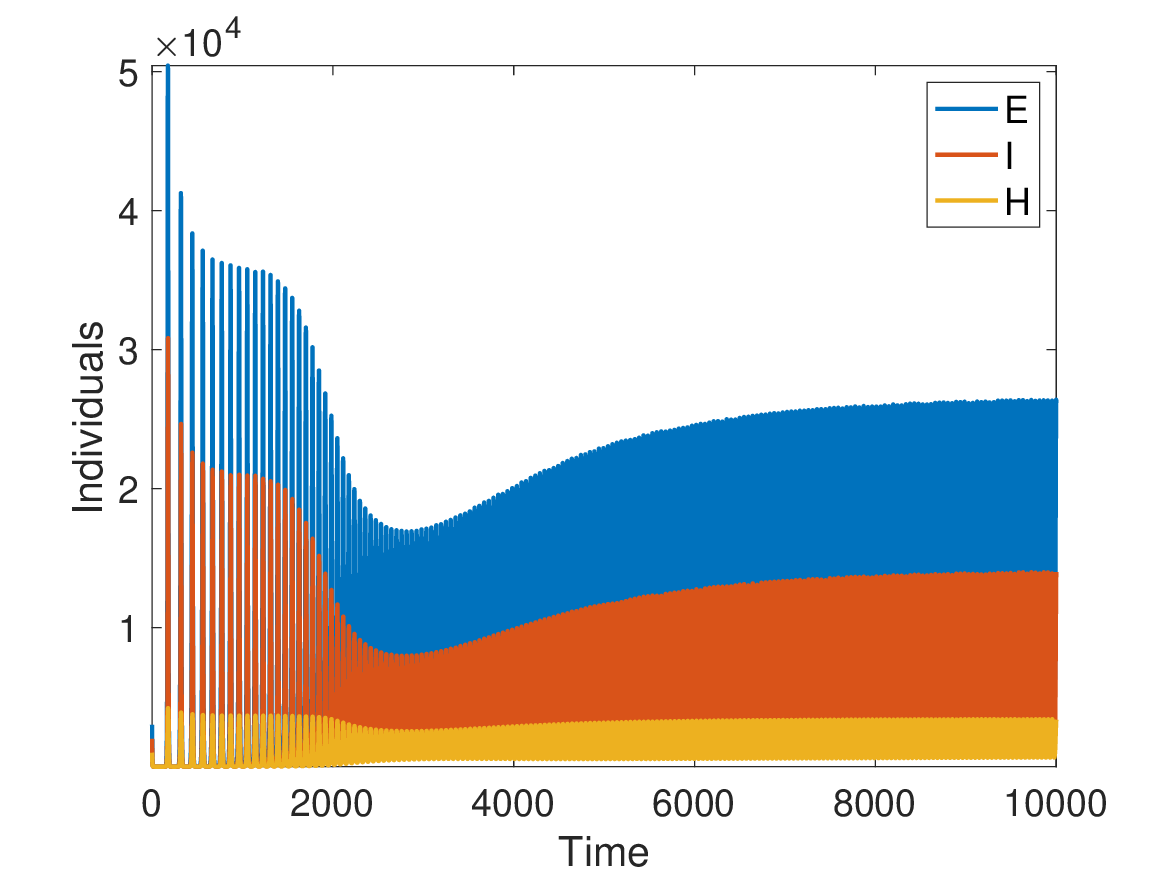}
		\caption{
			Dynamics of the delayed subsystem $F_{\Sigma_2}$ for
			(a) $\tau = 1$, (b) $\tau = 20$, (c) $\tau = 50$, and
			(d) $\tau = 60$, with $\beta=1.2\times10^{-5}$,
			$\Lambda=1400$, $\xi=0.001$, $K=11500$, $\mu_1=0.80$,
			$\omega=5000$, and $\lambda=0.93$. All other parameters are as in
			Table~\ref{param_val}, giving $\mathcal R_0\approx45.70$.
			The endemic equilibrium is stable for (a) and (c), whereas sustained
			oscillations are observed in (b) and (d).
		}
		\label{delaysub2}
	\end{figure}   
	
	We further investigate the effect of the time delay on subsystem
	$F_{\Sigma_2}$ when $\mathcal R_0>1$. For the parameter values used
	in Fig.~\ref{delaysub2}, we obtain $\mathcal R_0\approx45.70$,
	$I_1^*\approx5370.49$, and $B_1^*\approx15715.34$.
	The delayed subsystem is integrated using \texttt{dde23} over
	$t\in[0,10000]$ with constant initial history
	$(S,E,I,H,R,B)=(4000,3000,2000,1000,500,5000)$.
	In the absence of delay, all eigenvalues of the corresponding
	Jacobian matrix have negative real parts, with spectral abscissa
	approximately $-2.94275\times10^{-4}$, implying local asymptotic
	stability of the endemic equilibrium.
	
	Writing the characteristic equation as
	$P_c(s)+Q_c(s)e^{-s\tau}=0$, we compute the positive roots of
	$|P_c(i\sqrt\zeta)|^2-|Q_c(i\sqrt\zeta)|^2=0$ and evaluate the
	phase condition. The sign of the derivative of this polynomial
	at a simple positive root determines the crossing direction.
	Table~\ref{tab:critical_delay_check} lists the first three members
	of each critical-delay family. As the delay increases, the endemic
	equilibrium loses stability at $\tau\approx9.6833$, regains
	stability at $\tau\approx30.2884$, and subsequently loses stability
	again at $\tau\approx52.6685$. These analytical predictions agree
	with Fig.~\ref{delaysub2}, where convergence to the endemic
	equilibrium is observed for $\tau=1$ and $\tau=50$, whereas
	sustained oscillations are observed for $\tau=20$ and $\tau=60$.
	Thus, the healthcare-capacity response delay induces repeated
	local stability switching of the endemic equilibrium.
	
	Denote the destabilizing and stabilizing families by
	$\tau_+^{(k)}$ and $\tau_-^{(k)}$, respectively. Their spacings
	are approximately $42.9852$ and $45.1617$ weeks. The equilibrium
	is locally asymptotically stable for $0\leq\tau<\tau_+^{(0)}$
	and on $(\tau_-^{(k)},\tau_+^{(k+1)})$ for $k=0,\ldots,10$.
	The last such interval is approximately $(481.91,482.52)$ weeks.
	Beyond its upper endpoint, the number of destabilizing crossings
	exceeds that of stabilizing crossings, and the equilibrium remains
	unstable away from the critical delays.
	\begin{table}[H]
		\centering
		\caption{Critical frequencies, critical delays, and crossing directions for the parameter set in Figs.~\ref{delaysub2} and \ref{delaybif}. A positive (negative) sign denotes destabilizing (stabilizing) crossing.}
		\label{tab:critical_delay_check}
		\begin{tabular}{ccrrrr}
			\hline
			$\zeta_j$ & $\nu_j$ & $\tau_j^{(0)}$ & $\tau_j^{(1)}$ & $\tau_j^{(2)}$ & Sign\\
			\hline
			0.01935617 & 0.13912647 & 30.2884 & 75.4501 & 120.6118 & $-$\\
			0.02136592 & 0.14617085 & 9.6833 & 52.6685 & 95.6537 & $+$\\
			\hline
		\end{tabular}
	\end{table}

	\begin{figure}[H]
		\centering
		\includegraphics[scale=0.3]{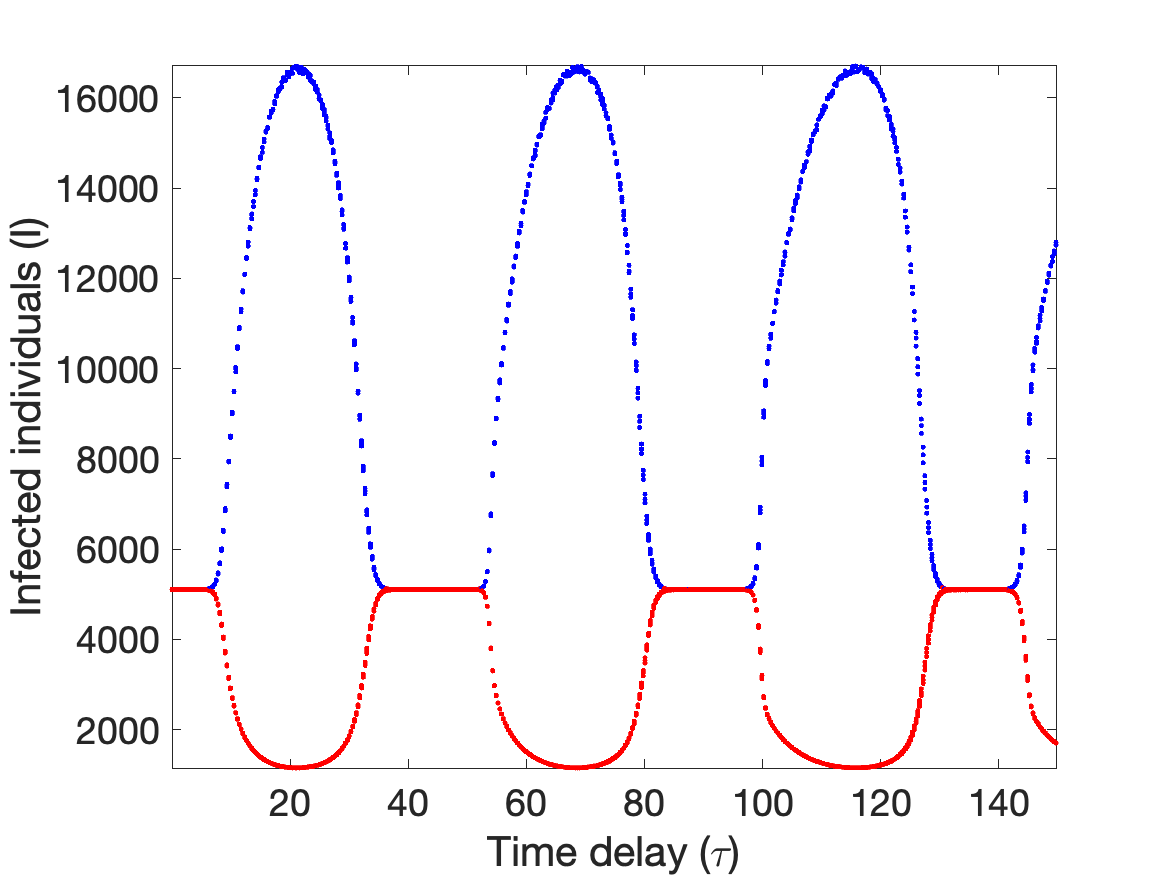}
		\caption{
			Bifurcation diagram illustrating multiple oscillatory endemic
			``bubbles'' in the delayed subsystem $F_{\Sigma_2}$ as the
			capacity-response delay $\tau$ varies. The blue and red curves
			represent, respectively, the local maxima and minima of $I(t)$
			along the post-transient trajectories. The envelopes illustrate changes in post-transient dynamics. The parameter values are
			$\beta=1.2\times10^{-5}$, $\Lambda=1400$, $\xi=0.001$,
			$K=11500$, $\mu_1=0.80$, $\omega=5000$, and $\lambda=0.93$, with all other
			parameters as in Table~\ref{param_val}.
		}
		\label{delaybif}
	\end{figure}
	
	To construct Fig.~\ref{delaybif}, we integrate the delayed system for each value of $\tau$, discard the initial transient, and record the post-transient local extrema of $I(t)$. The alternating steady and oscillatory windows are consistent with the critical delays and crossing directions in Table~\ref{tab:critical_delay_check}.
	
	\begin{figure}[H] 
		\centering
		(a) \includegraphics[scale=0.3]{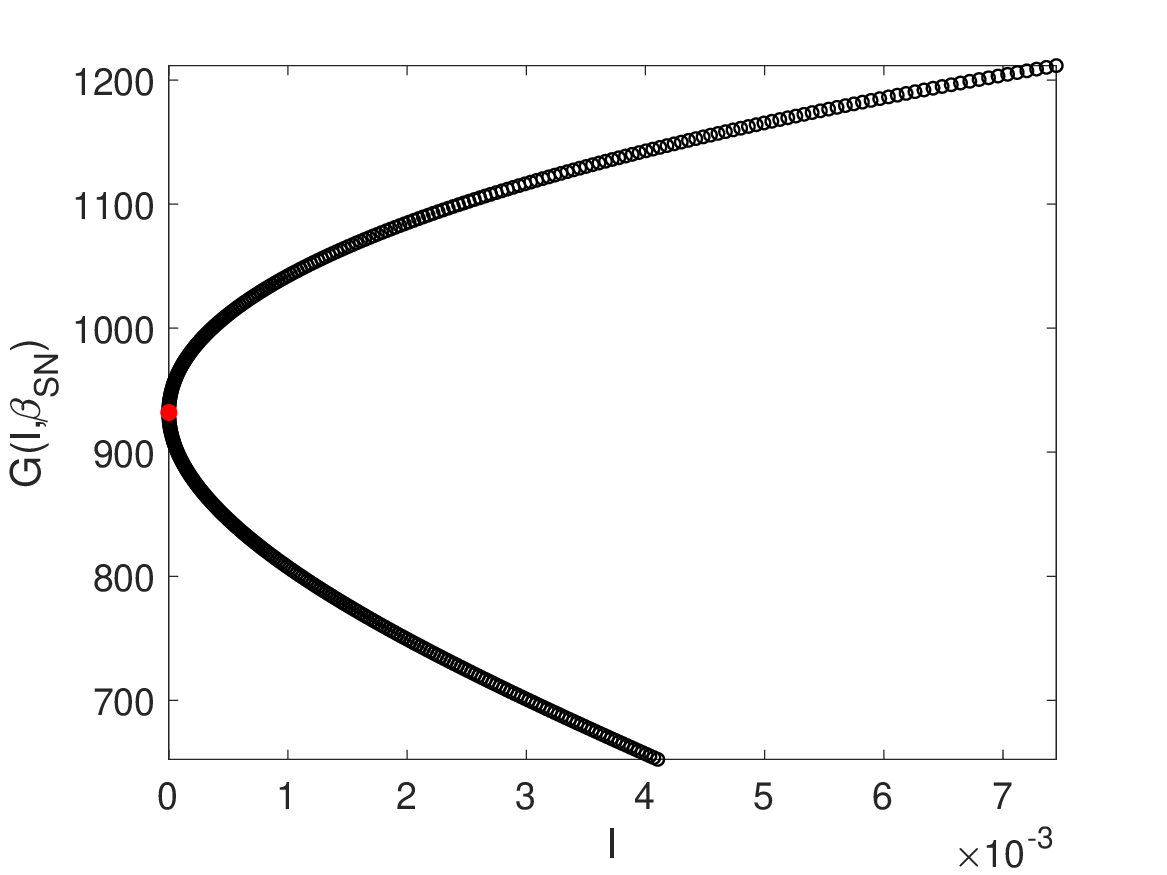} 
		(b) \includegraphics[scale=0.3]{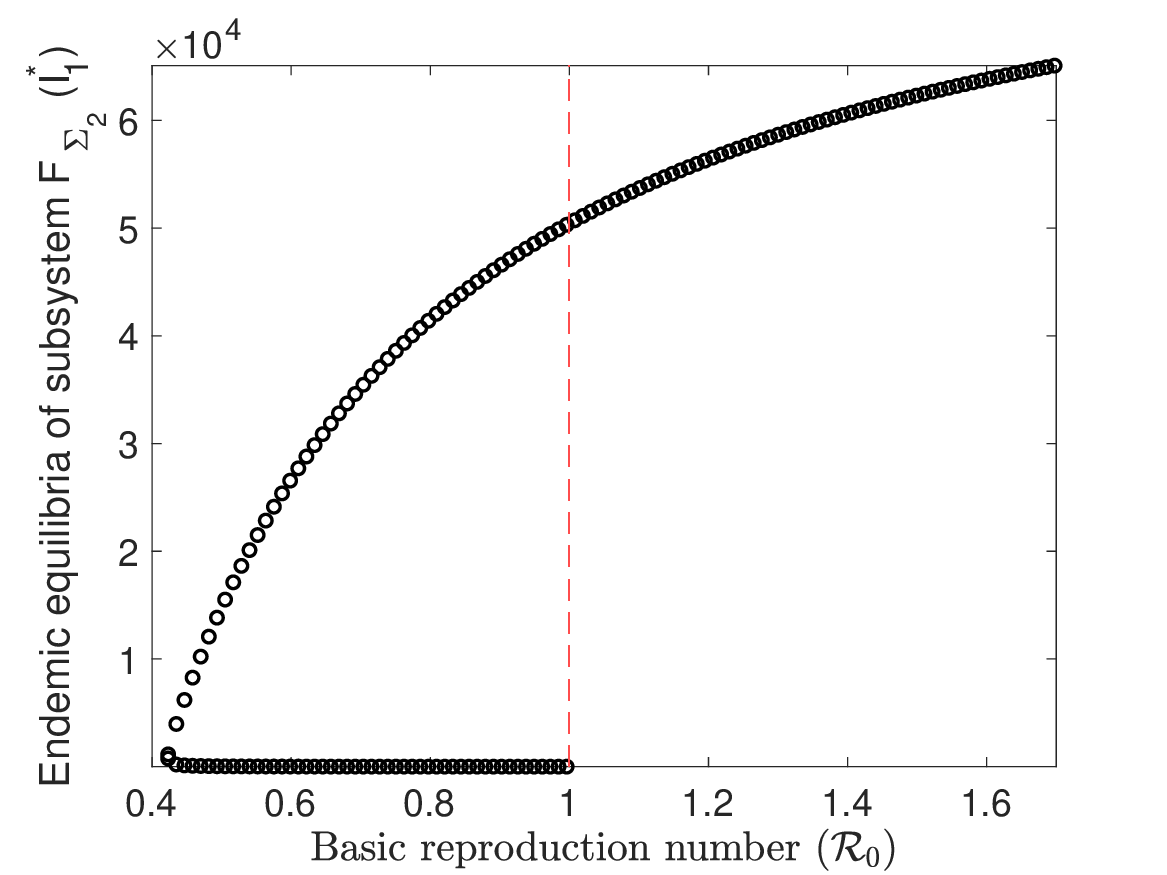} 
		\caption{
			(a) Local geometry of the algebraic equilibrium equation
			$G(I,\beta_{\rm SN})$ in the $(G,I)$-plane near the nondegenerate
			fold point (red dot) of the endemic-equilibrium branch of the smooth
			extension of subsystem $F_{\Sigma_2}$.
			(b) Numerical equilibrium diagram of the positive endemic branches
			of the smooth extension of $F_{\Sigma_2}$ with respect to
			$\mathcal R_0$. The fold occurs for $\mathcal R_0<1$, producing
			two positive endemic-equilibrium branches over an interval below
			the classical disease-invasion threshold. 
		} 
		\label{backbif}
	\end{figure} 
	
	Fig.~\ref{backbif}(a) illustrates the local geometry of the algebraic
	equilibrium equation \eqref{eq:G}, $G(I,\beta)=0$, near the fold point
	$(I_{\rm SN},\beta_{\rm SN})$ of the smooth extension of 
	$F_{\Sigma_2}$. The positive roots are computed numerically as $\beta$
	is varied, and the stability of the corresponding equilibrium branches
	is classified from the eigenvalues of the associated Jacobian matrix.
	For this illustration, we choose $\Lambda=1400$, $K=5$, $\lambda=0.7$,
	and the remaining parameters as in Table~\ref{param_val}, while varying
	$\beta$ over
	$ 
	1.0857\times10^{-7}<\beta<4.5\times10^{-7}.
	$ 
	Solving the equilibrium fold conditions yields
	\[
	(I_{\rm SN},\beta_{\rm SN})
	=
	(932.04,1.085214\times10^{-7}),
	\]
	at which
	$
	G(I_{\rm SN},\beta_{\rm SN})=0, 
	\partial_I G(I_{\rm SN},\beta_{\rm SN})=0,
	$
	with
	$
	\partial_{II}G(I_{\rm SN},\beta_{\rm SN})\neq0, 
	\partial_\beta G(I_{\rm SN},\beta_{\rm SN})\neq0.
	$
	Hence, the endemic-equilibrium equation undergoes a nondegenerate fold.
	Moreover,
	$ 
	\mathcal R_0(\beta_{\rm SN})\approx0.423<1,
	$ 
	so the fold occurs below the disease-invasion threshold
	$\mathcal R_0=1$. On one side of the fold, two positive branches are
	present, one locally stable and the other unstable.
	
	Fig.~\ref{backbif}(b) shows the corresponding equilibrium diagram in
	terms of $\mathcal R_0$. For this parameter set, no positive-capacity endemic
	equilibrium of the smooth extension is observed below the fold, two are observed for
	$\mathcal R_0(\beta_{\rm SN})<\mathcal R_0<1$, and a unique positive endemic equilibrium is
	observed for $\mathcal R_0>1$. These branches represent equilibria of
	the full Filippov system only when they are biologically feasible and
	satisfy the admissibility condition
	$ 
	I_1^*>I_0.
	$ 
	
	The numerical continuation therefore confirms a subthreshold fold; when the locally stable branch is admissible, it produces discontinuity-induced subthreshold bistability with $D_{02}$.
	
	For example, at $\beta=1.2\times10^{-7}$, with the other parameters
	as specified for Fig.~\ref{backbif}, we obtain
	$\mathcal R_0\approx0.46769$ and a positive-capacity endemic equilibrium
	with $I_1^*\approx9856.22>I_0=4500$. In the delay-free setting, the
	largest real part of the eigenvalues of the corresponding Jacobian is
	approximately $-1.19186\times10^{-3}<0$. Hence, this admissible locally
	asymptotically stable endemic equilibrium coexists with the locally
	asymptotically stable disease-free equilibrium $D_{02}$, providing a
	numerical example of discontinuity-induced subthreshold bistability in
	the full Filippov system.
	
	For this parameter set, the worst-case reproduction number at the fold
	satisfies $\mathcal R_0^{\max}(\beta_{\rm SN})\approx1.0225$, while at
	the illustrated bistability point $\beta=1.2\times10^{-7}$, we obtain
	$\mathcal R_0^{\max}\approx1.1306$. Thus, the observed bistability
	occurs in the intermediate regime $\mathcal R_0<1<\mathcal R_0^{\max}$,
	where local asymptotic stability of the disease-free equilibrium does
	not guarantee disease elimination across both hospitalization regimes.
	This observation is consistent with Theorem~\ref{thm:global_dfe_worstcase}.
	
	\subsection{Numerical simulations for the full Filippov system \eqref{ddemodel}} 
	We now present numerical simulations for the full Filippov system
	\eqref{ddemodel} in the absence of time delay, focusing on the repelling
	sliding region, pseudo-equilibria, boundary-equilibrium transitions, and
	crossing periodic dynamics near quadratic tangencies. We begin with the
	geometry of the repelling sliding region on the switching manifold
	$\Sigma=\{I=I_0\}$.
	
	\begin{figure}[H]
		\centering
		(a)\includegraphics[scale=0.3]{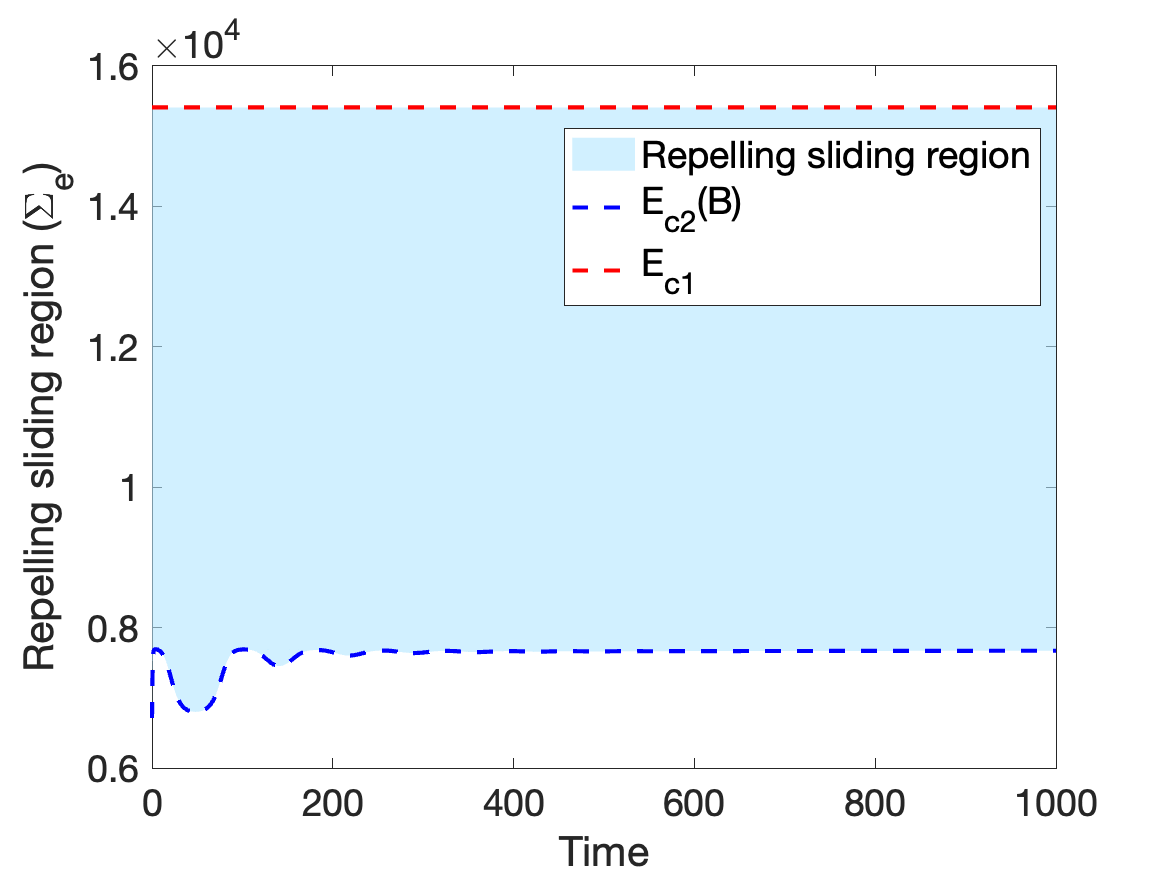}
		(b)	\includegraphics[scale=0.3]{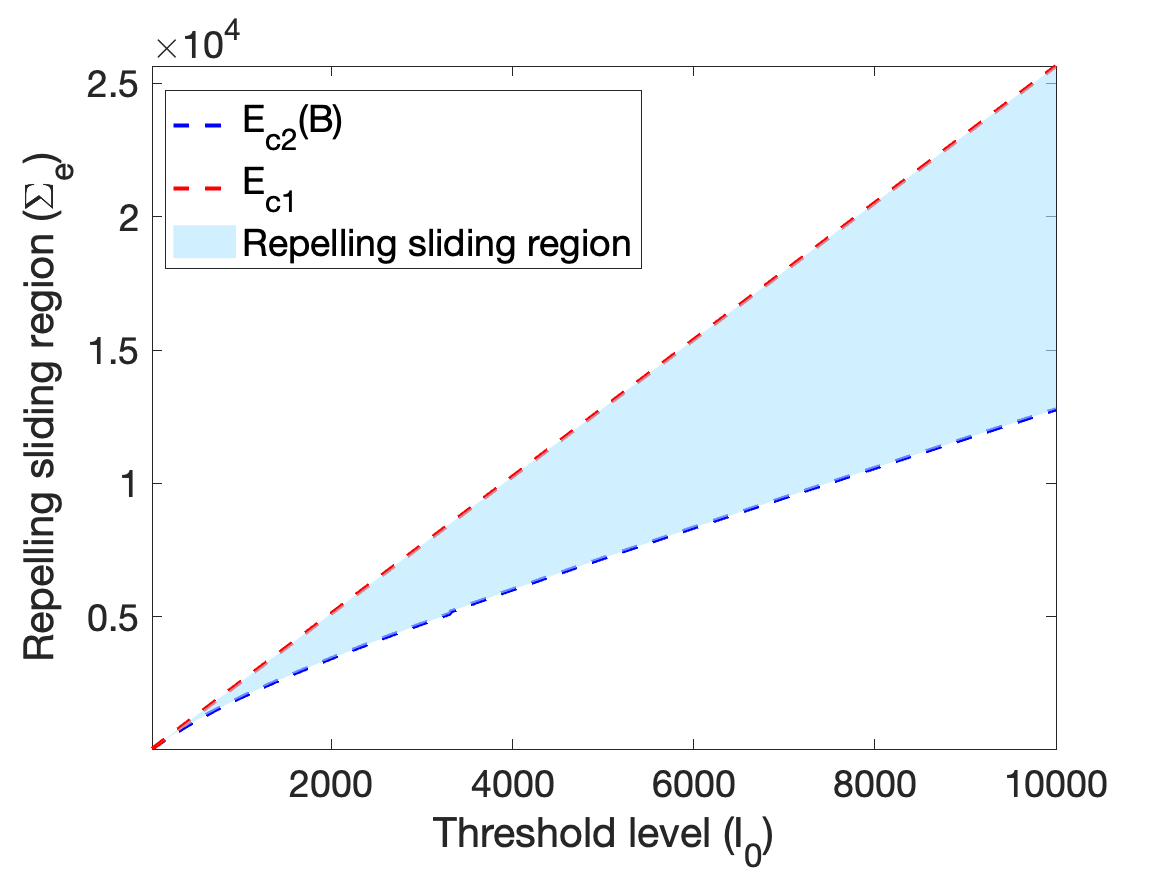}
		\caption{
			Repelling sliding region on the switching manifold $I=I_0$.
			(a) Time evolution of the state-dependent lower boundary $E_{c2}(B(t))$ (blue)
			and the upper boundary $E_{c1}$ (red dashed), which together delimit the
			repelling sliding interval on $\Sigma$.
			(b) Dependence of the repelling sliding region on the threshold level $I_0$,
			with the lower boundary $E_{c2}$ evaluated at the equilibrium healthcare capacity $B^s$.} 
		\label{fig:repelling_sliding_region}
	\end{figure}  
	
	Fig.~\ref{fig:repelling_sliding_region}(a) shows the evolution of the
	sliding boundaries. The lower boundary $E_{c2}(B(t))$ varies with the
	evolving healthcare capacity $B(t)$, whereas the upper boundary
	$E_{c1}$ is constant. Consequently, although the switching manifold is
	fixed, the width of the repelling sliding region changes during transient
	capacity dynamics and approaches a limiting configuration as $B(t)$
	converges. Thus, healthcare capacity influences the range of exposed
	individuals for which repelling sliding can occur near the hospitalization
	threshold. Fig.~\ref{fig:repelling_sliding_region}(b) illustrates the dependence of
	the sliding region on the threshold $I_0$. Both boundaries increase as
	$I_0$ increases, shifting the repelling sliding interval in the
	$E$-direction. Hence, a larger infection threshold requires a larger
	exposed population for sliding behavior to occur near $\Sigma$. For these simulations, we use
	$\beta=0.000021$, $\Lambda=1400$, $\mu=0.00028$, 
	$\lambda=0.93$, and $\mu_1=0.77$, with the remaining parameters as in
	Table~\ref{param_val}. The baseline value is $I_0=4500$, except in
	Fig.~\ref{fig:repelling_sliding_region}(b), where $I_0$ is varied. 
	
	\begin{figure}[H]
		\centering
		\includegraphics[scale=0.3]{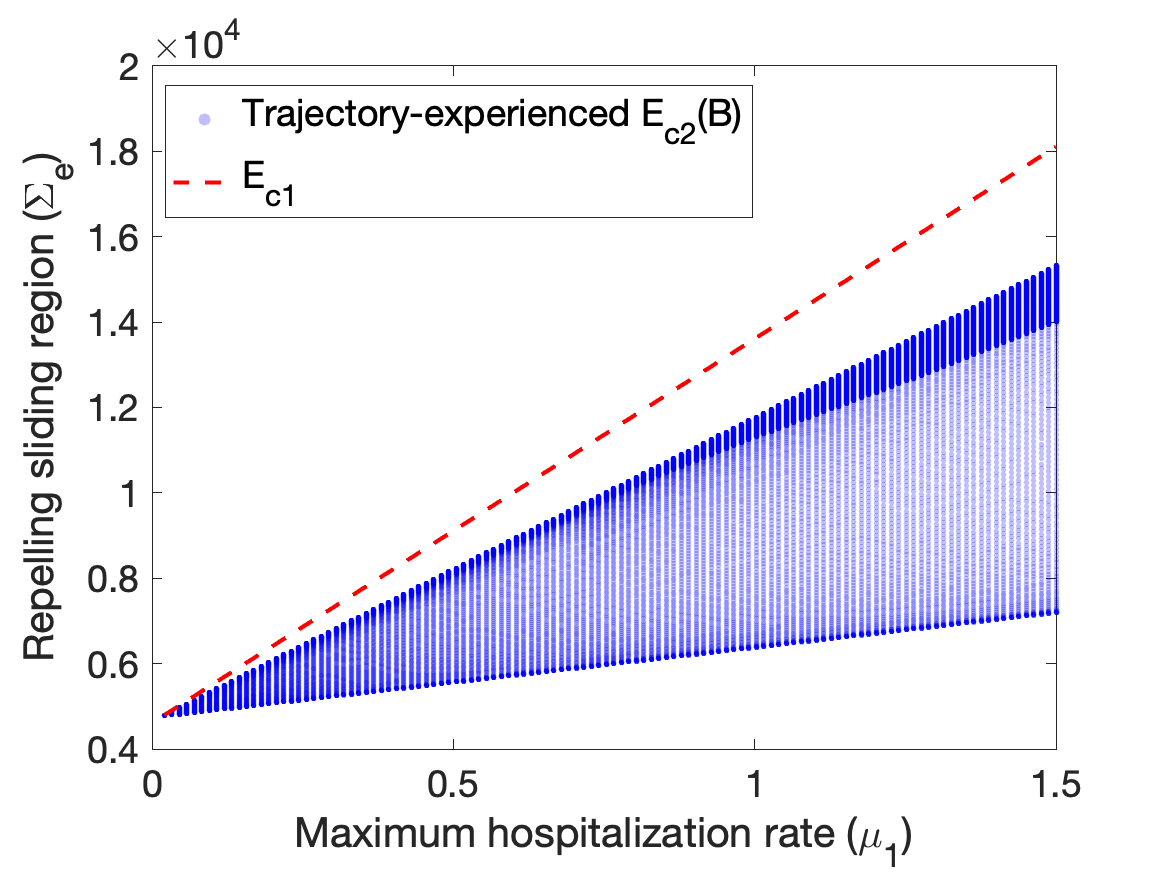} 
		\caption{
			Repelling sliding region on the switching manifold $I=I_0$ as a function of the
			maximum hospitalization rate $\mu_1$.
			Blue scatter dots denote trajectory-experienced values of the lower sliding boundary
			$E_{c2}(B(t))$ obtained from Filippov simulations; the solid blue curve shows
			$E_{c2}$ evaluated at the equilibrium healthcare  capacity $B^s$, and the red dashed
			curve represents the upper tangency boundary $E_{c1}(\mu_1)$.
		}
		\label{fig:repelling_sliding_region_mu1}
	\end{figure}  
	
	Fig.~\ref{fig:repelling_sliding_region_mu1} shows that increasing $\mu_1$ shifts and widens the repelling sliding region. The spread of the blue scatter points reflects the variation in $B(t)$ along simulated trajectories, and hence the range of lower-boundary values encountered near $I=I_0$.
	
	\begin{figure}[H]
		\centering
		\includegraphics[scale=0.3]{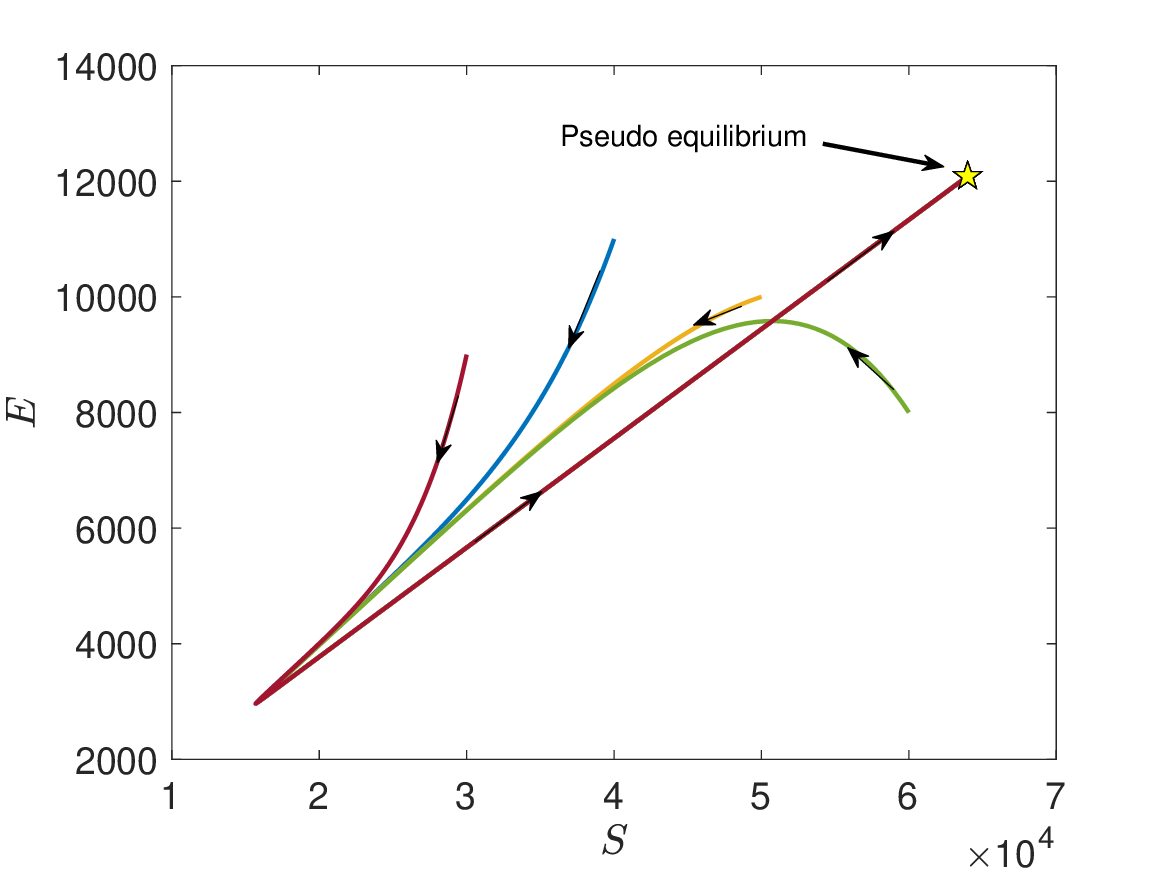}
		\includegraphics[scale=0.3]{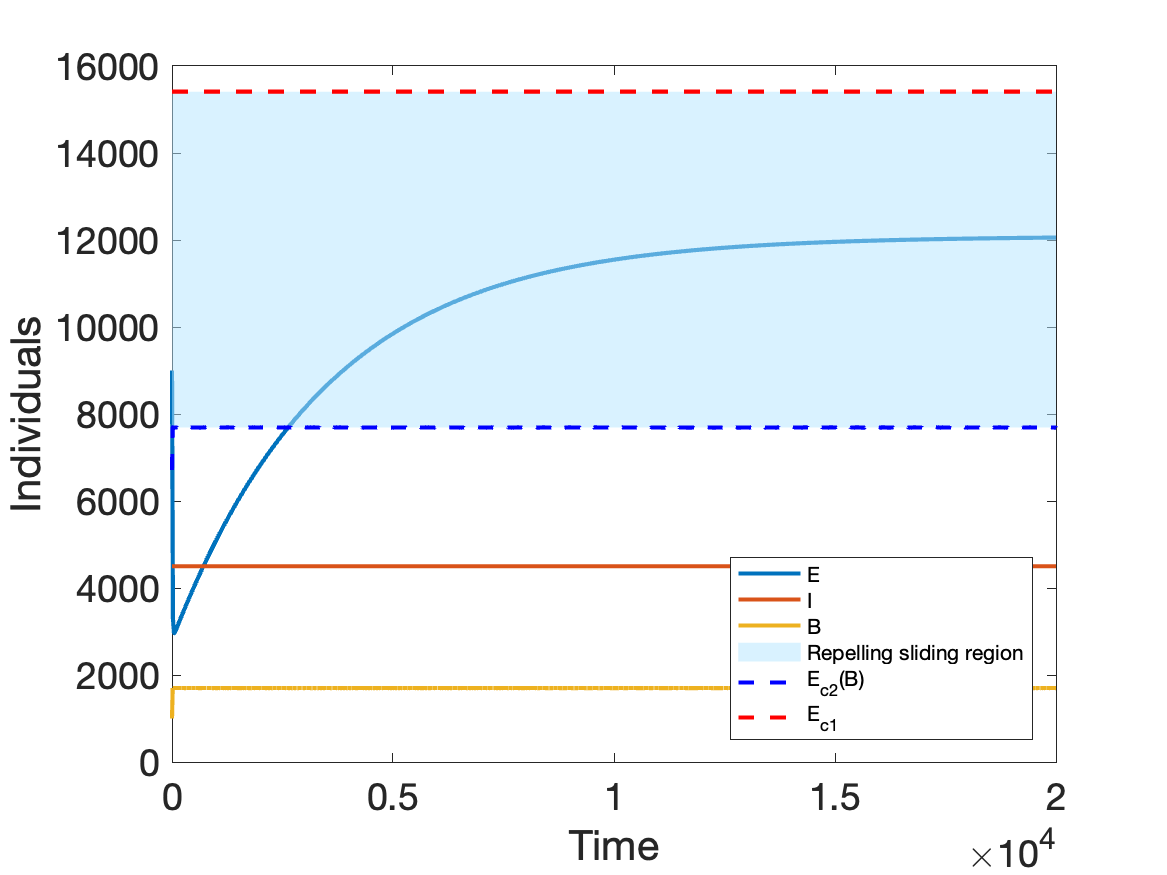}
		\caption{Two-dimensional phase portrait and time series illustrating convergence to the pseudo-equilibrium $D^{s}$ (yellow star) under the reduced sliding dynamics on the switching manifold. The left panel shows trajectories of the algebraically extended reduced tangential system in the $(S,E)$-plane for different initial conditions within the admissible repelling sliding region, while the right panel shows the corresponding time evolution of the state variables. Along these reduced dynamics, $I(t)=I_0$ and the remaining variables converge to $D^s$, consistently with Theorem~\ref{thm:pseudo_stability}. } 
		\label{fig:pseudo_eq} 
	\end{figure} 
	
	Fig.~\ref{fig:pseudo_eq} illustrates convergence to $D^s$ under the algebraically extended reduced tangential dynamics \eqref{Rep_sys}, consistently with Theorem~\ref{thm:pseudo_stability}. This convergence does not imply attraction in the full Filippov system. We use $\beta=0.000021$, $\xi=0.001$, $\lambda=0.93$, $K=1000$, $\mu_1=1.2$, and $I_0=4500$, with the remaining parameters as in Table~\ref{param_val}.
	
\begin{figure}[H]
		\centering
		(a) \includegraphics[scale=0.28]{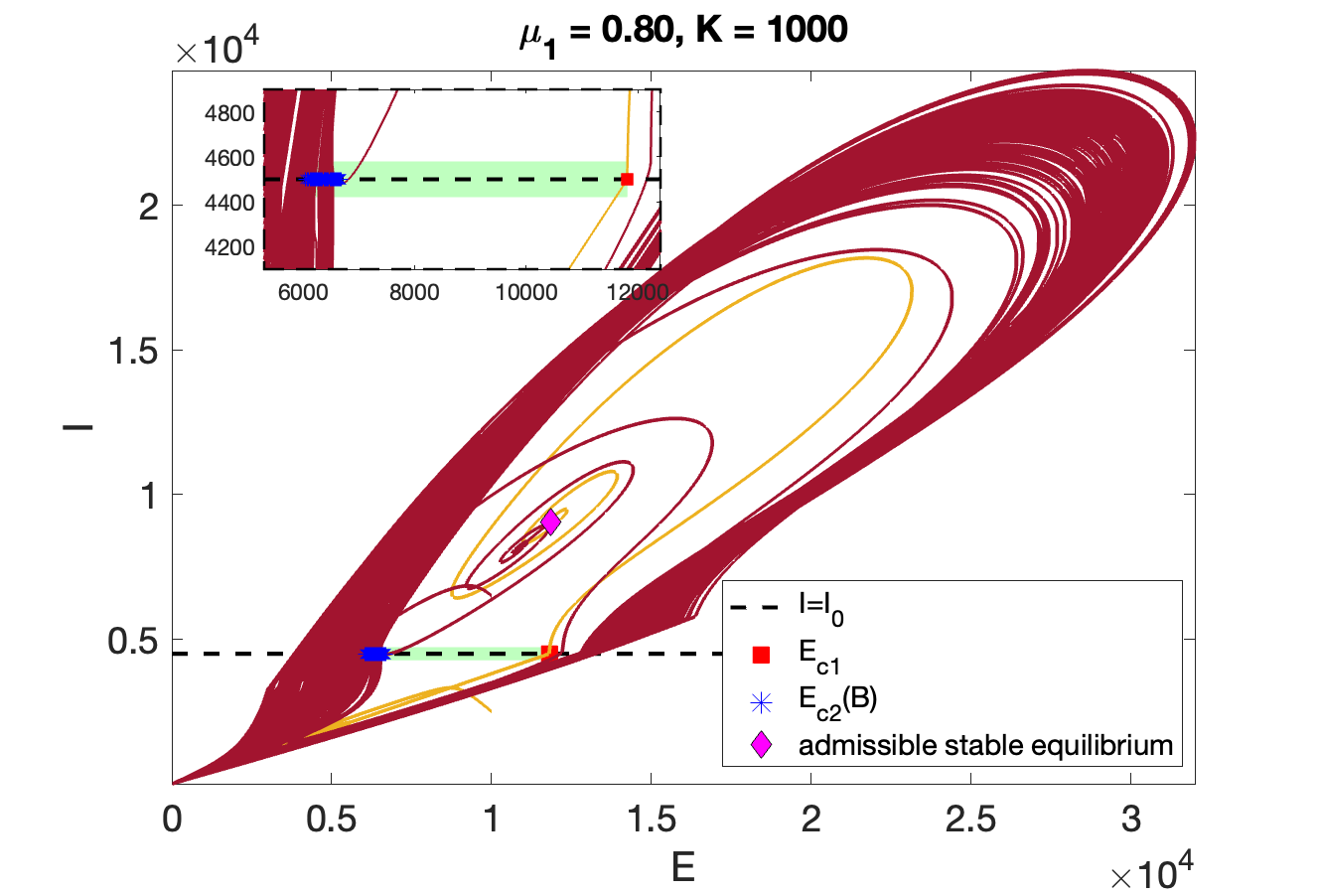}
		(b) \includegraphics[scale=0.28]{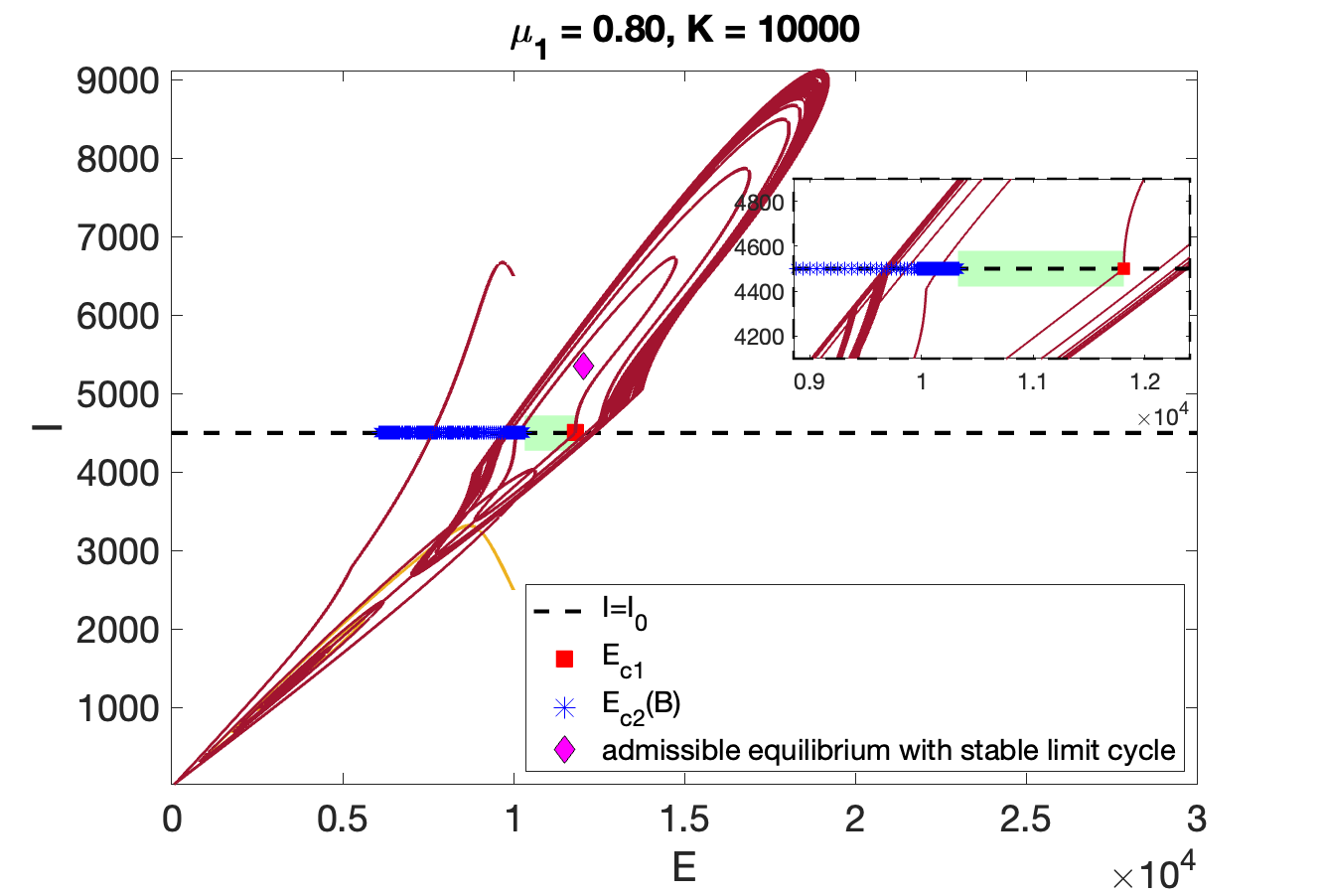}\\
		(c) \includegraphics[scale=0.28]{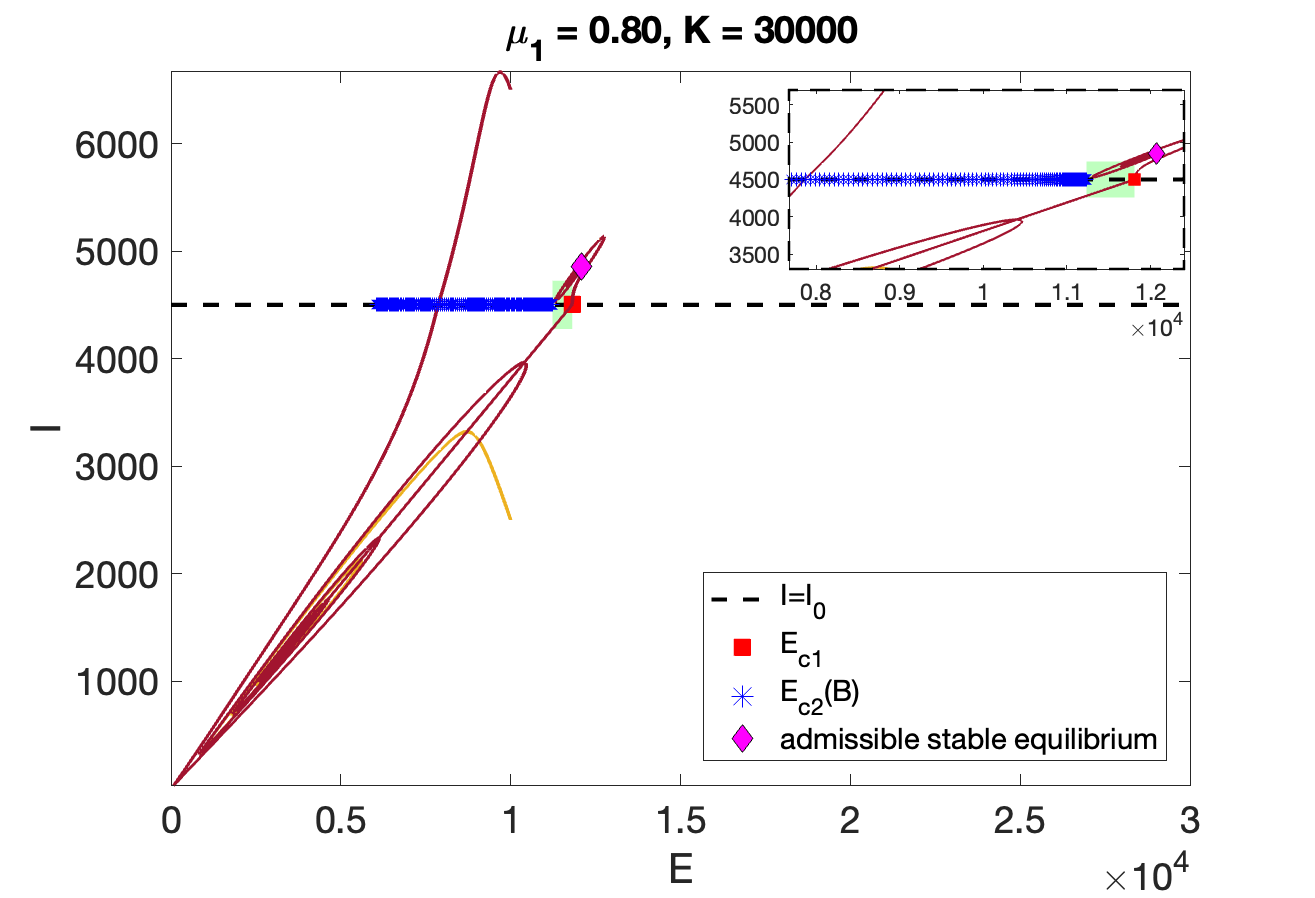}
		(d) \includegraphics[scale=0.28]{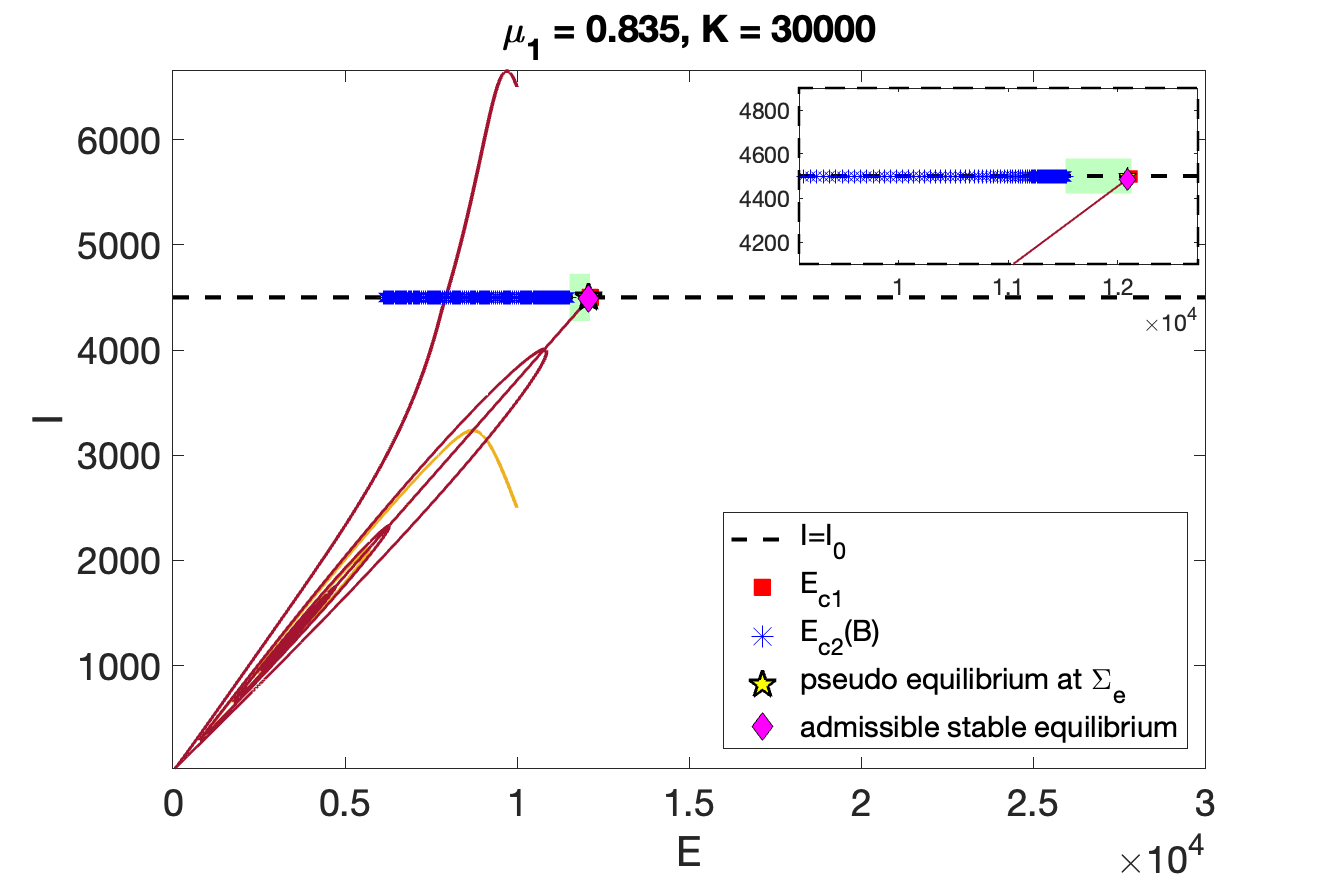} \\
		(e) \includegraphics[scale=0.28]{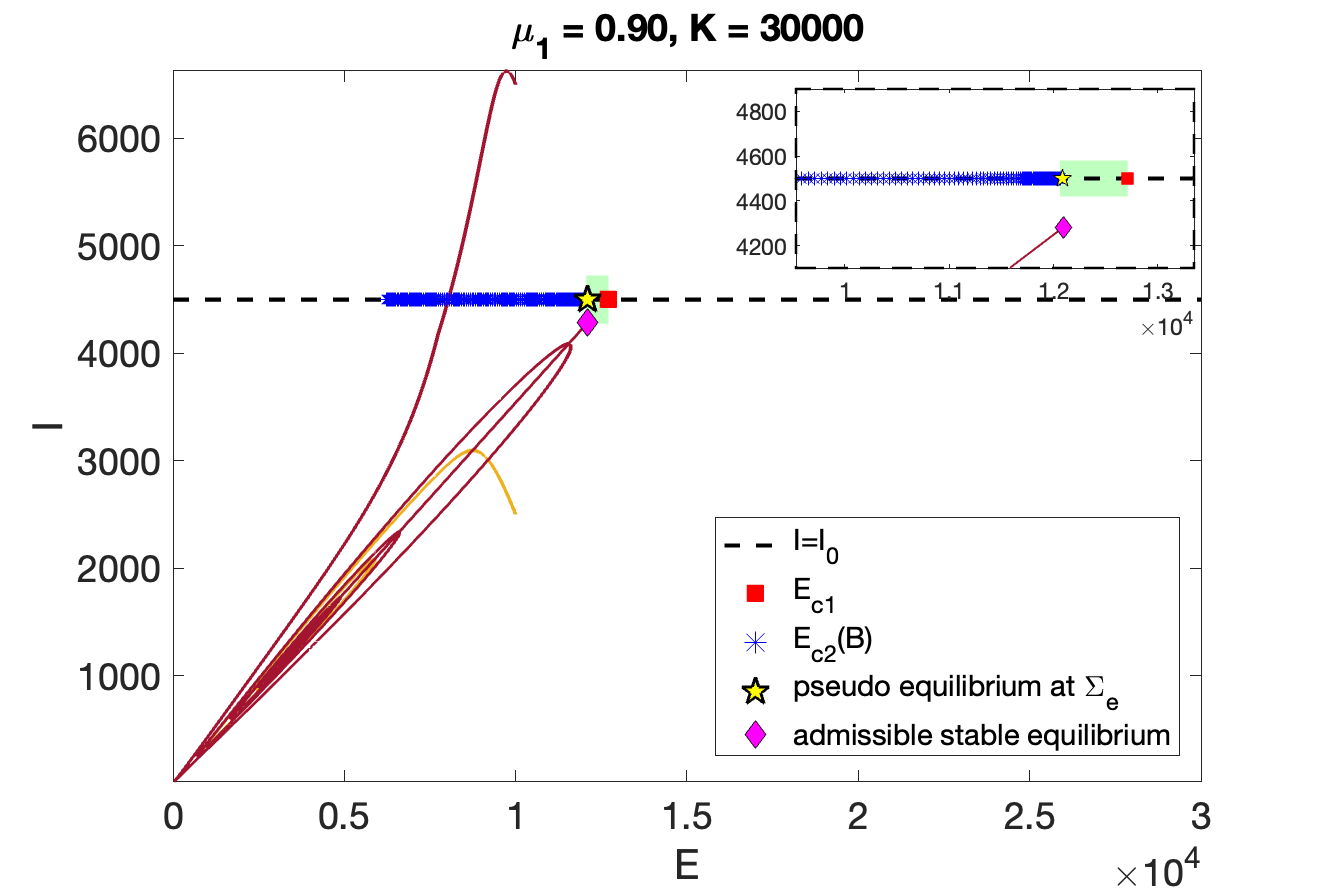}
		(f) \includegraphics[scale=0.28]{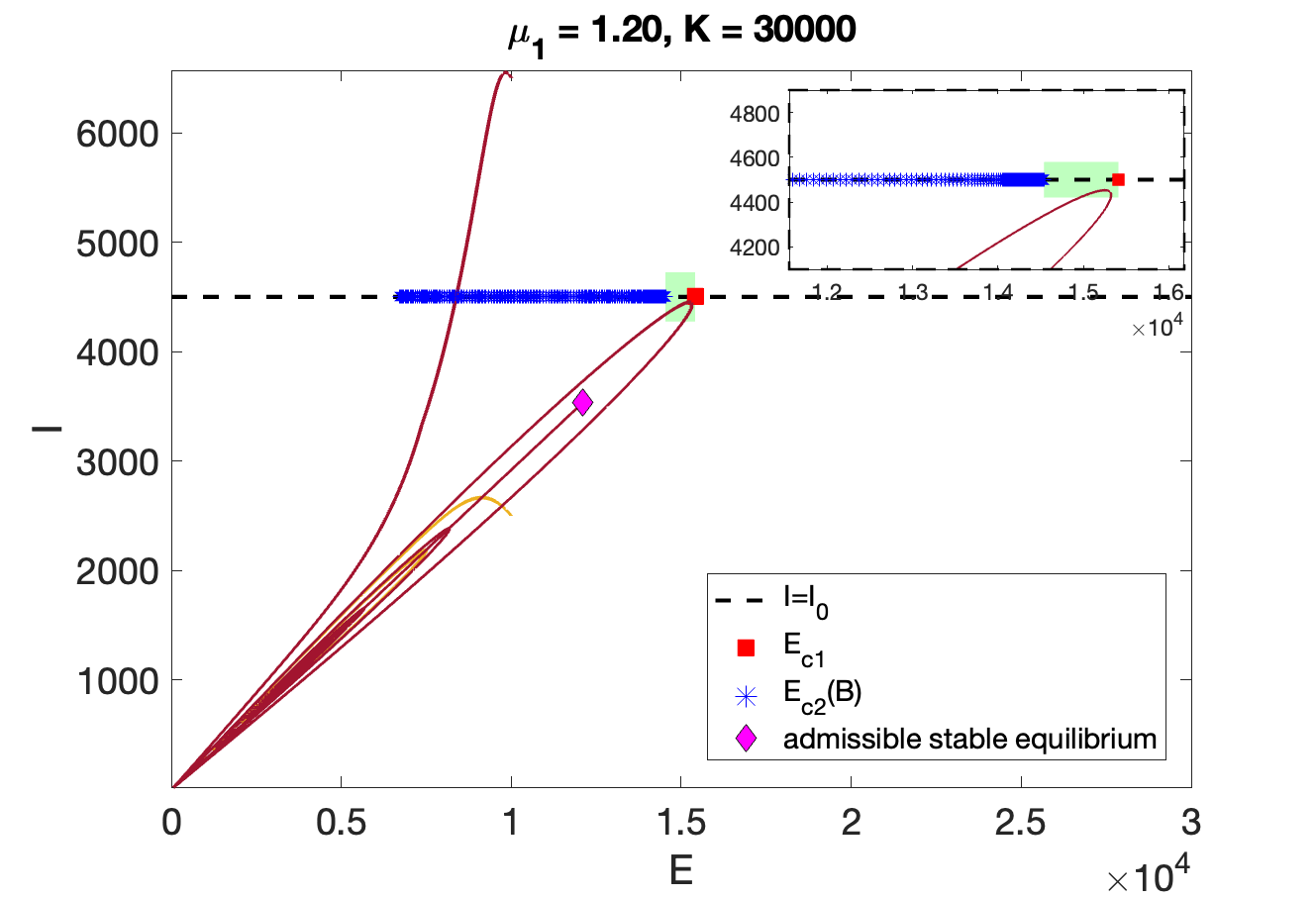} 
		\caption{
			Phase portraits in the $(E,I)$-plane for the Filippov system for different
			values of $\mu_1$ and $K$. Panels (a)-(c) correspond to $\mu_1=0.80$
			and increasing values of $K$. In (a), trajectories converge to a stable
			admissible equilibrium of $F_{\Sigma_2}$; in (b), a crossing periodic
			orbit intersects the switching manifold $I=I_0$ near the tangency
			boundaries without entering the interior of the repelling sliding region;
			and in (c), the admissible equilibrium of $F_{\Sigma_2}$ is again stable
			for larger $K$. Panels (d)-(f) correspond to $K=30000$ and increasing
			values of $\mu_1$. In (d) and (e), the admissibility of the equilibrium
			branches changes as they interact with the switching manifold, while the
			pseudo-equilibrium $D^s$ is admissible within $\Sigma_e$ over an
			intermediate parameter range. In (f), $D^s$ is no longer admissible and
			the trajectories converge to the admissible equilibrium of
			$F_{\Sigma_1}$.
		}
		\label{fig:phase_dde_system}
	\end{figure}

	We further analyze the qualitative dynamics induced by the switching
	manifold by varying key model parameters.  Fig.~\ref{fig:phase_dde_system}
	shows the $(E,I)$-phase portraits of the Filippov system
	\eqref{ddemodel} in the absence of time delay under variations of the
	baseline healthcare capacity $K$ and the maximum hospitalization rate
	$\mu_1$. Initial conditions are chosen on both sides of the switching
	manifold $\Sigma=\{I=I_0\}$ to capture dynamics governed by both
	subsystems $F_{\Sigma_1}$ and $F_{\Sigma_2}$. In each panel, the
	switching manifold is shown together with the upper tangency threshold
	$E_{c1}$ and the state-dependent lower threshold $E_{c2}(B(t))$, corresponding to the projections of the tangency sets $\Sigma_1^t$ and $\Sigma_2^t$ onto $\Sigma$. The trajectories were computed using \texttt{ode45} with relative and absolute tolerances of $10^{-7}$ and $10^{-9}$, respectively. The integration was advanced over short intervals, with the active vector field updated according to the trajectory position relative to $I=I_0$ and the Filippov conditions on $\Sigma$. \\ 
	
	\noindent
	\textbf{Crossing periodic dynamics near quadratic tangency:} For $\mu_1=0.80$ and $K=1000$
	(Fig.~\ref{fig:phase_dde_system}(a)), the admissible equilibrium of
	$F_{\Sigma_2}$ is stable and attracts trajectories initiated on both 
	sides of the switching manifold. These trajectories cross $\Sigma$
	outside the interior of the repelling sliding region and converge to a
	steady endemic state. For $K=10000$ (Fig.~\ref{fig:phase_dde_system}(b)), a stable crossing
	periodic orbit is observed. The orbit alternates between the dynamics of
	$F_{\Sigma_1}$ and $F_{\Sigma_2}$ and intersects the switching manifold
	near the tangency boundaries $\Sigma_1^t$ and $\Sigma_2^t$ without
	entering the interior of the repelling sliding region. At a generic
	quadratic tangency,
	$ 
	L=0, 
	\mathcal{L}_{F_{\Sigma_i}}L=0, 
	\mathcal{L}_{F_{\Sigma_i}}^2L\neq0, 
	i=1,2,
	$ 
	with the corresponding Lie-derivative conditions established in
	Supplementary Section~S8. The observed periodic orbit is therefore consistent
	with a nearby quadratic grazing transition. For larger $K$ (Fig.~\ref{fig:phase_dde_system}(c)), the equilibrium of
	$F_{\Sigma_2}$ is again stable, illustrating the nonmonotone effect of
	baseline healthcare capacity on the long-term dynamics. 
	
	\noindent
	\textbf{Boundary-equilibrium transition:}  Figs.~\ref{fig:phase_dde_system}(d)-(f) illustrate the effect of varying 
	the maximum hospitalization rate $\mu_1$ for a fixed baseline healthcare
	capacity $K=30000$. As $\mu_1$ increases, equilibrium branches of the
	smooth subsystems interact with the switching manifold
	$\Sigma=\{I=I_0\}$, resulting in changes in their admissibility within
	the full Filippov system. At a critical value of $\mu_1$, the endemic equilibrium of
	$F_{\Sigma_2}$ reaches $\Sigma$ and becomes a boundary equilibrium.
	As $\mu_1$ passes through this value, the equilibrium changes between
	admissible and virtual, consistently with a boundary-equilibrium
	transition. Over an intermediate parameter range, the pseudo-equilibrium $D^s$
	satisfies
	$ 
	E_{c2}(B^s)<E^s<E_{c1}
	$ 
	and is therefore admissible within the repelling sliding region
	$\Sigma_e$. Varying $\mu_1$ does not change the algebraic location of
	$D^s$; rather, it changes the tangency boundaries $E_{c1}$ and
	$E_{c2}(B^s)$ and hence its admissibility. Although $D^s$ organizes the
	reduced tangential dynamics on $\Sigma$, the repelling character of
	$\Sigma_e$ prevents it from being attracting in the full Filippov
	system. For larger $\mu_1$ (Fig.~\ref{fig:phase_dde_system}(f)), $D^s$ becomes
	inadmissible after reaching a boundary of $\Sigma_e$, and the displayed
	trajectories converge to the admissible equilibrium of
	$F_{\Sigma_1}$. Thus, variation in $\mu_1$ reorganizes the admissibility
	of both smooth and sliding equilibria and can alter the observed
	long-term state. 
	
	\begin{figure}[H]
		\centering
		\includegraphics[scale=0.3]{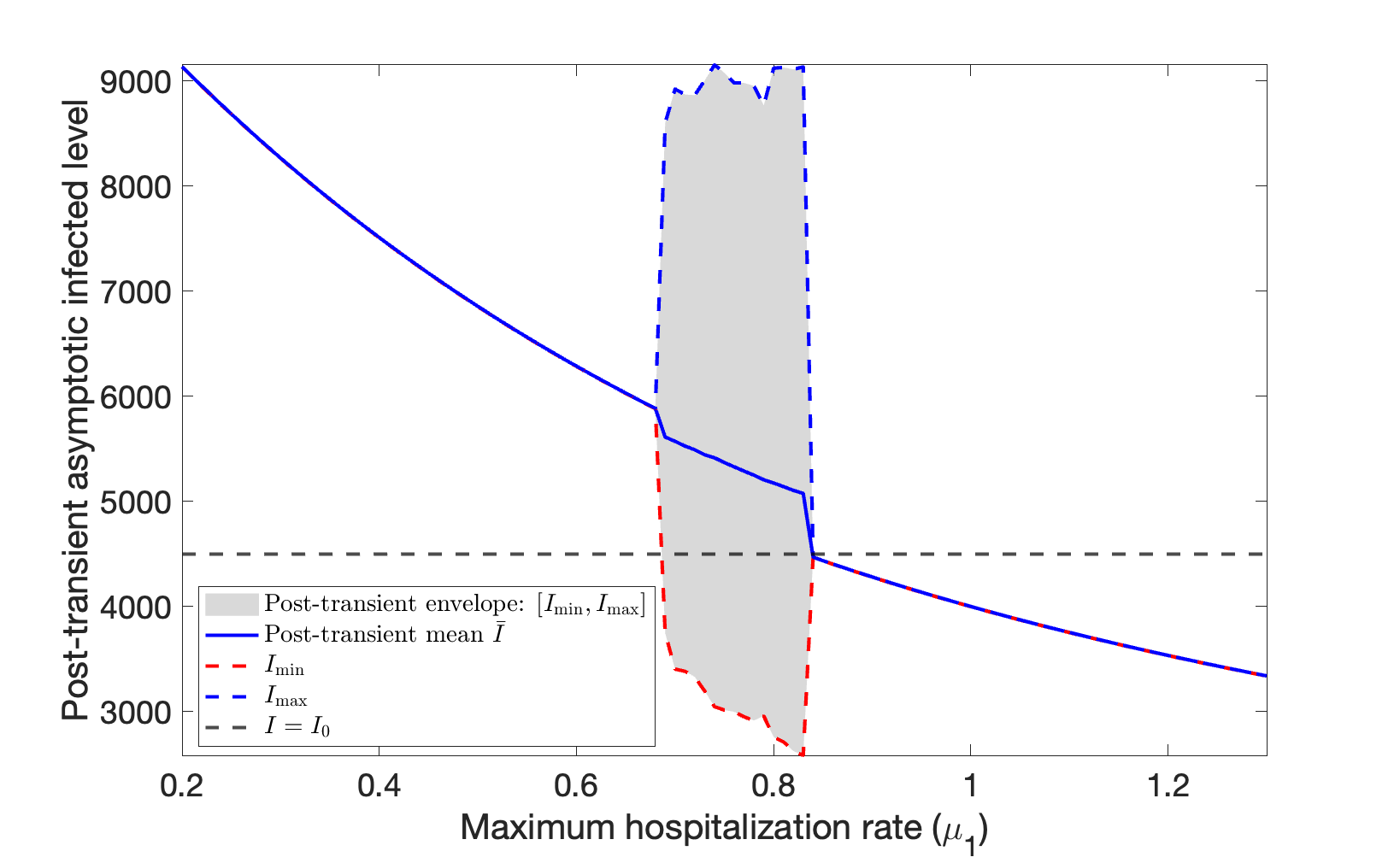}
		\caption{Post-transient infected levels as a function of the maximum hospitalization
			rate $\mu_1$.
			The solid blue curve shows the post-transient mean infected level $\bar I$,
			while the shaded region denotes the envelope $[I_{\min}, I_{\max}]$ of the
			long-time dynamics over multiple initial conditions.
			The crossing of $\bar I$ and the envelope through the switching threshold
			$I_0$ is consistent with a boundary-equilibrium transition in the
			Filippov system.} 
		\label{fig:boundary_collision_signature}
	\end{figure}  
	
	Fig.~\ref{fig:boundary_collision_signature} summarizes the long-time
	dynamics as $\mu_1$ varies. For this parameter sweep, the delay-free
	Filippov system was integrated using \texttt{ode45}.
	Simulations were continued to $t=35000$, with $t<25000$ discarded as transient. 
	For each value of $\mu_1$, trajectories were initialized on both sides
	of $I=I_0$, and the post-transient mean, minimum, and maximum values of
	$I(t)$ were recorded. The shaded region between the dashed curves
	represents the resulting infection envelope, with the blue and red
	dashed curves denoting $I_{\min}$ and $I_{\max}$, respectively, and the
	solid curve denoting the corresponding mean. For small $\mu_1$, the mean infected level remains above $I_0$,
	consistent with convergence to an admissible equilibrium of
	$F_{\Sigma_2}$. As $\mu_1$ increases, the envelope widens over a finite
	parameter interval, indicating sustained oscillations, and then narrows
	again as steady endemic behavior is restored. As $\mu_1$ increases, the post-transient mean crosses $I=I_0$, indicating a change in long-term infection levels relative to the switching manifold. The mean crossing alone does not identify a boundary-equilibrium bifurcation; the boundary-equilibrium condition is characterized analytically in Section~\ref{SLIDDYN}. For sufficiently large $\mu_1$, the
	simulated trajectories converge to infected levels below $I_0$. Throughout this parameter sweep, trajectories cross the switching
	manifold transversely and do not enter the repelling sliding region.
	Thus, the observed changes illustrate steady and oscillatory dynamics relative to the Filippov switching threshold. \\ 
	
\begin{figure}[H]
		\centering
		(a) \includegraphics[scale=0.3]{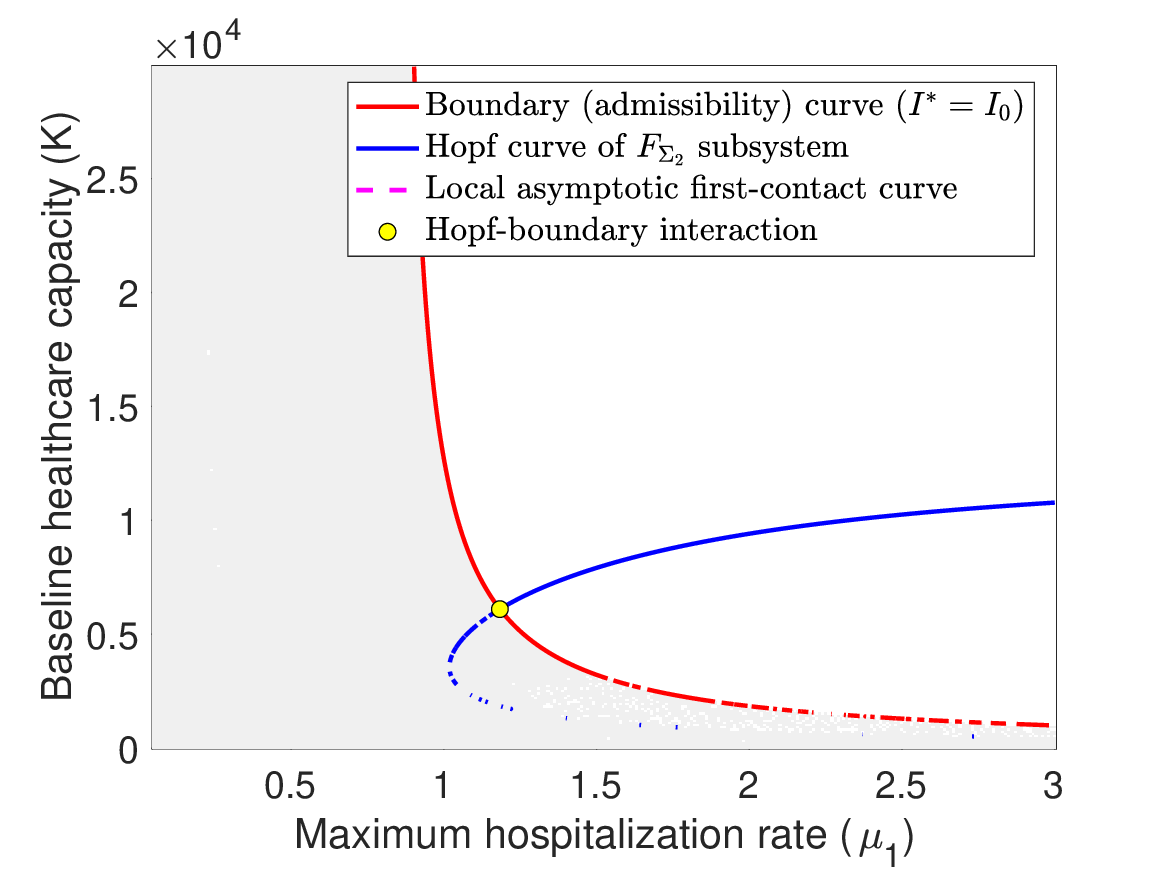} 
		(b)  \includegraphics[scale=0.3]{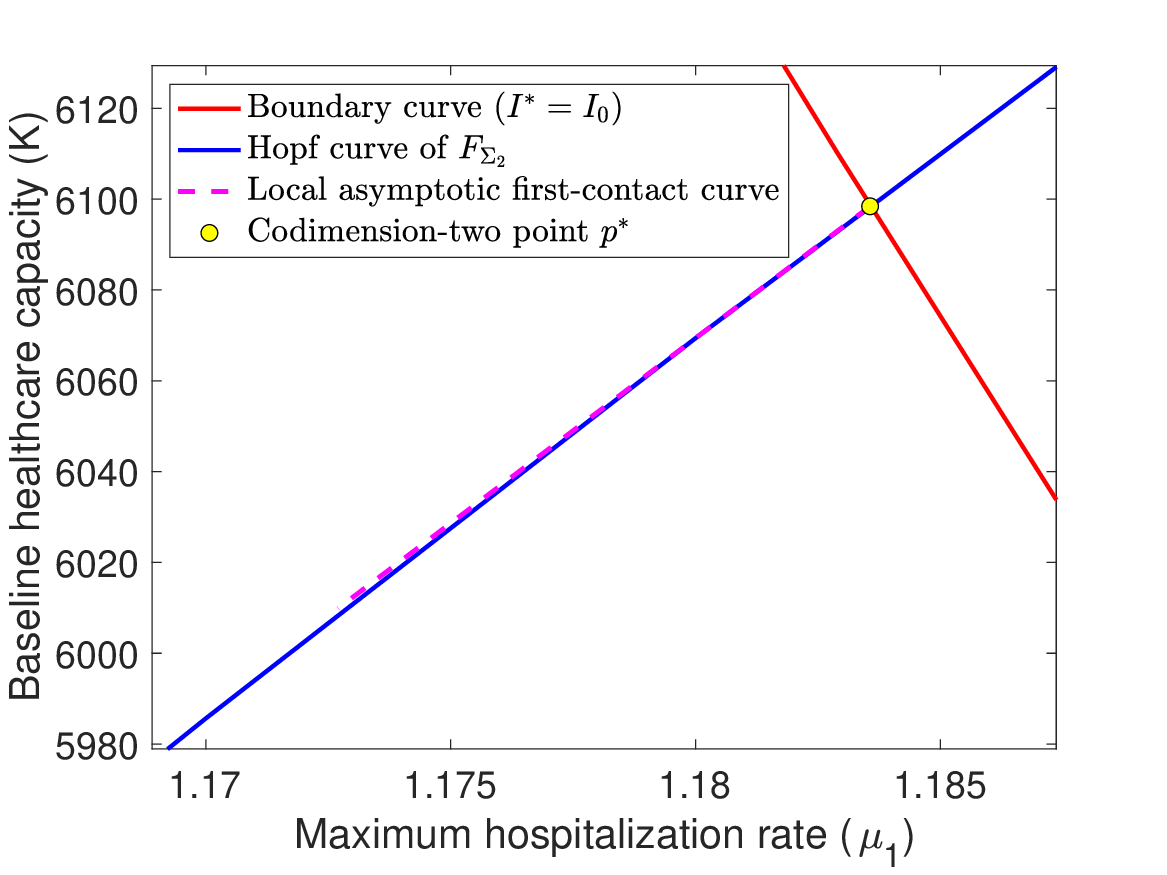} 
		\caption{
			Codimension-two Hopf-boundary interaction and local first-contact
			geometry in the $(\mu_1,K)$-plane for the Filippov system.
			In panel (a), the red curve represents the boundary-equilibrium
			(admissibility) locus $I_1^*(\mu_1,K)=I_0$, while the blue curve
			represents the numerically detected Hopf locus of the smooth extension
			of subsystem $F_{\Sigma_2}$.
			The magenta dashed curve denotes the local asymptotic first-contact
			locus of the Hopf-generated periodic orbit with the switching manifold,
			and the yellow marker denotes the transversal codimension-two
			Hopf-boundary interaction point.
			Panel (b) shows a local magnification near this point, illustrating
			the transversal intersection of the boundary-equilibrium and Hopf loci
			and the quadratic tangency of the first-contact curve to the Hopf locus,
			consistent with Proposition~\ref{prop:HB_grazing}.
		}
		\label{fig:hopf_boundary_scan}
	\end{figure}

	\noindent 
	\textbf{Codimension-two Hopf-boundary interaction:}
	Fig.~\ref{fig:hopf_boundary_scan} shows the two-parameter bifurcation
	structure in the $(\mu_1,K)$-plane together with a local magnification
	near the Hopf-boundary interaction point. For $\Lambda=1400$, $\xi=0.001$, $\beta=2.1\times10^{-5}$, $\lambda=0.93$, $\tau=0$, and the remaining fixed parameters from Table~\ref{param_val}, solving the boundary and Hopf conditions gives the interaction point and transversality determinants reported below. For each parameter pair
	$(\mu_1,K)$, we compute the endemic-equilibrium branch of the smooth
	extension of $F_{\Sigma_2}$ and evaluate the corresponding Jacobian
	matrix $J_{F_{\Sigma_2}}$. The endemic equilibria were computed using
	MATLAB's \texttt{fsolve} with multiple initial seeds and function and
	step tolerances of $10^{-12}$, while local stability and the Hopf locus
	were determined from the eigenvalues of the corresponding numerical
	Jacobian. The Hopf-boundary intersection was subsequently refined by
	simultaneously solving the boundary and Hopf conditions.
	
	The boundary-equilibrium locus, shown in red, is defined by
	$ 
	\Delta(\mu_1,K)
	=
	I_1^*(\mu_1,K)-I_0
	=
	0,
	$ 
	and separates parameter values for which the equilibrium is admissible
	for $F_{\Sigma_2}$,
	$
	I_1^*>I_0,
	$ 
	from those for which it is virtual in the full Filippov system. The
	Hopf locus, shown in blue, is obtained by tracking the critical
	complex-conjugate eigenvalue pair $s_H(\mu_1,K)$ and locating parameter
	values satisfying
	$ 
	\mathcal{H}(\mu_1,K)
	=
	\Re s_H(\mu_1,K)
	=
	0,
	$ 
	with nonzero imaginary part and all remaining eigenvalues having
	negative real parts. Their intersection, marked in yellow, identifies
	a parameter value at which the smooth Hopf and boundary-equilibrium
	conditions occur simultaneously. According to
	Theorem~\ref{thm:codim2_HB}, a transversal intersection of these two
	codimension-one loci constitutes a codimension-two Hopf-boundary
	interaction.
	
	To verify the conditions of Theorem~\ref{thm:codim2_HB}, we solve
	simultaneously
	$ 
	\Delta(\mu_1,K)=0, 
	\mathcal{H}(\mu_1,K)=0,
	$ 
	which yields the refined critical parameter pair
	$
	(\mu_1^*,K^*)\approx(1.18356,6098.42).
	$ 
	At this point,
	$ 
	I^*=4500=I_0,
	$ 
	and the corresponding Jacobian spectrum is
	\[
	\sigma(J_2)
	=
	\left\{
	-2.7385\times10^{-12}\pm0.17508623i,\,
	-2.9327\times10^{-4},\,
	-1.6846,\,
	-1.6957,\,
	-2.2303
	\right\}.
	\]
	Thus, within numerical precision, one simple complex-conjugate pair
	lies on the imaginary axis while all remaining eigenvalues have
	strictly negative real parts. The defining residuals are
	$ 
	\Delta(\mu_1^*,K^*)=0, 
	\mathcal{H}(\mu_1^*,K^*)
	=
	-2.7385\times10^{-12}.
	$ 
	
	The corresponding transversality determinant is
	\[
	\det D_{(\mu_1,K)}(\mathcal B,\mathcal{H})
	=
	-4.7905\times10^{-6}\neq0,
	\]
	showing that the gradients of the two defining conditions
	are linearly independent. Since $\mu_1$ and $K$ have
	substantially different numerical scales, we also introduce
	\[
	\kappa=10^{-4}K.
	\]
	In the rescaled coordinates,
	\[
	\det D_{(\mu_1,\kappa)}(\mathcal B,\mathcal{H})
	=
	-4.7905\times10^{-2}\neq0,
	\]
	confirming that transversality is preserved under rescaling.
	
	The first Lyapunov coefficient associated with the critical Hopf pair is
	$ 
	\ell_1(p^*)
	=
	-3.29403\times10^{-2}<0.
	$
	Hence, the Hopf bifurcation of the smooth extension of
	$F_{\Sigma_2}$ is nondegenerate and, under the standard Hopf
	convention, supercritical.
	
	We next verify numerically the additional hypotheses of
	Proposition~\ref{prop:HB_grazing}. The infectious component of the
	critical right eigenvector is
	$ 
	q_I
	=
	0.645074565
	$ 
	in the scaled coordinates, with corresponding physical magnitude
	$ 
	|q_I|_{\mathrm{phys}}
	=
	2.90284\times10^{3}.
	$ 
	Thus,
	$ 
	q_I\neq0,
	$
	verifying condition {\rm (G3)} and showing that the critical Hopf mode
	has a nonzero component transverse to the switching manifold. Numerical differentiation of
	$
	\Delta(p)=I^*(p)-I_0
	$ 
	at the critical point gives
	\[
	\nabla_p\Delta(p^*)
	=
	\left(
	-2.816178\times10^{3},
	-1.628837\times10^{-1}
	\right),
	\] 
	with
	$ 
	\|\nabla_p\Delta(p^*)\|
	=
	2.816178\times10^{3}\neq0.
	$ 
	Hence, condition {\rm (G2)} is satisfied. Since $\Delta=0$ and
	$\mathcal B=0$ define the same regular boundary-equilibrium locus,
	we obtain, at $p^*$,
	\[
	\begin{aligned}
		\det D_{(\mu_1,K)}(\Delta,\mathcal H)
		&=\frac{\partial\Delta}{\partial\mu_1}(p^*)
		\det D_{(\mu_1,K)}(\mathcal B,\mathcal H)\\
		&\approx 1.34909\times10^{-2}\neq0.
	\end{aligned}
	\]
	Consequently, $(\Delta,\mathcal H)$ form regular local coordinates
	near $p^*$.
	Using
	$
	\mathcal{H}=\Re s_H
	$ 
	as the local Hopf coordinate, the leading-order numerical
	first-contact relation is
	\[
	\mathcal{H}
	=
	1.71110\times10^{-10}\Delta^2
	+
	o(\Delta^2).
	\]
	Since the small-amplitude periodic orbit exists on the
	$\mathcal{H}>0$ side under this numerical orientation, whereas
	Proposition~\ref{prop:HB_grazing} orients the Hopf coordinate $\eta$
	so that the periodic orbit exists for $\eta<0$, setting
	$ 
	\eta=-\mathcal{H}
	$ 
	gives
	$ 
	\eta
	=
	-C_g\Delta^2
	+
	o(\Delta^2), 
	C_g>0,
	$ 
	in agreement with the analytical result.
	
	The magenta dashed curve in Fig.~\ref{fig:hopf_boundary_scan}
	represents this local asymptotic first-contact locus. Panel~(b)
	magnifies the neighborhood of $p^*$, showing the transversal
	intersection of the boundary-equilibrium and Hopf loci and the
	quadratic tangency of the first-contact curve to the Hopf locus.
	Thus, the numerical results support both the codimension-two
	Hopf-boundary interaction of Theorem~\ref{thm:codim2_HB} and the local
	first-contact geometry established in
	Proposition~\ref{prop:HB_grazing}.
	
\begin{figure}[H]
		\centering
		(a)\includegraphics[scale=0.27]{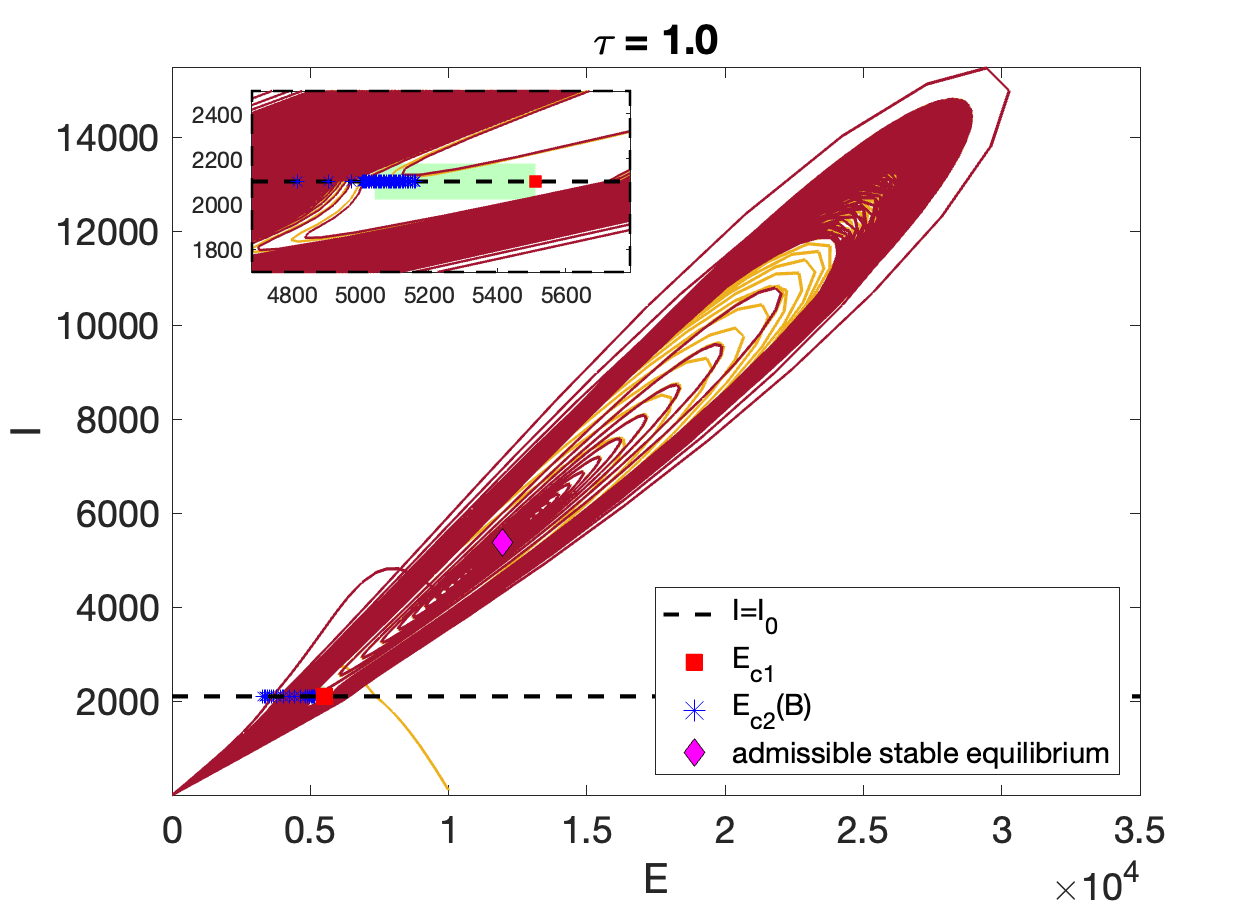}
		\includegraphics[scale=0.27]{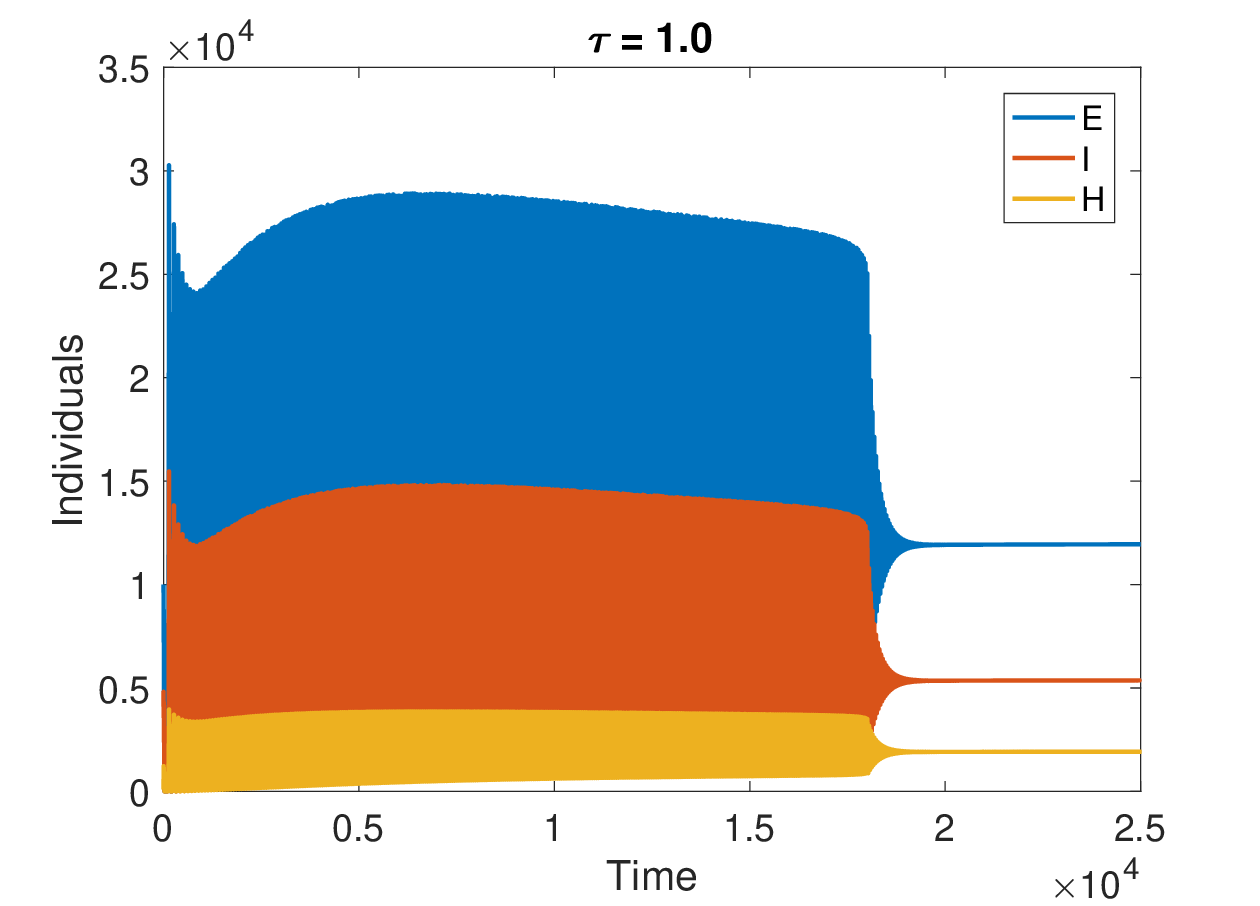} \\
		(b)\includegraphics[scale=0.27]{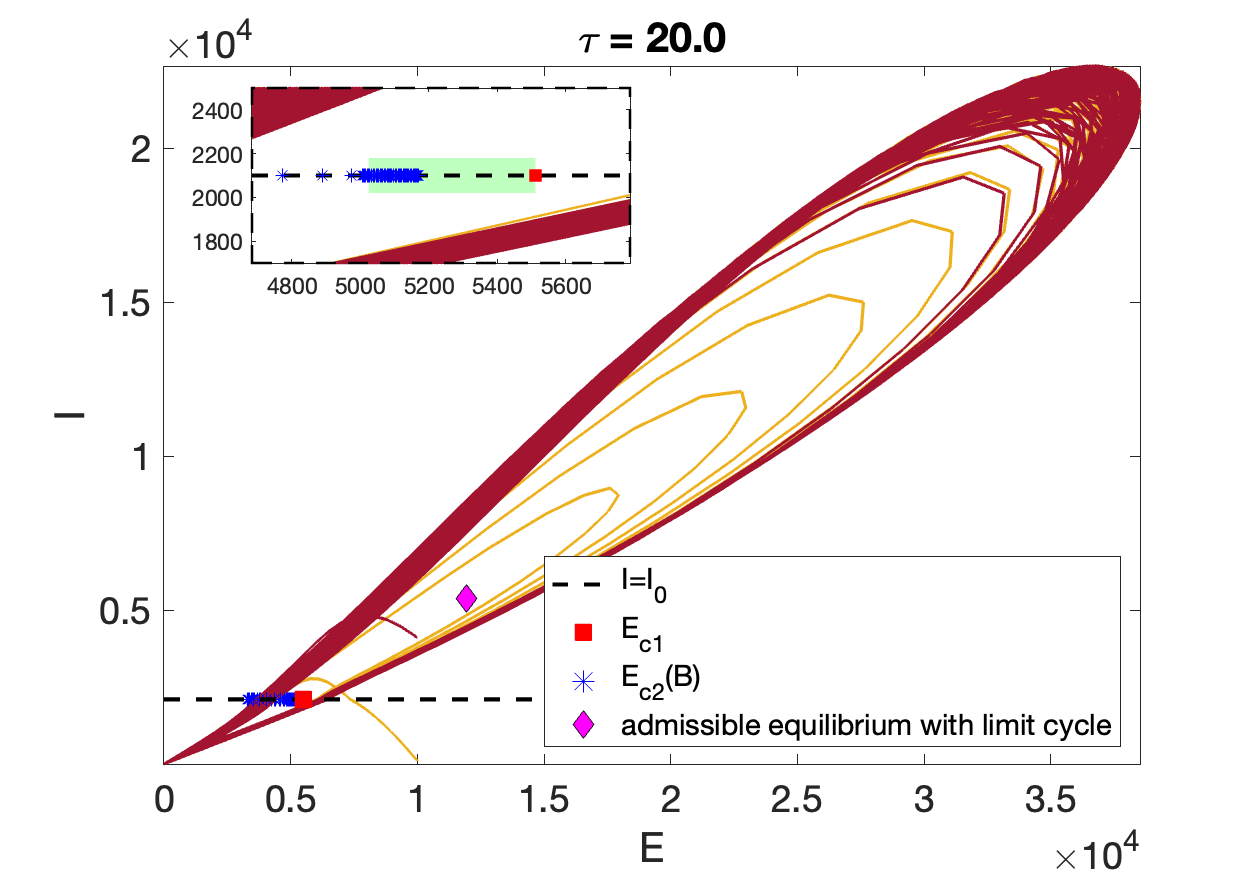}
		\includegraphics[scale=0.27]{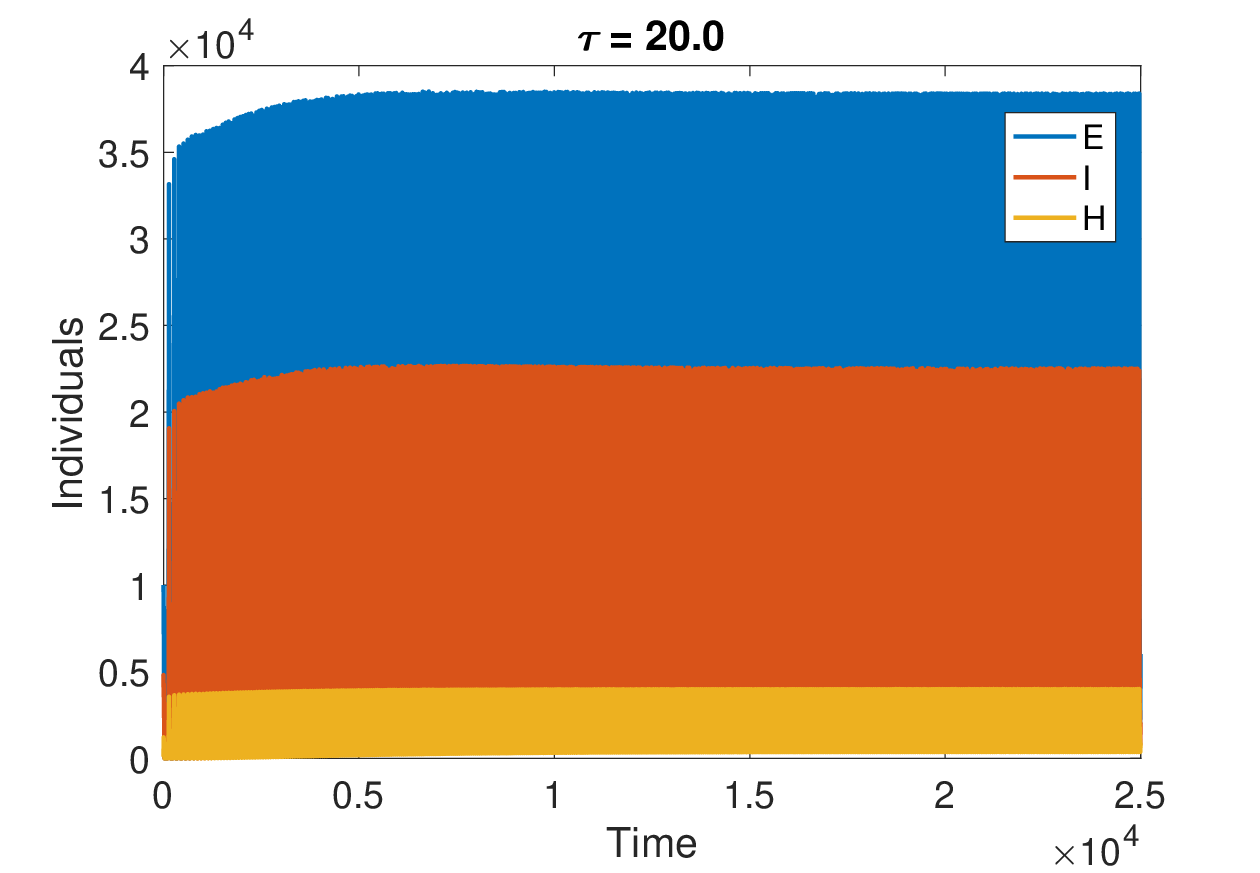} \\
		(c)\includegraphics[scale=0.27]{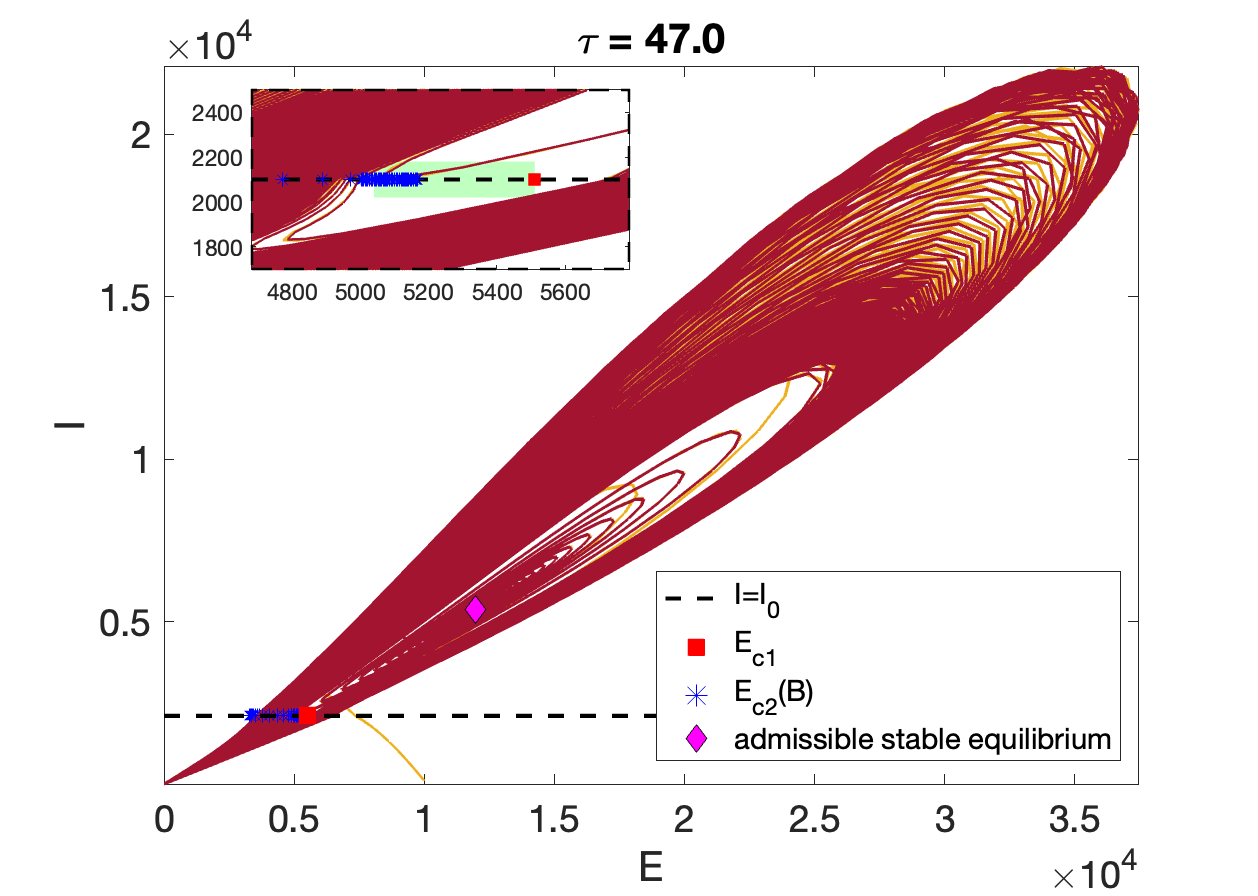}
		\includegraphics[scale=0.27]{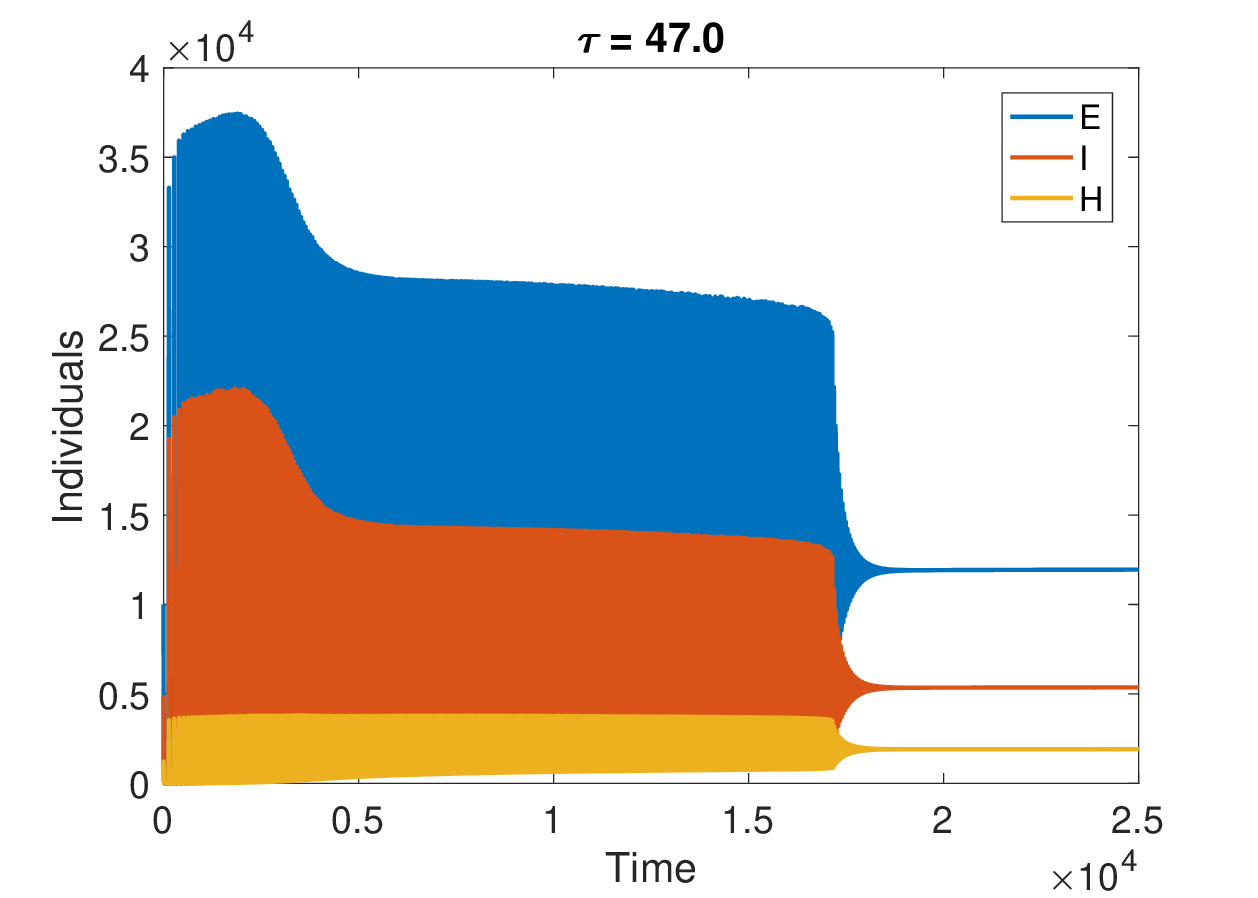}
		\caption{
			Phase portraits in the $(E,I)$-plane (left) and corresponding time series (right)
			for the Filippov system with time delay $\tau$.
			(a) For $\tau=1$, trajectories converge to a stable admissible equilibrium of
			$F_{\Sigma_2}$.
			(b) For $\tau=20$, a sustained periodic oscillation is observed, with
			a crossing cycle surrounding the repelling
			sliding mode region.
			(c) For $\tau=47$, the oscillations disappear and trajectories again converge to a
			stable admissible equilibrium of $F_{\Sigma_2}$.
			The dashed line denotes the switching manifold $I=I_0$, while $E_{c1}$ and
			$E_{c2}(B)$ indicate the upper and lower tangency thresholds, and the green region
			represents the repelling sliding region.
		} 
		\label{fig:Fig._vary_delay}
	\end{figure}

	Further, we present numerical simulation results for the Filippov system \eqref{ddemodel} with time delay $\tau$. The delayed Filippov system was integrated using \texttt{dde23}, with a constant prehistory equal to the prescribed initial state. The delayed infection level $I(t-\tau)$ was evaluated through the delayed-state argument of the solver, while the treatment regime at each time was selected according to the current infection level relative to $I_0$. The parameter values specified for these simulations are $\mu_1=0.80$, $K=11500$, $\omega=5000$, $\beta=0.000012$, and $I_0=2100$. For the parameter set of Figs.~\ref{delaysub2} and \ref{delaybif}, which additionally specifies $\Lambda=1400$, $\xi=0.001$, and $\lambda=0.93$, with the remaining parameters as in Table~\ref{param_val}, the smooth-subsystem characteristic equation is the one used in Table~\ref{tab:critical_delay_check}, because $I_0$ does not enter the smooth above-threshold vector field. Its endemic equilibrium satisfies $I_1^*\approx5370.49>I_0=2100$. For this parameter set, the equilibrium is locally asymptotically stable at $\tau=1$ and $\tau=47$, and unstable at $\tau=20$, in agreement with the qualitative sequence in Fig.~\ref{fig:Fig._vary_delay}. 
	
	Fig.~\ref{fig:Fig._vary_delay} illustrates the effect of $\tau$ on the long-term dynamics. For a small delay ($\tau=1$), trajectories converge to a stable admissible equilibrium of $F_{\Sigma_2}$. At $\tau=20$, sustained oscillations arise around the same equilibrium, with
	recurrent crossings of the switching
	manifold. The local characteristic-root analysis of $F_{\Sigma_2}$ concerns the admissible equilibrium and its nearby dynamics; once an oscillatory trajectory reaches the threshold, its subsequent evolution is governed by the full Filippov system. For $\tau=47$, the oscillations disappear and trajectories
	again converge to the admissible equilibrium of $F_{\Sigma_2}$. Thus,
	the numerical trajectories illustrate a nonmonotone dependence of the
	long-term dynamics on the healthcare-capacity response delay. 
	
	\noindent
	\textbf{Sequential Hopf and boundary-transition dynamics:}
	Motivated by the stable behavior in
	Fig.~\ref{fig:Fig._vary_delay}(a), we fix $\tau=1$ and vary
	$\mu_1$. We take $\Lambda=1400$, $\xi=0.001$,
	$\beta=1.2\times10^{-5}$, $K=11500$, $\lambda=0.93$,
	$\omega=5000$, and $I_0=2100$, with the remaining parameters
	as in Table~\ref{param_val}. The two constant initial
	histories, in the state order $(S,E,I,H,R,B)$, are
	$(60000,10000,100,100,50,1000)^\top$ and
	$(30000,10000,4100,100,50,1000)^\top$.
	The system is integrated using \texttt{dde23} on $[0,25000]$,
	with transients discarded up to $t=5000$.
	
\begin{figure}[H]
		\centering
		(a)\includegraphics[scale=0.3]{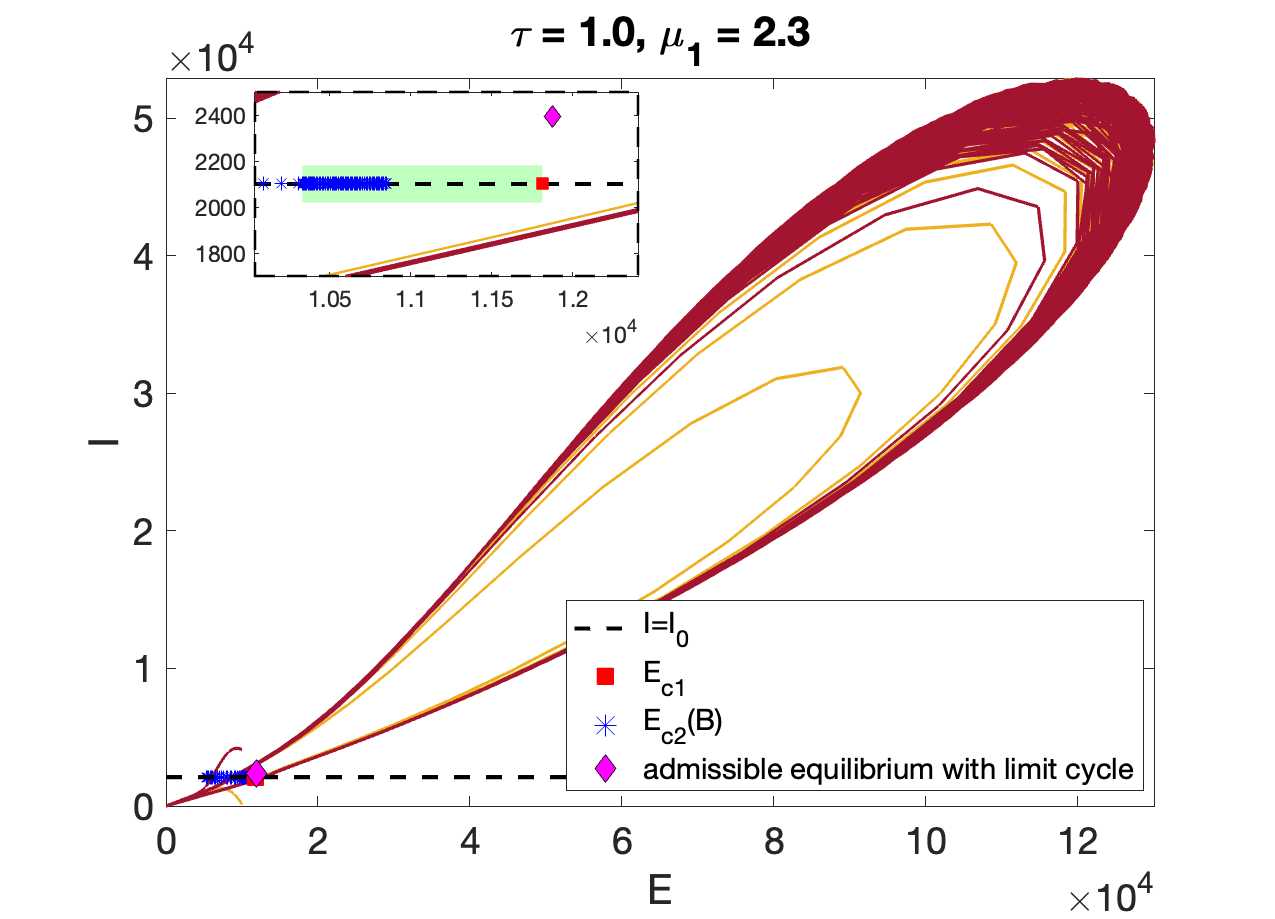}
		\includegraphics[scale=0.3]{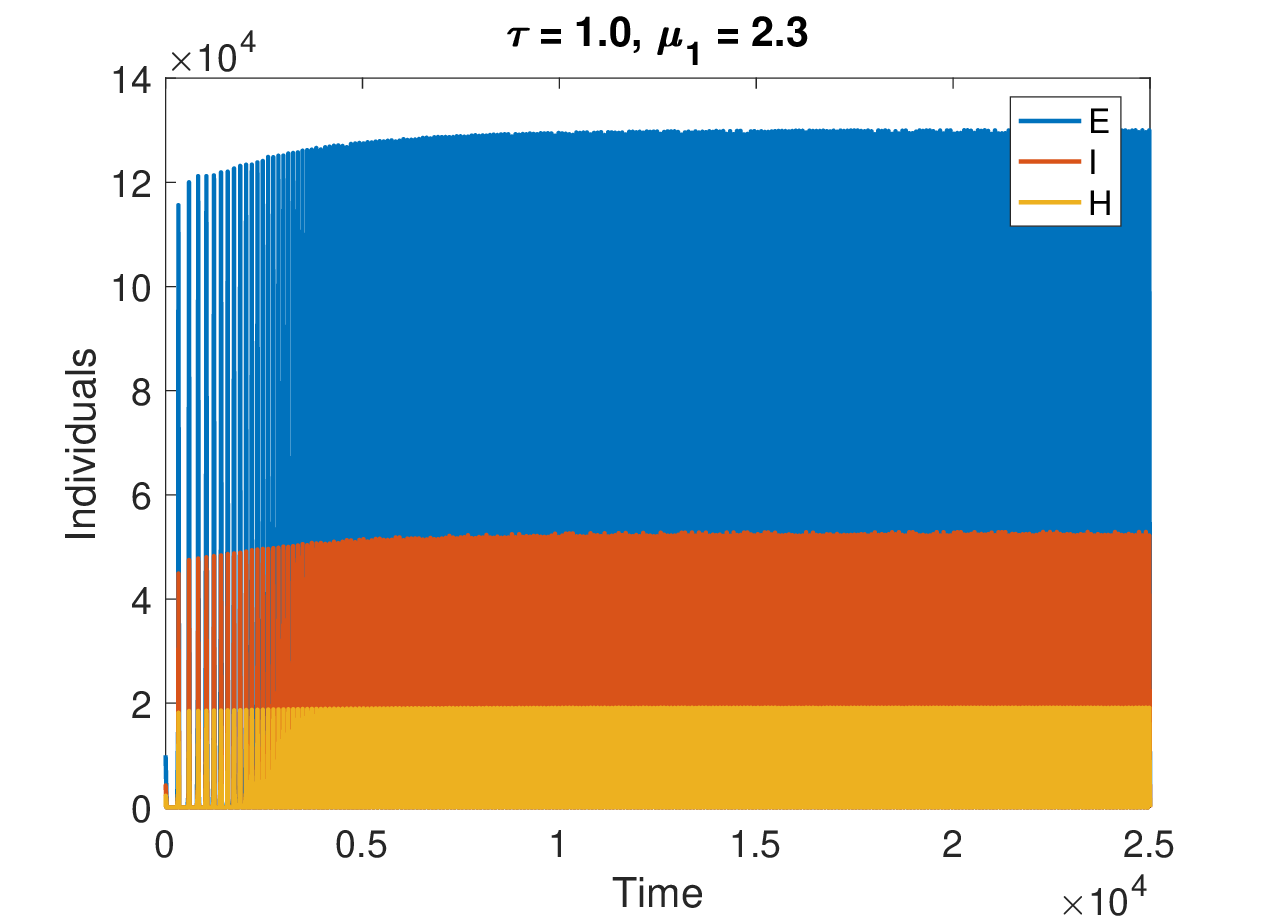}\\
		(b)\includegraphics[scale=0.3]{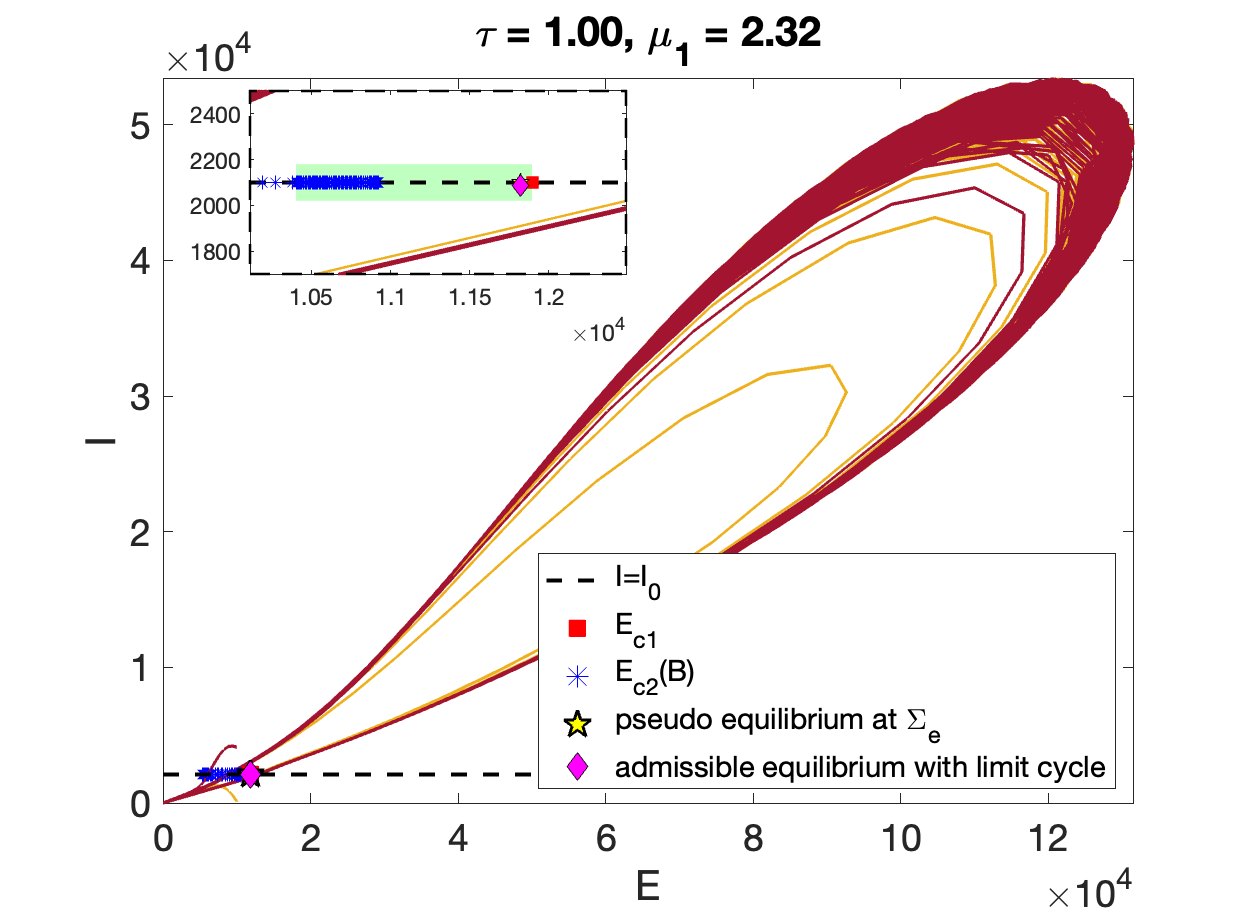}
		\includegraphics[scale=0.3]{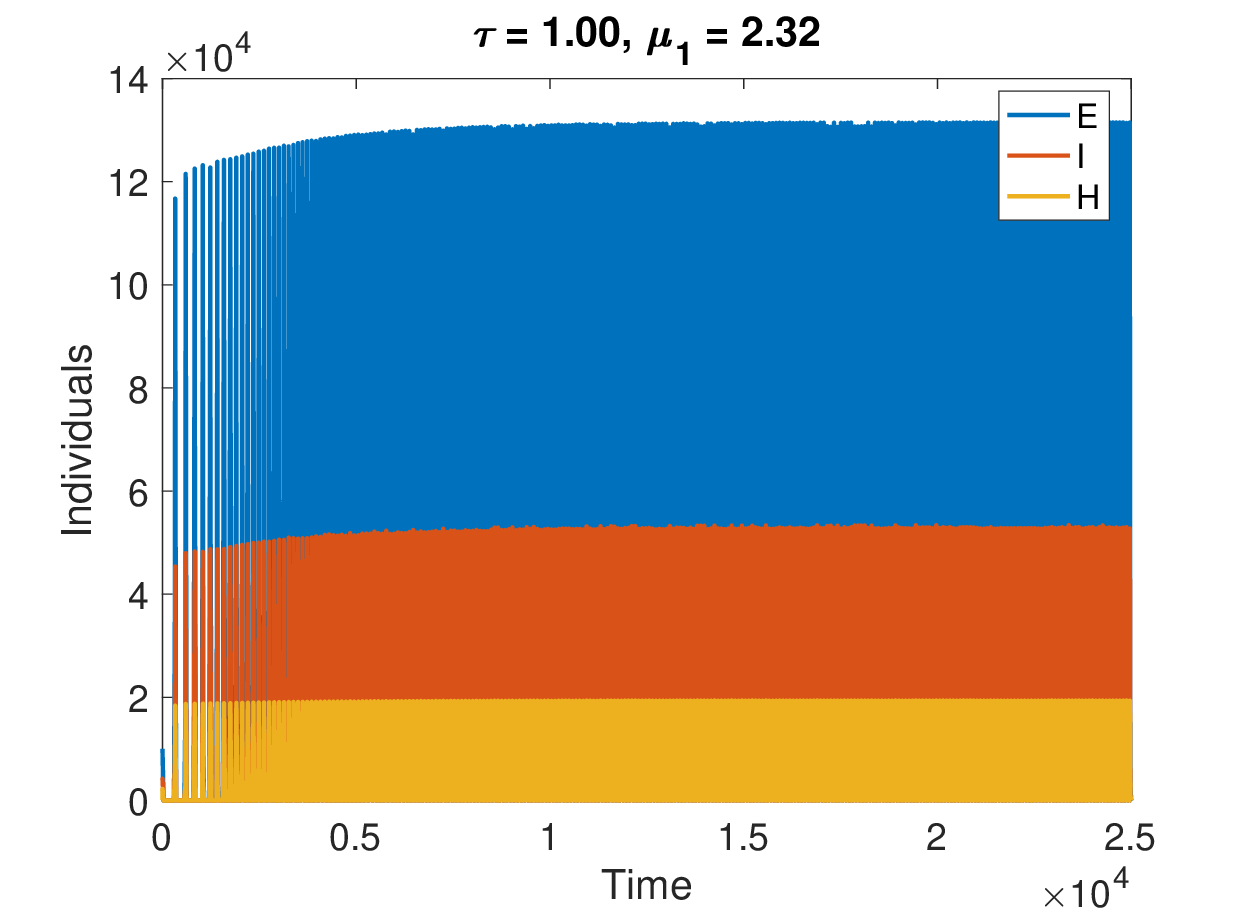}\\
		(c)\includegraphics[scale=0.3]{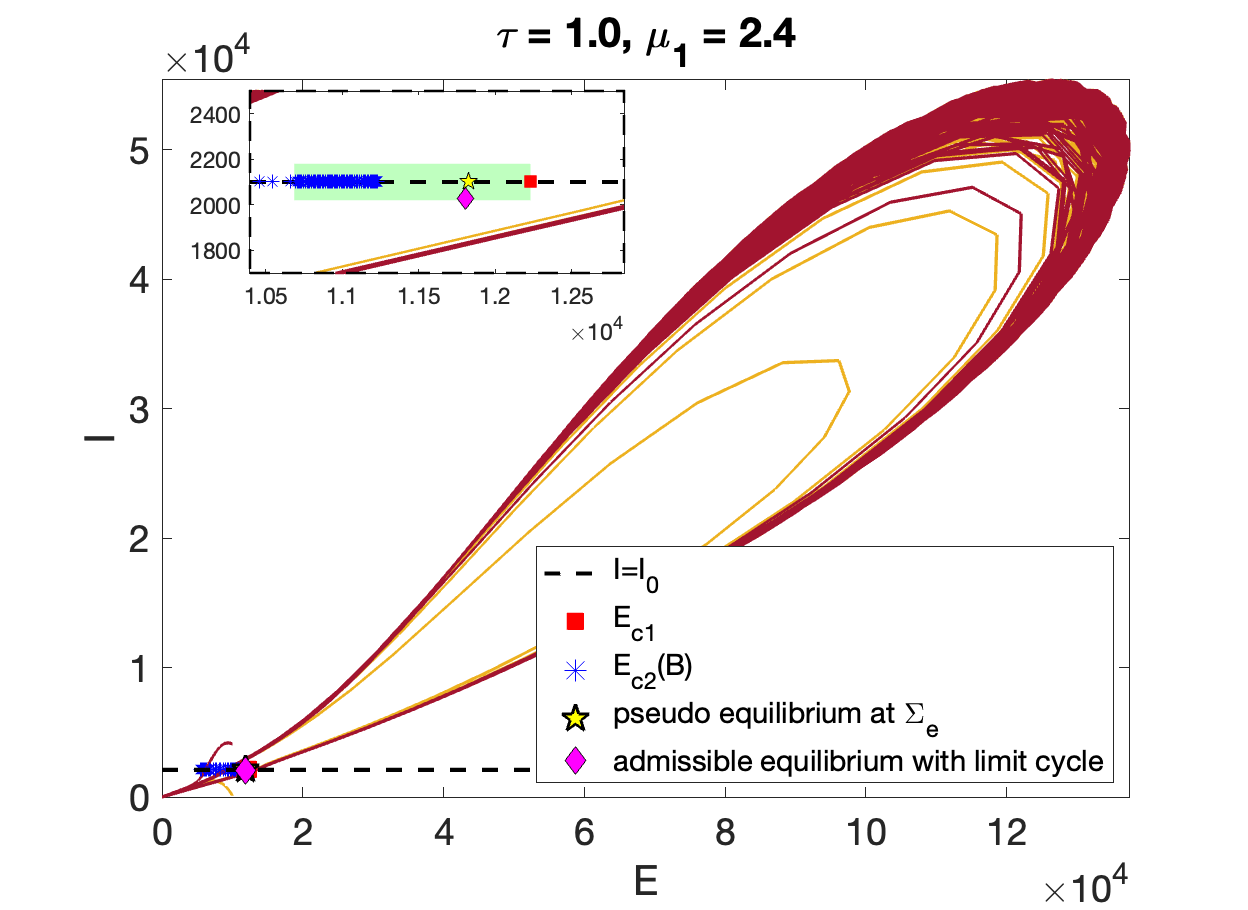}
		\includegraphics[scale=0.3]{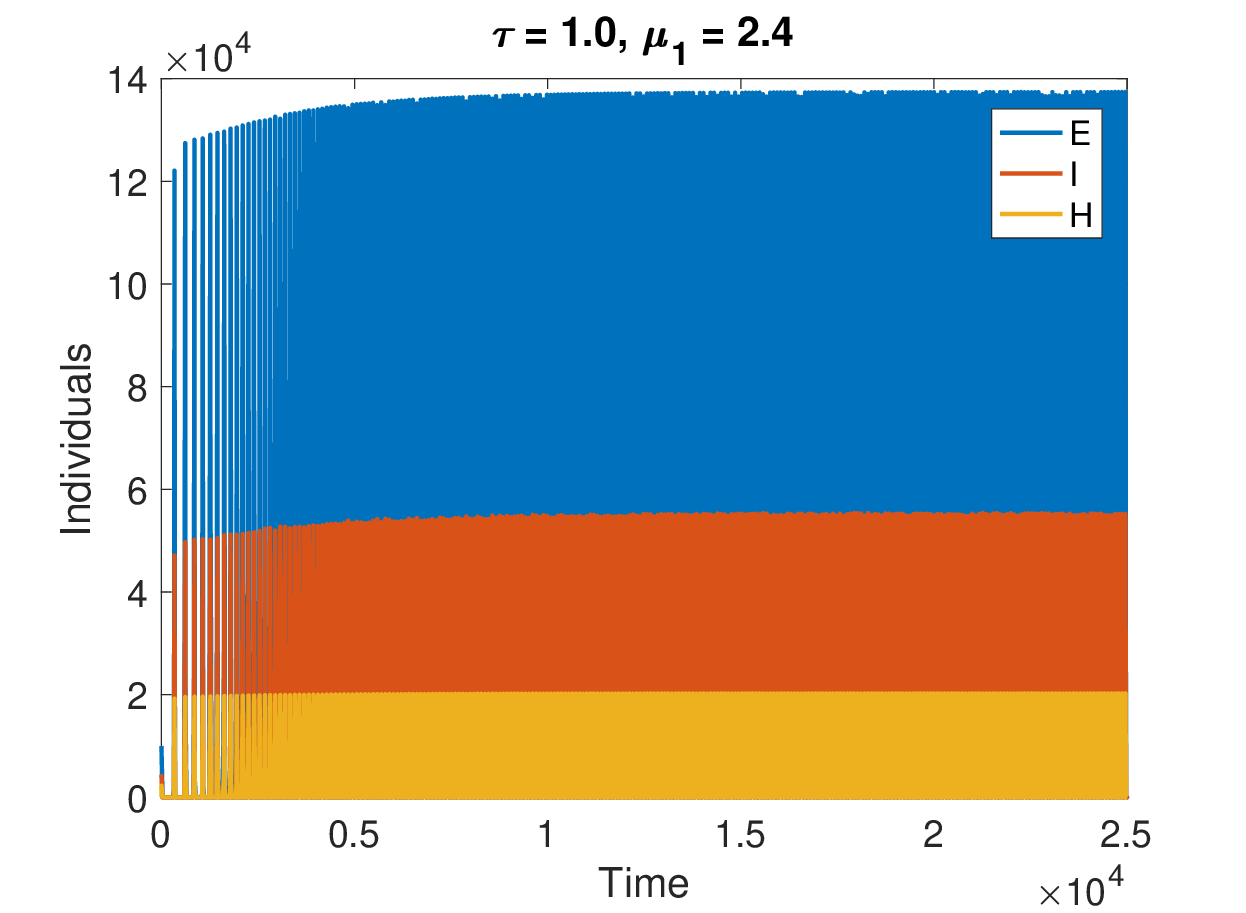}\\
		(d)\includegraphics[scale=0.3]{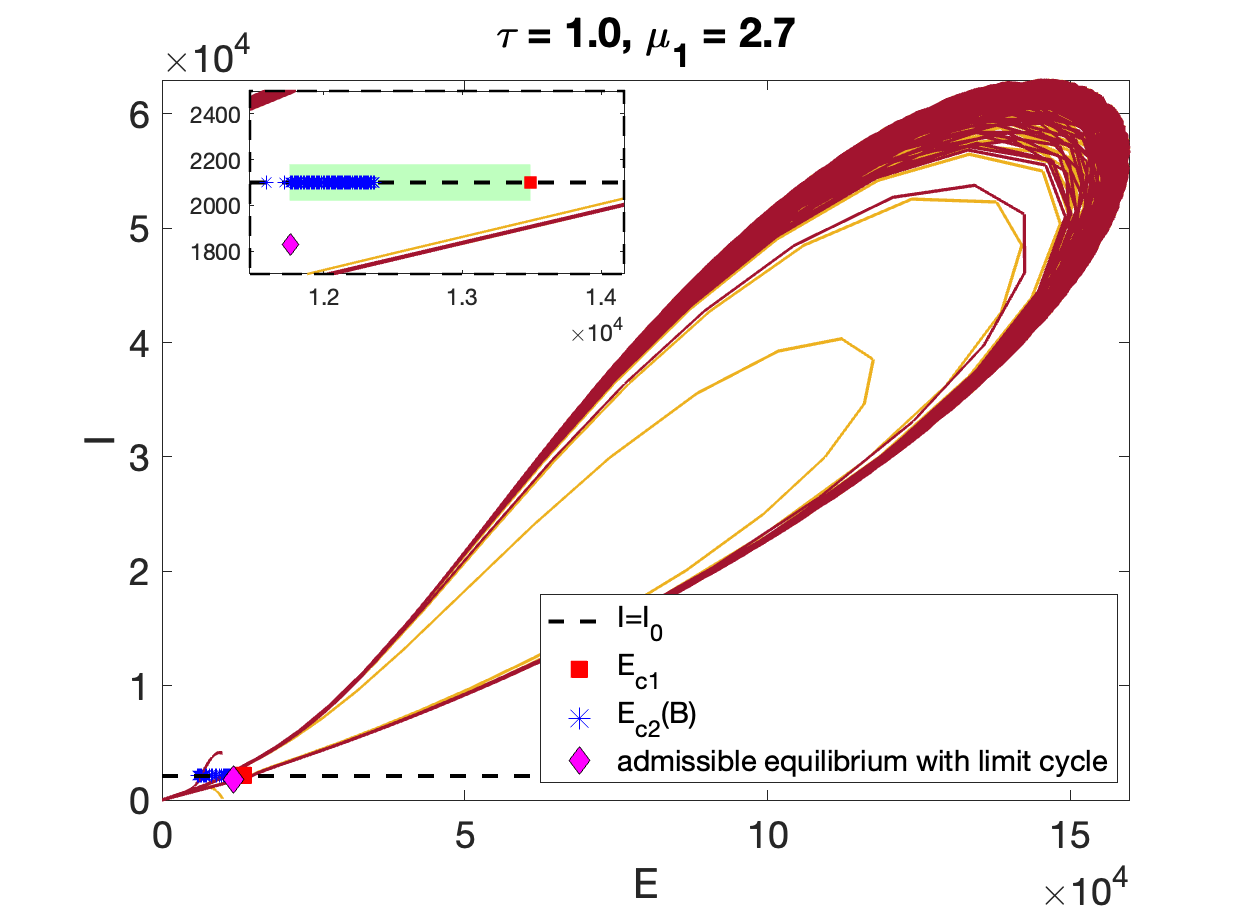}
		\includegraphics[scale=0.3]{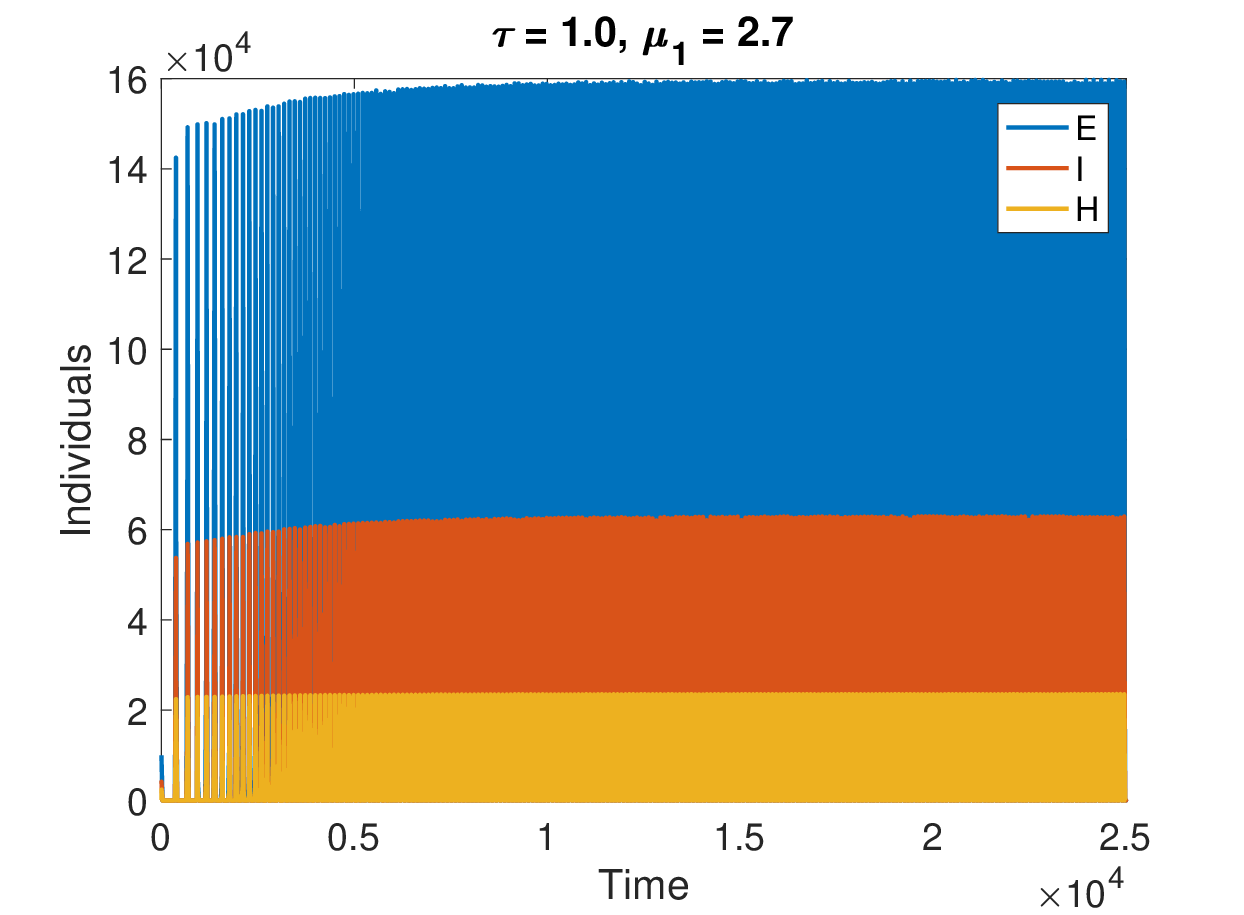}
		\caption{Phase portraits in the $(E,I)$-plane (left panels) and corresponding
			time series (right panels) of the Filippov system for varying maximum
			hospitalization rate $\mu_1$ with fixed time delay $\tau=1$.
			Here $\Lambda=1400$, $\xi=0.001$, $\beta=1.2\times10^{-5}$,
			$K=11500$, $\lambda=0.93$, $\omega=5000$, and $I_0=2100$;
			the remaining parameters are as in Table~\ref{param_val}.
			Panel (a) shows oscillatory dynamics around an admissible equilibrium
			of subsystem $F_{\Sigma_2}$ for $\mu_1=2.3$.
			Panels (b) and (c) illustrate a transition region in which the equilibria of $F_{\Sigma_2}$ and $F_{\Sigma_1}$ are both admissible, accompanied by the appearance of a
			pseudo-equilibrium on the switching manifold $\Sigma_e$, for $\mu_1=2.32$
			and $\mu_1=2.4$, respectively.
			Panel (d) shows that the pseudo-equilibrium loses admissibility and the
			dynamics exhibit sustained oscillations around an admissible equilibrium of
			$F_{\Sigma_1}$ for $\mu_1=2.7$.
		} 
		\label{fig:Fig._vary_mu1}
	\end{figure}

	For $\mu_1=2.3$, sustained oscillations occur around the
	admissible equilibrium of $F_{\Sigma_2}$
	(Fig.~\ref{fig:Fig._vary_mu1}(a)). This equilibrium is
	unstable for every $\tau\geq0$, having lost stability through
	a Hopf bifurcation of the undelayed subsystem at
	$\mu_1\approx0.975$. As $\mu_1$ increases
	(Figs.~\ref{fig:Fig._vary_mu1}(b) and (c)), both smooth
	endemic equilibria become admissible, together with a
	pseudo-equilibrium on the repelling sliding manifold
	$\Sigma_e$, while the oscillations persist. For larger
	$\mu_1$ (Fig.~\ref{fig:Fig._vary_mu1}(d)), the
	pseudo-equilibrium becomes inadmissible, and the oscillations
	occur around the admissible equilibrium of $F_{\Sigma_1}$.
	These observations are consistent with sequential Hopf and
	boundary-transition dynamics, whereby oscillations originating
	in the smooth above-threshold subsystem persist beyond the
	boundary-equilibrium transition through repeated crossings
	of the switching manifold. 
	
\begin{figure}[H]
		\centering
		\includegraphics[scale=0.3]{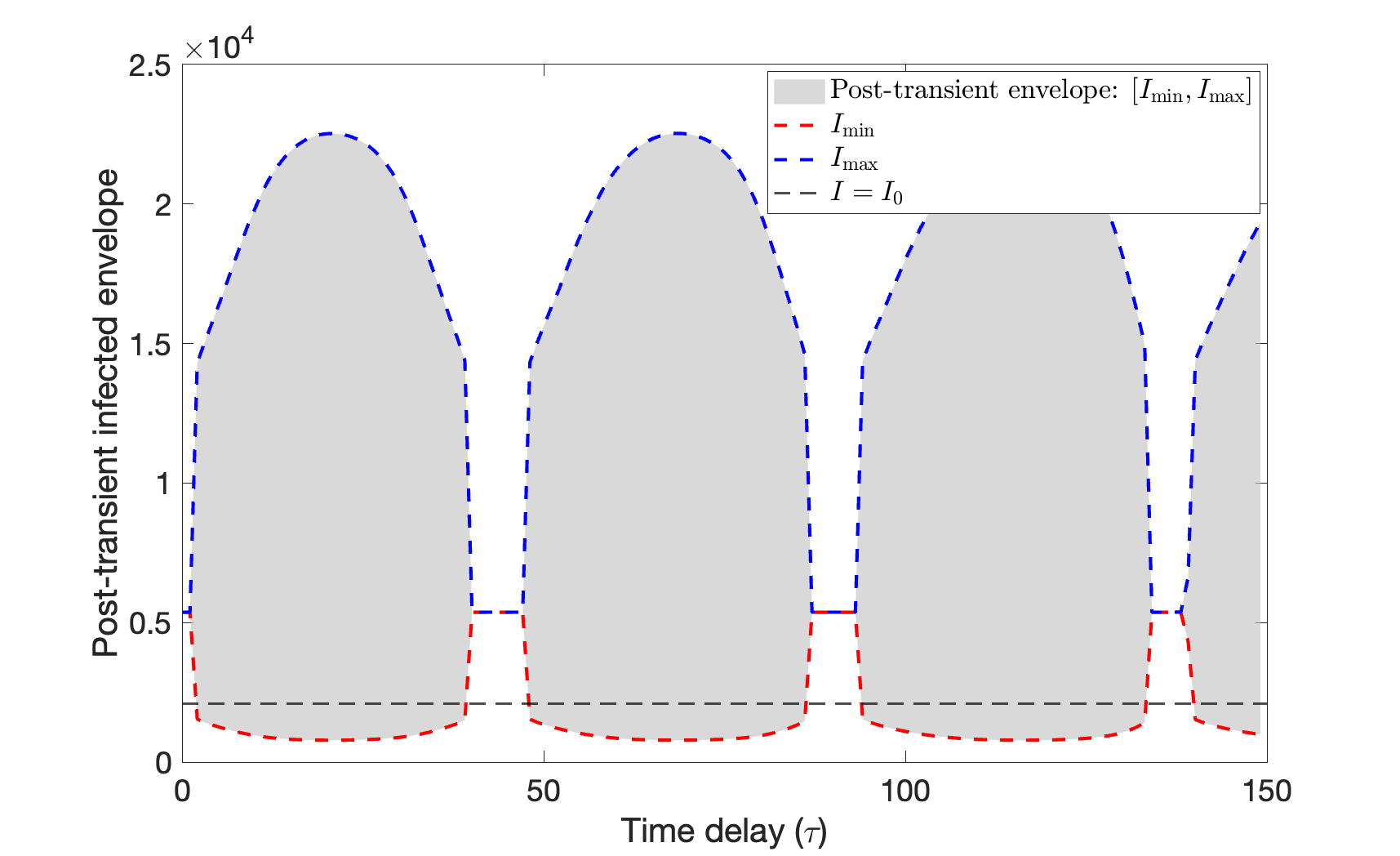}
		\caption{Post-transient infected envelope as a function of the time delay $\tau$. The dashed blue and red curves show the maximum  ($I_{\max}$) and minimum ($I_{\min}$) infected levels, while the shaded region denotes the envelope $[I_{\min}, I_{\max}]$ of the long-time dynamics over multiple initial conditions. The horizontal dashed line indicates the switching threshold $I=I_0$.}  
		\label{fig:non-smooth_hopf} 
	\end{figure}

	Fig.~\ref{fig:non-smooth_hopf} further summarizes the dependence of the
	long-term dynamics on the delay $\tau$. For each value of $\tau$,
	long-time simulations were performed and an initial transient was
	discarded before recording the post-transient minimum and maximum values
	of $I(t)$. Collapse of the envelope,
	$I_{\min}=I_{\max}$, corresponds to convergence to a steady infected
	level, whereas a nonzero separation between $I_{\min}$ and $I_{\max}$
	indicates persistent oscillations. The repeated widening and collapse of
	the envelope therefore illustrate delay-induced transitions between
	steady and oscillatory endemic dynamics. No loss of equilibrium
	admissibility or sliding-induced instability is observed over this
	parameter sweep, indicating that these oscillations arise from the
	delay-induced stability mechanism rather than from the repelling sliding
	dynamics.
	
	\section{Discussion and conclusion}
	\label{CONCLUSION}
	
	Sufficient availability of medical resources is widely recognized as a
	critical factor in controlling infectious disease outbreaks
	\cite{SHAN2014,SHAN2016,YU2023,WANG2018,QIN2016}. Motivated by this
	observation, a growing body of epidemic models has incorporated
	threshold-based intervention strategies to represent treatment responses
	under limited healthcare capacity, often leading to Filippov-type systems.
	However, comparatively little attention has been paid to the dynamical
	consequences of allowing healthcare resources to evolve with the epidemic
	state, particularly when their adjustment is delayed. In this work, we
	addressed this gap by developing a time-delayed Filippov epidemic model
	that couples threshold-based hospitalization with dynamically evolving
	effective healthcare capacity. This formulation allows the effects of
	non-smooth switching, adaptive healthcare resources, and delayed capacity
	response to be examined within a unified framework.
	
	We first analyzed the two smooth subsystems corresponding to infection
	levels below and above the hospitalization threshold. In the
	below-threshold subsystem $F_{\Sigma_1}$, the disease-free equilibrium
	$D_{02}$ is globally asymptotically stable in the smooth extension for
	$\mathcal R_0<1$ under the conditions of Theorem~\ref{thm:dfe_global}, while
	a forward transcritical bifurcation occurs at $\mathcal R_0=1$. For
	$\mathcal R_0>1$, a conditional global stability criterion for the positive-capacity endemic
	equilibrium, relative to infection-positive trajectories remaining below the threshold, is given in Theorem~\ref{thm_global_stab}. The healthcare-capacity response delay does not
	alter these equilibrium and stability properties because the capacity
	dynamics do not feed back into the epidemiological subsystem in this
	regime.
	
	Subsystem $F_{\Sigma_2}$ exhibits richer dynamics because hospitalization
	depends explicitly on both infection burden and healthcare capacity. In particular, the endemic-equilibrium
	equation of its smooth extension can undergo a nondegenerate fold at parameter values satisfying
	$\mathcal R_0<1$. Since the disease-free equilibrium at which
	$\mathcal R_0$ is defined belongs to the below-threshold regime, this
	phenomenon is more appropriately interpreted as a subthreshold fold
	rather than a classical backward bifurcation. When a locally stable
	positive endemic branch is biologically feasible and satisfies the
	admissibility condition $I_1^*>I_0$, it may coexist with the locally
	stable disease-free equilibrium $D_{02}$ for $\mathcal R_0<1$,
	producing discontinuity-induced subthreshold bistability in the full
	Filippov system. Thus, $\mathcal R_0<1$ ensures local stability of the
	positive-capacity disease-free equilibrium but does not necessarily exclude an admissible
	stable endemic state in the above-threshold regime. The subsystem $F_{\Sigma_2}$ also undergoes Hopf bifurcations in both the
	delay-free and delayed settings, with the latter giving rise to
	delay-induced stability switching.
	
	A central structural result of the model concerns the dynamics on the
	switching manifold $\Sigma=\{I=I_0\}$. The downward jump in the per
	capita hospitalization rate at the threshold structurally excludes
	attracting sliding: whenever a sliding region exists, it is necessarily
	repelling. Moreover, the lower tangency boundary depends explicitly on
	the evolving healthcare-capacity variable $B(t)$, so the width of the
	repelling sliding region varies with the current resource state. Thus,
	although the switching manifold itself is fixed, the local non-smooth
	geometry experienced by trajectories is state dependent. Epidemiologically,
	this implies that changes in effective healthcare capacity can alter the
	range of states near the hospitalization threshold for which repelling
	sliding behavior is possible.
	
	The sliding dynamics also admit pseudo-equilibria associated with zero
	and positive healthcare capacity. When admissible, the positive-capacity
	pseudo-equilibrium is globally asymptotically stable with respect to the
	algebraically extended reduced tangential dynamics for $B>0$, whereas the
	zero-capacity pseudo-equilibrium, when admissible, is unstable. Because the sliding region
	is transversely repelling, however, tangential attraction does not imply
	attraction in the full Filippov system. The positive-capacity
	pseudo-equilibrium therefore has a saddle-type organizing role, attracting
	the reduced tangential dynamics while transverse perturbations are
	directed away from the switching manifold. The analysis further
	identified tangency points, boundary equilibria, regular equilibria, and
	virtual equilibria, thereby clarifying how changes in parameter values
	reorganize the admissibility of smooth and sliding equilibria. 
	
	To obtain a regime-independent condition for disease elimination, we
	introduced the worst-case reproduction number
	$\mathcal R_0^{\max}$. The condition
	$\mathcal R_0^{\max}<1$ guarantees elimination across both
	hospitalization regimes and independently of the current healthcare
	capacity and the capacity-response delay. This condition is sufficient
	rather than necessary, and therefore complements the local
	disease-invasion threshold $\mathcal R_0$ by providing a stronger
	criterion for uniform elimination in the full Filippov system.
	
	An illustrative PINN-based calibration to cumulative reported COVID-19 cases from Kenya reproduces the overall cumulative-case trajectory. The calibrated parameters are effective rather than uniquely identifiable values, and the subsequent analyses explore broader parameter regimes to characterize the qualitative dynamics admitted by the model.
	
	The numerical analysis demonstrates that variations in baseline effective
	healthcare capacity $K$ can produce nonmonotone changes in the long-term
	dynamics. For intermediate values of $K$, sustained crossing cycles occur near quadratic tangencies, whereas other values of $K$ lead to
	convergence toward an admissible endemic equilibrium. These results show
	that increasing healthcare capacity does not necessarily produce a
	monotone change in endemic-equilibrium stability; instead, the feedback
	between infection, capacity, and hospitalization can shift the system
	between steady endemicity and recurrent oscillatory behavior.
	
	Variation in the maximum hospitalization rate $\mu_1$ produces an
	additional non-smooth reorganization. As $\mu_1$ changes, an endemic
	equilibrium of the smooth subsystem reaches the switching manifold and
	changes between admissible and virtual states, giving rise to a
	boundary-equilibrium transition. Over an intermediate parameter range, the positive-capacity
	pseudo-equilibrium is admissible within the repelling sliding region. Importantly, $\mu_1$ changes the tangency
	boundaries rather than the algebraic location of the pseudo-equilibrium,
	and therefore controls its admissibility within the sliding dynamics.
	For sufficiently large $\mu_1$ in the delay-free parameter sweep, the simulated trajectories converge to
	the admissible equilibrium of $F_{\Sigma_1}$, illustrating how changes in
	hospitalization intensity can reorganize the long-term Filippov dynamics. 
	
	A further contribution is the identification of a codimension-two
	Hopf-boundary interaction in the $(\mu_1,K)$-plane in the delay-free setting. The
	boundary-equilibrium locus and the Hopf locus of the smooth extension of
	$F_{\Sigma_2}$ intersect transversally at a parameter point where the
	endemic equilibrium simultaneously reaches the switching manifold and
	undergoes a nondegenerate supercritical Hopf bifurcation in the smooth
	extension. In addition,
	the critical Hopf mode has a nonzero component transverse to the switching
	manifold, and the associated small-amplitude periodic orbit possesses a
	local first-contact curve that emanates from the codimension-two point and
	is quadratically tangent to the Hopf locus. This establishes a local
	mechanism through which a smooth oscillatory instability approaches and
	interacts with the non-smooth threshold.
	
	From an epidemiological perspective, the Hopf-boundary interaction
	represents a critical balance between the maximum hospitalization rate
	and baseline effective healthcare capacity. Near this interaction,
	relatively small changes in healthcare response can simultaneously alter
	the stability of the endemic equilibrium, its admissibility with respect
	to the treatment threshold, and the interaction of emerging oscillations
	with the switching manifold. This local geometry provides a mechanism for emerging oscillations to
	reach the intervention threshold. 
	
	The healthcare-capacity response delay further enriches the dynamics.
	For the parameter sets considered, as the delay increases, the endemic
	equilibrium repeatedly loses and regains stability, producing alternating intervals of steady and
	oscillatory endemic behavior. Moreover, oscillations generated while an
	equilibrium is admissible in $F_{\Sigma_2}$ may persist after a subsequent
	boundary-equilibrium transition through repeated crossings of the
	switching manifold. These results show that delayed healthcare response
	and non-smooth switching can interact to produce recurrent epidemic
	dynamics that are not captured by either smooth subsystem alone.
	
	Overall, the analysis shows that coupling threshold-based hospitalization
	with evolving healthcare capacity and delayed response generates dynamical
	behavior in which equilibrium stability, admissibility, and sliding geometry
	are jointly influenced by resource dynamics and threshold switching. These results extend existing limited-resource and threshold-based
	epidemic models \cite{SHAN2014,WANG2018,MU2019,ZHAO2020,QIN2016} by
	showing that dynamically evolving healthcare resources can influence not
	only endemic infection levels but also the geometry and bifurcation
	structure associated with the intervention threshold.
	
	Despite the analytical and numerical insights obtained in this work, several limitations remain. The delay enters only through the
	healthcare-capacity dynamics, and the non-smooth analysis relies on the
	standard Filippov framework. Although the local first-contact geometry is established analytically
	in the delay-free setting and supported numerically, a complete continuation and stability analysis of
	post-contact crossing periodic orbits, particularly in the delayed
	setting, would require additional non-smooth DDE and discontinuity-map
	techniques. From a modeling perspective, the healthcare-capacity expansion
	mechanism and hospitalization function provide idealized representations
	of operational processes that may also depend on staffing, intensive-care
	availability, discharge rates, and triage policies. The use
	of a fixed infection threshold $I_0$ also excludes possible hysteresis or
	adaptive policy responses. Finally, the PINN calibration based on
	cumulative reported cases does not establish structural or practical
	identifiability of the calibrated parameters. These issues provide
	natural directions for further development of non-smooth epidemic models
	with dynamic healthcare resources and delayed response.

	\section*{Conflict of interest}
	The authors declare no conflict of interest. 
	
	\section*{Acknowledgments} 
	
	Sarita Bugalia acknowledges financial support from the Council of Scientific and Industrial Research (CSIR), India, under File No.~09/1131(0025)/2018-EMR-I.

	\section*{Data availability}
	
	The cumulative reported COVID-19 case data used in this study are
	publicly available through the World Health Organization COVID-19
	dashboard (\url{https://data.who.int/dashboards/covid19/data}).

\end{document}


\maketitle
	\section{Model and notation}\label{context:supp:model}
	This supplementary material provides supporting calculations, proofs, and numerical illustrations for the delayed Filippov epidemic model. 
	
	The human population is $N(t)=S(t)+E(t)+I(t)+H(t)+R(t)$, where $S$, $E$, $I$, $H$, and $R$ denote susceptible, exposed, infectious, hospitalized, and recovered individuals, respectively, and $B$ denotes effective healthcare capacity.
	
	The hospitalization law is given by
	\begin{equation*}
		T(I,B) =
		\left\{
		\begin{array}{ll}
			\mu_1, & I\leq I_0, \\[2mm]
			\displaystyle
			\mu_0+(\mu_1-\mu_0)\frac{B}{B+I},
			& I>I_0.
		\end{array}
		\right.
	\end{equation*}
	The delayed model is given by
	\begin{equation*}
		\begin{aligned}
			\frac{dS}{dt} &= \Lambda - \beta SI - \mu S + \xi R, \\
			\frac{dE}{dt} &= \beta SI - (\mu + a)E, \\
			\frac{dI}{dt} &= aE - (\mu + \delta_1 + \gamma_1)I - T(I,B)I, \\
			\frac{dH}{dt} &= T(I,B)I - (\mu + \delta_2 + \gamma_2)H, \\
			\frac{dR}{dt} &= \gamma_1 I + \gamma_2 H - (\mu + \xi)R, \\
			\frac{dB}{dt} &= rB\left(1-\frac{B}{K}\right) + \frac{\lambda I(t-\tau)B}{\omega + I(t-\tau)},
		\end{aligned}
	\end{equation*}
	The switching manifold is $\Sigma=\{I=I_0\}$, with smooth subsystems $F_{\Sigma_1}$ and $F_{\Sigma_2}$ corresponding to $I<I_0$ and $I>I_0$, respectively. At the switching manifold, solutions are interpreted in the Filippov sense, using the closed convex hull of the two one-sided vector fields. The history is continuous and nonnegative on $[-\tau,0]$.
	
	The disease-free equilibria are
	\[
	D_{01}=(\Lambda/\mu,0,0,0,0,0),\quad D_{02}=(\Lambda/\mu,0,0,0,0,K),
	\]
	and the basic reproduction number is
	\[
	\mathcal R_0=\frac{\beta\Lambda a}{\mu(\mu+a)(\mu+\delta_1+\gamma_1+\mu_1)}.
	\]
	The positive-capacity endemic equilibria of the two smooth subsystems are denoted by $D^*$ and $D_1^*$, respectively. The notation $D^*$ is also used locally for the above-threshold endemic equilibrium in stability calculations, as in the main manuscript. The algebraic pseudo-equilibria on $I=I_0$ are denoted by $D^{0s}$ and $D^s$, with zero and positive capacity, respectively.
	
	The biological interpretations and baseline parameter values are given in Tables~\ref{context:parameter_definitions} and \ref{context:parameter_values}, respectively.

	\begin{table}[H]
		\centering
		\caption{Biological interpretations of model parameters.}
		\label{context:parameter_definitions}
		\begin{tabular}{|p{2cm}|p{13cm}|}
			\hline
			Parameters & Biological interpretations \\
			\hline
			$\Lambda$ & Recruitment rate into the susceptible population \\
			$\beta$ & Disease transmission coefficient \\
			$1/\mu$ & Mean lifetime in the absence of disease-induced mortality \\
			$1/a$ & Mean latent period in the absence of natural mortality \\ 
			$\delta_1$ ($\delta_2$) & Disease-induced death rate of infectious (hospitalized) individuals \\
			$\gamma_1$ ($\gamma_2$) & Recovery rate of infectious (hospitalized) individuals \\
			$1/\xi$ & Mean duration of immunity in the absence of natural mortality \\ 
			$r$ & Intrinsic adjustment rate of effective healthcare capacity \\
			$\lambda$ & Maximal per capita infection-responsive expansion rate
			of effective healthcare capacity \\ 
			$K$ & Baseline effective healthcare capacity in the absence of
			infection-driven expansion \\ 
			$\omega$ & Half-saturation constant for infection-responsive capacity expansion \\ 
			$\mu_1$ ($\mu_0$) & Maximum (minimum) per capita hospitalization rate \\ 
			$I_0$ & Hospitalization threshold \\
			$\tau$ & Healthcare-capacity response delay to reported infection information \\  
			\hline
		\end{tabular}
	\end{table}

	\begin{table}[H]
		\centering
		\caption{Fixed and calibrated parameter values used in the 
			PINN-based fitting to cumulative reported COVID-19 cases for Kenya.}
		\label{context:parameter_values}
		\begin{tabular}{|p{2cm}|p{8cm}|p{2cm}|}
			\hline
			Parameters & Numerical values & References \\
			\hline
			
			$a$ & $1/2$ week$^{-1}$ & \cite{WANG2022} \\
			
			$\delta_1$ & $0.0119$ week$^{-1}$ & \cite{ZHOU2020} \\
			
			$\delta_2$ & $0.0077$ week$^{-1}$ & \cite{ZHOU2020} \\
			
			$\gamma_1$ & $1/2$ week$^{-1}$ & \cite{WANG2022,NIT2021} \\
			
			$\gamma_2$ & $1.6758$ week$^{-1}$ & \cite{ZHOU2020} \\
			
			$\xi$ & $0.0389$ week$^{-1}$ & \cite{BUG2023,CHOW2022} \\
			
			$r$ & $1.3139$ week$^{-1}$ & \cite{ZHOU2020} \\
			
			$\omega$ & $50$ individuals & Assumed \\
			
			$\mu_1$ & $0.77$ week$^{-1}$ & \cite{Bugalia2020} \\
			
			$\mu_0$ & $0.0182$ week$^{-1}$ & \cite{AJB2021} \\
			
			$\mu$ & $0.00028$ week$^{-1}$ & \cite{KenyaLifeExpectancy} \\ 
			
			$K$ & $1000$ individual-equivalent capacity units & Assumed \\
			
			$I_0$ & $4500$ individuals & Assumed \\
			
			$\Lambda$ & $1.4621\times10^{4}$ individuals week$^{-1}$ 
			& Computed \\
			
			$\beta$ & $1.2648\times10^{-8}$ individual$^{-1}$ week$^{-1}$ 
			& Calibrated \\
			
			$\lambda$ & $0.8583$ week$^{-1}$ & Calibrated \\
			
			$\tau$ & $0.9065$ weeks & Calibrated \\
			
			\hline
		\end{tabular}
	\end{table}

	\subsection{One-sided vector fields}
	Let $X=(S,E,I,H,R,B)$. The vector fields $F_1$ and $F_2$ are obtained from the governing equations above by substituting $T=\mu_1$ and $T=\mu_0+(\mu_1-\mu_0)B/(B+I)$, respectively. Switching depends on the current state, while delay enters only through the capacity equation.
	Then system (2.5) with (2.4) is interpreted as the Filippov differential inclusion 
	\begin{equation}\label{context:filisystem}
		\frac{dX}{dt}\in
		\mathcal{F}(X(t),X(t-\tau))
		=
		\begin{cases}
			F_1(X(t),X(t-\tau)), & I<I_0,\\[2mm]
			\operatorname{co}\left\{
			F_1(X(t),X(t-\tau)),
			F_2(X(t),X(t-\tau))
			\right\}, & I=I_0,\\[2mm]
			F_2(X(t),X(t-\tau)), & I>I_0,
		\end{cases}
	\end{equation}
	where $\operatorname{co}\{\cdot\}$ denotes the closed convex hull of the two one-sided vector fields. Thus, a Filippov solution is a continuous function $X(t)$ agreeing with the prescribed history on $[-\tau,0]$, locally absolutely continuous for $t\geq0$ on its interval of existence, and satisfying \eqref{context:filisystem} for almost every $t>0$ in that interval.
	
	The smooth regions are $\Sigma_1=\{X\in\mathbb R_+^6:I<I_0\}$ and $\Sigma_2=\{X\in\mathbb R_+^6:I>I_0\}$. We denote their vector fields by $F_{\Sigma_i}$, $i=1,2$, using their continuous one-sided extensions on $\Sigma$.
	
	\section{Positivity, boundedness and positive invariance}\label{Appndx_A}
	\noindent\textbf{Statement of Theorem~2.1 (main manuscript).} Let $(S(t),E(t),I(t),H(t),R(t),B(t))$ be any maximal Filippov solution of system (2.5) corresponding to a non-negative continuous initial history. Then the solution is defined for all $t\geq0$ and remains non-negative. Moreover, 
	\[
	N(t)\le
	\max\left\{N(0),\frac{\Lambda}{\mu}\right\},
	\]
	and
	\[
	B(t)\le
	\max\left\{
	B(0),
	K\left(1+\frac{\lambda}{r}\right)
	\right\},
	\quad t\ge0.
	\]
	If $B(0)>0$, then $B(t)>0$ for all $t\ge0$. In particular, the set
	\begin{equation}
		\Omega=\left\{(S,E,I,H,R,B)\in\mathbb{R}_+^6 \,\Big|\,
		0\le S+E+I+H+R\le\frac{\Lambda}{\mu},
		\; 0\le B\le K\left(1+\frac{\lambda}{r}\right)\right\}
	\end{equation}
	is positively invariant in the current-state space for every Filippov solution corresponding to a non-negative continuous initial history whose initial state belongs to $\Omega$.
	
	Let
	$
	X(t)=(S(t),E(t),I(t),H(t),R(t),B(t))
	$ 
	be any Filippov solution of system (2.5) corresponding to a non-negative continuous initial history. On the switching manifold $I=I_0$, the Filippov vector field belongs to the closed convex hull of the two one-sided vector fields $F_1$ and $F_2$. Therefore, it is sufficient to verify that both one-sided vector fields point into, or are tangent to, the non-negative orthant whenever one of the state variables reaches zero.
	
	For the epidemiological variables, we have
	
	\begin{align*}
		\frac{dS}{dt}\Big|_{S=0}
		&=\Lambda+\xi R>0,\\
		\frac{dE}{dt}\Big|_{E=0}
		&=\beta SI\ge0,\\
		\frac{dI}{dt}\Big|_{I=0}
		&=aE\ge0,\\
		\frac{dH}{dt}\Big|_{H=0}
		&=T(I,B)I\ge0,\\
		\frac{dR}{dt}\Big|_{R=0}
		&=\gamma_1I+\gamma_2H\ge0.
	\end{align*}
	Moreover, for the healthcare-capacity variable,
	\[
	\frac{dB}{dt}\Big|_{B=0}=0.
	\]
	These inequalities hold for both one-sided vector fields and hence also for every convex combination defining the Filippov regularization on $\Sigma$. Consequently, the non-negative orthant $\mathbb{R}_+^6$ is positively invariant, and every Filippov solution that starts with a non-negative initial history remains non-negative for all times for which it exists.
	
	We next establish boundedness of the total human population. Let
	$
	N(t)=S(t)+E(t)+I(t)+H(t)+R(t).
	$ 
	Adding the first five equations of system (2.5) gives
	\begin{equation}\label{eq:N_bound_app}
		\frac{dN}{dt}
		=
		\Lambda-\mu N-\delta_1I-\delta_2H
		\le
		\Lambda-\mu N.
	\end{equation}
	By the comparison principle,
	\[
	N(t)
	\le
	\frac{\Lambda}{\mu}
	+
	\left(
	N(0)-\frac{\Lambda}{\mu}
	\right)e^{-\mu t}.
	\]
	Therefore,
	\begin{equation}\label{eq:N_max_app}
		N(t)
		\le
		\max\left\{
		N(0),\frac{\Lambda}{\mu}
		\right\},
		\quad t\ge0.
	\end{equation}
	In particular, if
	\[
	N(0)\le\frac{\Lambda}{\mu},
	\]
	then
	\[
	N(t)\le\frac{\Lambda}{\mu}
	\quad \text{for all }t\ge0.
	\]
	
	We now consider the healthcare-capacity equation. Since the initial history of $I$ is non-negative and the non-negative orthant is  positively invariant, 
	$
	I(t-\tau)\ge0.
	$ 
	Hence,
	\[
	0
	\le
	\frac{I(t-\tau)}
	{\omega+I(t-\tau)}
	\le1.
	\]
	It follows that
	\begin{align}
		\frac{dB}{dt}
		&=
		rB\left(1-\frac{B}{K}\right)
		+
		\frac{\lambda I(t-\tau)B}
		{\omega+I(t-\tau)}
		\nonumber\\
		&\le
		(r+\lambda)B-\frac{r}{K}B^2.
		\label{eq:B_comparison_app}
	\end{align}
	Consider the comparison equation
	\[
	\frac{dy}{dt}
	=
	(r+\lambda)y-\frac{r}{K}y^2.
	\]
	Its positive equilibrium is
	\[
	y^*
	=
	\frac{K(r+\lambda)}{r}
	=
	K\left(1+\frac{\lambda}{r}\right).
	\]
	Therefore, by comparison,
	\begin{equation}\label{eq:B_bound_app}
		B(t)
		\le
		\max\left\{
		B(0),
		K\left(1+\frac{\lambda}{r}\right)
		\right\},
		\quad t\ge0.
	\end{equation}
	
	Furthermore, the $B$-equation can be written as
	\[
	\frac{dB}{dt}
	=
	B(t)
	\left[
	r\left(1-\frac{B(t)}{K}\right)
	+
	\frac{\lambda I(t-\tau)}
	{\omega+I(t-\tau)}
	\right].
	\]
	Therefore, whenever $B(0)>0$,
	\[
	B(t)
	=
	B(0)
	\exp\left\{
	\int_0^t
	\left[
	r\left(1-\frac{B(s)}{K}\right)
	+
	\frac{\lambda I(s-\tau)}
	{\omega+I(s-\tau)}
	\right]ds
	\right\}>0
	\]
	for all finite $t\ge0$. If $B(0)=0$, then $B(t)=0$ for all $t\ge0$. Thus, the hyperplane $\{B=0\}$ is invariant, while strictly positive healthcare capacity remains strictly positive. 
	
	It remains to verify the positive invariance of the set
	\[
	\Omega=
	\left\{
	(S,E,I,H,R,B)\in\mathbb{R}_+^6
	\,\Big|\,
	0\le N\le\frac{\Lambda}{\mu},
	\quad
	0\le B\le K\left(1+\frac{\lambda}{r}\right)
	\right\}.
	\]
	On the boundary
	$
	N=\frac{\Lambda}{\mu},
	$ 
	Eq. \eqref{eq:N_bound_app} gives
	\[
	\frac{dN}{dt}
	=
	-\delta_1I-\delta_2H
	\le0.
	\]
	Hence, the vector field cannot point outward through the population boundary. 
	
	Similarly, on the upper healthcare-capacity boundary 
	\[
	B=K\left(1+\frac{\lambda}{r}\right),
	\]
	using
	\[
	\frac{I(t-\tau)}{\omega+I(t-\tau)}\le1,
	\]
	we obtain
	\begin{align*}
		\frac{dB}{dt}
		&\le
		(r+\lambda)B-\frac{r}{K}B^2\\
		&=
		B\left[
		(r+\lambda)
		-
		\frac{r}{K}
		K\left(1+\frac{\lambda}{r}\right)
		\right]
		=0.
	\end{align*}
	At $B=0$, we have $dB/dt=0$. Together with the inward-pointing conditions on the coordinate boundaries established above, these inequalities show that $\Omega$ is positively invariant.
	
	Finally, estimates \eqref{eq:N_max_app} and \eqref{eq:B_bound_app} show that all components of any Filippov solution remain bounded on every finite time interval. Thus, no finite-time blow-up can occur, and any maximal Filippov solution can be continued for all $t\ge0$ under the adopted Filippov solution framework.
	
	Therefore, every Filippov solution corresponding to a non-negative initial history remains non-negative and bounded, and the set $\Omega$ is positively invariant whenever the initial state belongs to $\Omega$. This completes the proof.
	
	\section{Global stability of the disease-free equilibrium}\label{supp:dfeproof} 
	\noindent\textbf{Statement of Theorem~3.1 (main manuscript).} If $\mathcal R_0<1$, then the disease-free equilibrium $D_{02}$ is
	globally asymptotically stable for the smooth extension of subsystem $F_{\Sigma_1}$ with
	initial conditions in $\Omega$ satisfying $B(0)>0$.
	
	\begin{proof}
		Define
		\[
		L(E,I)=aE+(\mu+a)I.
		\]
		Then $L(E,I)\ge 0$, and $L(E,I)=0$ if and only if $E=I=0$. Differentiating $L$ with respect to $t$ along the solutions of $F_{\Sigma_1}$ yields 
		\begin{align*}
			\frac{dL}{dt}
			&= a\left(\beta S I-(\mu+a)E\right)
			+(\mu+a)\left(aE-(\mu+\delta_1+\gamma_1+\mu_1)I\right)\\
			&=\left(\beta a S-(\mu+a)(\mu+\delta_1+\gamma_1+\mu_1)\right)I .
		\end{align*}
		Since $S(t)\le \Lambda/\mu$ in the positively invariant region, we obtain
		\[
		\frac{dL}{dt}
		\le
		(\mu+a)(\mu+\delta_1+\gamma_1+\mu_1)(\mathcal R_0-1)I .
		\]
		Thus, if $\mathcal R_0<1$, then $\dot L\le 0$, and $\dot L=0$ only when $I=0$. On the largest invariant set contained in $\{I=0\}$, the equation $\dot I=aE-(\mu+\delta_1+\gamma_1+\mu_1)I $
		implies $E=0$. Hence, the largest invariant set contained in $\{\dot L=0\}$ satisfies $E=I=0$. By LaSalle's invariance principle \cite{LaSalle1976},  $(E(t),I(t))\to(0,0)$ as $t\to\infty$. 
		The remaining equations then imply
		$H(t)\to0, R(t)\to0,$ and $S(t)\to \frac{\Lambda}{\mu}.$  
		Since the delay appears only in the healthcare capacity equation and
		$I(t-\tau)\to0$ as $t\to\infty$, the healthcare-capacity  equation asymptotically reduces to the logistic equation
		\[
		\frac{dB}{dt} = rB\left(1-\frac{B}{K}\right). 
		\]
		For any $\varepsilon>0$, the infection-responsive per capita expansion term is eventually between $0$ and $\varepsilon$. Comparison with logistic equations with growth rates $r$ and $r+\varepsilon$ and common quadratic coefficient $r/K$ yields $K\leq\liminf B(t)\leq\limsup B(t)\leq K(1+\varepsilon/r)$ for $B(0)>0$. Letting $\varepsilon\to0^+$ gives $B(t)\to K$. Hence, all solutions in the positively invariant region with $B(0)>0$ converge to $D_{02}$. 
	\end{proof}
	\section{Forward transcritical bifurcation}\label{Appndx_B}
	\noindent\textbf{Statement of Theorem~3.2 (main manuscript).} The subsystem $F_{\Sigma_1}$ undergoes a forward transcritical bifurcation at the disease-free equilibrium $D_{02}$ when $\mathcal R_0 = 1$.
	
	We choose the transmission rate $\beta$ as the bifurcation parameter. Solving $\mathcal R_0 = 1$ yields the critical value
	\[
	\beta = \beta^*
	= \frac{\mu(\mu+a)(\mu+\delta_1+\gamma_1+\mu_1)}{\Lambda a}.
	\]
	
	Since $B$ does not feed back into the epidemiological subsystem of
	$F_{\Sigma_1}$, we apply the center manifold calculation to the autonomous
	five-dimensional subsystem $(S,E,I,H,R)$. At $\beta=\beta^*$ its Jacobian
	has a simple zero eigenvalue, while all remaining eigenvalues have negative
	real parts. Let $W=(w_1,w_2,w_3,w_4,w_5)^\top$ and
	$V=(v_1,v_2,v_3,v_4,v_5)$ denote the right and left eigenvectors
	associated with the zero eigenvalue, normalized so that $VW=1$.
	Choosing $w_3=1$ gives
	\begin{align*}
		w_1&=\frac{1}{\mu}\left[\frac{\xi}{\mu+\xi}
		\left(\gamma_1+\frac{\gamma_2\mu_1}{\mu+\delta_2+\gamma_2}\right)
		-\frac{(\mu+a)(\mu+\delta_1+\gamma_1+\mu_1)}{a}\right]<0,\\
		w_2&=\frac{\mu+\delta_1+\gamma_1+\mu_1}{a},\quad
		w_4=\frac{\mu_1}{\mu+\delta_2+\gamma_2},\\
		w_5&=\frac{1}{\mu+\xi}\left(\gamma_1+
		\frac{\gamma_2\mu_1}{\mu+\delta_2+\gamma_2}\right).
	\end{align*}
	Moreover, $v_1=v_4=v_5=0$, and
	\[
	v_2=\frac{a}{\mu+a}v_3,\quad
	v_3=\left(1+\frac{a w_2}{\mu+a}\right)^{-1}>0.
	\]
	Following Theorem 4.1 in \cite{CAS2004}, the bifurcation coefficients are
	\begin{align*}
		\mathfrak a&=\sum_{k,i,j=1}^{5}v_k w_i w_j
		\frac{\partial^2 f_k(D_{02},\beta^*)}{\partial x_i\partial x_j}
		=2\beta^*v_2w_1w_3<0,\\
		\mathfrak b&=\sum_{k,i=1}^{5}v_k w_i
		\frac{\partial^2 f_k(D_{02},\beta^*)}{\partial x_i\partial\beta}
		=v_2w_3\frac{\Lambda}{\mu}>0.
	\end{align*}
	Since $\mathfrak a<0$ and $\mathfrak b>0$, a forward transcritical
	bifurcation occurs in the epidemiological subsystem at $\beta=\beta^*$.
	The positive-capacity equilibrium branch lifts to the full smooth subsystem
	via $B^*=K(1+\lambda I^*/[r(\omega+I^*)])$. The additional capacity
	characteristic root is negative near $D_{02}$, so the local exchange of
	stability is retained in the delayed subsystem.
	
	\section{Characteristic equation coefficients and auxiliary expressions}\label{Appndx_C}
	The characteristic and equilibrium equations used in the following calculations are recalled here to specify the notation.
	
	\begin{equation}\label{context:chddee}
		(s-k_{66})
		\left(
		s^5+l_{11}s^4+l_{12}s^3+l_{13}s^2+l_{14}s+l_{15}
		\right)
		=0,
	\end{equation}
	\begin{equation}\label{context:eq12}
		G_1 I^3 + G_2 I^2 + G_3 I + G_4 = 0,
	\end{equation}
	\begin{equation} \label{context:chara_end} 
		C(s) = s^6 + f_{11} s^5 + f_{12} s^4 + f_{13} s^3 + f_{14} s^2 + f_{15} s +  f_{16}  = 0,
	\end{equation}  
	\begin{equation}
		\begin{aligned} \label{context:chdde2} 
			s^6 + g_{11}s^5 + g_{12}s^4 + g_{13}s^3 + g_{14}s^2 & + g_{15}s + g_{16} \\
			& + (g_{17}s^4 + g_{18}s^3 + g_{19}s^2 + g_{20}s + g_{21})e^{-s\tau} = 0, 
		\end{aligned} 
	\end{equation} 
	\begin{equation}\label{context:eq15}
		\nu^{12} + h_{11}\nu^{10} + h_{12}\nu^{8} + h_{13}\nu^6 + h_{14}\nu^4 + h_{15}\nu^2 + h_{16} = 0, 
	\end{equation}
	\begin{equation}\label{context:eq16}
		\zeta^{6} + h_{11}\zeta^{5} + h_{12}\zeta^4 + h_{13}\zeta^3 + h_{14}\zeta^2 + h_{15}\zeta + h_{16} = 0.
	\end{equation}
	
	\textbf{Coefficients of characteristic equation (3.4):}  
	
	\begin{align*}
		l_{11} =& -(k_{11}+k_{22}+k_{33}+k_{44}+k_{55}),\\
		l_{12} =& -k_{23} k_{32}+k_{33} k_{44}+\left(k_{33}+k_{44}\right) k_{55}+k_{22} \left(k_{33}+k_{44}+k_{55}\right), \\
		&+k_{11} \left(k_{22}+k_{33}+k_{44}+k_{55}\right),
	\end{align*} 
	\begin{align*}
		l_{13} =& -k_{13} k_{21} k_{32}+\left(k_{23} k_{32}-k_{22} k_{33}\right) k_{44}+(k_{23} k_{32}-k_{22} k_{33}\\
		&-\left(k_{22}+k_{33}\right) k_{44}) k_{55}+k_{11} (k_{23} k_{32}-k_{33} k_{44}-\left(k_{33}+k_{44}\right) k_{55}\\
		&-k_{22} \left(k_{33}+k_{44}+k_{55}\right)),\\
		l_{14} =& -k_{15} k_{21} k_{32} k_{53}+\left(k_{22} k_{33}-k_{23} k_{32}\right) k_{44} k_{55}+k_{13} k_{21} k_{32} \left(k_{44}+k_{55}\right)\\
		&+k_{11} \left(k_{22} k_{33} k_{44}+\left(k_{22} k_{33}+\left(k_{22}+k_{33}\right) k_{44}\right) k_{55}-k_{23} k_{32} \left(k_{44}+k_{55}\right)\right),\\
		l_{15} =& k_{15} k_{21} k_{32} \left(k_{44} k_{53}-k_{43} k_{54}\right)-\left(k_{13} k_{21} k_{32}+k_{11} \left(k_{22} k_{33}-k_{23} k_{32}\right)\right) k_{44} k_{55}. 
	\end{align*}
	with
	\begin{align*}
		k_{11} &= -\beta I^* - \mu, \quad
		k_{13} = -\beta S^*, \quad
		k_{15} = \xi, \quad
		k_{21} = \beta I^*, \\[2pt]
		k_{22} &= -(\mu + a), \quad
		k_{23} = \beta S^*, \quad
		k_{32} = a, \quad
		k_{33} = -(\mu + \delta_1 + \gamma_1 + \mu_1), \\[2pt]
		k_{43} &= \mu_1, \quad
		k_{44} = -(\mu + \delta_2 + \gamma_2), \quad
		k_{53} = \gamma_1, \quad
		k_{54} = \gamma_2, \\[2pt]
		k_{55} &= -(\mu + \xi), \quad
		k_{63} = \frac{\lambda \omega B^*}{(\omega + I^*)^2}, \quad
		k_{66} = r - \frac{2rB^*}{K} + \frac{\lambda I^*}{\omega + I^*}.
	\end{align*}
	
	\textbf{Coefficients of equation (3.5):} 
	\begin{align*}
		G_1 =& \beta r (\mu (a+\mu +\xi )(\gamma_1(\mu + \delta_2 + \gamma_2)+\mu_0 \gamma_2) + (a+\mu)(\mu + \xi)(\gamma_2(\mu + \delta_1) \\
		&+ (\mu + \delta_2)(\mu + \delta_1 + \mu_0))) >0,\\
		G_2 =& -a \beta  \left(\gamma _2+\delta _2+\mu \right) \left(\gamma _1 \xi  (K (\lambda +r)+r \omega )+\Lambda  r (\mu +\xi )\right) + (a+\mu ) (\mu \\
		&+\xi ) \left(\gamma _1+\delta _1+\mu \right) \left(\gamma _2+\delta _2+\mu \right) (\beta  K (\lambda +r)+r (\beta  \omega +\mu )) + (a+\mu ) (\mu \\
		&+\xi ) \left(\gamma _2+\delta _2+\mu \right) \left(\beta  K \mu _1 (\lambda +r)+\mu _0 r (\beta  \omega +\mu )\right) - a \beta  \gamma _2 \xi  \left(\mu _0 r \omega + K \mu _1 (\lambda +r)\right),\\
		G_3 =& -a \beta  \left(\gamma _2+\delta _2+\mu \right) \left(\gamma _1 K \xi  r \omega + \Lambda  (\mu +\xi ) (K (\lambda +r)+r \omega )\right) + (a+\mu ) (\mu +\xi ) (\gamma _1 \\
		&+\delta _1+\mu ) \left(\gamma _2+\delta _2+\mu \right) (r \omega  (\beta  K+\mu )+K \mu  (\lambda +r)) + (a+\mu ) (\mu +\xi ) (\gamma _2 \\
		&+\delta _2+\mu ) \left(K \mu _1 (\beta  r \omega +\mu  (\lambda +r))+\mu  \mu _0 r \omega \right)-a \beta  \gamma _2 K \mu _1 \xi  r \omega,\\
		G_4 =&  K r \omega \mu (\mu + \xi)(\mu + \delta_2 + \gamma_2)(a+\mu)(\mu + \delta_1 + \gamma_1 + \mu_1)(1- \mathcal R_0).  
	\end{align*}
	
	\textbf{Coefficients of characteristic equation (3.9):}  
	
	\begin{align*}
		f_{11} =& -(j_{11}+j_{22}+j_{33}+j_{44}+j_{55}+j_{66}),\\
		f_{12} =& -j_{23} j_{32}+j_{33} j_{44}+j_{33} j_{55}+j_{44} j_{55}-j_{36} j_{63}+\left(j_{33}+j_{44}+j_{55}\right) j_{66}\\
		&+j_{22} \left(j_{33}+j_{44}+j_{55}+j_{66}\right)+j_{11} \left(j_{22}+j_{33}+j_{44}+j_{55}+j_{66}\right),\\
		f_{13} =& -j_{13} j_{21} j_{32}-j_{33} j_{44} j_{55}-j_{22} \left(j_{33} j_{44}+\left(j_{33}+j_{44}\right) j_{55}\right)+j_{23} j_{32} \left(j_{44}+j_{55}\right)\\
		&+j_{36} \left(j_{22}+j_{44}+j_{55}\right) j_{63}+(j_{23} j_{32}-j_{33} j_{44}-\left(j_{33}+j_{44}\right) j_{55}\\
		&-j_{22} \left(j_{33}+j_{44}+j_{55}\right)) j_{66}+j_{11} (j_{23} j_{32}-j_{44} j_{55}+j_{36} j_{63}-\left(j_{44}+j_{55}\right) j_{66}\\
		&-j_{33} \left(j_{44}+j_{55}+j_{66}\right)-j_{22} \left(j_{33}+j_{44}+j_{55}+j_{66}\right)),\\
		f_{14} =& -j_{15} j_{21} j_{32} j_{53}+\left(j_{22} j_{33}-j_{23} j_{32}\right) j_{44} j_{55}-j_{36} \left(j_{44} j_{55}+j_{22} \left(j_{44}+j_{55}\right)\right) j_{63}\\
		&+\left(j_{33} j_{44} j_{55}-j_{23} j_{32} \left(j_{44}+j_{55}\right)+j_{22} \left(j_{33} j_{44}+\left(j_{33}+j_{44}\right) j_{55}\right)\right) j_{66}\\
		&+j_{13} j_{21} j_{32} \left(j_{44}+j_{55}+j_{66}\right)+j_{11} (j_{33} j_{44} j_{55}-j_{36} j_{63} j_{55}-j_{36} j_{44} j_{63}\\
		&+\left(j_{44} j_{55}+j_{33} \left(j_{44}+j_{55}\right)\right) j_{66}-j_{23} j_{32} \left(j_{44}+j_{55}+j_{66}\right)+j_{22} (j_{44} j_{55}\\
		&-j_{36} j_{63}+\left(j_{44}+j_{55}\right) j_{66}+j_{33} \left(j_{44}+j_{55}+j_{66}\right))),
	\end{align*}
	\begin{align*} 
		f_{15} =& j_{11} \left(j_{23} j_{32}-j_{22} j_{33}\right) j_{44} j_{55}+j_{36} \left(j_{22} j_{44} j_{55}+j_{11} \left(j_{44} j_{55}+j_{22} \left(j_{44}+j_{55}\right)\right)\right) j_{63}\\
		&+(\left(j_{23} j_{32}-j_{22} j_{33}\right) j_{44} j_{55}+j_{11} (-j_{22} j_{33} j_{44}-\left(j_{22} j_{33}+\left(j_{22}+j_{33}\right) j_{44}\right) j_{55}\\
		&+j_{23} j_{32} \left(j_{44}+j_{55}\right))) j_{66}+j_{15} j_{21} j_{32} \left(j_{53} \left(j_{44}+j_{66}\right)-j_{43} j_{54}\right)-j_{13} j_{21} j_{32} (j_{55} j_{66}\\
		&+j_{44} \left(j_{55}+j_{66}\right)),\\
		f_{16} =& j_{44} j_{55} \left(j_{13} j_{21} j_{32} j_{66}-j_{11} \left(j_{22} j_{36} j_{63}+\left(j_{23} j_{32}-j_{22} j_{33}\right) j_{66}\right)\right)\\
		&-j_{15} j_{21} j_{32} \left(j_{46} j_{54} j_{63}+\left(j_{44} j_{53}-j_{43} j_{54}\right) j_{66}\right),
	\end{align*}
	with 
	\begin{align*}
		j_{11} &= -\beta I^* - \mu, \quad j_{13} = -\beta S^*, \quad j_{15} = \xi, \quad j_{21} = \beta I^*, \quad j_{22} = -(\mu +a), \quad j_{23} = \beta S^*,\\
		j_{32} &= a, \quad j_{33} = -\left(\mu + \delta_1 + \gamma_1 + \mu_0 + (\mu_1 - \mu_0)\frac{B^{{*}^2}}{(B^* + I^*)^2} \right) , \quad j_{36} = -\frac{(\mu_1 - \mu_0)I^{{*}^2}}{(B^* + I^*)^2},\\
		j_{43} & = \mu_0 + (\mu_1 - \mu_0)\frac{B^{{*}^2}}{(B^* + I^*)^2}, \quad j_{44} = -(\mu + \delta_2 + \gamma_2), \quad j_{46} = \frac{(\mu_1 - \mu_0)I^{{*}^2}}{(B^* + I^*)^2}, \quad j_{53} = \gamma_1,\\
		j_{54} &= \gamma_2, \quad j_{55} = -(\mu + \xi), \quad j_{63} = \frac{\lambda \omega B^*}{(\omega + I^*)^2}, \quad j_{66} = r - \frac{2rB^*}{K}+\frac{\lambda I^*}{\omega + I^*}.
	\end{align*}
	
	\textbf{Coefficients of characteristic equation (3.10):}   
	\begin{align*}
		g_{11}  =& -(j_{11}+j_{22}+j_{33}+j_{44}+j_{55}+j_{66}),\\
		g_{12}  =& -j_{23} j_{32}+j_{33} j_{44}+j_{55} (j_{33} +j_{44}) +\left(j_{33}+j_{44}+j_{55}\right) j_{66}\\
		&+j_{22} \left(j_{33}+j_{44}+j_{55}+j_{66}\right)+j_{11} \left(j_{22}+j_{33}+j_{44}+j_{55}+j_{66}\right),\\
		g_{13} = & -j_{13} j_{21} j_{32}+j_{23} j_{32} (j_{44}+ j_{55}) -j_{22} j_{33} (j_{44}+ j_{55})- j_{44} j_{55} (j_{22} + j_{33}) \\
		&+(j_{23} j_{32}-j_{33} j_{44}-(j_{33}+j_{44}) j_{55}-j_{22} (j_{33}+j_{44}+j_{55})) j_{66}+j_{11} (j_{23} j_{32}\\
		&-j_{44} j_{55}-(j_{44}+j_{55}) j_{66}-j_{33} (j_{44}+j_{55}+j_{66})-j_{22} (j_{33}+j_{44}+j_{55}+j_{66})),\\
		g_{14} =& j_{22} j_{33} j_{44} j_{55} -j_{32} (j_{15} j_{21} j_{53}+j_{23} j_{44} j_{55}) +(j_{33} j_{44} j_{55}-j_{23} j_{32} (j_{44}+j_{55})\\
		&+j_{22} (j_{33} j_{44}+(j_{33}+j_{44}) j_{55})) j_{66}+j_{13} j_{21} j_{32} (j_{44}+j_{55}+j_{66})+j_{11} (j_{33} j_{44} j_{55} \\
		&+(j_{44} j_{55}+j_{33} (j_{44}+j_{55})) j_{66}-j_{23} j_{32} (j_{44}+j_{55}+j_{66})+j_{22} (j_{44} j_{55}+(j_{44}\\
		&+j_{55}) j_{66}+j_{33} (j_{44}+j_{55}+j_{66}))),\\
		g_{15} =& j_{11} j_{44} j_{55} (j_{23} j_{32}-j_{22} j_{33}) +((j_{23} j_{32}-j_{22} j_{33}) j_{44} j_{55}+j_{11} (-j_{22} j_{33} j_{44} \\
		&-(j_{22} j_{33}+(j_{22}+j_{33}) j_{44}) j_{55}+j_{23} j_{32} (j_{44}+j_{55}))) j_{66}+j_{15} j_{21} j_{32} (j_{53} (j_{44}\\
		&+j_{66})-j_{43} j_{54})-j_{13} j_{21} j_{32} (j_{55} j_{66}+j_{44} (j_{55}+j_{66})),\\
		g_{16}  = & (j_{15} j_{21} j_{32} (j_{43} j_{54}-j_{44} j_{53})+(j_{13} j_{21} j_{32}+j_{11} (j_{22} j_{33}-j_{23} j_{32})) j_{44} j_{55}) j_{66},\\
		g_{17}  =& -j_{36} j_{63},\\
		g_{18} =& j_{36} j_{63} (j_{11}+j_{22}+j_{44}+j_{55}),\\
		g_{19}  =& -j_{36} j_{63} (j_{44} j_{55}+j_{22} (j_{44}+j_{55})+j_{11} (j_{22}+j_{44}+j_{55})),\\
		g_{20} =& j_{36} j_{63} \left(j_{22} j_{44} j_{55}+j_{11} \left(j_{44} j_{55}+j_{22} \left(j_{44}+j_{55}\right)\right)\right),\\
		g_{21} =& -j_{63} \left(j_{15} j_{21} j_{32} j_{46} j_{54}+j_{11} j_{22} j_{36} j_{44} j_{55}\right).
	\end{align*}
	
	\textbf{Coefficients of equation (3.11):}    
	
	\begin{align*}
		h_{11} &= g_{11}^2 - 2g_{12},\quad 
		h_{12} = g_{12}^2-g_{17}^2 + 2 g_{14} - 2 g_{11} g_{13},\\
		h_{13} &= g_{13}^2 - g_{18}^2 + 2 g_{17} g_{19} -2g_{16} - 2 g_{12}g_{14} + 2 g_{11}g_{15},\\
		h_{14} &= g_{14}^2 - g_{19}^2 + 2 g_{18} g_{20} -2g_{17}g_{21} + 2 g_{12}g_{16} - 2 g_{13}g_{15},\\
		h_{15} &= g_{15}^2 - g_{20}^2 + 2 g_{19} g_{21} - 2 g_{14} g_{16}, \quad
		h_{16} = g_{16}^2 - g_{21}^2. 
	\end{align*}
	
	\begin{align*}
		M_1 &= 3 f_{13} \eta^2 - 5f_{11} \eta^4 - f_{15}, \quad
		M_2 = 4 f_{12} \eta^3 - 2 f_{14} \eta - 6 \eta^5, \\
		M_3 &= \eta^4 \frac{d f_{12}}{d \lambda} - \eta^2 \frac{d f_{14}}{d \lambda} + \frac{d f_{16}}{d \lambda}, \quad
		M_4 =  \eta^5 \frac{d f_{11}}{d \lambda} + \eta \frac{d f_{15}}{d \lambda} - \eta^3 \frac{d f_{13}}{d \lambda}.
	\end{align*}
	
	\begin{align*}
		A_{11} =& \left(5 g_{11} \nu^*{^4}-3 g_{13} \nu^*{^2}+g_{15}\right) \sin (\nu^*  \tau^* )+2 \nu^*  \left(-2 g_{12} \nu^*{^2}+g_{14}+3 \nu^*{^4}\right) \cos (\nu^*  \tau^* )\\
		& + g_{18} \nu^*{^3} \tau^* -4 g_{17} \nu^*{^3}-g_{20} \nu^*  \tau^* +2 g_{19} \nu^*,\\
		A_{12} =& 2 \nu^*  \left(-2 g_{12} \nu^*{^2}+g_{14}+3 \nu^*{^4}\right) \sin (\nu^*  \tau^* )-\left(5 g_{11} \nu^*{^4}-3 g_{13} \nu^*{^2}+g_{15}\right) \cos (\nu^*  \tau^* ) \\
		&+ g_{17} \nu^*{^4} \tau^* -g_{19} \nu^*{^2} \tau^* +3 g_{18} \nu^*{^2}+g_{21} \tau^* -g_{20},\\
		A_{13} =& g_{17} \nu^*{^5}-g_{19} \nu^*{^3}+g_{21} \nu^*,\\
		A_{14} =& g_{20} \nu^*{^2}-g_{18} (\nu^*)^4.
	\end{align*}
	
	\subsection{Next-generation matrix calculation}\label{supp:ngm}
	For the below-threshold subsystem, the disease-free susceptible population is $S^0=\Lambda/\mu$.
	For subsystem $F_{\Sigma_1}$, the infected compartments are $(E,I)$.  Writing the infected subsystem in the form $\frac{ dx}{dt} = \mathcal F(x) - \mathcal V(x)$ with $x=(E,I)^\top$, where $\mathcal F(x)$ and $\mathcal V(x)$ represent the new infection and transition terms as follows: 
	\[
	\mathcal F(E,I)=
	\begin{pmatrix}
		\beta S I\\[1mm]
		0
	\end{pmatrix} 
	\quad \text{and} \quad
	\mathcal V(E,I)=
	\begin{pmatrix}
		(\mu+a)E\\[1mm]
		(\mu+\delta_1+\gamma_1+\mu_1)I-aE
	\end{pmatrix},
	\]
	respectively. 
	Evaluating at the disease-free equilibrium $D_{02}$, where $S^0=\Lambda/\mu$,
	the Jacobians $F = D\mathcal F(D_{02})$ and $V = D\mathcal V(D_{02})$ are
	\[
	F=
	\begin{pmatrix}
		0 & \beta S^0\\
		0 & 0
	\end{pmatrix}
	\quad \text{and} \quad
	V=
	\begin{pmatrix}
		\mu+a & 0\\
		-a & \mu+\delta_1+\gamma_1+\mu_1
	\end{pmatrix}.
	\]
	The next-generation matrix is $\mathcal K = F V^{-1}$ and $\mathcal R_0=\rho(\mathcal K)$, where $\rho(\mathcal K)$ is the spectral radius of the matrix $\mathcal K$. 
	A direct calculation yields
	\[
	\mathcal R_0=\frac{\beta\Lambda a}{\mu(\mu+a)(\mu+\delta_1+\gamma_1+\mu_1)}.
	\]
	\subsection{Routh-Hurwitz determinants}\label{supp:hurwitz}
	For the fifth-degree polynomial in Eq.~\eqref{context:chddee}, we introduce the Routh-Hurwitz determinants
	\begin{equation*}
		\begin{aligned}
			H_1
			&=l_{11},
			\quad
			H_2
			=
			\begin{array}{|cc|}
				l_{11} & l_{13}\\
				1 & l_{12}
			\end{array},
			\quad
			H_3
			=
			\begin{array}{|ccc|}
				l_{11} & l_{13} & l_{15}\\
				1 & l_{12} & l_{14}\\
				0 & l_{11} & l_{13}
			\end{array},
			\\
			H_4
			&=
			\begin{array}{|cccc|}
				l_{11} & l_{13} & l_{15} & 0\\
				1 & l_{12} & l_{14} & 0\\
				0 & l_{11} & l_{13} & l_{15}\\
				0 & 1 & l_{12} & l_{14}
			\end{array},
			\quad
			H_5=l_{15}H_4.
		\end{aligned}
	\end{equation*}

	For the sixth-degree polynomial in Eq.~\eqref{context:chara_end}, we further define Routh-Hurwitz determinants 
	\begin{equation*}
		\begin{aligned}
			H_1 =& f_{11} , \quad H_2 = \begin{array}{|cc|}
				f_{11} & f_{13}  \\ 
				1 & f_{12}
			\end{array} , \quad
			H_3 = \begin{array}{|ccc|}
				f_{11} & f_{13}  & f_{15}  \\ 
				1 & f_{12}  & f_{14} \\ 
				0 & f_{11}  & f_{13} 
			\end{array} , \quad
			H_4 = \begin{array}{|cccc|}
				f_{11} & f_{13}  & f_{15} & 0 \\ 
				1 & f_{12}  & f_{14}  & f_{16} \\ 
				0 & f_{11}  & f_{13} & f_{15} \\
				0 & 1 & f_{12}  & f_{14} 
			\end{array} , \\
			H_5 =& \begin{array}{|ccccc|}
				f_{11} & f_{13}  & f_{15} & 0 & 0  \\ 
				1 & f_{12}  & f_{14}  & f_{16} & 0 \\ 
				0 & f_{11}  & f_{13} & f_{15} & 0 \\
				0 & 1 & f_{12}  & f_{14}  & f_{16} \\
				0 & 0 & f_{11}  & f_{13} & f_{15} 
			\end{array} , \quad  H_6 = f_{16} H_5.
		\end{aligned}
	\end{equation*}
	
	\subsection{Critical-delay phase calculations}\label{supp:delayphase}
	Substituting $s=i\nu$ into Eq.~\eqref{context:chdde2} and separating real and imaginary parts gives
	\begin{equation}\label{eq14} 
		\begin{aligned}
			(g_{20} \nu - g_{18}\nu^3) \sin(\nu \tau) + (g_{21} + g_{17}\nu^4 - g_{19} \nu^2) \cos(\nu \tau) &= \nu^6 - g_{12} \nu^4 + g_{14} \nu^2 - g_{16},\\
			(g_{20} \nu - g_{18}\nu^3) \cos(\nu \tau) - (g_{21} + g_{17}\nu^4 - g_{19} \nu^2) \sin(\nu \tau) &= g_{13} \nu^3 - g_{11} \nu^5 - g_{15} \nu. 
		\end{aligned}
	\end{equation}

	For each positive root $\zeta_j$ of Eq.~\eqref{context:eq16}, set $\nu_j=\sqrt{\zeta_j}$ and define
	\begin{align*}
		F_{1j}&=g_{20}\nu_j-g_{18}\nu_j^3,\\
		F_{2j}&=g_{21}+g_{17}\nu_j^4-g_{19}\nu_j^2,\\
		F_{3j}&=\nu_j^6-g_{12}\nu_j^4+g_{14}\nu_j^2-g_{16},\\
		F_{4j}&=g_{13}\nu_j^3-g_{11}\nu_j^5-g_{15}\nu_j.
	\end{align*}
	Assuming $F_{1j}^2+F_{2j}^2>0$, from equation \eqref{eq14},
	\begin{equation}\label{eq:cos_delay}
		\cos(\nu_j\tau)
		=
		\frac{F_{2j}F_{3j}+F_{1j}F_{4j}}
		{F_{1j}^2+F_{2j}^2},
	\end{equation}
	and
	\begin{equation}\label{eq:sin_delay}
		\sin(\nu_j\tau)
		=
		\frac{F_{1j}F_{3j}-F_{2j}F_{4j}}
		{F_{1j}^2+F_{2j}^2}.
	\end{equation}

	\paragraph{Crossing direction from the frequency polynomial.}
	Write the characteristic equation as $P_c(s)+Q_c(s)e^{-s\tau}=0$ and define $h(\zeta)=|P_c(i\sqrt\zeta)|^2-|Q_c(i\sqrt\zeta)|^2$. At a simple imaginary root $s=i\nu$ with $Q_c(i\nu)\ne0$, implicit differentiation and the characteristic equation yield
	\[
	\left(\frac{ds}{d\tau}\right)^{-1}=
	\frac{Q_c'(s)/Q_c(s)-P_c'(s)/P_c(s)}{s}-\frac{\tau}{s}.
	\]
	The last term is purely imaginary at $s=i\nu$. Differentiating $h$ with respect to $\zeta=\nu^2$ and using $|P_c(i\nu)|=|Q_c(i\nu)|$ gives
	\[
	\operatorname{Re}\left(\frac{ds}{d\tau}\right)^{-1}=\frac{h'(\nu^2)}{|Q_c(i\nu)|^2}.
	\]
	Thus the sign of $h'(\zeta)$ agrees with the sign of the real part of $ds/d\tau$, and with the crossing numerator used in the main manuscript. This provides an independent means of evaluating the critical-delay crossing directions.
	
	\section{Delay-free Hopf transversality calculation}\label{supp:hopftrans}
	For $\tau=0$, assume a smooth branch of above-threshold endemic equilibria $D(\lambda)$ near $\lambda=\lambda^*$, with characteristic polynomial
	\begin{equation}\label{context:chara_end_hopf}
	C(s,\lambda)=s^6+f_{11}s^5+f_{12}s^4+f_{13}s^3+f_{14}s^2+f_{15}s+f_{16},
	\end{equation}
	where the coefficients and Hurwitz determinants are defined in Supplementary Section~\ref{Appndx_C}. By Liu's criterion \cite{Liu1994}, the conditions at $\lambda^*$ are (i) $f_{16}>0$, $H_1,H_2,H_3,H_4>0$, $H_5=0$, and (ii) $dH_5/d\lambda\neq0$. Condition (i) yields a simple imaginary pair $\pm i\eta$ with all remaining roots in the open left half-plane; condition (ii) gives transversality. A nondegenerate Hopf bifurcation additionally requires a nonzero first Lyapunov coefficient. Differentiating the characteristic equation \eqref{context:chara_end_hopf} with respect to $\lambda$, we obtain
	\begin{equation*}
		\begin{aligned}
			(6 s^5 + 5 f_{11} s^4 + 4 f_{12} s^3 &+ 3 f_{13} s^2 + 2 f_{14} s + f_{15}) \frac{d s}{d \lambda} \\
			&+ s^5 \frac{d f_{11}}{d \lambda} + s^4 \frac{d f_{12}}{d \lambda} + s^3 \frac{d f_{13}}{d \lambda} + s^2 \frac{d f_{14}}{d \lambda} + s \frac{d f_{15}}{d \lambda} + \frac{d f_{16}}{d \lambda} = 0.
		\end{aligned} 
	\end{equation*}
	Assuming that the parameter derivative of the characteristic polynomial is nonzero at the critical root, we obtain
	\begin{equation*}
		\begin{aligned}
			\Big(\frac{d s}{d \lambda}\Big)^{-1} = - \frac{6 s^5 + 5 f_{11} s^4 + 4 f_{12} s^3 + 3 f_{13} s^2 + 2 f_{14} s + f_{15}}{s^5 \frac{d f_{11}}{d \lambda} + s^4 \frac{d f_{12}}{d \lambda} + s^3 \frac{d f_{13}}{d \lambda} + s^2 \frac{d f_{14}}{d \lambda} + s \frac{d f_{15}}{d \lambda} + \frac{d f_{16}}{d \lambda}}.
		\end{aligned}
	\end{equation*}
	Furthermore, we have
	\begin{equation*}
		\begin{aligned}
			sign \Bigg[\frac{d(Re(s))}{d \lambda}\Bigg]_{s=i\eta,\,\lambda=\lambda^*} &= sign \Bigg[Re \Bigg(\frac{d s}{d \lambda} \Bigg)^{-1} \Bigg]_{s=i\eta,\,\lambda=\lambda^*} \\
			& = sign [\Upsilon ], \\
		\end{aligned}
	\end{equation*}
	where, 
	\begin{equation*}
		\begin{aligned}
			\Upsilon & = Re \Bigg[\frac{(3 f_{13} \eta^2 - 5f_{11} \eta^4 - f_{15}) + i (4 f_{12} \eta^3 - 2 f_{14} \eta - 6 \eta^5)}{(\eta^4 \frac{d f_{12}}{d \lambda} - \eta^2 \frac{d f_{14}}{d \lambda} + \frac{d f_{16}}{d \lambda} ) + i ( \eta^5 \frac{d f_{11}}{d \lambda} + \eta \frac{d f_{15}}{d \lambda} - \eta^3 \frac{d f_{13}}{d \lambda} )} \Bigg] = \frac{M_1 M_3 + M_2 M_4}{M_3^2 + M_4^2 },
		\end{aligned}
	\end{equation*}
	where $M_{1}, M_{2}, M_{3}$ and $M_{4}$ are provided in Supplementary Section~\ref{Appndx_C}. If $M_1 M_3 + M_2 M_4 \neq 0,$ then \\  $sign \Big[\frac{d(Re(s))}{d \lambda}\Big]_{\lambda = \lambda^*} \neq 0 $ and the transversality condition (ii) holds. This gives the transversality calculation used in the main manuscript, provided $M_3^2+M_4^2>0$.
	
	\section{Lyapunov-based global stability analysis}\label{Appndx_D}
	\noindent\textbf{Auxiliary infection-integral result.} For $\mathcal R_0>1$ and nonnegative initial histories with $E(0)+I(0)>0$, we show that $\int_0^\infty I(t)\,dt=\infty$ in the smooth extension of $F_{\Sigma_1}$. Suppose, to the contrary, that $\int_0^\infty I(t)\,dt<\infty$. Boundedness of the epidemiological states and their derivatives implies uniform continuity of $I$, so Barbalat's lemma yields $I(t)\to0$. The stable linear equations for $E,H,R$ then give $E,H,R\to0$, and the susceptible equation gives $S\to\Lambda/\mu$. Define $L=aE+(\mu+a)I$. Along the below-threshold smooth extension,
	\[
	\dot L=[a\beta S-(\mu+a)(\mu+\delta_1+\gamma_1+\mu_1)]I.
	\]
	Because $\mathcal R_0>1$, the coefficient in brackets is eventually strictly positive. Moreover $L(t)>0$ for all finite $t$ by positivity of the infected subsystem and $E(0)+I(0)>0$. Thus $L$ is eventually nondecreasing and bounded below by a positive value, contradicting $L\to0$. This proves the auxiliary result.
	
	\noindent\textbf{Statement of Theorem~3.3 (main manuscript).} Suppose $\mathcal R_0>1$, $I^*<I_0$, and there exist $P=P^\top\succ0$ and $\eta>0$ such that $PQ(0)+Q(0)^\top P\preceq0$ and $PQ(I_0)+Q(I_0)^\top P\preceq-\eta I_5$. Then $D^*$ is globally asymptotically stable relative to nonnegative histories with $E(0)+I(0)>0$, $B(0)>0$ and trajectories remaining in $I<I_0$.
	
	\begin{proof}
		At $D^*$, the equilibrium relations are
		\begin{align}
			\beta S^*I^*+\mu S^*-\xi R^* & = \Lambda , \label{eq1}\\
			(\mu+a)E^* &= \beta S^*I^*, \label{eq2}\\
			(\mu+\delta_1+\gamma_1+\mu_1)I^* &= aE^*, \label{eq3}\\
			(\mu+\delta_2+\gamma_2)H^* &= \mu_1 I^*, \label{eq4}\\
			(\mu+\xi)R^* &= \gamma_1 I^*+\gamma_2 H^*. \label{eq5}
		\end{align}
		
		Define the deviation variables 
		\[
		x_1=S-S^*, ~ x_2=E-E^*, ~ x_3=I-I^*, ~ x_4=H-H^*, ~ x_5=R-R^*, ~ x=(x_1,x_2,x_3,x_4,x_5)^{\top}.
		\]
		
		Using the identity 
		\begin{equation}\label{eq:factor_SI}
			SI-S^*I^* = I\,x_1 + S^*x_3,
		\end{equation}
		(which follows from $SI=(S^*+x_1)(I^*+x_3)=S^*I^*+I^*x_1+S^*x_3+x_1x_3
		= S^*I^* + (I^*+x_3)x_1 + S^*x_3 = S^*I^* + Ix_1 + S^*x_3$),
		and subtracting \eqref{eq1}-\eqref{eq5} from equations of subsystem  $F_{\Sigma_1}$, the
		deviation dynamics can be written exactly in the state-dependent linear form 
		\begin{equation}\label{eq:sdlin}
			\dot x = Q(I(t))\,x,
		\end{equation}
		where $I(t)=I^*+x_3(t)$ and
		\[
		Q(I)=
		\begin{pmatrix}
			-(\mu+\beta I) & 0 & -\beta S^* & 0 & \xi\\
			\beta I & -(\mu+a) & \beta S^* & 0 & 0\\
			0 & a & -(\mu+\delta_1+\gamma_1+\mu_1) & 0 & 0\\
			0 & 0 & \mu_1 & -(\mu+\delta_2+\gamma_2) & 0\\
			0 & 0 & \gamma_1 & \gamma_2 & -(\mu+\xi)
		\end{pmatrix}.
		\]

		Since $Q(I)$ is affine in $I$, set $\theta=I/I_0$. Convexity gives
		\begin{equation}\label{eq:diagstab_assump}
			PQ(I)+Q(I)^\top P\preceq-\theta\eta I_5,\quad 0\leq I\leq I_0.
		\end{equation}
		For $V=x^\top Px$,
		\[
		\dot V\leq-\frac{\eta I(t)}{I_0}\|x\|^2\leq-cI(t)V,
		\quad c=\frac{\eta}{I_0\lambda_{\max}(P)}.
		\]
		Integration gives $V(t)\leq V(0)\exp[-c\int_0^t I(u)\,du]$. The auxiliary result above implies $V(t)\to0$ and hence $x(t)\to0$. For fixed $\tau$, $g(t)=\lambda I(t-\tau)/(\omega+I(t-\tau))\to g^*=\lambda I^*/(\omega+I^*)$. Comparing $\dot B=(r+g(t))B-rB^2/K$ with logistic equations having growth rates $r+g^*\pm\varepsilon$, for sufficiently small $\varepsilon>0$ and sufficiently large $t$, gives $B(t)\to K(1+g^*/r)=B^*$ when $B(0)>0$.
		
		Lyapunov stability follows from $V(t)\leq V(0)$ for the epidemiological variables. For histories sufficiently close to the equilibrium, the delayed infection remains uniformly close to $I^*$ (including the prescribed history). Logistic comparison with growth rates $r+g^*\pm\varepsilon$ bounds $B(t)$ between the initial capacity and the corresponding nearby positive equilibria. This yields stability in the history supremum norm. Since $I^*<I_0$, sufficiently small perturbations remain below the threshold. This completes the proof.
	\end{proof}
	
	\noindent\textbf{Nonvacuity of the matrix condition.} At the endemic equilibrium, $a\beta S^*=(\mu+a)(\mu+\delta_1+\gamma_1+\mu_1)$. Thus $Q(0)$ has a simple zero eigenvalue and four strictly stable eigenvalues. Let $w,v$ be its right and left null vectors, normalized by $v^\top w=1$, with $w_3=1$. As in Supplementary Section~\ref{Appndx_B}, $w_1<0$ and $v_2>0$. In a basis consisting of $w$ and a basis of $\ker v^\top$, $Q(0)$ is block diagonal with blocks $0$ and a Hurwitz matrix $Q_s$. Choose $P$ by transforming back the block diagonal matrix $\operatorname{diag}(1,P_s)$, where $P_s\succ0$ solves $P_sQ_s+Q_s^\top P_s=-I_4$. Then $Pw=v$ and $PQ(0)+Q(0)^\top P$ is negative semidefinite with kernel $\operatorname{span}\{w\}$. Write $Q(I)=Q(0)+\beta I A_1$, where the only nonzero entries of $A_1$ are $(A_1)_{11}=-1$ and $(A_1)_{21}=1$. Since
	\[
	w^\top(PA_1+A_1^\top P)w=2v_2w_1<0,
	\]
	the symmetric matrix $PQ(I)+Q(I)^\top P$ is strictly negative definite for sufficiently small positive $\beta I$, by the Schur complement in the decomposition into its null direction and its orthogonal complement. Hence, there exists $\varepsilon_0>0$ such that the revised matrix criterion is feasible whenever $0<\beta I_0<\varepsilon_0$. With the other parameters fixed, $Q(0)$ is independent of $\beta$, since $a\beta S^*=(\mu+a)(\mu+\delta_1+\gamma_1+\mu_1)$. Fix $I_0>0$ sufficiently small that $\beta^* I_0<\varepsilon_0$, where $\beta^*$ is the critical transmission rate corresponding to $\mathcal R_0=1$. Then choose $\beta>\beta^*$ sufficiently close to $\beta^*$ so that $\beta I_0<\varepsilon_0$ and $0<I^*<I_0$. Thus, the assumptions of the theorem can be satisfied. This establishes analytical nonvacuity of the criterion.
	
	\section{Lie derivatives and grazing classification}\label{Appndx_E}
	Let $L(X)=I-I_0$ be the switching function, with $\nabla L=(0,0,1,0,0,0)$. The normal components of the one-sided vector fields are
	\begin{align*}
		\mathcal{L}_1(X)
		&=
		\left\langle \nabla L,F_{\Sigma_1}\right\rangle
		=
		aE-(\mu+\delta_1+\gamma_1+\mu_1)I,\\
		\mathcal{L}_2(X)
		&=
		\left\langle \nabla L,F_{\Sigma_2}\right\rangle 
		=
		aE-
		\Big(
		\mu+\delta_1+\gamma_1+\mu_0
		+(\mu_1-\mu_0)\frac{B}{B+I}
		\Big)I.
	\end{align*}
	
	On $\Sigma$, their zeros define the critical values $E_{c1}$ and $E_{c2}(B)$ given in the main manuscript. For $0<\mu_0<\mu_1$, $I_0>0$ and finite $B\geq0$,
	\begin{equation}\label{context:Ec_order}
		E_{c2}(B)<E_{c1}.
	\end{equation}
	
	For the delay-free system, the first and second Lie derivatives of the switching function $L(x)$ along a smooth 
	vector field $F(x)$ are defined by
	\[
	\mathcal{L}_F L(x) = \nabla L(x)\cdot F(x),
	\]
	and
	\[
	\mathcal{L}_F^2 L(x)
	= \mathcal{L}_F(\mathcal{L}_F L)(x)
	= \nabla\bigl(\mathcal{L}_F L(x)\bigr)\cdot F(x),
	\]
	respectively. 
	
	Thus, $\mathcal L_F L=F_I$ and $\mathcal L_F^2 L=\nabla(F_I)\cdot F$, where $F_I$ is the infectious component of the corresponding one-sided vector field.
	Restricting to \(\Sigma\) (setting \(I = I_0\)) gives
	\begin{align}
		F_{\Sigma_1,I}\big|_{\Sigma} 
		&= aE - (\mu + \delta_1 + \gamma_1 + \mu_1) I_0, \label{eq:FminusA}\\[2pt]
		F_{\Sigma_2,I}\big|_{\Sigma} 
		&= aE - \bigl(\mu + \delta_1 + \gamma_1 + \mu_0 
		+ (\mu_1-\mu_0)\tfrac{B}{B+I_0}\bigr)\, I_0. \label{eq:FplusA}
	\end{align}
	
	\paragraph{Case $I<I_0$:}
	For this case, the treatment rate is constant, $T(I,B)\equiv\mu_1$, and
	\[
	F_{\Sigma_1,I} = aE - (\mu + \delta_1 + \gamma_1 + \mu_1) I.
	\]
	Thus
	\[
	\nabla F_{\Sigma_1,I}
	= (0,\; a,\; -(\mu + \delta_1 + \gamma_1 + \mu_1),\; 0,\; 0,\; 0).
	\]
	The second Lie derivative is
	\[
	\mathcal{L}^2_{F_{\Sigma_1}} L
	= \nabla F_{\Sigma_1,I} \cdot F_{\Sigma_1}
	= a\,\dot E_{\Sigma_1}
	- (\mu + \delta_1 + \gamma_1 + \mu_1)\,\dot I_{\Sigma_1},
	\]
	where $F_{\Sigma_1} = (\dot S_{\Sigma_1},\dot E_{\Sigma_1},\dot I_{\Sigma_1},
	\dot H_{\Sigma_1},\dot R_{\Sigma_1},\dot B_{\Sigma_1})$ is the one-sided 
	vector field for $I<I_0$.
	
	Evaluating on the switching manifold $\Sigma$ (i.e., setting $I=I_0$), we use
	\[
	\dot E_{\Sigma_1} 
	= \beta S I_0 - (\mu + a)E,
	\quad
	\dot I_{\Sigma_1} 
	= aE - (\mu + \delta_1 + \gamma_1 + \mu_1) I_0,
	\]
	to obtain
	\[
	\mathcal{L}^2_{F_{\Sigma_1}} L
	=
	a \bigl(\beta S I_0 - (\mu + a)E\bigr)
	- (\mu + \delta_1 + \gamma_1 + \mu_1)\bigl(aE - (\mu + \delta_1 + \gamma_1 + \mu_1)I_0\bigr).
	\]
	At the tangency curve
	\[
	\Sigma^t_1 = \{F_{\Sigma_1,I} = 0\}
	\quad \Longleftrightarrow \quad
	aE_{c1} = (\mu + \delta_1 + \gamma_1 + \mu_1) I_0,
	\]
	we obtain
	\[
	\mathcal{L}^2_{F_{\Sigma_1}}L\big|_{\Sigma^t_1}
	=
	a\bigl(\beta S I_0 - (\mu + a) E_{c1}\bigr),
	\]
	which is generically nonzero.
	
	\paragraph{Case $I>I_0$:} 
	For this case,
	\[
	F_{\Sigma_2,I}
	=
	aE - \Bigl(\mu + \delta_1 + \gamma_1 
	+ \mu_0 + (\mu_1-\mu_0)\tfrac{B}{B+I}\Bigr)I.
	\]
	The gradient of $F_{\Sigma_2,I}$ is given by
	\[
	\nabla F_{\Sigma_2,I}
	=
	\Bigl(
	0,\;
	a,\;
	-\Theta(I,B),\;
	0,\;
	0,\;
	-(\mu_1-\mu_0)\frac{I^2}{(B+I)^2}
	\Bigr),
	\]
	where
	\[
	\Theta(I,B)
	=
	\mu + \delta_1 + \gamma_1 + \mu_0
	+ (\mu_1-\mu_0)\frac{B}{B+I}
	- (\mu_1-\mu_0)\frac{IB}{(B+I)^2}.
	\]
	The second Lie derivative is therefore
	\[
	\mathcal{L}^2_{F_{\Sigma_2}}L
	=
	\nabla F_{\Sigma_2,I}\cdot F_{\Sigma_2}.
	\]
	
	Restricting to the switching manifold $I=I_0$ and evaluating on the tangency set
	$\Sigma^t_2=\{F_{\Sigma_2,I}=0\}$, we have $\dot I_{\Sigma_2}=0$, and hence the
	term involving $\Theta(I_0,B)\dot I_{\Sigma_2}$ vanishes. Using
	\[
	\dot E_{\Sigma_2} 
	= \beta S I_0 - (\mu + a)E,
	\quad
	\dot B_{\Sigma_2}
	= rB\Bigl(1-\tfrac{B}{K}\Bigr)
	+ \frac{\lambda I_0 B}{\omega+I_0},
	\]
	we obtain
	\[
	\mathcal{L}^2_{F_{\Sigma_2}}L
	=
	a\bigl(\beta S I_0 - (\mu + a)E\bigr)
	-
	(\mu_1-\mu_0)\frac{I_0^2}{(B+I_0)^2}
	\left[
	rB\Bigl(1-\tfrac{B}{K}\Bigr)
	+ \frac{\lambda I_0 B}{\omega+I_0}
	\right].
	\]
	At the tangency set
	\[
	\Sigma^t_2 = \{F_{\Sigma_2,I} = 0\}
	\quad \Longleftrightarrow \quad
	aE_{c2}(B) 
	= \Bigl(\mu + \delta_1 + \gamma_1 
	+ \mu_0 + (\mu_1-\mu_0)\tfrac{B}{B+I_0}\Bigr) I_0,
	\]
	we obtain
	\[
	\mathcal{L}^2_{F_{\Sigma_2}}L\Big|_{\Sigma^t_2}
	=
	a\bigl(\beta S I_0 - (\mu + a)E_{c2}(B)\bigr)
	-
	(\mu_1-\mu_0)\frac{I_0^2}{(B+I_0)^2}
	\left[
	rB\Bigl(1-\tfrac{B}{K}\Bigr)
	+ \frac{\lambda I_0 B}{\omega+I_0}
	\right],
	\]
	which is generically nonzero.
	
	Hence, both tangency sets  
	\[
	\Sigma^t_1 = \{I = I_0,\ E = E_{c1}\}, 
	\quad 
	\Sigma^t_2 = \{I = I_0,\ E = E_{c2}(B)\},
	\]
	correspond to quadratic grazing (fold) tangencies wherever the corresponding second Lie derivative is nonzero. Moreover, while the delay
	enters the system through the healthcare capacity  dynamics, it does not affect the
	location of the tangency sets, since it does not appear in $F_{\Sigma_1,I}$ or
	$F_{\Sigma_2,I}$. The location of the tangency sets is therefore unchanged by the delay, although the second derivative for $F_{\Sigma_2}$ depends on the delayed infection level through the healthcare-capacity equation. For positive delay, the corresponding expression is obtained by replacing the capacity expansion term $\lambda I_0 B/(\omega+I_0)$ by $\lambda I(t-\tau)B/(\omega+I(t-\tau))$. The quadratic grazing classification remains valid whenever the corresponding second derivative exists and is nonzero. 
	\section{PINN-based parameter calibration and computational methodology}\label{supp:pinn}
	We describe the computational methodology for the illustrative calibration to 236 weekly cumulative reported COVID-19 observations from Kenya, obtained from WHO \cite{WHO_Data}.
	
	For this purpose, we employ the Physics-Informed Neural Networks (PINNs)
	framework \cite{pinn1}. For numerical consistency, the epidemiological and healthcare-capacity state variables are internally normalized by a reference population $N$ corresponding to the Kenya data period, namely
	\[
	\bar S=\frac{S}{N},\quad
	\bar E=\frac{E}{N},\quad
	\bar I=\frac{I}{N},\quad
	\bar H=\frac{H}{N},\quad
	\bar R=\frac{R}{N},\quad
	\bar B=\frac{B}{N},\quad
	\bar C=\frac{C}{N}.
	\]
	Similarly, the quantities $K$, $\omega$, and $I_0$ are normalized by the same population scale. The PINN therefore operates on a dimensionally consistent normalized representation of the original model, while the predicted quantities are transformed back to their original population scale for interpretation and evaluation.
	
	The PINN is implemented as a fully connected feed-forward neural network that maps time $t\in\mathbb{R}$ to the seven-dimensional normalized output vector
	\[
	(\bar S(t),\bar E(t),\bar I(t),\bar H(t),\bar R(t),\bar B(t),\bar C(t)), 
	\]
	where $\bar C(t)$ represents the normalized cumulative reported cases. The cumulative case variable is incorporated through 
	\begin{equation}
		\frac{d\bar C}{dt}=a\bar E. 
	\end{equation}
	The architecture consists of three hidden layers with 64 neurons each, using the hyperbolic tangent (Tanh) activation function. The input layer has one node corresponding to time, and the output layer has seven nodes corresponding to the normalized model compartments and cumulative reported cases. Positivity of the predicted state variables is imposed through a smooth positive output transformation. 
	
	\medskip
	\noindent \textbf{Model constraints and loss function:}
	The network is trained using a composite loss function that combines the physics-informed residual loss, the data-fitting loss, and an initial-condition loss. The physics-informed component is derived from the governing differential equations describing the compartmental dynamics, including the delayed healthcare-capacity response. In particular, the infectious-class residual is constructed from the complete model equation
	\[
	\frac{dI}{dt}
	=
	aE-\left(\mu+\delta_1+\gamma_1+T(I,B)\right)I,
	\]
	so that hospitalization is consistently represented as a removal from the infectious compartment. The piecewise hospitalization term is evaluated branchwise at each training time, with $T=\mu_1$ for $I\leq I_0$ and the capacity-dependent expression used for $I>I_0$, according to the infection level predicted by the network. The total loss function is defined as 
	\[
	\mathcal{L}_{\text{total}}
	=
	\mathcal{L}_{\text{physics}}
	+
	\alpha_{\text{data}}\mathcal{L}_{\text{data}}
	+
	\alpha_{\text{IC}}\mathcal{L}_{\text{IC}},
	\]
	where $\mathcal{L}_{\text{physics}}$ penalizes the scaled residuals of the governing model equations using automatic differentiation, $\mathcal{L}_{\text{data}}$ penalizes the discrepancy between the predicted and observed cumulative cases, and $\mathcal{L}_{\text{IC}}$ penalizes deviations from the imposed initial constraints. In the present implementation, 
	$
	\alpha_{\text{data}}=0.8, 
	\alpha_{\text{IC}}=10. 
	$ 
	The initial-condition component constrains the predicted cumulative-case value at the initial observation and imposes the normalized population condition
	\[
	\bar S(0)+\bar E(0)+\bar I(0)+\bar H(0)+\bar R(0)=1.
	\]
	A delay-interpolation module is used to evaluate the delayed infection term $\bar I(t-\tau)$ appearing in the healthcare-capacity equation. For times preceding the first observation, a constant pre-history,
	$ 
	\bar I(\theta)=\bar I(0), \theta \leq 0, 
	$ 
	is assumed.
	
	\medskip
	\noindent \textbf{Parameter calibration:} The transmission rate $\beta$, infection-responsive expansion rate $\lambda$, and response delay $\tau$ are optimized jointly with the network weights using Adam. The remaining parameters are fixed, and time and rates are expressed on a weekly scale.

	\noindent \textbf{Demographic parameters:}
	For the demographic component, we use the reference population
	$N=52{,}217{,}334$ \cite{KENYA_Pop}. The natural mortality rate
	is chosen to correspond to an average life expectancy of approximately
	68 years, consistent with recent demographic estimates for Kenya
	\cite{KenyaLifeExpectancy}. Using the reciprocal-life-expectancy
	approximation on a weekly time scale gives
	$
	\mu \approx \frac{1}{68\times52}
	\approx 0.00028~\text{week}^{-1}.
	$ 
	The demographic recruitment rate is then specified through the
	population-balance relation
	$
	\Lambda=\mu N,
	$ 
	so that, after population normalization,
	$ 
	\frac{\Lambda}{N}=\mu.
	$
	Consequently,
	$
	\Lambda
	=
	\mu N
	=
	1.4621\times10^{4}
	$
	individuals/week.
	
	Under the normalized formulation, the transmission parameter used internally by the PINN is $\widetilde{\beta}=\beta N$, and the corresponding count-scale transmission coefficient is recovered after training as 
	$
	\beta=\frac{\widetilde{\beta}}{N}.
	$ 
	For the final calibration,
	$
	\widetilde{\beta}\approx0.66045~\text{week}^{-1},
	$ 
	which corresponds to the count-scale transmission coefficient reported. The healthcare-capacity expansion rate $\lambda$ is constrained to remain positive, while the response delay $\tau$ is restricted to the interval $0.05\leq\tau\leq12$ weeks. The fixed and calibrated parameter values are summarized in  Table~\ref{context:parameter_values}. 
	
	\noindent \textbf{Training and evaluation:} Training uses 236 weekly cumulative-case observations, 10{,}000 epochs, and Adam with learning rate $10^{-3}$. The best total loss is $2.474008\times10^{-4}$. The fit yields RMSE $5971.5250$, MAE $3387.1916$, NRMSE $0.023867$, and $R^2=0.997663$. Since only cumulative cases are observed, the calibrated parameters are effective values; structural and practical identifiability are not established. The numerical bifurcation parameter sets are not calibrated epidemiological estimates.
	\section{Supplementary numerical simulations for the smooth subsystems}\label{supp:sim}
	\subsection{Positive-delay simulations for subsystem $F_{\Sigma_1}$}\label{supp:sim1}
\begin{figure}[H]
		\centering
		(a)\includegraphics[scale=0.3]{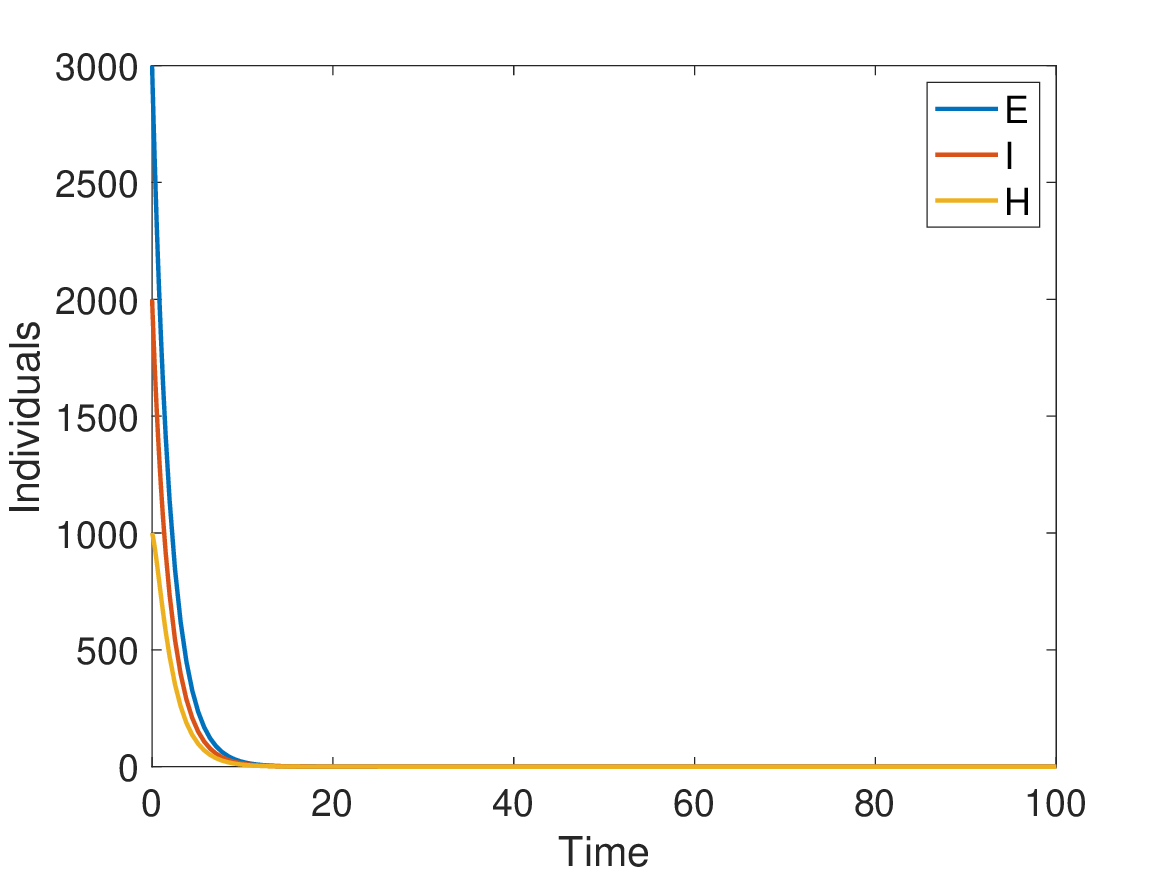} 
		(b)\includegraphics[scale=0.3]{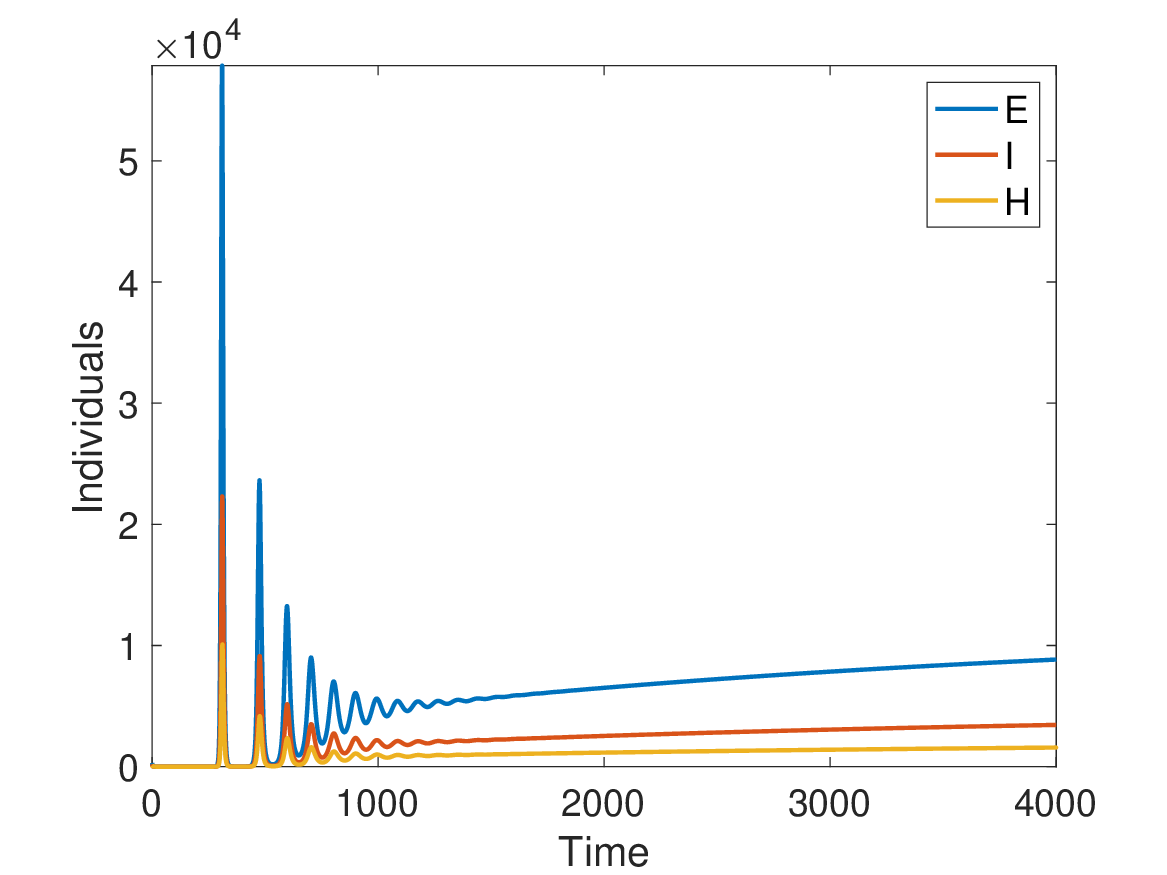}
		\caption{Dynamics of the delayed subsystem $F_{\Sigma_1}$ with time delay $\tau = 60$. 
			(a) The disease-free equilibrium ($D_{02}$) is asymptotically stable for 
			$\beta = 0.00000021$ and $\mathcal R_0 = 0.82 < 1$. 
			(b) The positive-capacity endemic equilibrium ($D^*$) is asymptotically stable for 
			$\beta = 0.000007$ and $\mathcal R_0 = 27.28 > 1$.}
		\label{delaysub1}
	\end{figure}

	Unless otherwise stated, the parameters are fixed as $\Lambda=1400$ and $\xi=0.001$, with all remaining values taken from Table~\ref{context:parameter_values}. Fig.~\ref{delaysub1} depicts the dynamics for a positive delay $\tau=60$. Here, $\tau$ represents the healthcare-capacity response delay in the number of infected individuals, so that healthcare capacity is adjusted according to the  previously reported infection level $I(t-\tau)$. The results indicate that the disease-free equilibrium $D_{02}$ remains asymptotically stable for $\mathcal R_0<1$, while the endemic equilibrium $D^*$ is asymptotically stable for the illustrated parameter set with $\mathcal R_0>1$. The corresponding stability behavior mirrors that of the undelayed case. Thus, the healthcare-capacity response delay does not alter the equilibrium or stability properties of subsystem $F_{\Sigma_1}$. From an epidemiological perspective, this suggests that, within this subsystem, healthcare-capacity response delays do not change the long-term outcome, i.e. the disease is eliminated when $\mathcal R_0<1$ and persists at a steady endemic level for the illustrated parameter set with $\mathcal R_0>1$. 
	
	\subsection{Disease-free equilibrium simulations for the smooth extension of subsystem $F_{\Sigma_2}$}\label{supp:sim2}
\begin{figure}[H]
		\centering
		(a)\includegraphics[scale=0.3]{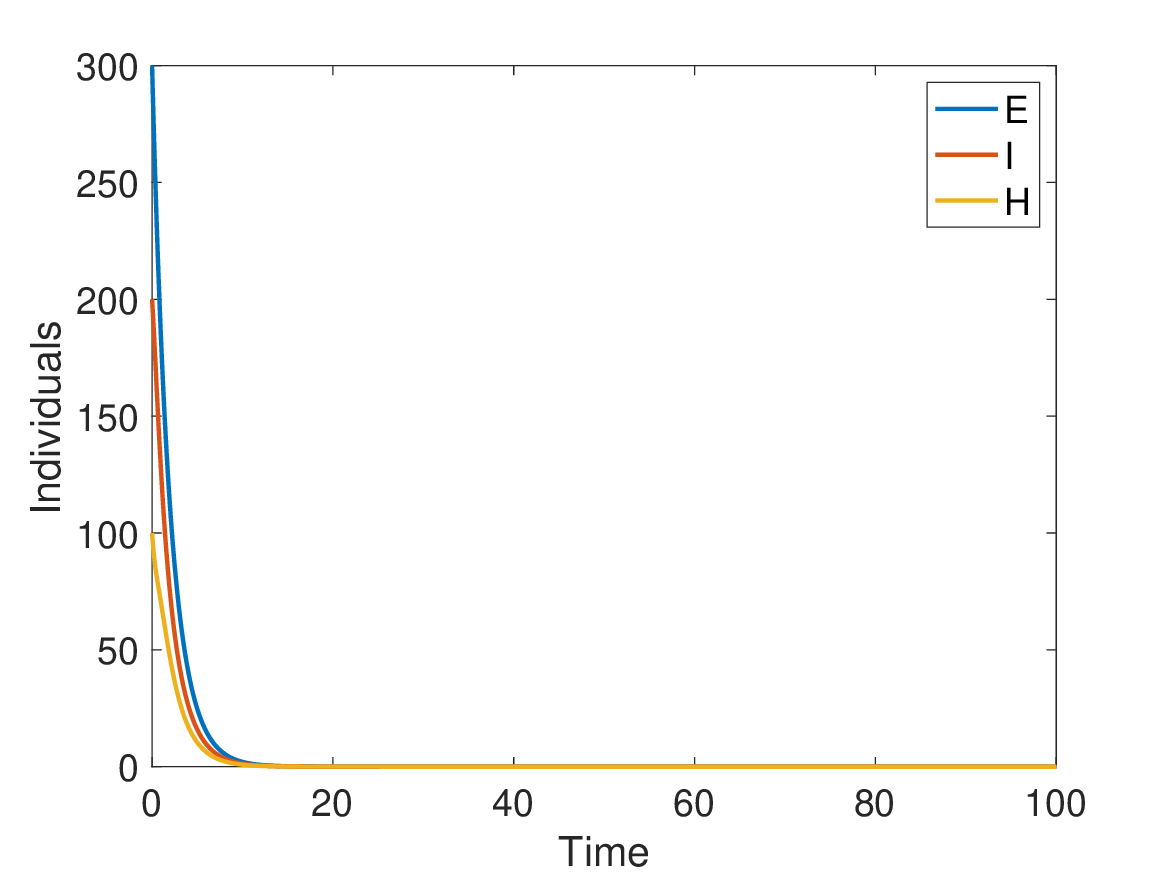} 
		(b)\includegraphics[scale=0.3]{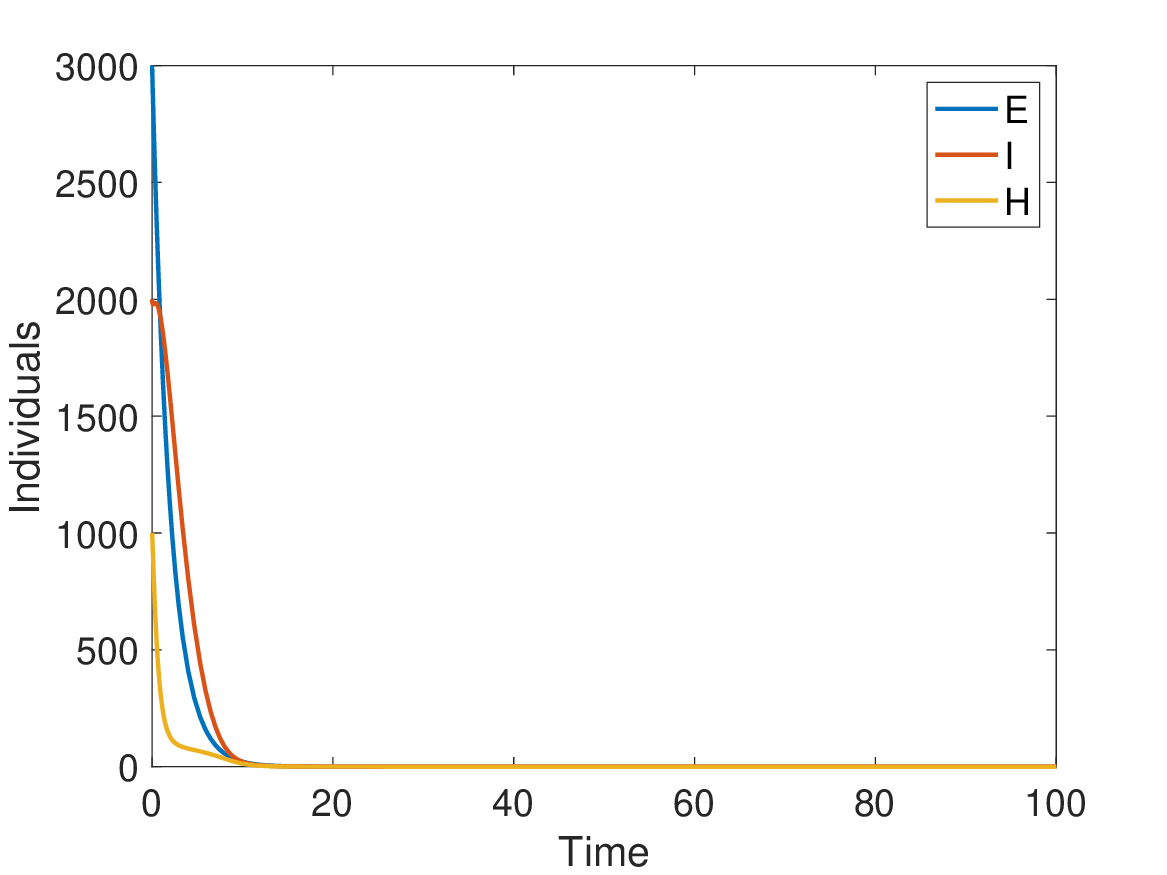} 
		\caption{Time evolution of the non-delayed (left panel) and delayed (right panel) subsystem $F_{\Sigma_2}$, illustrating asymptotic stability of the disease-free equilibrium ($D_{02}$) for $\mathcal R_0 = 0.82 < 1$.}
		\label{dfe_delaysub2}
	\end{figure}

	We now illustrate and validate the analytical stability and bifurcation results for the subsystem $F_{\Sigma_2}$. We begin by illustrating that the disease-free equilibrium $D_{02}$, which is virtual for subsystem $F_{\Sigma_2}$, is locally asymptotically stable in its smooth extension when $\mathcal R_0 < 1$, both in the absence and presence of time delay, as shown in Fig.~\ref{dfe_delaysub2}. Epidemiologically, this means that the disease is eliminated from the population, since the exposed, infectious, and hospitalized classes approach zero. The delay does not destabilize the disease-free equilibrium in this case because it acts only through the healthcare-capacity dynamics. For this simulation, we fixed $\Lambda = 1400$, $\xi = 0.001$, $\beta =  2.1\times10^{-7}$, and $\tau = 60$, while all remaining parameters were taken from Table~\ref{context:parameter_values}. 
	
	\section{Stability of sliding pseudo-equilibria}\label{supp:pseudo_stability}
	On the switching manifold $I=I_0$, the reduced tangential system is given by Eq.~(4.6) of the main text. Its admissibility condition is $E_{c2}(B^s)<E^s<E_{c1}$, where the two boundaries are defined in the main text. We restate the result before giving its proof.
	
	\noindent\textbf{Theorem~4.2 (main text).} Assume that $D^s$ is biologically feasible and satisfies the admissibility condition (4.7). Then $D^s$ is globally asymptotically stable for the algebraically extended reduced tangential system (4.6) in the region $B>0$. Consequently, every trajectory of the reduced sliding dynamics with $B(0)>0$ that remains in the admissible sliding region $\Sigma_e$ converges to $D^s$. If $D^{0s}$ is biologically feasible and admissible within $\Sigma_e$, then $D^{0s}$ is unstable with respect to the reduced sliding dynamics.
	
	\begin{proof}
		For $I=I_0$, the $(S,E,H,R)$-components of the reduced tangential system (4.6) form the affine system 
		\[
		\frac{dX}{dt}=AX+b,
		\quad
		X=(S,E,H,R)^{\top},
		\]
		where
		\[
		A=
		\begin{pmatrix}
			-(\beta I_0+\mu) & 0 & 0 & \xi \\[2pt]
			\beta I_0 & -(\mu+a) & 0 & 0 \\[2pt]
			0 & a & -(\mu+\delta_2+\gamma_2) & 0 \\[2pt]
			0 & 0 & \gamma_2 & -(\mu+\xi)
		\end{pmatrix},
		\]
		and
		\[
		b=
		\begin{pmatrix}
			\Lambda\\[2pt]
			0\\[2pt]
			-(\mu+\delta_1+\gamma_1)I_0\\[2pt]
			\gamma_1I_0
		\end{pmatrix}.
		\]
		The matrix $A$ is Metzler, and hence $-A$ is a $Z$-matrix. Because the only directed feedback cycle in $A$ involves all four variables, every proper principal minor of $-A$ is a product of positive diagonal factors and is therefore strictly positive. The determinant of $-A$ is
		\begin{align*}
			\det(-A)
			={}&
			(\beta I_0+\mu)(\mu+a)(\mu+\delta_2+\gamma_2)(\mu+\xi)
			-\xi a\gamma_2\beta I_0\\
			={}&
			\gamma_2\mu
			\Big[
			a(\beta I_0+\mu+\xi)
			+(\mu+\xi)(\beta I_0+\mu)
			\Big]\\
			&\quad
			+(a+\mu)(\delta_2+\mu)(\mu+\xi)(\beta I_0+\mu)
			>0.
		\end{align*}
		Thus, all principal minors of $-A$ are positive. Therefore, $-A$ is a nonsingular $M$-matrix, and consequently all eigenvalues of $A$ have strictly negative real parts. Hence, $A$ is Hurwitz.
		
		It follows that the affine subsystem has a unique equilibrium 
		$
		X^s=(S^s,E^s,H^s,R^s)^{\top},
		$
		and for every initial condition,
		\[
		X(t)-X^s
		=
		e^{At}\bigl(X(0)-X^s\bigr)
		\longrightarrow0
		\quad\text{as }t\to\infty.
		\]
		Thus, the equilibrium of the $(S,E,H,R)$-subsystem is globally asymptotically stable.
		
		The healthcare-capacity equation along the reduced tangential system is
		\[
		\frac{dB}{dt}
		=
		\left(
		r+\frac{\lambda I_0}{\omega+I_0}
		\right)B
		-\frac{r}{K}B^2.
		\]
		This scalar logistic equation has the equilibria
		\[
		B=0
		\quad\text{and}\quad
		B=B^s
		=
		\frac{K}{r}
		\left(
		r+\frac{\lambda I_0}{\omega+I_0}
		\right).
		\]
		For every $B(0)>0$, the positive equilibrium $B^s$ is globally asymptotically stable. Consequently,
		\[
		(S(t),E(t),H(t),R(t),B(t))
		\longrightarrow
		(S^s,E^s,H^s,R^s,B^s)
		\]
		for every initial condition with $B(0)>0$ in the algebraically extended reduced tangential system.
		
		Therefore, $D^s$ is globally asymptotically stable for the algebraically extended reduced system in the region $B>0$. When this vector field represents an actual Filippov sliding trajectory, the conclusion applies for as long as the trajectory remains in the admissible sliding region $\Sigma_e$. 
		
		For $D^{0s}$, linearization of the $B$-equation at $B=0$ gives
		\[
		\frac{db}{dt}
		=
		\left(
		r+\frac{\lambda I_0}{\omega+I_0}
		\right)b,
		\]
		and therefore the corresponding eigenvalue is
		\[
		s_B
		=
		r+\frac{\lambda I_0}{\omega+I_0}>0.
		\]
		Hence, whenever $D^{0s}$ is admissible, it is unstable with respect to the reduced sliding dynamics.
	\end{proof}
	
	\section{Effect of delay on sliding pseudo-equilibria}\label{supp:sliding_delay_stability}
	For histories on the switching manifold, $I(t)=I(t-\tau)=I_0$ and the delayed tangential system (Eq.~(4.5)) reduces to Eq.~(4.6). The stability result in Supplementary Section~\ref{supp:pseudo_stability} therefore applies to these histories.
	
	\noindent\textbf{Proposition~4.5 (main text).} For every $\tau\geq0$ and for histories contained in the switching
	manifold, the healthcare-capacity response delay does not alter the
	stability properties of the admissible pseudo-equilibria with respect to the
	algebraically extended reduced tangential dynamics. In particular, $D^{0s}$ remains unstable in
	the $B$-direction, whereas $D^{s}$ retains the global asymptotic
	stability for $B>0$ established in Theorem~4.2.
	
	\begin{proof}
		Along the reduced sliding dynamics, $I(t)=I(t-\tau)=I_0$. At $D^{0s}$, 
		the linearized $B$-equation is 
		\[
		\frac{db}{dt}
		=
		\left(
		r+\frac{\lambda I_0}{\omega+I_0}
		\right)b(t),
		\]
		whose coefficient is positive. Hence, $D^{0s}$ remains unstable in the
		$B$-direction.
		
		At $D^{s}$, the delayed perturbation of the infectious component
		vanishes for histories contained in the switching manifold, and the
		linearized $B$-equation reduces to
		\[
		\frac{db}{dt}
		=
		-\left(
		r+\frac{\lambda I_0}{\omega+I_0}
		\right)b(t).
		\]
		The remaining $(S,E,H,R)$-subsystem is identical to the delay-free
		reduced tangential system. Therefore, the stability result of
		Theorem~4.2 remains valid for every
		$\tau\geq0$.
	\end{proof}
	
	\section{Uniform disease elimination under a worst-case threshold}\label{supp:global_dfe_worstcase}
	For the full Filippov system, define
	\[
	\mathcal R_0^{\max}=\frac{\beta(\Lambda/\mu)a}{(\mu+a)(\mu+\delta_1+\gamma_1+\mu_0)}.
	\]
	The hospitalization rate is at least $\mu_0$ in both smooth regimes, and the same bound holds for the effective hospitalization rate under convexification on the switching manifold.
	
	\noindent\textbf{Theorem~4.6 (main text).} Suppose that $\mathcal R_0^{\max}<1$. Then, for every admissible
	nonnegative initial history and every delay $\tau\geq0$: 
	\begin{itemize}
		\item if $B(0)=0$, then $B(t)=0$ for all $t\geq0$, and the solution
		converges to
		$ 
		D_{01}=(\Lambda/\mu,0,0,0,0,0);
		$ 
		\item if $B(0)>0$, then $B(t)\to K$ and the solution converges to
		$
		D_{02}=(\Lambda/\mu,0,0,0,0,K).
		$ 
	\end{itemize}
	
	\begin{proof}
		Consider the Lyapunov function
		\[
		V(E,I)=E+\frac{\mu+a}{a}I,
		\]
		which is nonnegative in the biologically feasible region and vanishes
		only at $(E,I)=(0,0)$. Along any classical trajectory in either smooth
		region, $I<I_0$ or $I>I_0$, we obtain
		\begin{align*}
			\frac{dV}{dt}
			&=\frac{dE}{dt}
			+\frac{\mu+a}{a}\frac{dI}{dt}\\
			&=\bigl(\beta SI-(\mu+a)E\bigr)
			+\frac{\mu+a}{a}
			\bigl(aE-(\mu+\delta_1+\gamma_1+T(I,B))I\bigr)\\
			&=
			I\left( 
			\beta S
			-\frac{\mu+a}{a}
			\bigl(\mu+\delta_1+\gamma_1+T(I,B)\bigr)
			\right). 
		\end{align*}
		From $\dot N\leq\Lambda-\mu N$, we have
		$\limsup_{t\to\infty}S(t)\leq\Lambda/\mu$. Since
		$\mathcal R_0^{\max}<1$, choose $\varepsilon>0$ such that
		$\beta\varepsilon<\frac{(\mu+a)(\mu+\delta_1+\gamma_1+\mu_0)}{a}(1-\mathcal R_0^{\max})$.
		Then there exists $t_\varepsilon\geq0$ such that
		$S(t)\leq\Lambda/\mu+\varepsilon$ for $t\geq t_\varepsilon$.
		Also, $T(I,B)\geq\mu_0$ for all admissible states. Hence, for $t\geq t_\varepsilon$,
		\begin{align*}
			\frac{dV}{dt}
			&\leq
			I\left(
			\beta\left(\frac{\Lambda}{\mu}+\varepsilon\right)
			-\frac{\mu+a}{a}
			(\mu+\delta_1+\gamma_1+\mu_0)
			\right)\\
			&=
			\frac{(\mu+a)(\mu+\delta_1+\gamma_1+\mu_0)}{a}
			(\mathcal R_0^{\max}-1)I+\beta\varepsilon I.
		\end{align*}
		If $\mathcal R_0^{\max}<1$, there exists $k>0$ such that
		\[
		\frac{dV}{dt}\leq-kI
		\]
		for almost every $t\geq t_\varepsilon$.
		
		On $\Sigma$, the Filippov vector field belongs to the convex hull of the
		two one-sided vector fields. Since $\nabla V\cdot F$ is linear in $F$,
		the same estimate
		\[
		\frac{dV}{dt}\leq-kI
		\]
		holds along the Filippov convexification on the switching manifold.
		Therefore, $V(t)$ is non-increasing for $t\geq t_\varepsilon$ and bounded below, and
		\[
		\int_0^\infty I(t)\,dt<\infty.
		\]
		Because all state variables remain bounded by Theorem~2.1 and the one-sided
		vector fields are bounded on bounded subsets of the biologically feasible
		region, $I(t)$ is uniformly continuous along every Filippov solution.
		Hence, by Barbalat's lemma \cite{MIN2007},
		$
		I(t)\to0
		\quad\text{as }t\to\infty.
		$ 
		
		The exposed equation
		\[
		\frac{dE}{dt}
		=
		\beta SI-(\mu+a)E
		\]
		is a stable linear equation forced by the term $\beta SI$, which tends to
		zero. Consequently,
		$
		E(t)\to0.
		$ 
		Using $E(t)\to0$ and $I(t)\to0$ in the remaining epidemiological
		equations then gives
		$
		H(t)\to0, 
		R(t)\to0, 
		S(t)\to\frac{\Lambda}{\mu}.
		$
		
		Finally, since $\tau$ is fixed and $I(t)\to0$, we also have
		$
		I(t-\tau)\to0.
		$ 
		The healthcare-capacity equation is therefore asymptotically governed by
		\[
		\frac{dB}{dt}
		=
		rB\left(1-\frac{B}{K}\right).
		\]
		If $B(0)=0$, the multiplicative structure of the $B$-equation implies
		$
		B(t)=0~\text{for all}~t\geq0.
		$ 
		If $B(0)>0$, positivity of $B(t)$ is preserved. For any $\epsilon>0$,
		the per capita infection-responsive expansion term is eventually
		between $0$ and $\epsilon$. Comparison with logistic equations with
		growth rates $r$ and $r+\epsilon$ and the common quadratic coefficient
		$r/K$ gives $K\leq\liminf B(t)\leq\limsup B(t)\leq K(1+\epsilon/r)$.
		Letting $\epsilon\to0^+$ gives
		$ 
		B(t)\to K.
		$ 
		Therefore, the corresponding solutions converge to $D_{01}$ and
		$D_{02}$, respectively. This completes the proof. 
	\end{proof}